\documentclass[11pt,a4paper]{article}
\usepackage{amsmath,amsthm,amssymb,mathtools}
\usepackage{fontspec}

\usepackage[mathscr]{euscript}

\usepackage[a4paper,left=24mm,right=24mm,top=23mm,bottom=24mm,headheight=23pt,headsep=6mm]{geometry}

\usepackage{xcolor,titlesec,titling,tocloft,fancyhdr,needspace,aliascnt}
\usepackage{graphicx,booktabs,enumitem,float}
\setlist{itemsep=1.2pt,parsep=0pt,topsep=4pt}
\AtBeginDocument{\setlength{\abovedisplayskip}{8pt plus 2pt minus 2pt}%
  \setlength{\belowdisplayskip}{8pt plus 2pt minus 2pt}}
\definecolor{PaperInk}{HTML}{213747}
\definecolor{PaperAccent}{HTML}{2B6570}
\definecolor{PaperMuted}{HTML}{667882}
\definecolor{PaperRule}{HTML}{CBD7DC}
\usepackage[colorlinks=true,linkcolor=PaperAccent,citecolor=PaperAccent,urlcolor=PaperAccent,bookmarksnumbered=true,bookmarksopen=true,bookmarksopenlevel=1]{hyperref}
\hypersetup{pdftitle={Tidally forced traveling waves in a curved annular domain},pdfsubject={Direct supercriticality estimates, global continuation, crest stagnation and imperfect bifurcation}}
\allowdisplaybreaks[2]
\numberwithin{equation}{section}
\titleformat{\section}[hang]{\sffamily\Large\bfseries\color{PaperInk}}{\thesection}{.75em}{}[\vspace{.3ex}{\color{PaperRule}\titlerule}]
\titleformat{\subsection}[hang]{\sffamily\large\bfseries\color{PaperAccent}}{\thesubsection}{.7em}{}
\titlespacing*{\section}{0pt}{2.0\baselineskip}{.8\baselineskip}
\titlespacing*{\subsection}{0pt}{1.35\baselineskip}{.55\baselineskip}
\AddToHook{cmd/subsection/before}{\Needspace{6\baselineskip}}
\AddToHook{cmd/section/before}{\Needspace{12\baselineskip}}
\newtheoremstyle{paperplain}{6pt plus 2pt}{4pt plus 1pt}{\itshape}{} {\sffamily\bfseries\color{PaperInk}}{.}{.65em}{}
\newtheoremstyle{paperdefinition}{6pt plus 2pt}{4pt plus 1pt}{\normalfont}{} {\sffamily\bfseries\color{PaperInk}}{.}{.65em}{}
\theoremstyle{paperplain}
\newtheorem{theorem}{Theorem}[section]
\newaliascnt{lemma}{theorem}
\newtheorem{lemma}[lemma]{Lemma}
\aliascntresetthe{lemma}
\newaliascnt{proposition}{theorem}
\newtheorem{proposition}[proposition]{Proposition}
\aliascntresetthe{proposition}
\newaliascnt{corollary}{theorem}
\newtheorem{corollary}[corollary]{Corollary}
\aliascntresetthe{corollary}
\theoremstyle{paperdefinition}
\newaliascnt{definition}{theorem}
\newtheorem{definition}[definition]{Definition}
\aliascntresetthe{definition}
\newaliascnt{remark}{theorem}
\newtheorem{remark}[remark]{Remark}
\aliascntresetthe{remark}

\AddToHook{env/theorem/before}{\Needspace{5\baselineskip}}
\AddToHook{env/lemma/before}{\Needspace{4\baselineskip}}
\AddToHook{env/proposition/before}{\Needspace{5\baselineskip}}
\AddToHook{env/definition/before}{\Needspace{4\baselineskip}}
\renewcommand{\headrulewidth}{.35pt}
\renewcommand{\headrule}{\hbox to\headwidth{\color{PaperRule}\leaders\hrule height\headrulewidth\hfill}}

\date{}
\pretitle{\begin{flushleft}\sffamily\bfseries\color{PaperInk}\fontsize{19}{23}\selectfont}
\posttitle{\par\end{flushleft}\vspace{.8ex}{\color{PaperAccent}\hrule height 1pt}\vspace{1.2ex}}
\preauthor{}\postauthor{}\predate{}\postdate{}
\renewenvironment{abstract}{\par\noindent{\sffamily\bfseries\color{PaperAccent}Abstract}\par\medskip\small\noindent\ignorespaces}{\par\normalsize}

\AddToHook{cmd/maketitle/after}{\thispagestyle{empty}}

\AddToHook{cmd/appendix/before}{\clearpage}
\AddToHook{env/thebibliography/before}{%
  \Needspace{8\baselineskip}
  \markboth{References}{}%
  \phantomsection\addcontentsline{toc}{section}{References}}
\AddToHook{env/thebibliography/begin}{\small}

\usepackage[nameinlink,noabbrev]{cleveref}
\newaliascnt{example}{theorem}
\newtheorem{example}[example]{Example}
\aliascntresetthe{example}
\theoremstyle{paperplain}

\crefname{mainthm}{Theorem}{Theorems}
\Crefname{mainthm}{Theorem}{Theorems}
\AddToHook{env/mainthm/before}{\Needspace{6\baselineskip}}
\newcommand{\T}{\mathbb{T}}
\newcommand{\R}{\mathbb{R}}
\newcommand{\C}{\mathbb{C}}
\newcommand{\N}{\mathbb{N}}
\newcommand{\eps}{\varepsilon}
\newcommand{\p}{\partial}
\newcommand{\e}{\mathrm{e}}
\newcommand{\dd}{\,\mathrm{d}}
\newcommand{\Dop}{\mathcal{D}}
\newcommand{\Bop}{\mathcal{B}}
\newcommand{\Fop}{\mathcal{F}}
\newcommand{\Lop}{\mathcal{L}}
\newcommand{\Mop}{\mathcal{M}}
\newcommand{\Xs}{\mathcal{X}}
\newcommand{\Ys}{\mathcal{Y}}
\newcommand{\Zs}{\mathcal{Z}}
\newcommand{\Uo}{\mathcal{U}}
\newcommand{\cc}{\mathsf{c}}
\newcommand{\Ptide}{P}
\DeclareMathOperator{\spn}{span}
\DeclareMathOperator{\Fix}{Fix}
\DeclareMathOperator{\sgn}{sgn}

\begin{document}
\title{Tidally forced traveling waves in a curved annular domain:
direct estimates, global continuation and imperfect bifurcation}
\author{\begin{minipage}{\linewidth}\raggedright
\large Yiqiu Mao\textsuperscript{1} and Quan Wang\textsuperscript{2,*}\\[4pt]
\small \textsuperscript{1}School of Mathematics and Information Science, Guangzhou University, Guangzhou, 510000, China\\
\textsuperscript{2}College of Mathematics, Sichuan University, Chengdu, Sichuan, 610065, China\\
\textsuperscript{*}Corresponding author: \texttt{xihujunzi@scu.edu.cn}
\end{minipage}}
\hypersetup{pdfauthor={Yiqiu Mao and Quan Wang}}
\maketitle

\begin{abstract}
We study steady traveling waves in an annular fluid subject to inward gravity and a quadrupolar tidal potential. First, for the unforced problem with constant relative vorticity, we compute the Stokes coefficient and obtain explicit depth ranges in which the bifurcation is supercritical for all admissible fluxes and vorticities. We then establish global continuation of each branch. Along the unforced global curve, gravity, height and relative speed stay bounded, and the minimum magnitude of the azimuthal velocity relative to the satellite tends to zero as the continuation parameter tends to infinity. In the irrotational case we give conditions on depth and flux under which stagnation occurs at the crest. Under an additional slope bound, each subsequential limiting wave has a symmetric $120^\circ$ corner. Finally, we determine how weak tidal forcing changes the local bifurcation diagram. At an odd mode $n$ the two reflection-symmetric branches have critical parameters differing by $O(\eps^n)$. At the principal resonance, for $\eps>0$ small and nonzero cubic coefficient, a single fold at distance of order $\eps^{2/3}$ from resonance separates the regions with one and three nearby steady waves in the full periodic space. In the thin-layer limit we recover the resonance condition of Airy's canal theory.
\end{abstract}

\normalsize
\medskip
\begingroup\setcounter{tocdepth}{1}\tableofcontents\endgroup\clearpage
\section{Introduction}\label{sec:intro}

In an equatorial canal a moving tidal potential and a freely propagating long gravity wave can have the same speed. Airy's linear theory predicts a pole in the response at this speed and a reversal in the phase of the tide as the speed passes through resonance \cite{Airy1842,Lamb1932,Hendershott1981,Tokieda2013}. Our purpose is to determine how the weak tide continues through this resonance and how many steady profiles travel with the forcing at a given value of gravity. To this end, we consider the free-boundary Euler equations in a circular annulus. The rigid core has radius one, gravity acts radially inward, and a satellite supplies a degree-two potential. In the satellite frame this potential has a mean and a second angular harmonic. At fixed $\alpha$, varying the forcing holds the flux and Bernoulli constant fixed. When gravity varies at fixed $(p_0,\gamma,\kappa)$, equation~\eqref{eq:lambda} instead gives $\lambda\,'(\alpha)=e^\gamma$. In the fixed-gravity calculation, the mean component changes the laminar level. The second harmonic couples directly to the resonant wave.

A related question concerns the free waves on which the tide acts. For a horizontal bed and vertical gravity, Constantin and Strauss constructed connected families of periodic waves with vorticity by fixing the fluid domain in streamline coordinates \cite{CS2004}. Constantin, Strauss and V\u{a}rv\u{a}ruc\u{a} used conformal variables and analytic continuation for constant vorticity, allowing critical layers \cite{CSV2016}. Kozlov and Lokharu obtained limiting waves in finite depth by fixing the Bernoulli constant and varying the wavelength \cite{KozlovLokharu2023}. A further question concerns the shape of a wave near a stagnation point. Amick, Fraenkel and Toland proved the Stokes corner conjecture for the extreme irrotational wave \cite{AFT1982}. V\u{a}rv\u{a}ruc\u{a} and Weiss developed a blow-up analysis in the physical variables \cite{VarvarucaWeiss2011}. These problems and the roles of the different continuation parameters are surveyed in \cite{Haziot2022}.

Curvature changes both the laminar family and the dynamic boundary condition. Exact azimuthal Euler flows on a rotating sphere exhibit the interaction of curvature, rotation and the free surface \cite{CJ2016eq,HenryMartin2019}. For a two-dimensional annulus, Li and Wang \cite{LW,LWpreprint} bifurcated travelling waves by varying the radial gravity parameter. In their model the wave angular velocity is half the absolute vorticity, so the relative flow is irrotational. At the first angular mode they obtained a local pitchfork and, by degree theory, a global connected set subject to continuation alternatives. Their preprint also classified a shallow example as subcritical \cite[Example~6.2]{LWpreprint}. We retain their flux and Bernoulli convention for this comparison. We allow an independent satellite frequency, so that the relative vorticity $2\kappa$ can be nonzero, and ask what controls the direction and endpoint of each unforced branch and how the tide changes its local structure. References to specific equations and examples in that comparison are to the preprint.

\subsection*{The unforced branches}
We write $S_0=e^\gamma$ for the laminar surface radius, $w$ for its relative angular velocity and $r=w^2$. The laminar azimuthal velocity must be nonzero throughout the closed fluid layer. With $\alpha$ denoting dimensionless gravity, each angular mode $k\ge1$ has a simple critical value in the even subspace,
\[
 \alpha_k=S_0\{1-2\kappa w+r[1+k\coth(k\gamma)]\}.
\]
The local branch has $\alpha=\alpha_k+\alpha_k^{(2)}a^2+O(a^4)$, where $a$ is the $k$th cosine coefficient of the logarithmic surface radius. A single harmonic extension on the log-polar strip gives both this dispersion relation and the cubic solvability condition. To determine the direction of bifurcation, we must determine the sign of the cubic coefficient. Writing $\nu=\kappa/w$ makes the coefficient a positive factor times a quadratic in $r$. Its discriminant is itself quadratic in $\nu$, because its two highest powers cancel.

Section~\ref{sec:unforced} uses that cancellation to prove primary supercriticality at every flux for $\kappa=0$, $\gamma\le10.7959$, and at every admissible vorticity for $\gamma\le4.997$. In particular, the branch in Li and Wang's Example~6.2 turns towards increasing gravity. Explicit deeper examples turn in the opposite direction. For every mode, $k\gamma\le5/3$ suffices at zero vorticity and $k\gamma\le6/5$ suffices at arbitrary admissible vorticity. These bounds have direct proofs. The sharper five finite-mode bounds in Theorem~\ref{thm:kge2} use interval arithmetic.

Each local unforced branch continues along a curve, locally parametrised by analytic arcs, of waves with strictly decreasing streamlines between crest and trough. Gravity, height and speed stay bounded, while the minimum magnitude of the relative azimuthal velocity tends to zero. For a steep free surface, small azimuthal velocity alone does not imply stagnation, since the radial velocity may remain positive in magnitude. Moreover the centrifugal term makes the Bernoulli head a nonlinear function of height. Its derivative can change sign near the trough, where a direct maximum-principle argument for crest stagnation would fail. In Sections~\ref{sec:global} and~\ref{sec:stagnation} a conserved flow force controls the trough. After normalisation, ruling out the offending levels reduces to integrating a one-parameter family of profiles from its maximum to its first zero. This places stagnation at the crest for $\kappa=0$, $0<\gamma\le1.315$, $r\ge\gamma^2/(1+\gamma)^2$, with restricted-flux extensions up to $\gamma=1.4$. Under the additional slope bound of Corollary~\ref{cor:slope}(b), every subsequential limiting wave is smooth away from its crests and has a symmetric $120^\circ$ corner there.

For nonzero relative vorticity, differences of stream functions with the same vorticity are still harmonic. This permits comparison with laminar flows, while an inhomogeneous Cauchy--Riemann equation gives compactness away from stagnation. Section~\ref{sec:kappaglobal} proves global continuation for every admissible triple in $\alpha\in\R$. The curve stays in $\alpha>0$ if $|\kappa|\le1$ or $\kappa/w\le0$. Under the two head inequalities of Theorem~\ref{thm:kstag}, degeneration occurs on the crest line. For $\kappa/w\le0$ it occurs at the surface crest.

\subsection*{The forced branches}
The degree-two tide preserves the two axial reflections and the half-turn. Away from the even spectrum there is a unique nearby half-turn-invariant tidal equilibrium. At an odd mode $n$ each reflection subspace retains a pitchfork. At order $j$ in the forcing amplitude $\eps$, the Fourier index can change by at most $2j$. Coupling the modes $-n$ and $n$ therefore requires at least $n$ orders, giving $\alpha_n^+(\eps)-\alpha_n^-(\eps)=O(\eps^n)$. The leading splitting coefficient can vanish. At mode two the forcing already projects onto the kernel at first order. This is the setting of imperfect bifurcation \cite{GolubitskySchaeffer1979,GSS1988}, with the added issue that an even-subspace calculation alone would miss any wave whose crest is displaced from the forcing axes.

\begin{figure}[htbp]\centering
\includegraphics[width=.96\linewidth]{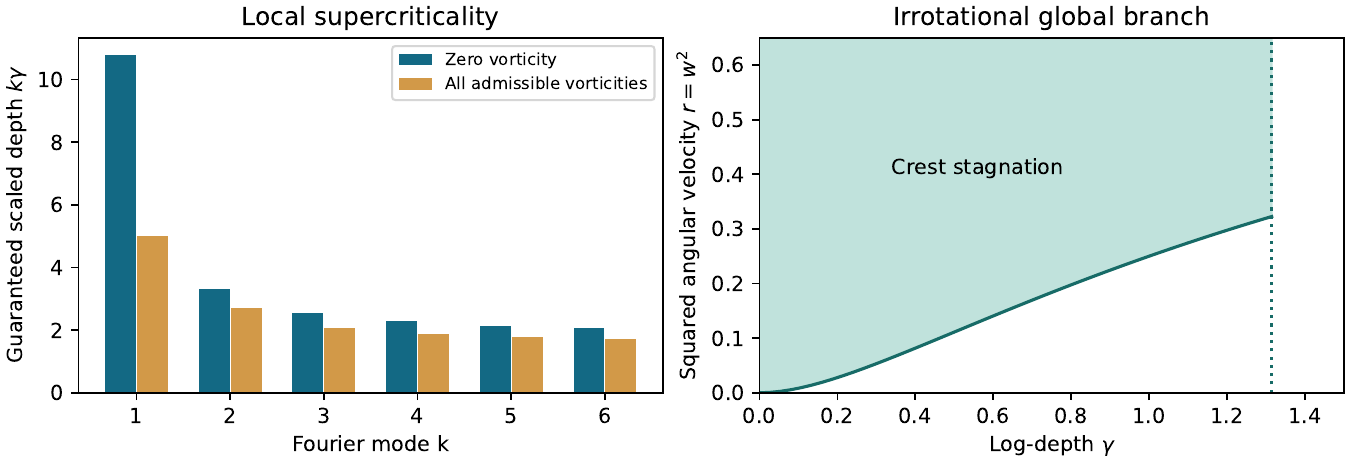}
\caption{Left: each bar covers the proved supercritical interval from zero (excluded) to the indicated scaled depth, uniformly in admissible flux. Right: the shaded portion of $0<\gamma\le1.315$, $r\ge(\gamma/(1+\gamma))^2$ guarantees crest stagnation for the unforced irrotational branch. The upper edge is a plotting limit. Both panels show sufficient regions; the corner conclusion requires the additional slope hypothesis.}\label{fig:ranges}
\end{figure}

Section~\ref{sec:forced} resolves the phase in the full two-dimensional kernel before deriving the scalar equation. If $\alpha_2^{(2)}\ne0$ and $\eps>0$ is small, there is one local fold at distance of order $\eps^{2/3}$ from $\alpha_2$, with amplitude of order $\eps^{1/3}$. On one side there is one nearby steady wave, on the other there are three. For $\alpha_2^{(2)}>0$ the three-wave side is the side of larger gravity. Section~\ref{sec:earth} relates this classification to the canal limit and illustrates it with computed profiles. Global continuation and the limiting-wave theorems concern zero forcing.
Figure~\ref{fig:ranges} compares the local sign estimates with a sufficient global crest criterion.
\section{The annular Euler problem and its spectrum}\label{sec:formulation}

\subsection{The physical problem}
Let $a>0$ be the radius of a rigid core. At time $t$ an incompressible inviscid fluid of
density $\rho$ occupies $\{x\in\R^2:\ a<|x|<\eta(t,\theta)\}$, where $(r,\theta)$ are polar
coordinates. The fluid has constant vorticity $2\omega_v$, and its velocity $v$ and pressure $P$
satisfy
\begin{equation}\label{eq:euler}
v_t+(v\cdot\nabla)v=-\frac1\rho\nabla P-g\,\mathbf e_r-\nabla\Phi,\qquad
\nabla\cdot v=0,\qquad \nabla\times v=2\omega_v,
\end{equation}
with the boundary conditions
\begin{equation}\label{eq:euler-bc}
\eta_t+\frac{v_\theta}{r}\eta_\theta=v_r\ \text{ and }\ P=p_{\rm atm}\ \text{ on } r=\eta(t,\theta),\qquad
v_r=0\ \text{ on } r=a .
\end{equation}

Here $\Phi$ is the tidal potential of a satellite of mass $M$ that revolves in the plane of
the flow on a circle of radius $d$ with angular velocity $\Omega_f>0$. Its Newtonian potential
is $-GM/|x-x_M(t)|$. In the frame of the planet's centre we remove the constant term and the
uniform field (the degree-one term, which is balanced by the planet's own orbital
acceleration), and truncate the multipole expansion at degree two:
\begin{equation}\label{eq:Phi}
\Phi(t,x)=-\frac{GM}{d^3}\,|x|^2P_2\big(\cos(\theta-\Omega_ft)\big),\qquad P_2(\xi)=\tfrac12(3\xi^2-1).
\end{equation}
The next multipole has an additional factor $r/d$, hence is uniformly smaller
by $a/d$ only when the dimensionless radius remains $O(1)$. Uniform validity
throughout the fluid requires $a\max S/d\ll1$.

We scale lengths by $a$, time by $\Omega_f^{-1}$, velocities by $a\Omega_f$, and pressure by
$\rho a^2\Omega_f^2$, and set
\[
\alpha=\frac{g}{a\Omega_f^2},\qquad \eps=\frac{GM}{d^3\Omega_f^2},\qquad \kappa=\frac{\omega_v}{\Omega_f}-1 .
\]
We look for traveling waves co-rotating with the satellite: in the variables
$\Theta=\theta-t$, $R=r$ all unknowns are independent of $t$. The relative velocity
$\mathbf w=(U,W)$, with radial component $U$ and azimuthal component $W=V-R$, has the constant
vorticity $2\omega_v/\Omega_f-2=2\kappa$. It is therefore given by a stream function
$\Psi(R,\Theta)$:
\begin{equation}\label{eq:stream}
W=-p_0\Psi_R,\qquad U=\frac{p_0}{R}\Psi_\Theta,\qquad \Delta\Psi=-\frac{2\kappa}{p_0}\ \text{ in } D_S:=\{1<R<S(\Theta)\},
\end{equation}
normalised by $\Psi=1$ on $R=1$ and $\Psi=0$ on $R=S(\Theta)$. Here $p_0\neq0$ is the
(constant) relative mass flux. We assume
\begin{equation}\label{eq:nostag}
W/p_0>0\quad\text{in }\overline{D_S},\qquad\text{i.e.}\quad \Psi_R<0 ,
\end{equation}
so that there are no stagnation points in the relative flow. For $\kappa=0$, Li and Wang
\cite{LWpreprint} take $p_0<0$, i.e.\ $W<0$.

At $\kappa=0$ the flow relative to the satellite is irrotational. For an ocean co-rotating with the Earth, the Earth--Moon frequencies instead give $\kappa\simeq26.4$ (Section~\ref{sec:earth}).
In the frame rotating with the satellite, the steady Euler equations read
$(\mathbf w\cdot\nabla)\mathbf w+2\mathbf e_z\times\mathbf w-\nabla\frac{R^2}2=-\nabla(P+\alpha R+\Phi_T)$, where
now $P$ is the scaled pressure and $\Phi_T(R,\Theta)=-\eps R^2P_2(\cos\Theta)$. Since
$(\mathbf w\cdot\nabla)\mathbf w=\nabla\frac12|\mathbf w|^2+2\kappa\,\mathbf e_z\times\mathbf w$ and
$\mathbf e_z\times\mathbf w=p_0\nabla\Psi$, we obtain
\[
\tfrac12|\mathbf w|^2+2(1+\kappa)p_0\Psi-\tfrac12R^2+P+\alpha R+\Phi_T=\mathrm{const}\qquad\text{in }\overline{D_S}.
\]
On the free surface $\Psi=0$ and $P=p_{\rm atm}$. Hence, for every $\kappa$, Bernoulli's law on
the free surface reads
\begin{equation}\label{eq:bern-phys}
p_0^2|\nabla\Psi|^2-R^2+2\alpha R-2\eps R^2P_2(\cos\Theta)=2\lambda\qquad\text{on }R=S(\Theta),
\end{equation}
for some constant $\lambda$. For $\kappa=0$ this is \cite[Lemma~2.1]{LWpreprint} with the tidal term
added. The problem is to find $(S,\Psi)$ satisfying \eqref{eq:stream}--\eqref{eq:bern-phys}.

\subsection{Log-polar variables and the Dubreil--Jacotin formulation}
Since $\Delta=R^{-2}(\p_s^2+\p_\Theta^2)$ with $s=\ln R$, the fluid domain becomes the
flat-bottomed periodic strip $\{0<s<\ln S(\Theta)\}$, and the interior equation in
\eqref{eq:stream} becomes
\begin{equation}\label{eq:flatpoisson}
\widehat\Delta\Psi:=\Psi_{ss}+\Psi_{\Theta\Theta}=-\frac{2\kappa}{p_0}\,\e^{2s}.
\end{equation}
The right-hand side depends on $s$ only. Hence $\p_\Theta\Psi$ is harmonic in $(\Theta,s)$ for
every $\kappa$, and $\Psi$ itself is harmonic when $\kappa=0$. We use the semi-hodograph variables
\[
q=\Theta,\qquad p=-\Psi(R,\Theta),\qquad H(q,p)=s=\ln R,
\]
which map $D_S$ onto the periodic strip $\Omega=\T\times(-1,0)$, $\T=\R/2\pi\mathbb Z$.

\begin{lemma}\label{lem:logpolar}
Let $\Psi$ be of class $C^2$ on the closure of a periodic strip $\{0<s<\eta(\Theta)\}$, with
$\Psi_s<0$, $\Psi=1$ on $s=0$ and $\Psi=0$ on $s=\eta(\Theta)$, and define $H$ on
$\overline\Omega$ by $\Psi(q,H(q,p))=-p$. Then $H_p=-1/\Psi_s>0$, $H_q=\Psi_\Theta H_p$,
\[
(\Psi_s^2+\Psi_\Theta^2)(q,H)=\frac{1+H_q^2}{H_p^2},\qquad
\Dop(H):=(1+H_q^2)H_{pp}-2H_qH_pH_{qp}+H_p^2H_{qq}=H_p^3\,(\widehat\Delta\Psi)(q,H).
\]
\end{lemma}

\begin{proof}
At fixed physical position the coordinate derivatives are
$\partial_s=H_p^{-1}\partial_p$ and $\partial_\Theta=\partial_q-(H_q/H_p)\partial_p$.
Applying them twice to $\Psi=-p$, we obtain the identities.
\end{proof}

By Lemma~\ref{lem:logpolar}, $|\nabla\Psi|^2=\e^{-2s}(\Psi_s^2+\Psi_\Theta^2)$ and
\eqref{eq:flatpoisson}, problem \eqref{eq:stream}--\eqref{eq:bern-phys} becomes
\begin{equation}\label{eq:DJ}
\Dop(H)+\frac{2\kappa}{p_0}\e^{2H}H_p^3=0\ \text{ in }\Omega,\qquad \Bop_\eps(\alpha;H)=0\ \text{ on }p=0,\qquad H=0\ \text{ on }p=-1,
\end{equation}
where
\[
\Bop_\eps(\alpha;H):=\frac{p_0^2(1+H_q^2)}{H_p^2}\,\e^{-2H}-\e^{2H}+2\alpha\e^{H}-2\eps\,\e^{2H}P_2(\cos q)-2\lambda .
\]
For $\kappa=0$ the interior operator coincides with the Dubreil--Jacotin operator for irrotational
flow over a flat bed \cite{CS2004}, and the curvature of the annulus survives only through the
exponential weights in $\Bop_\eps$. For $\kappa\neq0$ the vorticity adds the semilinear term
$\frac{2\kappa}{p_0}\e^{2H}H_p^3$, again with an exponential weight.

\begin{lemma}\label{lem:equiv}
Let $\beta\in(0,1)$. Solutions $(S,\Psi)$ of \eqref{eq:stream}--\eqref{eq:bern-phys} with
$S\in C^{3,\beta}(\T)$, $\Psi\in C^{3,\beta}(\overline{D_S})$ and \eqref{eq:nostag} correspond one-to-one
to solutions $H\in C^{3,\beta}_{\rm per}(\overline\Omega)$ of \eqref{eq:DJ} with $H_p>0$ in
$\overline\Omega$. The surface is recovered from $S(q)=\e^{H(q,0)}$.
\end{lemma}

\begin{proof}
The inverse function theorem on each radial fibre gives $p=-\Psi(q,s)$ and $S(q)=\exp H(q,0)$, with the asserted regularity. The two identities of Lemma~\ref{lem:logpolar} transport the interior and boundary equations in both directions. For comparison with \cite[(2.23)]{LWpreprint}, its height variable is $h=e^H-1$.
\end{proof}

\subsection{Laminar flows, parameters and function spaces}
For $\eps=0$, the $q$-independent solutions of \eqref{eq:DJ} with $H(0)=\gamma>0$ are the
laminar flows $S\equiv S_0=\e^{\gamma}$ with relative azimuthal velocity
\begin{equation}\label{eq:laminar}
W_0(R)=\kappa R-\frac{\mu}{R},\qquad p_0=\int_1^{S_0}W_0\dd R=\frac\kappa2\big(S_0^2-1\big)-\mu\gamma .
\end{equation}
Indeed, \eqref{eq:stream} gives
$\Psi_0(s)=p_0^{-1}\big[\tfrac\kappa2(\e^{2\gamma}-\e^{2s})-\mu(\gamma-s)\big]$, and the laminar
solution $H_0$ of \eqref{eq:DJ} is defined by $\Psi_0(H_0(p))=-p$. For $\kappa=0$,
$H_0(p)=\gamma(1+p)$ and $\mu=-p_0/\gamma$. We write
\[
W_s:=W_0(S_0),\qquad w:=\frac{W_s}{S_0},\qquad \cc:=W_s^2,\qquad r:=\cc\,\e^{-2\gamma}=w^2,\qquad T_k:=\coth(k\gamma).
\]
Thus $W_s$ is the relative velocity at the surface and $w$ the corresponding angular velocity,
in units of $\Omega_f$. Condition \eqref{eq:nostag} holds for the laminar flow if and only if
$W_0$ has no zero on $[1,S_0]$, that is,
\begin{equation}\label{eq:admissible}
w\big(w-\kappa(1-\e^{-2\gamma})\big)>0,\qquad\text{equivalently}\qquad w\ne0\ \text{ and }\ \frac\kappa w<\frac{1+T_1}2 .
\end{equation}
We call $(\gamma,\kappa,w)$ \emph{admissible} if \eqref{eq:admissible} holds. For $\kappa=0$
this only says $w\neq0$. If $(\gamma,\kappa,w)$ is admissible, then $H_{0,p}=-1/\Psi_0'(H_0)>0$ on
$[-1,0]$,
\begin{equation}\label{eq:h1}
H_{0,pp}=-\frac{2\kappa}{p_0}\e^{2H_0}H_{0,p}^3,\qquad
h_1:=H_{0,p}(0)=\frac{p_0}{S_0W_s},\qquad p_0^2=\cc\,h_1^2\e^{2\gamma}.
\end{equation}
For $\kappa=0$, $h_1=\gamma$ and $\cc=p_0^2/(\gamma^2\e^{2\gamma})=\mu^2/S_0^2$. The laminar flow solves
\eqref{eq:DJ} with $\eps=0$ if and only if
\begin{equation}\label{eq:lambda}
2\lambda=\cc+2\alpha\e^{\gamma}-\e^{2\gamma}.
\end{equation}
For $\kappa=0$ this is \cite[(4.2)--(4.3)]{LWpreprint}, with $H_0=\ln(1+H^{\rm LW})$.

\emph{Throughout, $p_0$, $\gamma>0$ and $\kappa$ are fixed and admissible, and
$\lambda=\lambda(\alpha)$ is given by \eqref{eq:lambda}, as in \cite{LWpreprint}.} Then $H_0$ solves
\eqref{eq:DJ} with $\eps=0$ for every $\alpha\in\R$. For fixed $\gamma$, the map
$(p_0,\kappa)\mapsto(\kappa,w)$ is a bijection, with
$w=\kappa\big(1-\frac{\sinh\gamma}{\gamma\e^{\gamma}}\big)+\frac{p_0}{\gamma\e^{2\gamma}}$, and we
state our results in terms of $(\gamma,\kappa,w)$. Problem \eqref{eq:DJ} depends on $p_0$ and
$\kappa$ only through $p_0^2$ and $\kappa/p_0$. It is therefore invariant under
$(p_0,\kappa,w)\mapsto(-p_0,-\kappa,-w)$, which reverses the relative flow, and all formulas
below depend on $(\kappa,w)$ only through $\kappa^2$, $\kappa w$ and $w^2$.

Substituting \eqref{eq:lambda} and \eqref{eq:h1} into $\Bop_\eps$, we obtain the form used below:
\begin{equation}\label{eq:Bform}
\Bop_\eps(\alpha;H)=\cc\Big[\frac{h_1^2(1+H_q^2)}{H_p^2}\e^{-2(H-\gamma)}-1\Big]
-\e^{2\gamma}\big(\e^{2(H-\gamma)}-1\big)+2\alpha\e^{\gamma}\big(\e^{H-\gamma}-1\big)-2\eps\,\e^{2H}P_2(\cos q).
\end{equation}

We fix $\beta\in(0,1)$ and let
\[
\Xs=\{H\in C^{3,\beta}_{\rm per}(\overline\Omega):\ H\text{ even in }q,\ H=0\text{ on }p=-1\},
\qquad
\Ys=C^{1,\beta}_{\rm per,e}(\overline\Omega)\times C^{2,\beta}_{\rm per,e}(\T),
\]
where ``per'' means $2\pi$-periodic in $q$ and ``e'' means even in $q$. Let
$\Uo=\{H\in\Xs:\ \min_{\overline\Omega}H_p>0\}$ and define
\begin{equation}\label{eq:Fdef}
\Fop:\R\times\R\times\Uo\to\Ys,\qquad \Fop(\alpha,\eps;H)=\Big(\Dop(H)+\frac{2\kappa}{p_0}\e^{2H}H_p^3,\ \Bop_\eps(\alpha;H)|_{p=0}\Big).
\end{equation}
The map $\Fop$ is real-analytic, and $H_0\in\Uo$. By Lemma~\ref{lem:equiv}, zeros of $\Fop$
correspond one-to-one to the solutions of the physical problem that are even in $\Theta$.

For $m\in\N$ let $\Xs_m\subset\Xs$ and $\Ys_m\subset\Ys$ denote the closed subspaces of
$2\pi/m$-periodic functions. We write $\Xs_\pi:=\Xs_2$ and $\Ys_\pi:=\Ys_2$.

\begin{remark}\label{rem:parameter-slice}
The data held fixed in the bifurcation problem are $(p_0,\gamma,\kappa)$,
while $\lambda=\lambda(\alpha)$ follows \eqref{eq:lambda}. No condition fixes
the fluid area $\tfrac12\int_{\mathbb T}(S(\Theta)^2-1)\,d\Theta$ along
the nonlinear branch. In particular the mean tidal response is allowed
to vary. A fixed-area formulation would require an additional constraint
and a corresponding rederivation of the reduced nonlinear coefficients.
\end{remark}

\subsection{Symmetries}
For a forcing profile $\Ptide\in C^\infty(\T)$, let $\Fop^{\Ptide}$ denote $\Fop$ with
$P_2(\cos q)$ replaced by $\Ptide(q)$. The two profiles relevant here are
\begin{equation}\label{eq:Ppm}
\Ptide^{\pm}(q)=\tfrac14\pm\tfrac34\cos2q,\qquad \Ptide^+(q)=P_2(\cos q).
\end{equation}
Let $(\tau_\vartheta H)(q,p)=H(q+\vartheta,p)$ and $(\iota H)(q,p)=H(-q,p)$. These act on the
spaces of $2\pi$-periodic functions without parity; in Lemma~\ref{lem:sym}(a) and (d),
$\Fop^{\Ptide}$ denotes the map given by the same formula on these larger spaces.

\begin{lemma}\label{lem:sym}
The following hold.
\begin{enumerate}[label=\textup{(\alph*)}]
\item For every profile $\Ptide$,
\[
\Fop^{\Ptide}(\alpha,\eps;\tau_\vartheta H)=\tau_\vartheta\Fop^{\tau_{-\vartheta}\Ptide}(\alpha,\eps;H),\qquad
\Fop^{\Ptide}(\alpha,\eps;\iota H)=\iota\Fop^{\iota\Ptide}(\alpha,\eps;H).
\]
\item For $\eps=0$, $\Fop$ commutes with all $\tau_\vartheta$ and with $\iota$ ($O(2)$-equivariance).
\item $\Fop^{\Ptide^+}$ commutes with $\rho:=\tau_\pi$ and with $\iota$; these generate the
Klein four-group $D_2$. The operator $\rho$ preserves $\Xs$ and $\Ys$.
\item $\Fop^{\Ptide^+}(\alpha,\eps;H)=0$ if and only if $\Fop^{\Ptide^-}(\alpha,\eps;\tau_{\pi/2}H)=0$.
Moreover, $H\in\Fix(\rho\iota)$, i.e.\ $H(\pi-q,p)=H(q,p)$, if and only if $\tau_{\pi/2}H$ is
even in $q$.
\end{enumerate}
\end{lemma}

\begin{proof}
The unforced terms are autonomous in $q$ and invariant under reflection. A change of $q$ therefore acts only on the forcing profile. The identities follow from $\tau_\pi P^+=\iota P^+=P^+$ and $\tau_{\pi/2}P^+=P^-$. Conjugation by the quarter-turn interchanges the reflection axes.
\end{proof}
We write $\Fop^\pm=\Fop^{P^\pm}$ and $\Fop=\Fop^+$.

\subsection{Linearisation and dispersion relation}\label{sec:linear}
Let $\Lop(\alpha):=D_H\Fop(\alpha,0;H_0)$. Since $H_{0,q}=0$, \eqref{eq:Fdef} and \eqref{eq:Bform} give
\begin{equation}\label{eq:L}
\Lop(\alpha)\psi=\Big(\psi_{pp}+H_{0,p}^2\psi_{qq}+\frac{2\kappa}{p_0}\e^{2H_0}H_{0,p}^2\big(3\psi_p+2H_{0,p}\psi\big),\ \
2b_0(\alpha)\psi-\frac{2\cc}{h_1}\psi_p\Big|_{p=0}\Big),
\end{equation}
where $b_0(\alpha):=\alpha\e^{\gamma}-\e^{2\gamma}-\cc$. For $\kappa=0$ we have $H_{0,p}=h_1=\gamma$, and the
interior operator is $\psi_{pp}+\gamma^2\psi_{qq}$. In particular
$\p_\alpha\Lop(\alpha)\psi=(0,\,2\e^{\gamma}\psi|_{p=0})$.

The following observation reduces the linear theory for every $\kappa$ to the irrotational one.
Since the vorticity is constant, a perturbation of the stream function at fixed
$(\Theta,s)$ is harmonic. In the variables $(q,p)$ this reads as follows.

\begin{lemma}\label{lem:conj}
For $\psi\in C^2(\overline\Omega)$ write $\psi=H_{0,p}\,u$. Then
\[
\Lop(\alpha)\psi=\Big(H_{0,p}^2\big[\p_p\big(H_{0,p}^{-1}u_p\big)+H_{0,p}u_{qq}\big],\ \
2h_1b(\alpha)u-2\cc\,u_p\big|_{p=0}\Big),\qquad b(\alpha):=b_0(\alpha)+2\kappa\e^{\gamma}W_s .
\]
In the variable $s=H_0(p)\in[0,\gamma]$ the interior component equals $H_{0,p}^3(u_{ss}+u_{qq})$.
Moreover $u(q,H_0^{-1}(s))$ is, to first order, the perturbation of the stream function $\Psi$ at
fixed $(\Theta,s)=(q,s)$.
\end{lemma}

\begin{proof}
For $H_t=H_0+t\psi$, the identity $\Psi_t(q,H_t)=-p$ gives the Eulerian variation $u=\psi/H_{0,p}$. Differentiating the identity between the height residual and $H_p^3(\widehat\Delta\Psi+2\kappa e^{2s}/p_0)$ at the laminar solution, we obtain $H_{0,p}^3\widehat\Delta u$. At the surface insert $\psi_p=H_{0,pp}u+h_1u_p$ into \eqref{eq:L} and use \eqref{eq:h1}.
\end{proof}
\begin{lemma}\label{lem:fredholm}
For every $\alpha\in\R$, $\Lop(\alpha):\Xs\to\Ys$ is a Fredholm operator of index zero. The
same holds for its restriction $\Xs_m\to\Ys_m$, $m\in\N$.
\end{lemma}

\begin{proof}
Take $\Mop\psi=(\psi_{pp}+H_{0,p}^2\psi_{qq},-2\cc\psi_p/h_1|_{p=0})$ with zero bottom trace.
The energy $\int(\psi_p^2+H_{0,p}^2\psi_q^2)$ is coercive by the bottom condition. The weak mixed problem and boundary Schauder estimates on the two disjoint boundary circles give an isomorphism $\Xs\to\Ys$ \cite{ADN1959,GT2001}. Uniqueness preserves parity and subperiod. The lower-order interior and trace terms in $\Lop-\Mop$ gain one derivative relative to the target spaces and hence are compact.
\end{proof}

\begin{proposition}\label{prop:dispersion}
For $k\ge1$ let $\sigma_k(p):=H_{0,p}(p)\sinh\big(kH_0(p)\big)$, let
$\sigma_0(p):=H_{0,p}(p)H_0(p)$, and let
\[
\alpha_k:=S_0\big[1+r(1+kT_k)-2\kappa w\big]\quad(k\ge1),\qquad
\alpha_0:=S_0\big[1+r(1+\gamma^{-1})-2\kappa w\big].
\]
Then $\ker\Lop(\alpha)\ne\{0\}$ if and only if $\alpha=\alpha_k$ for some $k\ge0$. In that case
\[
\ker\Lop(\alpha_k)=\spn\{\varphi_k\},\qquad \varphi_k(q,p):=\sigma_k(p)\cos kq .
\]
Moreover $\alpha_0<\alpha_1<\alpha_2<\cdots$, $\alpha_k\to\infty$, and the differences
$\alpha_k-\alpha_j=S_0r(kT_k-jT_j)$, with $0\cdot T_0:=\gamma^{-1}$, do not depend on $\kappa$ at fixed $(\gamma,w)$. In terms of the relative surface
velocity,
\[
\alpha_k=S_0-2\kappa W_s+\frac{W_s^2}{S_0}\Big(1+k\coth(k\ln S_0)\Big).
\]
For $\kappa=0$, $\sigma_k=\gamma s_k$ with $s_k(p):=\sinh(k\gamma(1+p))$ for $k\ge1$,
$\sigma_0=\gamma^2(1+p)$, and $W_s^2=\mu^2/S_0^2$ with $\mu=-p_0/\gamma$.
\end{proposition}

\begin{proof}
Let $\psi\in\ker\Lop(\alpha)$, write $\psi=H_{0,p}u$ and let $u_m(p)$ be the $m$-th Fourier
cosine coefficient of $u$. Differentiating under the integral sign and integrating by parts in
$q$, Lemma~\ref{lem:conj} gives, in the variable $s=H_0(p)$,
\[
\p_s^2u_m=m^2u_m\ \text{ on }(0,\gamma),\qquad u_m=0\ \text{ at }s=0,\qquad b(\alpha)u_m=\cc\,\p_su_m\ \text{ at }s=\gamma,
\]
where we used $u_p=H_{0,p}\p_su$ and $H_{0,p}(0)=h_1$. Hence $u_m=A_m\sinh(ms)$ for $m\ge1$ and
$u_0=A_0s$. A non-zero $A_m$ is possible only if $b(\alpha)=\cc\,mT_m$ when $m\ge1$, or
$b(\alpha)=\cc/\gamma$ when $m=0$. Since $\cc\e^{-\gamma}=rS_0$ and $W_s=wS_0$, these conditions are
exactly $\alpha=\alpha_m$.

The function $x\mapsto x\coth x$ tends to $1$ as $x\to0$ and is strictly increasing on
$(0,\infty)$, since its derivative is $(\sinh x\cosh x-x)/\sinh^2x>0$. Hence
$1/\gamma<T_1<2T_2<\cdots$. This gives the ordering of the $\alpha_k$, which are therefore
pairwise distinct. Consequently at most one Fourier mode can lie in the kernel.
\end{proof}
\begin{figure}[htbp]\centering
\includegraphics[width=.91\linewidth]{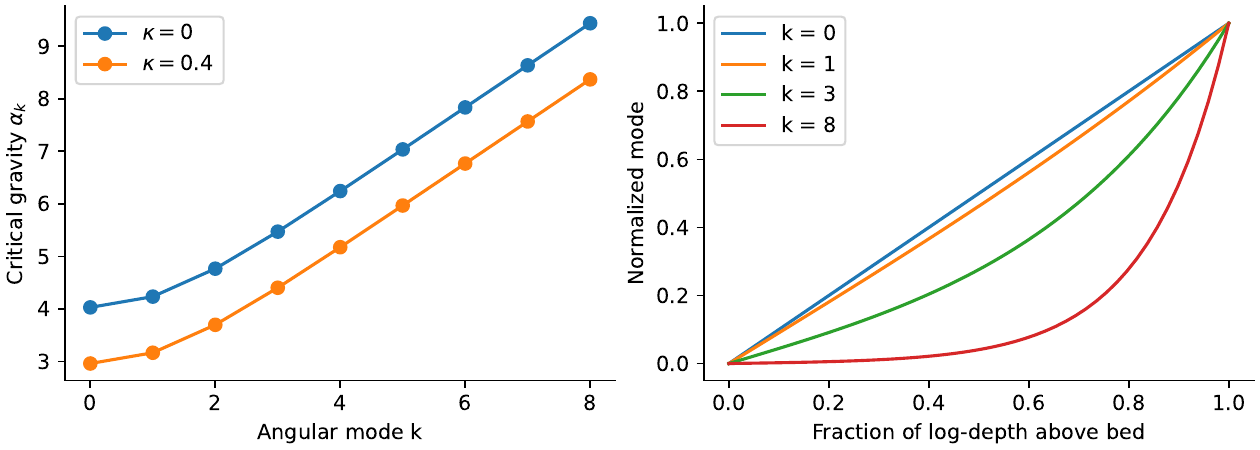}
\caption{Dispersion relation of Proposition~\ref{prop:dispersion} at $\gamma=0.8$, $w=0.6$. Left: changing $\kappa$ translates all critical values by the same amount. Only the integer modes are spectral points. Right: at $\kappa=0$, the vertical kernel profiles normalized by their surface values are $\sinh(k\gamma z)/\sinh(k\gamma)$, with $z=1+p$ and mean mode $z$. Higher modes concentrate near the surface.}\label{fig:atlas-spectrum}
\end{figure}

Admissibility gives $2\kappa/w<1+T_1\le1+kT_k$, hence $\alpha_k>S_0$ for $k\ge1$. For $\kappa=0$, $\alpha_1$ is the critical value of \cite[Theorem~2.1]{LWpreprint}. For fixed integers $j,k\ge0$, expansion of the strip multiplier gives
\[
 \alpha_k-\alpha_j=S_0r(k^2-j^2)\gamma/3+O(r\gamma^3).
\]
Thus absolute gaps close at bounded $r$. At fixed nonzero flux, $r\sim p_0^2/\gamma^2$ and only the relative gaps close.

\begin{lemma}\label{lem:range}
For $k\ge1$ let $v_k(q,p):=\sinh\big(kH_0(p)\big)\cos kq$, so that $\varphi_k=H_{0,p}v_k$, and define
the bounded linear functional on $\Ys$
\begin{equation}\label{eq:ell}
\ell_k(f,g):=\iint_\Omega f\,v_k\,\frac{h_1^2}{H_{0,p}^2}\dd q\dd p+\frac{h_1}{2\cc}\int_\T g(q)\,v_k(q,0)\dd q .
\end{equation}
Then:
\begin{enumerate}[label=\textup{(\alph*)}]
\item $\mathrm{R}(\Lop(\alpha_k))=\ker\ell_k$;
\item $\ell_k\big(\p_\alpha\Lop(\alpha_k)\hat\varphi_k\big)=\pi h_1\cc^{-1}\e^{\gamma}\sinh(k\gamma)>0$, where
$\hat\varphi_k:=\varphi_k/(h_1\sinh(k\gamma))$, so that $\hat\varphi_k(q,0)=\cos kq$;
\item if $\alpha\notin\{\alpha_j\}_{j\ge0}$, then $\Lop(\alpha):\Xs\to\Ys$ is an isomorphism;
\item if $\alpha\notin\{\alpha_{jm}\}_{j\ge0}$, then $\Lop(\alpha):\Xs_m\to\Ys_m$ is an isomorphism.
\end{enumerate}
For $\kappa=0$, $v_k=s_k\cos kq$ and $\ell_k(f,g)=\iint fv_k+\frac{\gamma}{2\cc}\int gv_k(\cdot,0)$.
\end{lemma}

\begin{proof}
(a) Let $(f,g)=\Lop(\alpha_k)\psi$ and $u=\psi/H_{0,p}$. By Lemma~\ref{lem:conj},
$f/H_{0,p}^2=\p_p(H_{0,p}^{-1}u_p)+H_{0,p}u_{qq}$, and $v_k$ satisfies the homogeneous equation. Green's
identity for this symmetric operator, with $u=v_k=0$ on $p=-1$ and periodicity in $q$, gives $\iint_\Omega\frac{f\,v_k}{H_{0,p}^2}=\frac1{h_1}\int_\T(u_pv_k-uv_{k,p})\big|_{p=0}\dd q .$ At $p=0$ we have $\cc u_p=h_1b(\alpha_k)u-\frac12g$ and $\cc v_{k,p}=h_1b(\alpha_k)v_k$, the latter
because $v_k$ solves the homogeneous problem. Hence we obtain $\ell_k(f,g)=0$, i.e.\
$\mathrm R(\Lop(\alpha_k))\subset\ker\ell_k$. Moreover $\ell_k\neq0$, since
$\ell_k(0,\cos kq)=\pi h_1\sinh(k\gamma)/(2\cc)$. By Lemma~\ref{lem:fredholm} and
Proposition~\ref{prop:dispersion} the range has codimension one, so equality holds.

(b) This follows from $\p_\alpha\Lop\hat\varphi_k=(0,2\e^{\gamma}\cos kq)$.

(c) and (d) follow from Lemma~\ref{lem:fredholm} and Proposition~\ref{prop:dispersion}: an
injective Fredholm operator of index zero is an isomorphism.
\end{proof}

\section{Local bifurcation and the sign of the Stokes coefficient}\label{sec:unforced}
\subsection{The local curve and its coefficient}

For a function $H$ on $\overline\Omega$ we define its \emph{$k$-th surface amplitude}
\[
\mathsf a_k[H]:=\frac1\pi\int_\T H(q,0)\cos kq\dd q .
\]
Since $H(q,0)=\ln S(q)$ and $q=\Theta$, $\mathsf a_k[H]$ is the $\cos k\Theta$-coefficient of the
logarithm of the physical free surface. Note that $\mathsf a_k[\hat\varphi_k]=1$.

\begin{theorem}\label{thm:local}
Let $\eps=0$ and $k\ge1$. There exist $\bar s>0$ and a real-analytic curve
\[
(-\bar s,\bar s)\ni s\mapsto(\alpha_k(s),H_k(s))\in\R\times(\Xs_k\cap\Uo)
\]
with the following properties.
\begin{enumerate}[label=\textup{(\roman*)}]
\item $\Fop(\alpha_k(s),0;H_k(s))=0$, $(\alpha_k(0),H_k(0))=(\alpha_k,H_0)$,
$H_k(s)=H_0+s\hat\varphi_k+o(s)$, and $\mathsf a_k[H_k(s)-H_0]=s$.
\item Every solution of $\Fop(\alpha,0;H)=0$ in a suitable neighbourhood $\mathcal W$ of
$(\alpha_k,H_0)$ in $\R\times\Xs$ either has $H=H_0$ or lies on the curve.
\item $\alpha_k(-s)=\alpha_k(s)$ and $H_k(-s)=\tau_{\pi/k}H_k(s)$. The two arms are related by translation. The pitchfork is nondegenerate when $\alpha_k^{(2)}\ne0$, and
\[
\alpha_k(s)=\alpha_k+\alpha_k^{(2)}s^2+O(s^4).
\]
\end{enumerate}
\end{theorem}

\begin{proof}
Proposition~\ref{prop:dispersion} and Lemma~\ref{lem:range} verify the simple-kernel and transversality hypotheses of analytic local bifurcation \cite{CR1971,BuffoniToland2003}. We choose the closed complement $\Zs_k=\ker\mathsf a_k$, since $\mathsf a_k[\hat\varphi_k]=1$. This gives the amplitude parametrisation and uniqueness in $\Xs$. Application in $\Xs_k$ gives the same curve by uniqueness. There $\tau_{\pi/k}$ preserves evenness, commutes with the equation and changes the sign of $\mathsf a_k$. Uniqueness therefore gives (iii).
\end{proof}

\begin{theorem}\label{thm:direction}
For every $k\ge1$ and every admissible $(\gamma,\kappa,w)$, with $T=T_k=\coth(k\gamma)$ and $r=w^2$,
\begin{equation}\label{eq:alpha2}
\alpha_k^{(2)}=\frac{\e^{\gamma}\,P_k(\kappa,w)}{8kr\,(k\gamma T-1)},
\end{equation}
where
\begin{equation}\label{eq:Pk}
P_k(\kappa,w)=(1-4\kappa^2)\big[a_0(1-4\kappa^2)-12X_k\,\kappa w\big]+a_1r+b_1\kappa wr+b_2\kappa^2r+a_2r^2
\end{equation}
and
\begin{align*}
a_0&=T-k\gamma\,(T^2-2),\\
a_1&=-\big(6k^2T^3+10kT^2+6k^2T+6T+9k\big)\\
&\quad+k\gamma\big(6k^2T^4+10kT^3-6k^2T^2+6T^2-11kT+4k^2-12\big),\\
a_2&=9k^4T^5+30k^3T^4+6k^4T^3+43k^2T^3+30k^3T^2+30kT^2+k^4T+37k^2T+9T+15k\\
&\quad-k\gamma\big(9k^4T^6+30k^3T^5-12k^4T^4+43k^2T^4-30k^3T^3+30kT^3\\
&\qquad\qquad+13k^4T^2-49k^2T^2+9T^2+20k^3T-45kT-2k^4+12k^2-18\big),\\
X_k&=k\gamma\,(kT^3+T^2-kT-2)-(kT^2+T+k),\\
b_1&=2k\gamma\big(18k^3T^5+48k^2T^4-12k^3T^3+48kT^3-40k^2T^2+18T^2\\
&\qquad\qquad+10k^3T-63kT+12k^2-36\big)\\
&\quad-2\big(18k^3T^4+48k^2T^3+24k^3T^2+48kT^2+56k^2T+18T-2k^3+33k\big),\\
b_2&=4\big(15k^2T^3+28kT^2+21k^2T+15T+24k\big)\\
&\quad-4k\gamma\big(15k^2T^4+28kT^3-9k^2T^2+15T^2-32kT+4k^2-30\big).
\end{align*}
Since $k\gamma T>1$, the bifurcation at $\alpha_k$ is supercritical if $P_k(\kappa,w)>0$ and
subcritical if $P_k(\kappa,w)<0$. For $\kappa=0$, $P_k(0,w)=a_0+a_1r+a_2r^2$.
\end{theorem}

\begin{proof}
We write $v=\hat\varphi_k$, $F_j=D_H^j\Fop(\alpha_k,0;H_0)$ and let $z\in\Zs_k$ solve $\Lop(\alpha_k)z=-F_2[v,v]/2$. This is solvable by Lemma~\ref{lem:range} since the forcing has only modes $0,2k$. At cubic order application of $\ell_k$ eliminates the unknown third-order height and yields $\alpha_k^{(2)}=-\frac{\ell_k(F_2[v,z]+F_3[v,v,v]/6)} {\ell_k(\partial_\alpha\Lop(\alpha_k)v)}.$ The trace calculation below evaluates this quotient and gives \eqref{eq:alpha2}. A smooth change of amplitude with nonzero derivative preserves the sign of the coefficient.
\end{proof}

\subsection{Calculation by harmonic extension}\label{app:computations}\label{app:flat}
We write $\psi:=p_0\Psi=\psi_0(s)+\phi(\Theta,s)$, where
$\psi_0(s)=\frac\kappa2(\e^{2\gamma}-\e^{2s})-\mu(\gamma-s)$ is the laminar stream function. Then
$\phi$ is harmonic in the fluid, $\phi=0$ on $s=0$, and on the free surface $s=\eta(\Theta)$
\begin{subequations}\label{eq:flatbc}
\begin{gather}
\psi_0(\eta)+\phi(\Theta,\eta)=0,\label{eq:flatK}\\
\e^{-2\eta}\big[(\psi_0'(\eta)+\phi_s)^2+\phi_\Theta^2\big]-\e^{2\eta}+2\alpha\e^{\eta}-2\eps\e^{2\eta}\Ptide(\Theta)-\big(\cc+2\alpha\e^{\gamma}-\e^{2\gamma}\big)=0,\label{eq:flatB}
\end{gather}
\end{subequations}
by \eqref{eq:bern-phys}, \eqref{eq:lambda} and $|\nabla\Psi|^2=\e^{-2s}(\Psi_s^2+\Psi_\Theta^2)$. At $s=\gamma$,
$\psi_0'=-\e^{\gamma}W_s$ and $\psi_0^{(j)}=-2^{j-1}\kappa\e^{2\gamma}$ for $j\ge2$.

We write $\eta=\gamma+\sum_{j\ge1}t^j\eta_j$ and $\phi=\sum_{j\ge1}t^j\phi_j$. For an analytic solution curve the coefficients of $\phi$ can be determined on the fixed strip $0\le s\le\gamma$. Indeed, differentiating $\Psi(q,H(q,p))=-p$ successively determines the Eulerian coefficient at each order from the corresponding height coefficient and the earlier derivatives. The right-hand side of \eqref{eq:flatpoisson} is independent of the curve parameter, so each coefficient is harmonic, with zero bottom trace. Induction in the height equation gives finitely many Fourier modes at each order. Each therefore extends to $s=\gamma$ by the strip harmonic extension, and Taylor expansion of \eqref{eq:flatbc} there gives the same coefficients as differentiation in fixed height variables.
For the curve of Theorem~\ref{thm:local} we use $t=\mathsf a_k[H-H_0]$, so $\eta_1=\cos k\Theta$ and
$\eta_j$ has no $\cos k\Theta$-component for $j\ge2$, and $\alpha=\alpha_k+t^2\alpha_k^{(2)}+O(t^4)$.
For a mode $\cos(mk\Theta)$ with trace $A$ at $s=\gamma$, the harmonic extension vanishing at $s=0$
has $\p_s^{2i}=(mk)^{2i}A$ and $\p_s^{2i+1}=(mk)^{2i+1}T_{(m)}A$ at $s=\gamma$, where
$T_{(m)}=\coth(mk\gamma)$, $\coth(2x)=\frac{1+\coth^2x}{2\coth x}$ and
$\coth(3x)=\frac{\coth^3x+3\coth x}{3\coth^2x+1}$; for $m=0$ the extension is linear in $s$.

\emph{First order.} \eqref{eq:flatK} gives $\phi_1=S_0W_s\frac{\sinh(ks)}{\sinh(k\gamma)}\cos k\Theta$, and
the $\cos k\Theta$-component of \eqref{eq:flatB} vanishes exactly for $\alpha=\alpha_k$
(Proposition~\ref{prop:dispersion}).

\begin{lemma}\label{lem:flat2}
Let $T=T_k$ and $r=w^2$. At second order, $\eta_2=E_0+E_{2k}\cos2k\Theta$ with
\begin{gather*}
E_0=\frac{\gamma(1-4\kappa^2)+2\kappa w\,(2k\gamma T+3\gamma+1)-r\big(k^2\gamma T^2+5k\gamma T+2kT-(k^2-3)\gamma\big)}{4r\,(k\gamma T-1)},\\
E_{2k}=\frac{r\,(3k^2T^3+5kT^2-k^2T+3T)-T(1-4\kappa^2)-2\kappa w\,(3kT^2+3T+k)}{4kr}.
\end{gather*}
\end{lemma}

\begin{proof}
For $m=0,2k$, let $R_m$ be the $m$-th cosine coefficient (the mean when $m=0$) of the quadratic Bernoulli residual with $\eta=\gamma+t\cos k\Theta$ and the harmonic Dirichlet correction supplied by \eqref{eq:flatK}. The terms involving $\eta_2$ are the linearised residual: $R_m+2S_0(\alpha_k-\alpha_m)E_m=0$. To evaluate $R_m$, use $b_1=-\psi_0'(\gamma)\cos k\Theta$ and $b_2=-\psi_0''(\gamma)\cos^2(k\Theta)/2-\cos(k\Theta)\partial_s\mathcal E b_1$ as the two harmonic traces, where $\mathcal E$ is the strip extension and the derivative is at $s=\gamma$. Inserting these in \eqref{eq:flatB}, we obtain the displayed coefficients. The division in each mode is justified by Proposition~\ref{prop:dispersion}.
\end{proof}

For convenience, the cubic projection can be evaluated without solving a third-order interior equation. Let $c=\cos k\Theta$, $d=E_0+E_{2k}\cos2k\Theta$ and $u=tc+t^2d$. For $j=1,2,3$ let $b_j$ be the surface trace of the harmonic term and let $\mathcal E b$ denote the harmonic extension of $b$ with zero bottom trace. All derivatives in the following recursion are evaluated at $s=\gamma$:
\[
 b_j=-[t^j]\left\{\psi_0(\gamma+u)+
 \sum_{i=1}^{j-1}t^i\sum_{l=0}^{j-i}\frac{u^l}{l!}\,\partial_s^l\mathcal E b_i\right\}.
\]
Here $[t^j]$ extracts the coefficient of $t^j$. The omitted third-order surface displacement has zero $k$-th mode by the amplitude normalisation, and its linear response is diagonal. We set $\phi^{[3]}=\sum_{i=1}^3t^i\mathcal E b_i$. Since the stream function is zero on $s=\gamma+u$, its tangential derivative gives $\phi_\Theta=-u_\Theta(\psi_0'+\phi_s)$ there, to the required order. If $[\cos k\Theta]$ denotes the cosine coefficient, the cubic answer is
\[
 -\frac{[\cos k\Theta][t^3]}{2S_0}\!
 \left\{e^{-2(\gamma+u)}(1+u_\Theta^2)
     [\psi_0'(\gamma+u)+\phi_s^{[3]}(\Theta,\gamma+u)]^2
       -e^{2(\gamma+u)}+2\alpha_k e^{\gamma+u}\right\}.
\]
Substitution of $E_0,E_{2k}$ and of the multipliers above gives \eqref{eq:alpha2}--\eqref{eq:Pk}. The denominator $2S_0$ is the $k$-th coefficient of the $t^3$ term from varying $\alpha$ on the slice \eqref{eq:lambda}.

For later use, $a_0>0$. In fact, with $x=k\gamma$, $a_0\sinh^2x=(\sinh(2x)+x\cosh(2x)-3x)/2$. This vanishes at zero and its derivative is $3(\cosh(2x)-1)/2+x\sinh(2x)>0$ for $x>0$.
\begin{theorem}\label{thm:k1}
\leavevmode
\begin{enumerate}[label=\textup{(\roman*)}]
\item Let $\kappa=0$. For all $p_0\ne0$ and all $\gamma\in(0,10.7959]$ one has $\alpha_1^{(2)}>0$, so the
bifurcation at $\alpha_1=\alpha_c$ is supercritical. For $(\gamma,r)=(12,\frac25)$ one has
$\alpha_1^{(2)}<0$.
\item Let $0<\gamma\le4.997$. Then $\alpha_1^{(2)}>0$ for every admissible $(\gamma,\kappa,w)$.
\item The triple $(\gamma,\kappa,w)=(6,-\frac{12}{25},\frac1{50})$ is admissible, and
$\alpha_1^{(2)}<0$ there.
\end{enumerate}
\end{theorem}

\begin{proof}
For $\kappa=0$, Theorem~\ref{thm:k1-D} proves
$a_1^2-4a_0a_2<0$ throughout the stated range by explicit inequalities.
Since $a_0>0$, the quadratic $a_0+a_1r+a_2r^2$ is positive.
For arbitrary admissible vorticity, Theorem~\ref{thm:k1-vorticity}
proves the same positivity after writing $\nu=\kappa/w$ and showing that
$A_1(\nu)^2-4a_0A_2(\nu)<0$ on $\nu<\nu_c$.
The two subcritical examples follow from the negative-coefficient
identities of Theorem~\ref{thm:k1-witnesses}.
All three arguments are contained in Section~\ref{app:primarydirect}.
\end{proof}

\begin{theorem}\label{thm:uniform-modes}
For every integer $k\ge1$ the following sufficient conditions hold:
\begin{enumerate}[label=\textup{(\roman*)}]
\item If $\kappa=0$ and $k\gamma\le5/3$, then $\alpha_k^{(2)}>0$ for every $w\ne0$.
\item If $k\gamma\le6/5$, then $\alpha_k^{(2)}>0$ for every admissible $(\gamma,\kappa,w)$.
\end{enumerate}
\end{theorem}
\begin{proof}
For $k=1$, both conclusions follow from Theorem~\ref{thm:k1}.
For $k\ge2$, the discriminant factorisations and direct coefficient bounds
are proved in Theorems~\ref{thm:mode-zero} and~\ref{thm:mode-full} of
Section~\ref{app:allmodesdirect}.
\end{proof}

\begin{theorem}\label{thm:kge2}
Let $2\le k\le6$, and let
\begin{center}
\begin{tabular}{lccccc}\toprule
$k$ & 2 & 3 & 4 & 5 & 6\\\midrule
$\gamma_k^{0}$ & 1.655 & 0.8517 & 0.5695 & 0.4275 & 0.3423\\
$\gamma_k^{\flat}$ & 1.357 & 0.693 & 0.468 & 0.354 & 0.286\\\bottomrule
\end{tabular}
\end{center}
\begin{enumerate}[label=\textup{(\roman*)}]
\item If $\kappa=0$ and $0<\gamma\le\gamma_k^0$, then $\alpha_k^{(2)}>0$ for all $p_0\ne0$.
\item If $0<\gamma\le\gamma_k^\flat$, then $\alpha_k^{(2)}>0$ for every admissible $(\gamma,\kappa,w)$.
\item For $\kappa=0$ and $\gamma=\gamma_k^0+10^{-3}$, and for $\gamma=\gamma_k^\flat+2\times10^{-3}$, there are admissible
$(\gamma,\kappa,w)$ with $\alpha_k^{(2)}<0$.
\end{enumerate}
\end{theorem}

\begin{proof}
With $\nu=\kappa/w$, the polynomial has the form
$P_k=a_0+A_1(\nu)r+A_2(\nu)r^2$. Its discriminant is quadratic in $\nu$:
the terms of degrees four and three cancel identically. By
Proposition~\ref{prop:capkge2}, the relevant discriminant is negative
throughout each stated depth range. Since $a_0>0$,
$A_2>A_1^2/(4a_0)\ge0$ and hence $P_k>0$ for every $r$.
For (iii), the rational parameters and strictly negative interval
enclosures in Section~\ref{app:cap} give the stated admissible witnesses.
Theorem~\ref{thm:direction} then determines the sign of $\alpha_k^{(2)}$.
\end{proof}

For comparison with \cite[Examples 6.1--6.2]{LWpreprint}, put $Y=e^\gamma$ and $T=\coth\gamma$. At the first critical point the convention of \cite[(4.3)]{LWpreprint} gives
\[
 r=\frac{2\lambda/Y^2-1}{3+2T},\qquad p_0^2=\gamma^2Y^4r,\qquad
 \alpha_1=Y[1+r(1+T)].
\]
Substitution in \eqref{eq:alpha2} gives the following values; $\lambda$ denotes its value at the critical point.
\begin{center}\begin{tabular}{ccrrr}\toprule
$\gamma$&$\lambda$&$p_0^2$&$\alpha_1$&$\alpha_1^{(2)}$\\\midrule
0.3&1.15&0.00794366757&1.508927834&22.50601846\\
0.2&1.40&0.00594402292&1.716147572&2015.14799993\\\bottomrule
\end{tabular}\end{center}
A change of amplitude $a=bt+O(t^2)$ multiplies the quadratic coefficient by $b^2>0$ when $b\ne0$. There is a separate discrepancy in the second-order correction printed in arXiv v1 of \cite{LWpreprint}. Equation (2.23) there requires a zero trace at $p=-1$. In its (4.30), with $y=e^{-2\gamma}$ and the normalization $\widehat C$ used there, the trace gives
\[
 \widehat C^2[\cos(2q)]h(q,-1)
 =\frac{33-22y-6y^2-5y^4}{18}
 =\frac{(1-y)(5y^3+5y^2+11y+33)}{18}>0.
\]
Here brackets extract the indicated Fourier coefficient. At $\gamma=0.3$ this quantity is approximately $1.03696605$. Adding a first-mode kernel element preserves the bed trace. Our computation uses the homogeneous bed condition at each order; its sign also follows from Theorem~\ref{thm:k1}. This comparison concerns the cited preprint version.

\subsection{The primary discriminant}
\label{sec:k1-direct}\label{app:primarydirect}

We set $x=\gamma$ and $t=\coth x$, with the angular velocity parameter $w\ne0$ of Section~\ref{sec:formulation}.
We write $r=w^2$, $\nu=\kappa/w$ and $\nu_c=(1+t)/2$; admissibility is the condition $\nu<\nu_c$.
The primary Stokes coefficient is
\begin{equation}\label{eq:k1-alpha}
\alpha_1^{(2)}=\frac{e^x P_1(\kappa,w)}{8r(xt-1)},\qquad
P_1=a_0+A_1(\nu)r+A_2(\nu)r^2,
\end{equation}
where
\begin{align*}
A_1(\nu)&=a_1-12X\nu-8a_0\nu^2,\\
A_2(\nu)&=a_2+b_1\nu+b_2\nu^2+48X\nu^3+16a_0\nu^4,
\end{align*}
with $X=X_1$ and the coefficients of \eqref{eq:Pk} evaluated at $k=1$. Decimal endpoints in sign estimates denote exact rationals.

\begin{lemma}\label{lem:k1-aux}
For $x>0$, define
\[
\omega(x)=\frac{x\coth x-1}{x^2},\qquad
u(x)=\coth x-1=\frac{2}{e^{2x}-1}.
\]
Then
\begin{equation}\label{eq:k1-omega}
0<\omega<\frac13\quad(x>0),\qquad
\frac14<\omega<\frac13\quad(0<x\le2).
\end{equation}
Both $xu(x)$ and $x^2u(x)$ decrease for $x\ge1$. Moreover, $a_0>0$ and $xt-1>0$ for every $x>0$.
\end{lemma}

\begin{proof}
We set $h=x\cosh x-\sinh x$. Since $h'=x\sinh x>0$ and $[(1+x^2/3)\sinh x-x\cosh x]'=xh/3$, the upper and lower bounds $0<\omega<1/3$ follow. For $g=x\sinh x\cosh x+x^2-2\sinh^2x$ we have $g(0)=g'(0)=0$, $g''=4\sinh x\,h>0$ and $\omega'=-g/(x^3\sinh^2x)<0$. Hence $\omega(x)\ge\omega(2)>1/4$ for $x\le2$. For $j=1,2$ the logarithmic derivative of $x^ju(x)$ is $j/x-2/(1-e^{-2x})<0$ on $x\ge1$. Positivity of $a_0$ follows from the calculation preceding Theorem~\ref{thm:k1}.
\end{proof}

We record the elementary bounds
\begin{equation}\label{eq:k1-exp}
e^2>7,\qquad e^4>50,\qquad e^{9.8}>10000.
\end{equation}
These follow by truncating the positive exponential series at orders four and seven for $e^2$ and $e^4$, and using $e^{9.8}=(e^4)^2e^{1.8}$.
\begin{theorem}\label{thm:k1-D}
Let $D=a_1^2-4a_0a_2$. Then
\[
D(x)<0\qquad(0<x\le10.7959).
\]
Consequently, if $\kappa=0$ and $0<\gamma\le10.7959$, the primary bifurcation is supercritical for every $w\ne0$.
\end{theorem}

\begin{proof}
Substituting $t=x^{-1}+x\omega$ into the displayed coefficients, we obtain
\begin{equation}\label{eq:k1-Domega}
\begin{aligned}
D={}&48\omega^2-96\omega+x(108\omega^2-216\omega)\\
&+x^2(192\omega^3-152\omega^2-144\omega)\\
&+x^3(324\omega^3-276\omega^2-24\omega)\\
&+x^4(288\omega^4-16\omega^3-63\omega^2)\\
&+x^5(324\omega^4-60\omega^3)\\
&+x^6(192\omega^5+40\omega^4)+108\omega^5x^7+48\omega^6x^8.
\end{aligned}
\end{equation}

\emph{The range $0<x\le2$.}
In \eqref{eq:k1-Domega} the terms through degree four sum to at most $-\omega(80+180x+144x^2)$. Dividing the remaining terms by $\omega x^2$ and using $x\le2$, $\omega<1/3$, we obtain the upper bound $128/3+1664/27+128/3+1024/81<160$. Thus
\begin{equation}\label{eq:k1-Dsmall}
D\le\omega(-80-180x+16x^2)
\le\omega(-80-148x)<0.
\end{equation}

\emph{The range $x\ge2$.}
Substitution of $t=1+u$ gives the identity
\begin{equation}\label{eq:k1-Du}
D=(x-1)(49x-529)+\sum_{j=1}^{6}c_j(x)u^j,
\end{equation}
where
\[
\begin{array}{c|l}
j&c_j(x)\\\hline
1&658x^2-2082x+944\\
2&1797x^2-2880x+748\\
3&2140x^2-2000x+300\\
4&1300x^2-696x+48\\
5&396x^2-96x\\
6&48x^2.
\end{array}
\]
For $x\ge2$, the table gives $c_j/x^2\le M_j$, where $M=(894,1984,2215,1312,396,48)$. Since $u<1/24$ by \eqref{eq:k1-exp},
\begin{equation}\label{eq:k1-error}
\sum_{j=1}^6 c_j u^j\le x^2u\sum_{j=1}^6 M_j24^{1-j}<1000x^2u.
\end{equation}
For $2\le x\le10$, the convex quadratic $(x-1)(49x-529)$ is at most $-351$:
its endpoint values are $-431$ and $-351$. By Lemma~\ref{lem:k1-aux}, $D<-351+\frac{8000}{e^4-1}<-351+\frac{8000}{49}<0.$ For $10\le x\le107959/10000$, one has $u\le\frac2{e^{20}-1}<10^{-8},\qquad x^2<121,$ because $e^{20}=(e^2)^{10}>7^{10}>200000001$.
The quadratic in \eqref{eq:k1-Du} increases on this interval, so $(x-1)(49x-529)\le-\frac{881631}{10^8}.$ Therefore $D<-\frac{881631}{10^8}+\frac{121}{100000} =-\frac{760631}{10^8}<0.$ This proves the asserted range. For $\kappa=0$, one has
$P_1=a_0+a_1r+a_2r^2$. Since $a_0>0$ and $D<0$,
$a_2>a_1^2/(4a_0)\ge0$, and hence $P_1>0$ for all real $r$.
Equation~\eqref{eq:k1-alpha} gives $\alpha_1^{(2)}>0$.
\end{proof}

Define
\begin{align*}
Q_2={}&xt^3+xt^2-t^2+xt-t-2x-3,\\
Q_D={}&12xt^5+20xt^4-12t^4-20xt^3-20t^3\\
&\quad-36xt^2-4t^2+9xt-4t+16x-9,\\
Q_S={}&6xt^4+7xt^3-6t^3-4xt^2-7t^2-5xt-8t-4x-9.
\end{align*}

\begin{lemma}\label{lem:k1-Qsigns}
The following strict inequalities hold:
\[
\begin{array}{c|c}
\text{factor}&\text{range}\\\hline
Q_2<0&0<x\le4.997\\
Q_D<0&0<x\le5\\
Q_S<0&0<x<\infty.
\end{array}
\]
\end{lemma}

\begin{proof}
\emph{The factor $Q_2$.}
The two substitutions used above give
\begin{align}
Q_2={}&\omega^3x^4+\omega^2x^3+(2\omega^2+\omega)x^2
+(\omega-2)x+\omega-2,\label{eq:k1-Q2omega}\\
Q_2={}&x-5+(6x-3)u+(4x-1)u^2+xu^3.\label{eq:k1-Q2u}
\end{align}
For $0<x\le2$, \eqref{eq:k1-omega} gives
\begin{align*}
Q_2&\le-\frac53-\frac53x+\frac59x^2+\frac19x^3+\frac1{27}x^4\\
&\le-\frac53-\frac53x+\frac{25}{27}x^2
\le-\frac53+\frac5{27}x
\le-\frac{35}{27}<0.
\end{align*}
For $x\ge2$, one has $u<1/24$, so
\begin{equation}\label{eq:k1-Q2bound}
Q_2<x-5+7xu.
\end{equation}
The monotonicity of $xu$, together with \eqref{eq:k1-exp}, yields
\[
\begin{array}{c|c}
\text{range of }x&\text{upper bound for }Q_2\\\hline
{[2,4]}&-1+28/49<0\\
{[4,4.9]}&-1/10+56/2499<0\\
{[4.9,4.99]}&-1/100+70/9999<0.
\end{array}
\]
For the last row it is enough to use $e^{2x}>10000$ and $x\le5$ directly.
To reach $4.997$, we use the additional elementary bound
\begin{equation}\label{eq:k1-exp-refined}
e^{9.98}>21000.
\end{equation}
Indeed, $e^2>147/20$, since its Taylor sum through degree six is $331/45>147/20$.
The geometric-series comparison gives $e^{1/50}<50/49$. Consequently, $e^{9.98}>\frac{49}{50}\left(\frac{147}{20}\right)^5 =\frac{3363432789843}{160000000}>21000.$ For $4.99\le x\le4.997$, we therefore have $u<2/20999<1/10000$.
Using \eqref{eq:k1-Q2u},
\begin{align*}
Q_2&\le x-5+xu(6+4u+u^2)\\
&<-\frac3{1000}+\frac{6001}{2099900}<0.
\end{align*}
Here we used $x\le5$ and $6+4u+u^2<6001/1000$.

\emph{The factor $Q_D$.}
The small-$x$ identity is
\begin{equation}\label{eq:k1-QDomega}
\begin{aligned}
x^2Q_D={}&12\omega-24+x(20\omega-40)
+x^2(48\omega^2-68\omega)\\
&+x^3(60\omega^2-76\omega+16)\\
&+x^4(72\omega^3-64\omega^2+9\omega)\\
&+x^5(60\omega^3-36\omega^2)+x^6(48\omega^4-20\omega^3)\\
&+20\omega^4x^7+12\omega^5x^8.
\end{aligned}
\end{equation}
For $0<x\le2$, Lemma~\ref{lem:k1-aux} gives $1/4<\omega<1/3$.
The constant and linear terms are at most $-20$ and $-100x/3$.
The coefficients of $x^2,x^4,x^5,x^6$ are negative; for the coefficient of $x^4$, use $72\omega^2-64\omega+9 \le72(1/4)^2-64(1/4)+9=-\frac52,$ since this quadratic decreases on $[1/4,1/3]$.
The coefficient of $x^3$ is at most $16-56\omega\le2$.
It follows that
\begin{align*}
x^2Q_D&\le-20-\frac{100}{3}x+2x^3
+\frac{20}{81}x^7+\frac4{81}x^8\\
&\le-20-\frac{260}{81}x<0.
\end{align*}
For $2\le x\le5$, the other identity is
\begin{equation}\label{eq:k1-QDu}
\begin{aligned}
Q_D={}&x-49+(17x-120)u+(144x-136)u^2\\
&+(180x-68)u^3+(80x-12)u^4+12xu^5.
\end{aligned}
\end{equation}
The coefficient of $u$ is negative. Since $u<1/24$, $Q_D\le-44+\frac{720}{24^2}+\frac{900}{24^3} +\frac{400}{24^4}+\frac{60}{24^5}<-42<0.$ \emph{The factor $Q_S$.}
For small $x$ we use
\begin{equation}\label{eq:k1-QSomega}
\begin{aligned}
xQ_S={}&6\omega-12+x(7\omega-14)
+x^2(18\omega^2-16\omega-4)\\
&+x^3(14\omega^2-5\omega)+x^4(18\omega^3-4\omega^2)\\
&+7\omega^3x^5+6\omega^4x^6.
\end{aligned}
\end{equation}
If $0<x\le2$, the coefficients of $x^2$ and $x^3$ are negative, and
\begin{align*}
xQ_S&\le-10-\frac{35}{3}x+\frac29x^4+\frac7{27}x^5+\frac2{27}x^6\\
&\le-10-\frac{91}{27}x<0.
\end{align*}
For every $x\ge2$, use
\begin{equation}\label{eq:k1-QSu}
Q_S=-30+(32x-40)u+(53x-25)u^2+(31x-6)u^3+6xu^4.
\end{equation}
The bounds $u<1/24$ and
$32+53/24+31/24^2+6/24^3<35$ give $Q_S<-30+35xu\le-30+\frac{140}{49}<0.$ Thus the last inequality holds on the whole positive half-line.
\end{proof}

\begin{theorem}\label{thm:k1-vorticity}
For $0<\gamma\le4.997$, the primary bifurcation is supercritical for every admissible pair $(\kappa,w)$.
\end{theorem}

The coefficient estimates proving this assertion are given in Appendix~\ref{app:details-k1-vorticity}.

\begin{theorem}\label{thm:k1-witnesses}
The primary bifurcation is subcritical at
\[
(\gamma,\kappa,r)=(12,0,2/5)
\]
and at the admissible triple
\[
(\gamma,\kappa,w)=(6,-12/25,1/50).
\]
\end{theorem}

\begin{proof}
For the first parameter choice, put $u=\coth12-1>0$.
Direct substitution into the coefficients gives
\begin{equation}\label{eq:k1-witness-zero}
\begin{aligned}
25P_1(0,\sqrt{2/5})=-\bigl(&45+5951u+19264u^2+23936u^3\\
&+14148u^4+3996u^5+432u^6\bigr)<-45.
\end{aligned}
\end{equation}
For the second parameter choice, put $u=\coth6-1>0$.
The corresponding identity is
\begin{equation}\label{eq:k1-witness-vorticity}
\begin{aligned}
6250000P_1(-12/25,1/50)=-\bigl(&3942+180263u+479407u^2\\
&+463373u^3+158061u^4\\
&+5679u^5+54u^6\bigr)<0.
\end{aligned}
\end{equation}
Admissibility in the second case follows from $\nu=-24<\nu_c$;
in the first case it follows from $\nu=0<\nu_c$.
Every coefficient inside each parenthesis is positive. Thus $P_1<0$ in both cases,
and \eqref{eq:k1-alpha} gives $\alpha_1^{(2)}<0$.
\end{proof}

\subsection{Uniform estimates in the Fourier mode}
\label{sec:all-modes-direct}\label{app:allmodesdirect}

In this subsection we derive sufficient conditions valid for every Fourier mode. The proof uses factorisations of the discriminant and applies directly to all modes. The sharper finite-mode bounds of Theorem~\ref{thm:kge2} are treated separately in Section~\ref{app:cap}.

Use $r=w^2$, $\nu=\kappa/w$ and the polynomial $P_k$ of Theorem~\ref{thm:direction}. With $x=k\gamma$, $T=\coth x$ we have $a_0=T-x(T^2-2)>0$ and
\begin{equation}\label{eq:all-mode-prefactor}
 \alpha_k^{(2)}=\frac{e^\gamma P_k}{8kr(xT-1)}.
\end{equation}
The admissibility condition is
$2\nu<1+\coth\gamma$. We write
\[
 \Delta_k(\nu)=A_1(\nu)^2-4a_0A_2(\nu),\qquad
 D_k=\Delta_k(0)=a_1^2-4a_0a_2.
\]

\begin{lemma}\label{lem:mode-hyperbolic}
For $x>0$, set
\[
 h=x\coth x-1,\qquad d=\operatorname{csch}^2x=T^2-1,
 \qquad z=hd.
\]
Then
\begin{equation}\label{eq:mode-hyperbolic-bounds}
 0<h<\frac{x^2}{3},\qquad
 h>\frac{x^2}{3+x^2},\qquad
 0<z<\frac13,\qquad Th<x.
\end{equation}
In addition,
\begin{equation}\label{eq:mode-d-bounds}
 0<x\le1\ \Longrightarrow\ d>\frac23,
 \qquad
 0<x\le\frac32\ \Longrightarrow\ d>\frac15.
\end{equation}
\end{lemma}

\begin{proof}
The upper bound for $h$ is \eqref{eq:k1-omega}. For the lower bound, $x(3+x^2)\cosh x-(3+2x^2)\sinh x$ vanishes at zero and has positive derivative $x(x\cosh x-\sinh x+x^2\sinh x)$. Likewise $\sinh^3x+3\sinh x-3x\cosh x$ has derivative $3\sinh x(\sinh x\cosh x-x)>0$, giving $z<1/3$. The same inequality $\sinh(2x)>2x$ gives $Th<x$.
For $j\ge1$, $(2j+1)!\ge6\cdot20^{j-1}$, hence $\sinh x<x+x^3/[6(1-x^2/20)]$ for $0<x<\sqrt{20}$. At $1$ and $3/2$ the right side is respectively below $\sqrt{3/2}$ and $\sqrt5$, proving \eqref{eq:mode-d-bounds}.
\end{proof}

\begin{lemma}\label{lem:mode-D-factor}
For every $k\ge1$, the discriminant at zero vorticity has the factorisation
\begin{equation}\label{eq:mode-D-factor}
 D_k=khR_0,
\end{equation}
where
\begin{align}
 R_0={}&16k^3\left[d(3z-6)+2z+3h-2\right]\notag\\
 &+36k^2\left[T(3z-6)+2x\right]\notag\\
 &+k(104z-23h-208)+48(Th-2x).
 \label{eq:mode-R0}
\end{align}
\end{lemma}

\begin{proof}
Substitute $h=xT-1$, $d=T^2-1$ and $z=(xT-1)(T^2-1)$ into
\eqref{eq:mode-R0} and collect powers of $k$. Expansion gives
$a_1^2-4a_0a_2=khR_0$. This is a polynomial identity in $(k,x,T)$;
no sign information from a symbolic or numerical computation is needed.
\end{proof}

\begin{theorem}
\label{thm:mode-zero}
For every integer $k\ge1$, every nonzero flux, and $\kappa=0$,
\[
 0<k\gamma\le\frac53\quad\Longrightarrow\quad\alpha_k^{(2)}>0.
\]
\end{theorem}

\begin{proof}
Suppose that $0<x\le3/2$. Lemma~\ref{lem:mode-hyperbolic} gives
\[
\begin{split}
 d(3z-6)+2z+3h-2
 &<-5d+3h-\frac43\\
 &<-1+\frac94-\frac43=-\frac1{12},
\end{split}
\]
and $T(3z-6)+2x<-5T+2x<-5+3=-2.$ The remaining two terms in \eqref{eq:mode-R0} are strictly negative because $104z-23h-208<\frac{104}{3}-208=-\frac{520}{3}, \qquad Th-2x<-x.$ These estimates prove $R_0<0$ for $0<x\le3/2$.

We extend the same conclusion to $3/2\le x\le5/3$. The function $h$ is strictly
increasing, whereas $z$ is strictly decreasing. For the latter assertion, write $z'(x)=-\frac{B(x)}{\sinh^4x},\qquad B(x)=x(\cosh2x+2)-3\sinh x\cosh x.$ Here $B(0)=0$ and $B'(x)=4\sinh x\,(x\cosh x-\sinh x)>0.$ The elementary Taylor estimate $e^3>\sum_{j=0}^{8}\frac{3^j}{j!} =\frac{89641}{4480}>20$ gives
\[
 T(3/2)<\frac{21}{19},\qquad
 h(3/2)<\frac{25}{38},\qquad
 d(3/2)<\frac{80}{361},\qquad
 z(x)\le z(3/2)<\frac{1000}{6859}<\frac16.
\]
Also $e^{10/3}=e^3e^{1/3}>20(1+1/3)>26$, so $h(x)\le h(5/3)<\frac53\frac{27}{25}-1=\frac45.$ The geometric Taylor bound used in Lemma~\ref{lem:mode-hyperbolic} gives $\sinh\frac53<\frac53+ \frac{(5/3)^3}{6(1-(5/3)^2/20)} =\frac{715}{279}<\sqrt{\frac{20}{3}},$ whence $d(x)>3/20$ throughout this interval. Consequently, $d(3z-6)+2z+3h-2 <-\frac{11}{2}\frac3{20}+\frac13+\frac{12}{5}-2 =-\frac{11}{120}<0,$ and $T(3z-6)+2x<-\frac{11}{2}+\frac{10}{3}<0.$ The last two terms of \eqref{eq:mode-R0} remain strictly negative by
\eqref{eq:mode-hyperbolic-bounds}. Thus $R_0<0$ and $D_k<0$ for
$0<x\le5/3$. Since $a_0>0$, it follows that
$a_2>a_1^2/(4a_0)\ge0$, and the quadratic
$P_k=a_0+a_1r+a_2r^2$ is positive for every real $r$.
The prefactor in \eqref{eq:all-mode-prefactor} is positive, proving the claim.
\end{proof}

The useful variable for the full discriminant is
\begin{equation}\label{eq:mode-b}
 b=kT+1-2\nu.
\end{equation}
It is positive throughout the physical admissibility range. Indeed,
$x\mapsto x\coth x$ is strictly increasing, since its derivative is
$(\sinh x\cosh x-x)/\sinh^2x>0$. Therefore
\[
 \coth\gamma\le k\coth(k\gamma)=kT,
 \qquad
 2\nu<1+\coth\gamma\le1+kT.
\]

\begin{lemma}\label{lem:mode-full-factor}
For every $k\ge1$,
\begin{equation}\label{eq:mode-full-factor}
 \Delta_k(\nu)=kh\left(12q_2b^2+Lb+C\right),
\end{equation}
where
\begin{align}
 q_2={}&k(2h+z-2)+Th-2x,\label{eq:mode-q2}\\
 L={}&k^2\left[T(-16h+24z-48)+32x\right]\notag\\
 &+k(8h+28z-56)+24(Th-2x),\label{eq:mode-L}\\
 C={}&k^3\left[d(12z-24)-12z+8h+8\right]\notag\\
 &+k^2\left[T(-28h+20z-40)+40x\right]\notag\\
 &+k(16z-7h-32)+12(Th-2x).
 \label{eq:mode-C}
\end{align}
\end{lemma}

\begin{proof}
The coefficients of $\nu^4$ and $\nu^3$ in $A_1^2-4a_0A_2$ cancel.
Substitute $2\nu=kT+1-b$ in the remaining quadratic. Collection of the
coefficients of $b^2$, $b$ and $1$, followed by
$h=xT-1$, $d=T^2-1$ and $z=hd$, gives the displayed identity.
\end{proof}

\begin{theorem}
\label{thm:mode-full}
For every integer $k\ge1$ and every admissible triple $(\gamma,\kappa,w)$,
\[
 0<k\gamma\le\frac65\quad\Longrightarrow\quad\alpha_k^{(2)}>0.
\]
\end{theorem}

The coefficient estimates proving this assertion are given in Appendix~\ref{app:details-mode-full}.

\begin{remark}
These sufficient conditions cover every Fourier mode. For $2\le k\le6$,
the finite interval certificates in Section~\ref{app:cap} extend the
depth ranges towards the numerical transition values. No finite-mode
certificate is needed for the uniform criteria proved here.
\end{remark}

\subsection{Sharpening the bounds in five modes}\label{app:cap}
It remains to determine the sign of the discriminant on the compact intervals specified below.
\begin{proposition}\label{prop:capkge2}
For $2\le k\le6$, the discriminant $D_k$ is negative for
$0<\gamma\le\gamma_k^0$. For every $0<\gamma\le\gamma_k^\flat$,
\[
 \Delta_k(\nu)<0\qquad
 \text{for every }\nu<\nu_c=\tfrac12(1+\coth\gamma).
\]
Here $\gamma_k^0$ and $\gamma_k^\flat$ are the exact terminating
decimals in Theorem~\ref{thm:kge2}.
\end{proposition}

\begin{proof}
Theorems~\ref{thm:mode-zero} and~\ref{thm:mode-full} give the assertions,
respectively, when $k\gamma\le5/3$ and $k\gamma\le6/5$.
On the remaining intervals $I_k^0=[5/(3k),\gamma_k^0],\qquad I_k^\flat=[6/(5k),\gamma_k^\flat],$ write $\Delta_k(\nu)=d_2\nu^2+d_1\nu+D_k$.
We certify $D_k<0$ on $I_k^0$, and the three inequalities $d_2<0,\qquad \Delta_k(\nu_c)<0,\qquad \Delta_k'(\nu_c)>0 \quad\text{on }I_k^\flat.$ Together these inequalities imply $\Delta_k(\nu)<0$ on $\nu<\nu_c$ by concavity
and monotonicity. No assertion that these three separate inequalities
follow from the uniform small-depth theorem is needed.

The evaluation uses the exact factorisations
\eqref{eq:mode-D-factor} and~\eqref{eq:mode-full-factor}.
With $h=k\gamma\coth(k\gamma)-1$ and the coefficients $q_2,L,C,R_0$
defined there, put $b_c=k\coth(k\gamma)-\coth\gamma$. Then
\begin{align*}
 D_k&=khR_0,\qquad d_2=48khq_2,\\
 \Delta_k(\nu_c)&=kh(12q_2b_c^2+Lb_c+C),\\
 \Delta_k'(\nu_c)&=-2kh(24q_2b_c+L).
\end{align*}
These forms avoid the large cancellations in direct expansion of $P_k$.

At 160-bit precision the Arb verifier \cite{Arb} bisects rational intervals until each ball enclosure has the prescribed strict sign. Sorting the accepted intervals and checking their rational endpoints verifies both outer endpoints and equality of every pair of adjacent endpoints. All twenty tests terminate, using $217$ accepted intervals in total. The source supplement includes the verifier and the exact partitions with exact rational lower and upper bounds for every ball. A separate reader verifies the saved signs, ordering, outer endpoints and adjacency using rational arithmetic.
\end{proof}

The following parameters are exact terminating decimals. For each row,
the verifier evaluates $P_k=a_0+A_1(\nu)r+A_2(\nu)r^2$ at $200$-bit
precision, proves $r>0$ and $\nu<\nu_c$, and obtains a strictly negative
ball. The values of $P_k$ in the table are rounded centres for readability;
the full certified enclosures are included in the certificate file.
The triple $(\gamma,\nu\sqrt r,\sqrt r)$ is therefore admissible and
subcritical in each case.
\begin{center}\small
\setlength{\tabcolsep}{4pt}
\begin{tabular}{c|ccc|cccc}\toprule
& \multicolumn{3}{c|}{$\kappa=0$} & \multicolumn{4}{c}{general $\kappa$}\\
$k$ & $\gamma$ & $r$ & $P_k$ & $\gamma$ & $\nu$ & $r$ & $P_k$\\\midrule
2 & $1.656$ & $0.132296$ & $-0.0035871$ & $1.359$ & $1.070670418$ & $1.57085244$ & $-0.249317$\\
3 & $0.8527$ & $0.070501$ & $-0.0085586$ & $0.695$ & $1.331690455$ & $0.45945318$ & $-0.159498$\\
4 & $0.5705$ & $0.043859$ & $-0.0119594$ & $0.470$ & $1.641032277$ & $0.25955723$ & $-0.223669$\\
5 & $0.4285$ & $0.029968$ & $-0.0155842$ & $0.356$ & $1.963331379$ & $0.16855206$ & $-0.201703$\\
6 & $0.3433$ & $0.021821$ & $-0.0209378$ & $0.288$ & $2.28384677$ & $0.1343368$ & $-0.482917$\\\bottomrule
\end{tabular}
\end{center}

For each of the ten parameter choices, the accompanying program reconstructs the partition and checks both admissibility and strict negativity.

\section{Global continuation in the irrotational case}\label{sec:global}

In Sections~\ref{sec:BT}--\ref{sec:corner} we take $\eps=0$, $\kappa=0$ and, as in \cite{LWpreprint},
$p_0<0$ (for $\kappa=0$ the problem depends on $p_0$ only through $p_0^2$); $k\ge1$ is fixed.
Section~\ref{sec:kappaglobal} extends Sections~\ref{sec:BT}--\ref{sec:globalthm} to every admissible
$(\gamma,\kappa,w)$, and Section~\ref{sec:kstag} extends the case $m=1$ of Section~\ref{sec:slope}, and with
it part of Section~\ref{sec:stagnation}, to $\kappa\ne0$. We write $\varkappa$ for the curvature of the free surface,
$\Fop(\alpha,H):=\Fop(\alpha,0;H)$
and $\Bop_0$ for $\Bop_\eps$ with $\eps=0$. By Lemma~\ref{lem:sym}(b), $\Fop$ maps the open set
\[
\mathcal V_k:=\R\times(\Xs_k\cap\Uo)\subset\R\times\Xs_k
\]
into $\Ys_k$. For $x>0$ we call $x(1+p)$ the \emph{laminar flow of depth $x$}, and we write
\[
\alpha_0(x):=\e^{x}+p_0^2\,\e^{-3x}\,\frac{1+x}{x^3},\qquad
\alpha_k(x):=\e^{x}+p_0^2\,\e^{-3x}\,\frac{1+k\coth(kx)}{x^2}
\]
for its critical values. These are the numbers of Proposition~\ref{prop:dispersion} with
$(\gamma,\cc)$ replaced by $(x,\,p_0^2\e^{-2x}/x^2)$; in particular $\alpha_k(\gamma)=\alpha_k$ and
$\alpha_0(\gamma)=\alpha_0$.

We first establish the continuation and estimates for zero relative vorticity. The proof for general vorticity in Section~\ref{sec:kappaglobal} uses the same nodal argument and modifies the velocity and flow-force estimates.
\subsection{An analytic global bifurcation theorem}\label{sec:BT}
We use the following result of Dancer \cite{Dancer1973} and of Buffoni and Toland
\cite[Theorem~9.1.1]{BuffoniToland2003}, in the form of \cite[Theorem~6]{CSV2016}.

\begin{theorem}\label{thm:BT}
Let $X$, $Y$ be Banach spaces, $\mathcal V\subset\R\times X$ open, and
$F:\mathcal V\to Y$ real-analytic. Assume:
\begin{enumerate}[label=\textup{(H\arabic*)}]
\item $(\lambda,0)\in\mathcal V$ and $F(\lambda,0)=0$ for all $\lambda\in\R$;
\item for some $\lambda_*$, $\ker D_xF(\lambda_*,0)=\spn\{\xi_*\}$ is one-dimensional,
$\mathrm R(D_xF(\lambda_*,0))$ is closed of codimension one, and
$D_\lambda D_xF(\lambda_*,0)\xi_*\notin\mathrm R(D_xF(\lambda_*,0))$;
\item $D_xF(\lambda,x)$ is Fredholm of index zero whenever $(\lambda,x)\in\mathcal V$ and
$F(\lambda,x)=0$;
\item there are closed bounded sets $Q_j\subset\mathcal V$, $j\in\N$, with
$\mathcal V=\bigcup_jQ_j$, such that $F^{-1}(0)\cap Q_j$ is compact for every $j$.
\end{enumerate}
Then there is a continuous curve $\mathcal K=\{(\lambda(s),x(s)):s\in\R\}\subset\mathcal V$ of
solutions of $F=0$ such that:
\begin{enumerate}[label=\textup{(C\arabic*)}]
\item $(\lambda(0),x(0))=(\lambda_*,0)$;
\item $x(s)=s\xi_*+o(s)$ in $X$ as $s\to0$;
\item there are a neighbourhood $W$ of $(\lambda_*,0)$ and $\epsilon>0$ such that
$\{(\lambda,x)\in W:\ x\ne0,\ F(\lambda,x)=0\}=\{(\lambda(s),x(s)):\ 0<|s|<\epsilon\}$;
\item The curve is the maximal continuation of the local bifurcation arc by
distinguished real-analytic arcs. In particular, for every $s_0\in\R$ there are
an interval $I$ containing $s_0$ and an increasing homeomorphism
$\sigma:(-1,1)\to I$, with $\sigma(0)=s_0$, such that
$t\mapsto(\lambda(\sigma(t)),x(\sigma(t)))$ is injective and real-analytic.
\item one of the following holds:
\textup{(I)} for every $j$ there is $s_j>0$ such that $(\lambda(s),x(s))\notin Q_j$ for $|s|>s_j$;
\textup{(II)} there is $T>0$ such that $(\lambda(s+T),x(s+T))=(\lambda(s),x(s))$ for all $s\in\R$.
\end{enumerate}
The maximally continued route, with the orientation fixed by \textup{(C2)},
is unique up to an increasing reparametrisation.
\end{theorem}

The distinguished-arc convention in \cite[Theorem~6 and Remarks~7--8]{CSV2016} fixes the outgoing germ at a singular point by the incoming injective analytic uniformisation. The no-return argument below uses uniqueness of this oriented route and its local injectivity.

We apply Theorem~\ref{thm:BT} with $X=\Xs_k$, $Y=\Ys_k$, $\lambda=\alpha$, $x=H-H_0$,
$\mathcal V=\{(\alpha,x):(\alpha,H_0+x)\in\mathcal V_k\}$, $\lambda_*=\alpha_k$ and
$\xi_*=\hat\varphi_k$. Hypothesis (H1) holds because $H_0$ is a solution for every $\alpha$.
Hypothesis (H2) is Lemma~\ref{lem:range}: in $\Xs_k$ the kernel of $\Lop(\alpha_k)$ is still
spanned by $\hat\varphi_k$, since $\alpha_{jk}\ne\alpha_k$ for $j\ne1$. We verify (H3) and (H4)
next.

\subsection{Fredholm property and compactness}
For $(\alpha,H)\in\mathcal V_k$ we have
$D_H\Fop(\alpha,H)\psi=(\mathcal A_H\psi,\ \mathcal B_{\alpha,H}\psi)$, where
\begin{align*}
\mathcal A_H\psi&=(1+H_q^2)\psi_{pp}-2H_qH_p\psi_{qp}+H_p^2\psi_{qq}
+2(H_qH_{pp}-H_pH_{qp})\psi_q+2(H_pH_{qq}-H_qH_{qp})\psi_p,\\
\mathcal B_{\alpha,H}\psi&=\big(b_p\psi_p+b_q\psi_q+b_0\psi\big)\big|_{p=0},
\end{align*}
with
\begin{gather*}
b_p=-\frac{2\cc\gamma^2(1+H_q^2)}{H_p^3}\e^{-2(H-\gamma)},\qquad
b_q=\frac{2\cc\gamma^2H_q}{H_p^2}\e^{-2(H-\gamma)},\\
b_0=-\frac{2\cc\gamma^2(1+H_q^2)}{H_p^2}\e^{-2(H-\gamma)}-2\e^{2H}+2\alpha\e^{H}.
\end{gather*}
At $H=H_0$ this is \eqref{eq:L}. The operator $\mathcal A_H$ has no zeroth-order term, because
$\Dop$ depends only on the derivatives of $H$. Its principal symbol is positive definite, with
determinant $H_p^2>0$. The boundary operator is oblique, since $b_p<0$.

\begin{lemma}\label{lem:fredholm-global}
For every $(\alpha,H)\in\mathcal V_k$, $D_H\Fop(\alpha,H):\Xs_k\to\Ys_k$ is Fredholm of index
zero.
\end{lemma}

\begin{proof}
The coefficients of $\mathcal A_H$ belong to $C^{1,\beta}(\overline\Omega)$ and those of
$\mathcal B_{\alpha,H}$ to $C^{2,\beta}(\T)$. Their parities are such that $\Xs_k$ is mapped into
$\Ys_k$. For instance, $H_qH_p$ is odd in $q$ and multiplies $\psi_{qp}$, which is odd when $\psi$
is even. The smallest eigenvalue of the principal symbol is at least
$H_p^2/(1+H_q^2+H_p^2)$, which is bounded below on $\overline\Omega$, and $b_p$ is bounded away
from zero on $\T$. The Schauder theory for the mixed problem with a Dirichlet condition on
$p=-1$ and an oblique condition on $p=0$ \cite{ADN1959,GT2001} gives
\begin{equation}\label{eq:schauder}
\|\psi\|_{C^{3,\beta}}\le C\big(\|\mathcal A_H\psi\|_{C^{1,\beta}}+\|\mathcal B_{\alpha,H}\psi\|_{C^{2,\beta}(\T)}+\|\psi\|_{C^0}\big),
\end{equation}
where $C$ depends only on the ellipticity and obliqueness constants and on the norms of the
coefficients. Since $\Xs_k\hookrightarrow C^0(\overline\Omega)$ is compact, \eqref{eq:schauder}
implies that $D_H\Fop(\alpha,H)$ has a finite-dimensional kernel and a closed range.

Now let $\Mop$ be the isomorphism of Lemma~\ref{lem:fredholm}, restricted to $\Xs_k$, and set
$T_t:=(1-t)D_H\Fop(\alpha,H)+t\Mop$ for $t\in[0,1]$. Each $T_t$ is of the same form: its
principal symbol is a convex combination of positive definite symbols, and the coefficient of
$\psi_p$ in its boundary operator is $(1-t)b_p-2t\cc/\gamma<0$. Hence \eqref{eq:schauder} holds
for $T_t$ with a constant independent of $t$, and every $T_t$ is semi-Fredholm. The map
$t\mapsto T_t$ is continuous in operator norm, and the index of semi-Fredholm operators is
locally constant \cite[Section~IV.5]{Kato1966}. Therefore
$\operatorname{ind}D_H\Fop(\alpha,H)=\operatorname{ind}\Mop=0$.
\end{proof}

\begin{lemma}\label{lem:proper}
For $j\in\N$ let
$Q_j:=\{(\alpha,H)\in\R\times\Xs_k:\ |\alpha|\le j,\ \|H\|_{C^{3,\beta}}\le j,\ \min_{\overline\Omega}H_p\ge1/j\}$.
Then $Q_j\subset\mathcal V_k$ is closed and bounded, $\bigcup_jQ_j=\mathcal V_k$, and
$\Fop^{-1}(0)\cap Q_j$ is compact in $\R\times\Xs_k$.
\end{lemma}

\begin{proof}
Only the compactness needs proof. Let $(\alpha_n,H_n)\in\Fop^{-1}(0)\cap Q_j$. After passing to
a subsequence, $\alpha_n\to\bar\alpha$ and $H_n\to\bar H$ in $C^3(\overline\Omega)$, by the
Arzel\`a--Ascoli theorem. The limit satisfies $(\bar\alpha,\bar H)\in Q_j$ and
$\Fop(\bar\alpha,\bar H)=0$. To obtain convergence in $C^{3,\beta}$, write $0=\Fop(\alpha_n,H_n)-\Fop(\alpha_m,H_m)=\mathcal L_{nm}(H_n-H_m)+(\alpha_n-\alpha_m)\,e_{nm},$ where $\mathcal L_{nm}:=\int_0^1D_H\Fop(\alpha^t,H^t)\dd t$ and
$e_{nm}:=\int_0^1\p_\alpha\Fop(\alpha^t,H^t)\dd t$, with
$(\alpha^t,H^t):=(\alpha_m,H_m)+t(\alpha_n-\alpha_m,H_n-H_m)\in Q_j$. The operator
$\mathcal L_{nm}$ has the structure described above, with ellipticity and obliqueness constants
and coefficient norms that depend only on $j$. Thus \eqref{eq:schauder} gives
\[
\|H_n-H_m\|_{C^{3,\beta}}\le C(j)\big(|\alpha_n-\alpha_m|\,\|e_{nm}\|_{\Ys}+\|H_n-H_m\|_{C^0}\big)\longrightarrow0
\qquad(n,m\to\infty).
\]
Hence $H_n\to\bar H$ in $\Xs_k$.
\end{proof}

\subsection{Nodal properties}

\begin{definition}\label{def:nodal}
Let $\mathcal N_k$ be the set of $H\in\Xs_k$ such that
\begin{enumerate}[label=\textup{(N\arabic*)}]
\item $H_q<0$ on $(0,\pi/k)\times(-1,0]$;
\item $H_{qp}<0$ on $(0,\pi/k)\times\{-1\}$;
\item $H_{qq}(0,p)<0<H_{qq}(\pi/k,p)$ for $p\in(-1,0]$;
\item $H_{qqp}(0,-1)<0<H_{qqp}(\pi/k,-1)$.
\end{enumerate}
\end{definition}

For $H\in\Xs_k$ the function $H_q$ vanishes on the lines $q=0$, $q=\pi/k$ (by evenness and
$2\pi/k$-periodicity) and on $p=-1$. Condition (N1) says that the free surface
$q\mapsto H(q,0)=\ln S(q)$, and every streamline $q\mapsto H(q,p)$, decreases strictly from a
crest at $q=0$ to a trough at $q=\pi/k$. The set $\tau_{\pi/k}\mathcal N_k$ consists of the
profiles with a trough at $q=0$, and $\mathcal N_k\cap\tau_{\pi/k}\mathcal N_k=\emptyset$.
Conditions (N1)--(N4) involve only $H_q$, so $\mathcal N_k$ is invariant under addition of
functions of $p$ and under multiplication by positive numbers.

\begin{lemma}\label{lem:Nopen}
$\mathcal N_k$ is open in $\Xs_k$; in fact it is open in the $C^3(\overline\Omega)$-topology.
\end{lemma}

\begin{proof}
Let $H\in\mathcal N_k$, $w:=H_q$, and let $G\in\Xs_k$ with $v:=G_q$ satisfy
$\|v-w\|_{C^2}<\delta$. Note that $v$ vanishes on $q=0$, $q=\pi/k$ and $p=-1$.

Near the corner $(0,-1)$: by (N4) there is $\eta>0$ with $w_{qp}<0$ on
$[0,\eta]\times[-1,-1+\eta]$, and the same holds for $v_{qp}$ if $\delta$ is small. Since
$v=0$ on $q=0$ and on $p=-1$,
\begin{gather*}
v(q,p)=\int_0^q\!\!\int_{-1}^pv_{qp}\dd p'\dd q'<0,\\
v_q(0,p)=\int_{-1}^pv_{qp}(0,p')\dd p'<0,\qquad
v_p(q,-1)=\int_0^qv_{qp}(q',-1)\dd q'<0
\end{gather*}
for $0<q\le\eta$ and $-1<p\le-1+\eta$. The corner $(\pi/k,-1)$ is treated in the same way.

Near the side $q=0$, away from the bottom: by (N3), $w_q(0,p)\le-c_1<0$ for
$p\in[-1+\eta,0]$, so $v_q<-c_1/2$ on $[0,\eta']\times[-1+\eta,0]$ for some $\eta'>0$ if
$\delta$ is small. Hence $v(q,p)=\int_0^qv_q<0$ there. The side $q=\pi/k$ is similar. Near the
bottom, away from the corners, (N2) gives $w_p\le-c_2<0$ on
$[\eta',\pi/k-\eta']\times\{-1\}$. Therefore $v_p<0$ on $[\eta',\pi/k-\eta']\times[-1,-1+\eta'']$
and $v<0$ there for $p>-1$. On the remaining compact set
$[\eta',\pi/k-\eta']\times[-1+\eta'',0]$, which contains the free surface away from the corners,
$w\le-c_3<0$ by (N1), so $v<0$. Together with (N4) for $G$, which holds for small $\delta$,
this shows $G\in\mathcal N_k$.
\end{proof}

\begin{lemma}\label{lem:Nlocal}
Let $s\mapsto(\alpha_k(s),H_k(s))$ be the curve of Theorem~\ref{thm:local}. There is
$\epsilon_1>0$ such that $H_k(s)\in\mathcal N_k$ for $0<s<\epsilon_1$ and
$H_k(s)\in\tau_{\pi/k}\mathcal N_k$ for $-\epsilon_1<s<0$.
\end{lemma}

\begin{proof}
We write $\hat\varphi_k=\sigma(p)\cos kq$, where $\sigma=s_k/\sinh(k\gamma)$ satisfies $\sigma>0$
on $(-1,0]$, $\sigma(-1)=0$ and $\sigma'(-1)>0$. Then
$(\hat\varphi_k)_q=-k\sigma\sin kq$, $(\hat\varphi_k)_{qp}(q,-1)=-k\sigma'(-1)\sin kq$,
$(\hat\varphi_k)_{qq}=-k^2\sigma\cos kq$ and $(\hat\varphi_k)_{qqp}=-k^2\sigma'\cos kq$. These
satisfy (N1)--(N4), so $\hat\varphi_k\in\mathcal N_k$. By analyticity,
$H_k(s)=H_0+s\big(\hat\varphi_k+r(s)\big)$ with $r(s)\to0$ in $\Xs_k$. By
Lemma~\ref{lem:Nopen}, $\hat\varphi_k+r(s)\in\mathcal N_k$ for small $|s|$. By the invariance
properties of $\mathcal N_k$, $H_k(s)\in\mathcal N_k$ for small $s>0$. For $s<0$ use
$H_k(s)=\tau_{\pi/k}H_k(-s)$.
\end{proof}

\begin{lemma}\label{lem:Nclosed}
Let $(\alpha_n,H_n)\in\mathcal V_k$ be solutions of $\Fop=0$ with $H_n\in\mathcal N_k$, and
suppose that $(\alpha_n,H_n)\to(\alpha,H)$ in $\R\times\Xs_k$ with $(\alpha,H)\in\mathcal V_k$.
Then either $H\in\mathcal N_k$, or $H=\bar x(1+p)$ is a laminar flow of some depth $\bar x>0$.
\end{lemma}

\begin{proof}
Let $w:=H_q\in C^{2,\beta}(\overline\Omega)$. Neither $\Dop$ nor $\Bop_0$ depends explicitly
on $q$. Differentiating $\Fop(\alpha,H)=0$ with respect to $q$ therefore, we obtain
\begin{equation}\label{eq:wq}
\mathcal A_Hw=0\ \text{ in }\Omega,\qquad b_pw_p+b_qw_q+b_0w=0\ \text{ on }p=0,\qquad w=0\ \text{ on }p=-1.
\end{equation}
Moreover $w$ is odd with respect to $q=0$ and to $q=\pi/k$. Passing to the limit, $w\le0$ on
$\overline{\Omega_k^+}$, where $\Omega_k^+:=(0,\pi/k)\times(-1,0)$, and (N2)--(N4) hold with
non-strict inequalities.

The operator $\mathcal A_H$ is uniformly elliptic with bounded coefficients and has no
zeroth-order term. If $w$ vanishes at a point of $\Omega_k^+$, the strong maximum principle of \cite[Theorem~3.5]{GT2001} gives $w\equiv0$ in $\Omega_k^+$, hence everywhere by oddness and
periodicity. Then $H$ is independent of $q$, $\Dop(H)=H_{pp}=0$ and $H(-1)=0$, so
$H=\bar x(1+p)$ with $\bar x=H_p>0$.

Suppose now that $w<0$ in $\Omega_k^+$; we show (N1)--(N4).

(N1) Let $q_0\in(0,\pi/k)$ with $w(q_0,0)=0$. This is a maximum of $w$ on
$\overline{\Omega_k^+}$, so $w_q(q_0,0)=0$, and the Hopf lemma \cite[Lemma~3.4]{GT2001} gives
$w_p(q_0,0)>0$. The boundary condition in \eqref{eq:wq} then reads $b_pw_p=0$, which contradicts
$b_p\ne0$.

(N2) and (N3) on $(-1,0)$: at points of the bottom and of the lines $q=0$, $q=\pi/k$ we have $w=0$
and $w<0$ inside $\Omega_k^+$. The Hopf lemma gives $w_p<0$ on the bottom,
$w_q(0,p)<0$ and $w_q(\pi/k,p)>0$.

(N3) at the top corners. Suppose that $w_q(0,0)=0$. Since $w(0,\cdot)\equiv0$ and $w$ is odd in $q$, $w=w_p=w_q=w_{pp}=w_{qq}=0\qquad\text{at }(0,0).$ The boundary condition in \eqref{eq:wq} holds identically in $q$ as an equation in
$C^{1,\beta}(\T)$. Differentiating it in $q$ at $q=0$ leaves only $b_pw_{qp}=0$, so
$w_{qp}(0,0)=0$ as well. Thus all derivatives of $w$ of order at most two vanish at the corner.
The mixed coefficient $-H_qH_p$ of $\mathcal A_H$ vanishes on $q=0$ and is therefore $O(q)$.
This is the condition on the coefficients required in Serrin's corner lemma
\cite{Serrin1971}, which then implies that $\p_\nu w<0$ or
$\p_\nu^2w<0$ at $(0,0)$ for every direction $\nu$ entering $\Omega_k^+$ non-tangentially, which is
a contradiction. Hence $w_q(0,0)<0$. The corner $(\pi/k,0)$ is treated in the same way.

(N4) At $(0,-1)$, $w$ vanishes on both sides of the corner. Hence all its first and second
derivatives vanish there, except possibly $w_{qp}(0,-1)=H_{qqp}(0,-1)\le0$. If this were zero,
Serrin's lemma would again give a contradiction. The corner $(\pi/k,-1)$ is similar.
\end{proof}

\begin{lemma}\label{lem:laminarlimit}
Let $(\alpha_n,H_n)\in\mathcal V_k$ be solutions with $H_n\in\mathcal N_k$ that converge in
$\R\times\Xs_k$ to $(\bar\alpha,\bar x(1+p))$ for some $\bar x>0$. Then
$\bar\alpha=\alpha_k(\bar x)$.
\end{lemma}

\begin{proof}
Let $w_n:=\p_qH_n\not\equiv0$ and $v_n:=w_n/\|w_n\|_{C^0}$. Each
$v_n$ solves \eqref{eq:wq} with $(\alpha,H)$ replaced by $(\alpha_n,H_n)$. The Schauder estimate
for this problem at the level $C^{2,\beta}$, whose constants are uniform in $n$ because $H_n$
converges in $C^{3,\beta}$, gives $\|v_n\|_{C^{2,\beta}}\le C$. A subsequence converges in
$C^2$ to some $\bar v$ with $\|\bar v\|_{C^0}=1$. The limit is odd and $2\pi/k$-periodic in
$q$, satisfies $\bar v\le0$ in $\Omega_k^+$, and solves the problem linearised at the laminar
flow $\bar H=\bar x(1+p)$.

At $\bar H$ we have $\mathcal A_{\bar H}=\p_p^2+\bar x^2\p_q^2$. With $\cc\gamma^2\e^{2\gamma}=p_0^2$
the boundary coefficients become $b_p=-2\bar\cc/\bar x$ and
$b_0=2(\bar\alpha\e^{\bar x}-\e^{2\bar x}-\bar\cc)$, where $\bar\cc:=p_0^2\e^{-2\bar x}/\bar x^2$.
This is the operator \eqref{eq:L} with $(\gamma,\cc)$ replaced by $(\bar x,\bar\cc)$. As in the
proof of Proposition~\ref{prop:dispersion}, the kernel in the class of odd $2\pi/k$-periodic
functions is spanned by $\sinh(jk\bar x(1+p))\sin(jkq)$ for the single $j\ge1$, if any, with
$\bar\alpha=\alpha_{jk}(\bar x)$. Since $\bar v\le0$ in $\Omega_k^+$ and $\bar v\not\equiv0$, its
coefficient of $\sin kq$, namely $\frac{2k}\pi\int_0^{\pi/k}\bar v(q,p)\sin kq\dd q$, is negative
for some $p$. Hence $j=1$ and $\bar\alpha=\alpha_k(\bar x)$.
\end{proof}

\begin{lemma}\label{lem:noreturn}
Let $\mathcal K=\{(\alpha(s),H(s)):s\in\R\}$ be the curve given by Theorem~\ref{thm:BT}, and
suppose that $H(s)\in\mathcal N_k$ for $0<s<s_1$. Then $(\alpha(s_1),H(s_1))\ne(\alpha_k,H_0)$.
\end{lemma}

\begin{proof}
Let $\Gamma(\sigma)=(\alpha_k(\sigma),H_k(\sigma))$ denote the local branch,
normalised by $\mathsf a_k[H_k(\sigma)-H_0]=\sigma$. Its positive half lies in
$\mathcal N_k$, and its negative half in $\tau_{\pi/k}\mathcal N_k$.
Suppose that $\mathcal K(s_1)=(\alpha_k,H_0)$. For $s<s_1$ sufficiently close
to $s_1$, local uniqueness and $H(s)\in\mathcal N_k$ imply
$\mathcal K(s)\in\Gamma((0,\epsilon))$.

Use the injective analytic uniformisation of the route near $s_1$, with
parameter $t=0$ at $s_1$ and with $t$ increasing with $s$. Its amplitude is
analytic, positive for small $t<0$, and zero at $0$. It is nonzero for small
$t>0$ by the identity theorem. Local uniqueness places the corresponding
points on $\Gamma$. Their amplitude cannot be positive: continuity would then
give the same positive amplitude, hence the same point of $\Gamma$, on both
sides of $0$, which contradicts local injectivity. Thus the route passes from the
positive to the negative half of $\Gamma$ at $s_1$.

Consequently the reversed route $s\mapsto\mathcal K(s_1-s)$ has at $s=0$
the same oriented analytic germ as $s\mapsto\mathcal K(s)$. An increasing
change of parameter near $0$ normalises its amplitude as in \textup{(C2)};
this change can be extended to an increasing homeomorphism of $\R$.
The uniqueness of the oriented maximal route therefore gives an increasing
homeomorphism $h:\R\to\R$ such that $\mathcal K(s_1-s)=\mathcal K(h(s))\qquad(s\in\R).$ The function $s_1-s-h(s)$ is strictly decreasing from $+\infty$ to $-\infty$.
Let $m$ be its zero and set $\mu=s_1-m=h(m)$. For every small $u>0$, $\mathcal K(\mu-u) =\mathcal K\bigl(\mu+h(m+u)-h(m)\bigr).$ The two parameters lie on opposite sides of $\mu$ and both tend to $\mu$.
This contradicts local injectivity of the route at $\mu$.
\end{proof}

\subsection{Conjugate flows}\label{sec:comparison}
Let $(\alpha,H)\in\mathcal V_k$ solve $\Fop(\alpha,H)=0$. We return to the flat variables
$(\Theta,s)$, $s=\ln R$, of Section~\ref{sec:formulation}, and write
\begin{gather*}
\eta(\Theta):=H(\Theta,0)=\ln S(\Theta),\qquad D:=\{(\Theta,s):\ 0<s<\eta(\Theta)\},\\
\eta_{\max}:=\max\eta,\qquad \eta_{\min}:=\min\eta .
\end{gather*}
Since $H_p>0$, the map $(q,p)\mapsto(q,H(q,p))$ is a $C^{3,\beta}$-diffeomorphism of
$\overline\Omega$ onto $\overline D$, with inverse $(\Theta,s)\mapsto(\Theta,-\Psi(\Theta,s))$.
The stream function $\Psi\in C^{3,\beta}(\overline D)$ is harmonic in $D$, equals $1$ on $s=0$
and $0$ on $s=\eta(\Theta)$, and differentiation of $H(\Theta,-\Psi(\Theta,s))=s$ gives
\begin{equation}\label{eq:flat}
\Psi_s=-\frac1{H_p},\qquad \Psi_\Theta=\frac{H_q}{H_p},\qquad
|\widehat\nabla\Psi|^2:=\Psi_s^2+\Psi_\Theta^2=\frac{1+H_q^2}{H_p^2}
\end{equation}
at corresponding points, where $\widehat\nabla:=(\p_\Theta,\p_s)$. The physical gradient is
$|\nabla\Psi|=\e^{-s}|\widehat\nabla\Psi|$. On the free surface $\Psi_\Theta=-\eta'\Psi_s$, because
$\Psi(\Theta,\eta(\Theta))=0$. For $\eps=0$, Bernoulli's law \eqref{eq:bern-phys} becomes
\begin{equation}\label{eq:bern-flat}
|\widehat\nabla\Psi|^2=G_\alpha(\eta)\quad\text{on }s=\eta(\Theta),\qquad
G_\alpha(x):=\frac{\e^{2x}}{p_0^2}\big(2\lambda(\alpha)+\e^{2x}-2\alpha\e^{x}\big).
\end{equation}
The laminar flow $x(1+p)$ of depth $x>0$ has $\Psi=1-s/x$ and $|\widehat\nabla\Psi|=1/x$. It is a
solution if and only if $G_\alpha(x)=x^{-2}$, that is, if and only if $g_\alpha(x)=0$, where
\begin{equation}\label{eq:galpha}
g_\alpha(x):=\frac{p_0^2\e^{-2x}}{x^2}+2\alpha\e^{x}-\e^{2x}-2\lambda(\alpha)
=-p_0^2\e^{-2x}\Big(G_\alpha(x)-\frac1{x^2}\Big).
\end{equation}
A direct computation gives $g_\alpha'(x)=2\e^{x}\big(\alpha-\alpha_0(x)\big)$.
Since $H_0$ is a solution for every $\alpha$, $g_\alpha(\gamma)=0$ for all $\alpha$. The other
zeros of $g_\alpha$ are the depths of the \emph{conjugate flows} of $H_0$, the laminar flows with
the same flux and the same Bernoulli constant \cite{BenjaminLighthill1954}.

\begin{lemma}\label{lem:galpha}
\leavevmode
\begin{enumerate}[label=\textup{(\alph*)}]
\item $\alpha_0$ is strictly convex on $(0,\infty)$, and $\alpha_0(x)\to\infty$ as $x\to0$ and as
$x\to\infty$.
\item For $\alpha>\alpha_0$ the function $g_\alpha$ has exactly three zeros
$x_1(\alpha)<\gamma<x_3(\alpha)$. They are simple; $g_\alpha>0$ on
$(0,x_1(\alpha))\cup(\gamma,x_3(\alpha))$ and $g_\alpha<0$ on
$(x_1(\alpha),\gamma)\cup(x_3(\alpha),\infty)$. Moreover, $x_1$ is strictly decreasing and $x_3$
is strictly increasing on $(\alpha_0,\infty)$.
\end{enumerate}
\end{lemma}

\begin{proof}
Differentiation gives $\alpha_0'(x)=e^x-p_0^2e^{-3x}(3x^{-2}+5x^{-3}+3x^{-4})$. The positive expression subtracted from $e^x$ is strictly decreasing. Thus $\alpha_0$ is strictly convex, with the stated endpoint limits. When $\alpha>\alpha_0(\gamma)$ its sublevel interval $(a,b)$ contains $\gamma$. The sign of $g_\alpha'$ is successively $-,+,-$, with $g_\alpha(0+)=+\infty$, $g_\alpha(\gamma)=0$, $g_\alpha(+\infty)=-\infty$. This gives the three simple roots. At the outer roots $g_\alpha'<0$, while $\partial_\alpha g_\alpha(x)=2(e^x-e^\gamma)$, which gives their monotonicity.
\end{proof}

\begin{lemma}\label{lem:compare}
Let $(\alpha,H)\in\mathcal V_k$ be a solution of $\Fop(\alpha,H)=0$ that is not laminar. Then:
\begin{enumerate}[label=\textup{(\alph*)}]
\item $|\widehat\nabla\Psi|<1/\eta_{\max}$ at every point $(\Theta,\eta_{\max})$ of the free surface;
\item $|\widehat\nabla\Psi|>1/\eta_{\min}$ at every point $(\Theta,\eta_{\min})$ of the free surface;
\item $1/\eta_{\max}<|\Psi_s|<1/\eta_{\min}$ on the bottom $s=0$.
\end{enumerate}
\end{lemma}

\begin{proof}
The harmonic function $U_x=\Psi-(1-s/x)$ vanishes at the bed and has top trace $\eta/x-1$. For $x=\eta_{\max}$ the trace is nonpositive and nonconstant, hence $U_x<0$ inside. Hopf's lemma gives $(U_x)_s>0$ at a crest and $(U_x)_s<0$ at the bed. At a crest $\Psi_\Theta=0$, so these are precisely the required upper surface and lower bed bounds. Taking $x=\eta_{\min}$ reverses the signs and, we obtain the other two bounds.
\end{proof}

\begin{proposition}\label{prop:heights}
Let $(\alpha,H)\in\mathcal V_k$ be a solution of $\Fop(\alpha,H)=0$ that is not laminar. Then
$g_\alpha(\eta_{\max})>0>g_\alpha(\eta_{\min})$ and $\alpha\ne\alpha_0$. If $\alpha>\alpha_0$, then
\begin{equation}\label{eq:heights}
x_1(\alpha)<\eta_{\min}<\gamma<\eta_{\max}<x_3(\alpha).
\end{equation}
\end{proposition}

\begin{proof}
The crest and trough estimates in Lemma~\ref{lem:compare}, inserted in \eqref{eq:galpha}, give
$g_\alpha(\eta_{\min})<0<g_\alpha(\eta_{\max})$. Its limits at $0$ and $\infty$ now force a zero in each of $(0,\eta_{\min})$, $(\eta_{\min},\eta_{\max})$, $(\eta_{\max},\infty)$. Strict convexity of $\alpha_0(x)$ gives at most two zeros of $g_\alpha'(x)=2e^x(\alpha-\alpha_0(x))$. Rolle's theorem then makes the three zeros simple and exhaustive. Since $\gamma$ is one of them, $\alpha\ne\alpha_0(\gamma)$. If $\alpha>\alpha_0(\gamma)$, only the middle zero has positive derivative, so it equals $\gamma$.
\end{proof}

\subsection{A priori bounds}

\begin{lemma}\label{lem:velbound}
Let $M>\alpha_0$. There is $C_1=C_1(M)$ such that every non-laminar solution
$(\alpha,H)\in\mathcal V_k$ of $\Fop(\alpha,H)=0$ with $\alpha_0<\alpha\le M$ satisfies
\begin{gather*}
x_1(M)<\eta_{\min}<\gamma<\eta_{\max}<x_3(M),\qquad |\widehat\nabla\Psi|\le C_1\ \text{ in }\overline D,\\
H_p\ge1/C_1\ \text{ and }\ |H_q|\le C_1H_p\ \text{ in }\overline\Omega .
\end{gather*}
Moreover, $\Psi_s$ attains its maximum over $\overline D$ only on the free surface.
\end{lemma}

\begin{proof}
The bounds on the heights follow from Proposition~\ref{prop:heights} and the monotonicity of
$x_1$ and $x_3$ (Lemma~\ref{lem:galpha}). The function $F:=\Psi_s+\mathrm i\Psi_\Theta$ is
holomorphic in $\zeta=\Theta+\mathrm is$, since the Cauchy--Riemann equations read
$\Psi_{s\Theta}=\Psi_{\Theta s}$ and $\Psi_{ss}=-\Psi_{\Theta\Theta}$. Hence
$|\widehat\nabla\Psi|^2=|F|^2$ is subharmonic, and it attains its maximum over $\overline D$ on
$\p D$. On the bottom, $|\widehat\nabla\Psi|<1/\eta_{\min}<1/x_1(M)$ by
Lemma~\ref{lem:compare}(c). On the free surface, $|\widehat\nabla\Psi|^2=G_\alpha(\eta)$ by
\eqref{eq:bern-flat}, and $G_\alpha(x)$ is bounded for $\alpha\in(\alpha_0,M]$ and
$x\in[x_1(M),x_3(M)]$. This gives $C_1$, and the bounds on $H_p$ and $H_q$ follow from
\eqref{eq:flat}. Finally, $\Psi_s=\operatorname{Re}F$ is harmonic, so its maximum over
$\overline D$ is attained on $\p D$. By Lemma~\ref{lem:compare}, $\Psi_s<-1/\eta_{\max}$ on the
bottom, whereas $\Psi_s>-1/\eta_{\max}$ at a crest.
\end{proof}

Thus, as long as $\alpha$ stays bounded, the relative velocity stays bounded, and nothing
degenerates at the rigid core.

\begin{lemma}\label{lem:apriori}
Let $M>\alpha_0$. There is $C_2=C_2(M)$ such that every non-laminar solution
$(\alpha,H)\in\mathcal V_k$ of $\Fop(\alpha,H)=0$ with $\alpha_0<\alpha\le M$ and
$\max_{\overline\Omega}H_p\le M$ satisfies $\|H\|_{C^{3,\beta}(\overline\Omega)}\le C_2$ and
$\min_{\overline\Omega}H_p\ge1/C_2$.
\end{lemma}

\begin{proof}
In this proof $C$ denotes constants that depend only on $M$ (and on $p_0$, $\gamma$, $k$,
$\beta$). Let $F=\Psi_s+\mathrm i\Psi_\Theta$ be as in the proof of Lemma~\ref{lem:velbound}. By
that lemma and \eqref{eq:flat},
\[
M^{-1}\le-\Psi_s\le|F|\le C_1\qquad\text{in }\overline D .
\]
Since $\operatorname{Re}(-F)>0$, the function $\mathfrak f:=\log(-F)=\mathfrak q+\mathrm i\chi$,
with $\mathfrak q=\ln|\widehat\nabla\Psi|$ and $\chi\in(-\pi/2,\pi/2)$, is holomorphic in $D$ and
of class $C^{2,\beta}$ up to $\p D$. Moreover $|\mathfrak q|\le C$ and $|\chi|\le\chi_0<\pi/2$,
where $\cos\chi_0=1/(MC_1)$. By \eqref{eq:flat},
\begin{equation}\label{eq:qchi}
H_q=-\tan\chi,\qquad H_p=\frac{\e^{-\mathfrak q}}{\cos\chi}
\end{equation}
at corresponding points. In particular $|\eta'|\le L:=\tan\chi_0$.

\emph{Reflection.} On the bottom $\Psi_\Theta=0$, so $\chi=0$ there. By the Schwarz reflection
principle, $\mathfrak f$ extends to a holomorphic function on
$\widetilde D:=\{(\Theta,s):|s|<\eta(\Theta)\}$, with $\mathfrak q$ even and $\chi$ odd in $s$.
By \eqref{eq:bern-flat},
\begin{equation}\label{eq:qbdry}
\mathfrak q=\tfrac12\ln G_\alpha\big(\eta(\Theta)\big)\qquad\text{on }s=\pm\eta(\Theta),
\end{equation}
and $G_\alpha(\eta(\Theta))=|\widehat\nabla\Psi|^2\ge M^{-2}$ there.

\emph{Hölder bounds.} The domain $\widetilde D$ is bounded by the Lipschitz graphs
$s=\pm\eta(\Theta)$, with $\eta\ge x_1(M)$ and $|\eta'|\le L$. It therefore satisfies a uniform
exterior cone condition, and by \eqref{eq:qbdry} the boundary values of $\mathfrak q$ are
Lipschitz in $\Theta$ with constant $C$. The global Hölder estimate \cite[Theorem~8.29]{GT2001},
applied to the harmonic function $\mathfrak q$ on $\widetilde D\cap\{|\Theta|<3\pi/k\}$, together
with periodicity, gives $\|\mathfrak q\|_{C^{\beta_0}(\overline{\widetilde D})}\le C$ for some
$\beta_0\in(0,1)$ that depends only on $M$. The interior gradient estimate for harmonic functions
\cite[Theorem~2.10]{GT2001} then gives
\[
|\widehat\nabla\chi(z)|=|\widehat\nabla\mathfrak q(z)|\le C\,d(z)^{\beta_0-1},\qquad
d(z):=\operatorname{dist}(z,\p\widetilde D).
\]
This yields a Hölder bound for $\chi$. Let $z_j=(\Theta_j,s_j)\in\widetilde D$, $j=1,2$, with
$r:=|z_1-z_2|\le x_1(M)/(4+4L)$, and put $a:=(2+2L)r$ and $\delta(z):=\eta(\Theta)-|s|$, so that
$d\ge\delta/\sqrt{1+L^2}$. Move $z_j$ vertically towards $s=0$ by $\min\{a,|s_j|\}$, to a point
$z_j'=(\Theta_j,s_j')$. Along these vertical segments $\delta$ increases at unit rate, so the
integral of $|\widehat\nabla\chi|$ over each of them is at most $C\int_0^at^{\beta_0-1}\dd t\le
Cr^{\beta_0}$. On the segment $[z_1',z_2']$ we have $|s|\le\max_j|s_j'|$, and
$\eta(\Theta)\ge|s_j|-Lr$ for $j=1,2$. Hence $\delta\ge a-Lr\ge r$ there if $\max_j|s_j'|>0$, and
$\delta\ge x_1(M)\ge r$ otherwise. The length of $[z_1',z_2']$ is at most $r+2a$, so the
integral of $|\widehat\nabla\chi|$ over it is also at most $Cr^{\beta_0}$. Thus
$|\chi(z_1)-\chi(z_2)|\le C|z_1-z_2|^{\beta_0}$; for $|z_1-z_2|>x_1(M)/(4+4L)$ this follows from
$|\chi|<\pi/2$.

\emph{Bootstrap.} On the free surface, \eqref{eq:qchi} gives
$\eta'(\Theta)=-\tan\chi(\Theta,\eta(\Theta))$, so $\|\eta\|_{C^{1,\beta_0}}\le C$. Suppose that
$\|\eta\|_{C^{j,\beta_0}}\le C$ for some $j\ge1$. Then $\widetilde D$ is a domain of class
$C^{j,\beta_0}$, and by \eqref{eq:qbdry} the boundary values of $\mathfrak q$ have
$C^{j,\beta_0}$-norm at most $C$. The boundary Schauder estimates for the Dirichlet problem
(\cite[Corollary~8.36]{GT2001} for $j=1$; \cite{ADN1959} and \cite[Section~6.2]{GT2001} for
$j\ge2$) are local, so they apply to the periodic domain $\widetilde D$ and give
$\|\mathfrak q\|_{C^{j,\beta_0}(\overline{\widetilde D})}\le C$. By the Cauchy--Riemann equations
$\widehat\nabla\chi=(-\mathfrak q_s,\mathfrak q_\Theta)$, so $\|\chi\|_{C^{j,\beta_0}}\le C$, and
then $\|\eta\|_{C^{j+1,\beta_0}}\le C$. After finitely many steps,
$\|\mathfrak q\|_{C^{4,\beta_0}}+\|\chi\|_{C^{4,\beta_0}}\le C$. Since
$\Psi_s=-\e^{\mathfrak q}\cos\chi$ and $\Psi_\Theta=-\e^{\mathfrak q}\sin\chi$, this gives
$\|\Psi\|_{C^{5,\beta_0}(\overline D)}\le C$. Finally, $H$ is obtained by solving
$\Psi(q,H(q,p))=-p$, where $\Psi_s\le-1/M$; hence $\|H\|_{C^4(\overline\Omega)}\le C$, which bounds
$\|H\|_{C^{3,\beta}(\overline\Omega)}$. The lower bound for $H_p$ is part of
Lemma~\ref{lem:velbound}.
\end{proof}

\subsection{A bound for \texorpdfstring{$\alpha$}{alpha}}
We now show that $\alpha$ is bounded on the set of all waves with monotone profiles. Write
$Y:=\e^\gamma$ and $y:=\e^x$. By \eqref{eq:lambda},
\begin{equation}\label{eq:Bquad}
G_\alpha(x)=\frac{y^2}{p_0^2}B_\alpha(y),\qquad
B_\alpha(y):=y^2-2\alpha y+\cc+2\alpha Y-Y^2=\cc+(Y-y)(2\alpha-y-Y),
\end{equation}
and a short computation gives
\begin{equation}\label{eq:gstar}
G_\alpha'(x)=\frac{2y^2}{p_0^2}\big(2y^2-3\alpha y+\cc+2\alpha Y-Y^2\big),\qquad
g_*(x):=-\frac{G_\alpha'(x)}{2G_\alpha(x)}=\frac{y(\alpha-y)}{B_\alpha(y)}-1 .
\end{equation}
By \eqref{eq:bern-flat}, $B_\alpha(\e^{\eta})=p_0^2\e^{-2\eta}|\widehat\nabla\Psi|^2>0$ on the free
surface. If $(\alpha-Y)^2\ge\cc$, the quadratic $B_\alpha$ has the zeros
$z_\mp(\alpha):=\alpha\mp\big((\alpha-Y)^2-\cc\big)^{1/2}$, which are the levels $\e^x$ at which
the Bernoulli head $G_\alpha$ vanishes.

\begin{lemma}\label{lem:surfspeed}
Let $(\alpha,H)\in\mathcal V_k$ be a non-laminar solution of $\Fop(\alpha,H)=0$ with
$\alpha>\alpha_0$ and $(\alpha-\e^\gamma)^2\ge\cc$. Then
$\e^{\eta}<z_-(\alpha)\le\e^\gamma+\cc/(\alpha-\e^\gamma)$ on the free surface, and
$|\widehat\nabla\Psi|^2>G_\alpha(\eta_{\max})$ at every point of the free surface where
$\eta<\eta_{\max}$. In particular, the relative speed on the free surface is smallest exactly at the
crests.
\end{lemma}

\begin{proof}
$B_\alpha$ is negative between $z_-$ and $z_+$, and $Y<z_-$ because $B_\alpha(Y)=\cc>0$ and
$Y<\alpha$. Since $\eta$ is continuous, attains the value $\eta_{\min}<\gamma$, and
$B_\alpha(\e^\eta)>0$, the free surface stays below $\ln z_-$. The bound for $z_-$ follows from
$(a^2-\cc)^{1/2}\ge a-\cc/a$ with $a:=\alpha-Y\ge\cc^{1/2}$.

Let $\sigma_\alpha(y):=2y^2-3\alpha y+\cc+2\alpha Y-Y^2$, so that $G_\alpha'$ has the sign of
$\sigma_\alpha(\e^x)$ by \eqref{eq:gstar}. The quadratic $\sigma_\alpha$ is convex,
$\sigma_\alpha(Y)=\cc+Y^2-\alpha Y<0$ because $\alpha>\alpha_0>Y+\cc/Y$, and
$\sigma_\alpha(z_-)=z_-(z_--\alpha)\le0$ because
$\sigma_\alpha-B_\alpha=y^2-\alpha y$ and $B_\alpha(z_-)=0$.
The endpoint can be a double root; convexity gives strict negativity on
$[Y,z_-)$. Hence $G_\alpha$ is strictly
decreasing on $[\gamma,\ln z_-]$, while on $[\eta_{\min},\gamma]$ it is first increasing (possibly on
an empty interval) and then decreasing. Thus $G_\alpha(x)>\min\{G_\alpha(\eta_{\min}),G_\alpha(\eta_{\max})\}$
for $\eta_{\min}<x<\eta_{\max}$. Finally, by \eqref{eq:bern-flat} and Lemma~\ref{lem:compare},
$G_\alpha(\eta_{\max})<\eta_{\max}^{-2}<\eta_{\min}^{-2}<G_\alpha(\eta_{\min})$.
\end{proof}

The condition $(\alpha-\e^\gamma)^2\ge\cc$ holds for all $\alpha>\alpha_0$ if
$r\ge\gamma^2/(1+\gamma)^2$, since $\alpha_0-\e^\gamma=\cc\,\e^{-\gamma}(1+1/\gamma)$.

\begin{lemma}\label{lem:flowforce}
Let $(\alpha,H)\in\mathcal V_k$ solve $\Fop(\alpha,H)=0$. Then
\[
\mathcal S(\Theta):=\int_0^{\eta(\Theta)}\big(\Psi_s^2-\Psi_\Theta^2\big)(\Theta,s)\dd s
+\int_\gamma^{\eta(\Theta)}G_\alpha(x)\dd x
\]
is independent of $\Theta$. If $H\in\mathcal N_k$ and $\alpha>\alpha_0$, then
\begin{equation}\label{eq:ffineq}
\frac1{\eta_{\max}}+\int_{\eta_{\min}}^{\gamma}G_\alpha(x)\dd x\le\eta_{\min}\max_{[\eta_{\min},\eta_{\max}]}G_\alpha .
\end{equation}
\end{lemma}

\begin{proof}
Since $\Psi_{\Theta\Theta}=-\Psi_{ss}$, we have
$\p_\Theta(\Psi_s^2-\Psi_\Theta^2)=2\p_s(\Psi_\Theta\Psi_s)$, and $\Psi_\Theta=0$ on $s=0$.
Using $\Psi_\Theta=-\eta'\Psi_s$ on the free surface and \eqref{eq:bern-flat}, we get
\[
\mathcal S'=\eta'\big(\Psi_s^2-\Psi_\Theta^2\big)+2\Psi_\Theta\Psi_s+\eta'G_\alpha(\eta)
=-\eta'\,\Psi_s^2(1+\eta'^2)+\eta'G_\alpha(\eta)=0,
\]
where all functions are evaluated at $s=\eta(\Theta)$.

Let now $H\in\mathcal N_k$, so that the crest is at $\Theta=0$, the trough at $\Theta=\pi/k$, and
$\Psi_\Theta=0$ on both lines by symmetry. Every level in $[\eta_{\min},\eta_{\max}]$ is attained on
the free surface, where $G_\alpha(\eta)=|\widehat\nabla\Psi|^2>0$; hence $G_\alpha>0$ on this
interval. On the crest line $\int_0^{\eta_{\max}}|\Psi_s|\dd s=1$, so the Cauchy--Schwarz inequality
gives $\mathcal S(0)\ge1/\eta_{\max}$. On the trough line,
$\Psi_s^2\le\max_{\overline D}|\widehat\nabla\Psi|^2$. This maximum is attained on $\p D$, because
$|\widehat\nabla\Psi|^2$ is subharmonic (see the proof of Lemma~\ref{lem:velbound}). By
Lemma~\ref{lem:compare}, it is not attained on the bottom, where $|\widehat\nabla\Psi|<1/\eta_{\min}$,
while $|\widehat\nabla\Psi|>1/\eta_{\min}$ at the trough. Hence it is at most
$\max_{[\eta_{\min},\eta_{\max}]}G_\alpha$, and
$\mathcal S(\pi/k)\le\eta_{\min}\max G_\alpha-\int_{\eta_{\min}}^\gamma G_\alpha$. Since
$\mathcal S(0)=\mathcal S(\pi/k)$, \eqref{eq:ffineq} follows.
\end{proof}

\begin{lemma}\label{lem:robin}
Let $(\alpha,H)\in\mathcal V_k$ solve $\Fop(\alpha,H)=0$ with $H\in\mathcal N_k$. Let
$D_+:=D\cap\{0<\Theta<\pi/k\}$, let $\Gamma_+$ be the part of the free surface with
$0<\Theta<\pi/k$, let $\chi:=-\arctan\eta'\in[0,\pi/2)$ be the angle between the free surface and
the circles $R=\mathrm{const}$, and let $\varkappa$ be the curvature of $\Gamma_+$, traversed from the
trough to the crest, so that $\varkappa>0$ where $\eta''<0$. Then $v:=-\Psi_\Theta$ is positive in
$D_+\cup\Gamma_+$, vanishes on the rest of $\p D_+$, and
\begin{equation}\label{eq:robin}
\p_nv=\big(g_*(\eta)\cos\chi-\varkappa\big)\,v\qquad\text{on }\Gamma_+,
\end{equation}
where $n$ is the outer normal. Consequently, if $\varphi(\Theta,s)=f(\Theta)\zeta(s)$ with
$f\in C^1([0,\pi/k])$, $f(0)=f(\pi/k)=0$, and $\zeta$ Lipschitz with $\zeta(0)=0$, then
$\varphi_\Gamma(\Theta):=\varphi(\Theta,\eta(\Theta))$ satisfies
\begin{equation}\label{eq:eigineq}
\int_0^{\pi/k}g_*(\eta(\Theta))\,\varphi_\Gamma(\Theta)^2\dd\Theta
\le\int_{D_+}|\widehat\nabla\varphi|^2\dd\Theta\dd s+\frac\pi2\int_0^{\pi/k}\big|(\varphi_\Gamma^2)'(\Theta)\big|\dd\Theta .
\end{equation}
\end{lemma}

\begin{proof}
By (N1) and \eqref{eq:flat}, $\Psi_\Theta=H_q/H_p<0$ in $D_+\cup\Gamma_+$. Moreover, $\Psi_\Theta$
vanishes on $s=0$ and, by symmetry, on $\Theta=0$ and $\Theta=\pi/k$.

\emph{The Robin condition.} Write $Q:=|\widehat\nabla\Psi|$. Differentiating $\Psi(\Theta,\eta(\Theta))=0$, we obtain
$v=\eta'\Psi_s=-\eta'Q\cos\chi$ on $\Gamma_+$, by \eqref{eq:qchi}. Differentiating \eqref{eq:bern-flat}
gives $\p_\Theta(Q^2)+\eta'\p_s(Q^2)=G_\alpha'(\eta)\eta'$. Here
$\p_\Theta(Q^2)=-2\widehat\nabla\Psi\cdot\widehat\nabla v=2Q\,\p_nv$, because $\widehat\nabla\Psi=-Q\,n$
on $\Gamma_+$. Eliminating $\eta'$, we get
\begin{equation}\label{eq:robin1}
\p_nv=\frac{\p_s(Q^2)-G_\alpha'(\eta)}{2Q^2\cos\chi}\,v .
\end{equation}
To compute $\p_s(Q^2)$ we use the holomorphic function $\mathfrak f=\mathfrak q+\mathrm i\chi=\log(-\Psi_s-\mathrm i\Psi_\Theta)$
from the proof of Lemma~\ref{lem:apriori}; it is defined for every solution. On $\Gamma_+$ the unit
tangent in the direction of the flow is $t=(-\cos\chi,\sin\chi)$, and $\nu:=(-\sin\chi,-\cos\chi)$ is
the inner normal; $\nu$ is obtained from $t$ by a rotation through $+\pi/2$. In the frame $(t,\nu)$ the
Cauchy--Riemann equations give $\p_\nu\mathfrak q=-\p_t\chi=\varkappa$, while \eqref{eq:bern-flat} gives
$\p_t\mathfrak q=\p_t\tfrac12\ln G_\alpha(\eta)=-g_*(\eta)\sin\chi$. Since $(0,1)=\sin\chi\,t-\cos\chi\,\nu$,
\[
\p_s\mathfrak q=-g_*(\eta)\sin^2\chi-\varkappa\cos\chi,\qquad
\p_s(Q^2)=2Q^2\p_s\mathfrak q=G_\alpha'(\eta)\sin^2\chi-2Q^2\varkappa\cos\chi .
\]
Inserting this into \eqref{eq:robin1}, we obtain \eqref{eq:robin}.

\emph{The inequality.} By (N2)--(N4), $v$ vanishes to first order on $s=0$ and on $\Theta=0,\pi/k$,
and to second order at the bottom corners; hence
$v\ge c\,\min\{\Theta,\pi/k-\Theta\}\min\{s,1\}$ in $D_+$ for some $c>0$, whereas
$|\varphi|\le C\min\{\Theta,\pi/k-\Theta\}\min\{s,1\}$. For small $\epsilon>0$ let
$D_\epsilon:=D_+\cap\{\epsilon<\Theta<\pi/k-\epsilon,\ s>\epsilon\}$ and write $\varphi=v\psi$ there.
Since $v$ is harmonic, Green's formula gives
\[
\int_{D_\epsilon}|\widehat\nabla\varphi|^2=\int_{D_\epsilon}v^2|\widehat\nabla\psi|^2
+\int_{\p D_\epsilon}\frac{\varphi^2}{v}\,\p_nv
\ge\int_{\p D_\epsilon}\frac{\varphi^2}{v}\,\p_nv .
\]
The integrals over the segments $\Theta=\epsilon$, $\Theta=\pi/k-\epsilon$ and $s=\epsilon$ tend to
zero as $\epsilon\to0$, because $\varphi^2/v=O(\epsilon)$ there. On the remaining part of
$\p D_\epsilon$, which lies on $\Gamma_+$, we use \eqref{eq:robin}. Letting $\epsilon\to0$ we obtain
$\int_{D_+}|\widehat\nabla\varphi|^2\ge\int_{\Gamma_+}(g_*\cos\chi-\varkappa)\varphi^2\dd\ell$. Since
$\dd\ell=\dd\Theta/\cos\chi$, the first term on the right equals
$\int_0^{\pi/k}g_*(\eta)\varphi_\Gamma^2\dd\Theta$. Since $\varkappa=-\p_t\chi$ and $\chi$ vanishes at
the crest and at the trough, an integration by parts gives
$\big|\int_{\Gamma_+}\varkappa\varphi^2\dd\ell\big|=\big|\int_{\Gamma_+}\chi\,\p_t(\varphi^2)\dd\ell\big|
\le\frac\pi2\int_0^{\pi/k}|(\varphi_\Gamma^2)'|\dd\Theta$.
\end{proof}

In other words, the Robin problem \eqref{eq:robin} in $D_+$, with a Dirichlet condition on the rest
of $\p D_+$, has the positive solution $v$, so its quadratic form is non-negative. For the laminar
flow of depth $x$ the coefficient in \eqref{eq:robin} is the constant $-\tfrac12x^2G_\alpha'(x)$,
and the problem has the positive solution $\sinh(ks)\sin(k\Theta)$ exactly when
$\alpha=\alpha_k(x)$; this is the dispersion relation of Proposition~\ref{prop:dispersion}.

\begin{proposition}\label{prop:alphabound}
There is $A_*=A_*(p_0,\gamma,k)$ such that $\alpha\le A_*$ for every solution
$(\alpha,H)\in\mathcal V_k$ of $\Fop(\alpha,H)=0$ with $H\in\mathcal N_k$ and $\alpha>\alpha_0$.
\end{proposition}

\begin{proof}
Put
\begin{gather*}
t_0:=\frac{3(Y-1)^2(Y+2)}{16Y^3},\qquad J(t):=\frac1{p_0^2}\int_t^\gamma\e^{2x}(Y-\e^x)\dd x,\\
a_1:=\frac{8Y^3}{27p_0^2},\qquad a_0:=\frac{\cc\,(Y^2+(Y+1)^2)}{p_0^2},
\end{gather*}
$K:=16p_0^2/t_0^2$ and $E_0:=(\pi/k)\big(k^2(\gamma+1)+1/t_0\big)$, and note that
$J(0)=(Y-1)^2(Y+2)/(6p_0^2)>0$. We may assume
$\alpha>\alpha_1:=\max\{Y+\cc^{1/2},\,Y+\cc,\,2Y+2,\,2K/Y\}$. Then Lemma~\ref{lem:surfspeed} applies
and gives $\e^\eta<z_-\le Y+1$ on the free surface; in particular $\eta_{\max}\le\gamma+1$, and
$\alpha-y\ge\alpha/2$ for $y=\e^\eta$. For $1\le y\le Y$ we have
\begin{equation}\label{eq:Bbounds}
\alpha(Y-y)\le B_\alpha(y)\le\cc+2\alpha(Y-y),
\end{equation}
the lower bound because $2\alpha-y-Y\ge\alpha$. For $Y\le y\le z_-$ we have $0<B_\alpha(y)\le\cc$.

\emph{Deep troughs: $\eta_{\min}\le t_0$.} By \eqref{eq:Bquad} and \eqref{eq:Bbounds},
$\int_{\eta_{\min}}^\gamma G_\alpha\ge\alpha J(\eta_{\min})$ and
$\max_{[\eta_{\min},\eta_{\max}]}G_\alpha\le a_1\alpha+a_0$; here we used
$y^2(Y-y)\le4Y^3/27$ for $1\le y\le Y$. Then \eqref{eq:ffineq} gives
$\alpha\big(J(\eta_{\min})-a_1\eta_{\min}\big)\le a_0t_0$. Moreover
$J(t)-a_1t\ge J(0)-4Y^3t/(9p_0^2)\ge J(0)/2$ for $t\le t_0$. Hence $\alpha\le2a_0t_0/J(0)$.

\emph{Shallow troughs: $\eta_{\min}>t_0$.} The divergence theorem in $D_+$, together with
$\Psi_\Theta=0$ on $\Theta=0,\pi/k$ and Lemma~\ref{lem:compare}(c), gives
\[
\int_{\Gamma_+}|\widehat\nabla\Psi|\dd\ell=\int_0^{\pi/k}|\Psi_s(\Theta,0)|\dd\Theta<\frac{\pi}{k\eta_{\min}}<\frac\pi{kt_0}.
\]
Where $\eta<\gamma$, \eqref{eq:bern-flat} and \eqref{eq:Bbounds} give
$|\widehat\nabla\Psi|^2\ge\alpha(Y-\e^\eta)/p_0^2$. Since $\dd\ell\ge\dd\Theta$, it follows that
$\int_{\{\eta<\gamma\}}(Y-\e^{\eta})^{1/2}\dd\Theta\le\pi|p_0|/(kt_0\alpha^{1/2})$. Let
$\mathsf S:=\{\Theta\in[0,\pi/k]:\e^{\eta(\Theta)}\ge Y-K/\alpha\}$. On its complement
$(Y-\e^\eta)^{1/2}>(K/\alpha)^{1/2}$, so the complement has length at most
$\pi|p_0|/(kt_0K^{1/2})=\pi/(4k)$. On $\mathsf S$ we have $y\ge Y/2$, $\alpha-y\ge\alpha/2$ and
$B_\alpha(y)\le\cc+2K$, so that \eqref{eq:gstar} gives $g_*\ge Y\alpha/(4(\cc+2K))-1$; on the whole free
surface $g_*>-1$. Now apply \eqref{eq:eigineq} with $f(\Theta)=\sin k\Theta$ and
$\zeta(s)=\min\{s/t_0,1\}$. Then $\varphi_\Gamma=\sin k\Theta$,
$\int_0^{\pi/k}|(\varphi_\Gamma^2)'|\dd\Theta=2$ and $\int_{D_+}|\widehat\nabla\varphi|^2\le E_0$, so that
\[
\frac{\pi}{4k}\Big(\frac{Y\alpha}{4(\cc+2K)}-2\Big)\le\int_0^{\pi/k}g_*(\eta)\sin^2(k\Theta)\dd\Theta\le E_0+\pi .
\]
This bounds $\alpha$, and we may take
$A_*:=\max\big\{\alpha_1,\ 2a_0t_0/J(0),\ \tfrac{4(\cc+2K)}{Y}\big(2+\tfrac{4k}{\pi}(E_0+\pi)\big)\big\}$.
\end{proof}

\subsection{The global theorem}\label{sec:globalthm}

\begin{theorem}\label{thm:global}
Let $k\ge1$, $p_0<0$ and $\gamma>0$, and let
$\mathcal K_k=\{(\alpha(s),H(s)):s\in\R\}$ be the curve given by Theorem~\ref{thm:BT} at
$(\alpha_k,H_0)$. Then $\mathcal K_k\subset\mathcal V_k$ consists of solutions of
$\Fop(\alpha,H)=0$, and the following hold.
\begin{enumerate}[label=\textup{(\roman*)}]
\item \emph{(Local structure and symmetry.)} Near $s=0$, $\mathcal K_k$ coincides with the curve
of Theorem~\ref{thm:local}. It has a local real-analytic reparametrisation around each of its
points, and it is invariant under $(\alpha,H)\mapsto(\alpha,\tau_{\pi/k}H)$.
\item \emph{(Monotone profiles.)} $H(s)\in\mathcal N_k$ for every $s>0$. Thus, for $s>0$, the
free surface and all streamlines decrease strictly from a crest at $\Theta=0$ to a trough at
$\Theta=\pi/k$; in particular $H(s)$ is not laminar.
\item \emph{(Bounds.)} For every $s>0$, $\alpha_0<\alpha(s)\le A_*$ and
\[
x_1(\alpha(s))<\eta_{\min}(s)<\gamma<\eta_{\max}(s)<x_3(\alpha(s)),
\]
where $\eta(s)=\ln S=H(s)|_{p=0}$, $x_1(\alpha)<\gamma<x_3(\alpha)$ are the depths of the
conjugate flows (Lemma~\ref{lem:galpha}), and $A_*$ is the constant of
Proposition~\ref{prop:alphabound}. Moreover, the relative velocity is bounded along
$\{(\alpha(s),H(s)):s>0\}$.
\item \emph{(Global alternative.)} As $s\to\infty$,
\[
\max_{\overline\Omega}\p_pH(s)\longrightarrow\infty .
\]
\item \emph{(No loop.)} $\mathcal K_k$ is not a closed loop.
\end{enumerate}
\end{theorem}

\begin{proof}
Lemmas~\ref{lem:fredholm-global} and~\ref{lem:proper} verify (H3) and (H4), and (H1), (H2) were
checked above, so Theorem~\ref{thm:BT} applies.

(i) The first two statements follow from (C2)--(C4) and Theorem~\ref{thm:local}(ii). The
reflected curve $s\mapsto(\alpha(-s),\tau_{\pi/k}H(-s))$ also has (C1)--(C5). Indeed,
$\tau_{\pi/k}$ preserves $\mathcal V_k$, the solution set and the sets $Q_j$, and
$\tau_{\pi/k}(H_0-s\hat\varphi_k)=H_0+s\hat\varphi_k$. By uniqueness it is a reparametrisation of
$\mathcal K_k$, so $\mathcal K_k$ is invariant as a set.

(ii), (iii) By (C2), (C3) and Lemma~\ref{lem:Nlocal}, $H(s)\in\mathcal N_k$ for small $s>0$. Let
$s_*\in(0,\infty]$ be the supremum of all $s_1>0$ such that $H(s)\in\mathcal N_k$ for $0<s<s_1$.
For $0<s<s_*$ the function $H(s)$ depends on $q$, so it is not laminar, and
Proposition~\ref{prop:heights} gives $\alpha(s)\ne\alpha_0$. Since $\alpha$ is continuous and
$\alpha(s)\to\alpha_k>\alpha_0$ as $s\downarrow0$, we get $\alpha(s)>\alpha_0$ on $(0,s_*)$, and
then \eqref{eq:heights} holds for $H(s)$, $0<s<s_*$.

Suppose that $s_*<\infty$. By Lemma~\ref{lem:Nopen}, $H(s_*)\notin\mathcal N_k$; otherwise
$H(s)\in\mathcal N_k$ for $s$ slightly larger than $s_*$. By Lemma~\ref{lem:Nclosed}, $H(s_*)$ is
a laminar flow $\bar x(1+p)$. As $s\uparrow s_*$, both $\eta_{\min}(s)$ and $\eta_{\max}(s)$ tend
to $\bar x$, so \eqref{eq:heights} forces $\bar x=\gamma$, that is, $H(s_*)=H_0$. By
Lemma~\ref{lem:laminarlimit}, $\alpha(s_*)=\alpha_k(\gamma)=\alpha_k$. This contradicts
Lemma~\ref{lem:noreturn}. Hence we obtain $s_*=\infty$, which proves (ii). Now (iii) follows from
Proposition~\ref{prop:alphabound} and Lemma~\ref{lem:velbound} with $M:=A_*$.

(v) If (C5)(II) held with period $T$, then by (C2), (C3) and Lemma~\ref{lem:Nlocal},
$H(T-u)=H(-u)\in\tau_{\pi/k}\mathcal N_k$ for small $u>0$. This contradicts (ii).

(iv) By (v), alternative (C5)(I) holds. Suppose that (iv) fails. Then there are $M\ge A_*$ and
$s_n\to\infty$ with $\max\p_pH(s_n)\le M$. By (iii) and Lemma~\ref{lem:apriori},
$\|H(s_n)\|_{C^{3,\beta}}\le C_2(M)$ and $\min\p_pH(s_n)\ge1/C_2(M)$, and
$0<\alpha_0<\alpha(s_n)\le M$. Hence $(\alpha(s_n),H(s_n))\in Q_j$ for all $n$ if
$j\ge\max\{M,C_2(M)\}$, which contradicts (C5)(I).
\end{proof}

The analytic curve in Theorem~\ref{thm:global} can be a proper subset of the connected component if a secondary branch attaches.
\begin{remark}\label{rem:globalinterp}
At corresponding points $|V-R|=|p_0|e^{-s}/H_p$. By Lemma~\ref{lem:velbound}, the minimum of $-\Psi_s$ is attained on the free surface, where $-\Psi_s=|\widehat\nabla\Psi|\cos\chi$ and $\chi=-\arctan\eta'$. Thus surface speed can approach zero or the surface can approach a radial tangent. Lemma~\ref{lem:surfspeed} identifies the speed minimum with the crest when $(\alpha-e^\gamma)^2\ge\cc$.
\end{remark}
\section{Crest stagnation and limiting waves}\label{sec:stagnation}
\subsection{A criterion for stagnation at the crest}
By Remark~\ref{rem:globalinterp}, Theorem~\ref{thm:global}(iv) leaves two possibilities along a
sequence $s_n\to\infty$: the relative speed tends to zero somewhere on the free surface, or the free
surface becomes radial somewhere. We now show that, under two conditions on the Bernoulli function
$G_\alpha$, the maximum of $\p_pH$ is attained at the crest along the whole curve, so that (iv)
means stagnation at the crest. The proof propagates along $\mathcal K_k$ one more nodal property,
the monotonicity of $H_p$ in $q$.

Let $(\alpha,H)\in\mathcal V_k$ solve $\Fop(\alpha,H)=0$ with $H\in\mathcal N_k$ and
$\alpha>\alpha_0$, and use the notation of Lemma~\ref{lem:robin}. Let $\overline\Gamma_+$ be the closed arc of the free surface
from the crest $(0,\eta_{\max})$ to the trough $(\pi/k,\eta_{\min})$. At its end points $\chi=0$,
and by (N3) the curvature there satisfies $\varkappa_c:=\varkappa(0)>0$ at the crest and
$\varkappa_t:=\varkappa(\pi/k)<0$ at the trough. Let
\begin{equation}\label{eq:rho}
\varrho:=g_*(\eta)\cos\chi-\varkappa\qquad\text{on }\overline\Gamma_+
\end{equation}
be the coefficient of the Robin condition \eqref{eq:robin}, and let
\[
\mathfrak a:=\p_\Theta\mathfrak q=\p_\Theta\ln|\widehat\nabla\Psi|,\qquad
\mathfrak b:=\p_\Theta\chi,
\]
so that $\mathfrak f'=\mathfrak a+\mathrm i\mathfrak b$, where the prime denotes the complex
derivative with respect to $\zeta=\Theta+\mathrm is$. The bootstrap in the proof of
Lemma~\ref{lem:apriori} gives $\mathfrak f\in C^{4,\beta_0}(\overline D)$ and $\eta\in C^{5,\beta_0}$,
so $\mathfrak a,\mathfrak b\in C^{3,\beta_0}(\overline D)$ and $\varkappa\in C^{3,\beta_0}$. We consider
the condition
\begin{enumerate}[label=\textup{(N5)}]
\item $H_{qp}<0$ on $(0,\pi/k)\times(-1,0]$, and $H_{qqp}(0,p)<0<H_{qqp}(\pi/k,p)$ for
$p\in(-1,0]$.
\end{enumerate}
Together with (N2) and (N4), it gives $H_{qp}$ on the closed half-cell the sign structure that
(N1)--(N4) give $H_q$. Since $-\Psi_s=1/H_p$ by \eqref{eq:flat}, (N5) says that $-\Psi_s$ decreases
strictly along every streamline from the trough line to the crest line; equivalently, every stream
tube of the flow in the flat variables widens from the trough to the crest. By
Theorem~\ref{thm:global}(ii),(iii), every $H(s)$ with $s>0$ satisfies the standing assumptions.

\begin{lemma}\label{lem:frakA}
Let $(\alpha,H)\in\mathcal V_k$ solve $\Fop(\alpha,H)=0$ with $H\in\mathcal N_k$ and $\alpha>\alpha_0$.
\begin{enumerate}[label=\textup{(\alph*)}]
\item At corresponding points, $\mathfrak a=-H_{qp}/\big(H_p(1+H_q^2)\big)$.
\item $\mathfrak a$ is harmonic in $D$, odd with respect to $\Theta=0$ and to $\Theta=\pi/k$, and
positive on the bottom $(0,\pi/k)\times\{0\}$.
\item On $\overline\Gamma_+$, $\mathfrak a=\varrho\sin\chi$ and
$\mathfrak b=g_*(\eta)\sin^2\chi+\varkappa\cos\chi$. On $\Gamma_+$, $\p_n\mathfrak a=\p_t\mathfrak b$.
\item At the crest and at the trough, $\mathfrak a_\Theta=\varkappa\varrho$.
\item $\varrho>0$ on $\overline\Gamma_+$ if and only if $H$ satisfies \textup{(N5)}. In this case
$H_p$ attains its maximum over $\overline\Omega$ exactly at the crests $(2\pi j/k,0)$,
$j\in\mathbb Z$, and
\[
\max_{\overline\Omega}H_p=\frac1{|\widehat\nabla\Psi|(0,\eta_{\max})}=G_\alpha(\eta_{\max})^{-1/2}.
\]
\end{enumerate}
\end{lemma}

\begin{proof}
(a) By \eqref{eq:flat}, $\mathfrak q=\frac12\ln(1+H_q^2)-\ln H_p$ at corresponding points. Since
$p=-\Psi(\Theta,s)$ and $\Psi_\Theta=H_q/H_p$, we have $\p_\Theta=\p_q-(H_q/H_p)\p_p$, and therefore
\[
\mathfrak a=\frac{H_qH_{qq}}{1+H_q^2}-\frac{H_{qp}}{H_p}
-\frac{H_q}{H_p}\Big(\frac{H_qH_{qp}}{1+H_q^2}-\frac{H_{pp}}{H_p}\Big).
\]
We eliminate $H_{qq}$ with $\Dop(H)=0$, that is, $H_p^2H_{qq}=2H_qH_pH_{qp}-(1+H_q^2)H_{pp}$. The
terms with $H_{pp}$ cancel, and (a) follows.

(b) $\mathfrak a=\operatorname{Re}\mathfrak f'$ is harmonic. Since $H$ is even in $q$ and
$2\pi/k$-periodic, $|\widehat\nabla\Psi|$ is even with respect to $\Theta=0$ and $\Theta=\pi/k$, so
$\mathfrak a$ is odd. On the bottom $p=-1$ and $H_q=0$, so $\mathfrak a=-H_{qp}/H_p>0$ by (N2).

(c) In the proof of Lemma~\ref{lem:robin} we found $\p_t\mathfrak q=-g_*(\eta)\sin\chi$ and
$\p_\nu\mathfrak q=-\p_t\chi=\varkappa$; the Cauchy--Riemann equations also give
$\p_\nu\chi=\p_t\mathfrak q$. Since $(1,0)=-\cos\chi\,t-\sin\chi\,\nu$,
\[
\mathfrak a=-\cos\chi\,\p_t\mathfrak q-\sin\chi\,\p_\nu\mathfrak q=\varrho\sin\chi,\qquad
\mathfrak b=-\cos\chi\,\p_t\chi-\sin\chi\,\p_\nu\chi=\varkappa\cos\chi+g_*(\eta)\sin^2\chi
\]
on $\Gamma_+$, and by continuity on $\overline\Gamma_+$. In the positively oriented frame $(t,\nu)$
the Cauchy--Riemann equations for $\mathfrak a+\mathrm i\mathfrak b$ give
$\p_\nu\mathfrak a=-\p_t\mathfrak b$; since $n=-\nu$, $\p_n\mathfrak a=\p_t\mathfrak b$.

(d) At the crest and at the trough $\chi=0$ and $t=(-1,0)$. By (c) and $\p_t\chi=-\varkappa$,
$\p_t\mathfrak a=\sin\chi\,\p_t\varrho+\varrho\cos\chi\,\p_t\chi=-\varkappa\varrho$ there, so
$\mathfrak a_\Theta=-\p_t\mathfrak a=\varkappa\varrho$.

(e) Let $\varrho>0$ on $\overline\Gamma_+$. Then $\mathfrak a>0$ on $\Gamma_+$ by (c), since
$\sin\chi>0$ there; $\mathfrak a=0$ on the lines $\Theta=0$ and $\Theta=\pi/k$; and $\mathfrak a>0$ on
the bottom by (b). The strong maximum principle gives $\mathfrak a>0$ in $D_+$, and the Hopf lemma
gives $\mathfrak a_\Theta>0$ on $\{0\}\times(0,\eta_{\max})$ and $\mathfrak a_\Theta<0$ on
$\{\pi/k\}\times(0,\eta_{\min})$. By (d), $\mathfrak a_\Theta=\varkappa_c\varrho>0$ at the crest and
$\mathfrak a_\Theta=\varkappa_t\varrho<0$ at the trough. On the line $q=0$ we have $H_q=H_{qp}=0$, so
(a) gives $\p_q[\mathfrak a(q,H(q,p))]=-H_{qqp}/H_p$ there, and the left side equals
$\mathfrak a_\Theta$ because $H_q=0$; the same holds on $q=\pi/k$. By (a), these sign conditions on
$\mathfrak a$ and $\mathfrak a_\Theta$ are (N5). Conversely, (N5) and (a) give $\mathfrak a>0$ on
$\Gamma_+$, hence $\varrho>0$ there, and $\varrho=\mathfrak a_\Theta/\varkappa>0$ at the crest and at
the trough by (d), (N3) and (N5) at $p=0$.

If (N5) holds, then $H_p(\cdot,p)$ is strictly decreasing on $[0,\pi/k]$ for every
$p\in(-1,0]$, and by (N2) also for $p=-1$. By evenness and periodicity, $\max_qH_p(q,p)$ is
attained exactly at $q\in(2\pi/k)\mathbb Z$. By Lemma~\ref{lem:velbound}, $H_p=-1/\Psi_s$ attains its
maximum over $\overline\Omega$ only on $p=0$. Hence the maximum of $H_p$ is attained exactly at the
crests, where $\Psi_\Theta=0$ and, by \eqref{eq:flat} and \eqref{eq:bern-flat},
$H_p=1/|\widehat\nabla\Psi|=G_\alpha(\eta_{\max})^{-1/2}$.
\end{proof}

The flow force of Lemma~\ref{lem:flowforce}, evaluated on the two symmetry lines, gives an identity
for the azimuthal velocities there, and under (N5) in the non-strict form it controls the trough.

\begin{lemma}\label{lem:ffsym}
Let $(\alpha,H)\in\mathcal V_k$ solve $\Fop(\alpha,H)=0$ with $H\in\mathcal N_k$ and $\alpha>\alpha_0$.
\begin{enumerate}[label=\textup{(\alph*)}]
\item
$\displaystyle\int_{-1}^0\Big(\frac1{H_p(\pi/k,p)}-\frac1{H_p(0,p)}\Big)\dd p=\int_{\eta_{\min}}^{\eta_{\max}}G_\alpha(x)\dd x .$
\item If $\mathfrak a\ge0$ in $D_+$, that is, $H_{qp}\le0$ on $(0,\pi/k)\times(-1,0]$, then
$-\Psi_s$ is non-decreasing along the trough line from the bottom to the trough, and
\begin{equation}\label{eq:FF3}
\int_{\eta_{\min}}^{\eta_{\max}}G_\alpha(x)\dd x<G_\alpha(\eta_{\min})^{1/2}\,\frac{\eta_{\max}-\eta_{\min}}{\eta_{\max}} .
\end{equation}
\item \emph{(Trough identity.)} Let $t:=\eta_{\min}$ and
$K_t(y):=\dfrac{\pi}{t^2}\,\dfrac{\e^{-\pi|y|/t}}{(1-\e^{-\pi|y|/t})^2}$ for $y\ne0$. Then
\begin{equation}\label{eq:troughid}
G_\alpha(t)^{1/2}\,t-1=t\int_{\R}K_t(y)\,\Psi\Big(\frac\pi k+y,t\Big)\dd y
\ \ge\ t\,K_t\Big(\frac\pi k\Big)\int_0^{2\pi/k}\Psi(\Theta,t)\dd\Theta .
\end{equation}
\end{enumerate}
\end{lemma}

\begin{proof}
(a) On the lines $\Theta=0$ and $\Theta=\pi/k$ we have $\Psi_\Theta=0$, and the substitution
$s=H(\Theta,p)$ together with \eqref{eq:flat} gives $\int_0^{\eta(\Theta)}\Psi_s^2\dd s=\int_{-1}^0H_p(\Theta,p)^{-1}\dd p$.
Hence (a) is the identity $\mathcal S(0)=\mathcal S(\pi/k)$ of Lemma~\ref{lem:flowforce}.

(b) Write $u:=-\Psi_s=1/H_p$. On the trough line, $\mathfrak a=0$ and $\mathfrak a\ge0$ for
$\Theta<\pi/k$, so $\mathfrak a_\Theta\le0$ there. By the Cauchy--Riemann equations,
$\mathfrak b_s=\mathfrak a_\Theta\le0$ and $\p_s\mathfrak q=-\mathfrak b$. Since $\chi=0$ on the bottom,
$\mathfrak b=0$ at $(\pi/k,0)$, so $\mathfrak b\le0$ and $\p_s\mathfrak q\ge0$ on the trough line. There
$\Psi_\Theta=0$ and $u=\e^{\mathfrak q}$, so $u$ is non-decreasing in $s$, from
$a:=u(\pi/k,0)>1/\eta_{\max}$ (Lemma~\ref{lem:compare}(c)) to $b:=G_\alpha(\eta_{\min})^{1/2}$ at the
trough, and $b\,\eta_{\min}>1$ by Lemma~\ref{lem:compare}(b). Since $(u-a)(b-u)\ge0$ and
$\int_0^{\eta_{\min}}u\dd s=1$ (the flux),
\[
\int_0^{\eta_{\min}}u^2\dd s\le(a+b)\int_0^{\eta_{\min}}u\dd s-ab\,\eta_{\min}=b-a(b\eta_{\min}-1)
<b-\frac{b\eta_{\min}-1}{\eta_{\max}} .
\]
On the crest line the Cauchy--Schwarz inequality gives $\int_0^{\eta_{\max}}u^2\dd s\ge1/\eta_{\max}$,
again because the flux is $1$. Inserting both bounds into (a), we obtain \eqref{eq:FF3}.

(c) The strip $\Sigma:=\R\times(0,t)$ lies in $D$, because $\eta\ge t$. In $\Sigma$, $\Psi$ is
harmonic and $2\pi/k$-periodic, $\Psi=1$ on $s=0$, and $f:=\Psi(\cdot,t)$ is of class $C^{2,\beta}$,
non-negative, with $f(\pi/k)=0$, since the line $s=t$ touches the free surface at the troughs.
Expanding $f=\sum_jf_j\e^{\mathrm ijk\Theta}$, we get
$\Psi(\Theta,s)=1+(f_0-1)s/t+\sum_{j\ne0}f_j\e^{\mathrm ijk\Theta}\sinh(|j|ks)/\sinh(|j|kt)$, and hence $-\Psi_s(\Theta,t)=\frac{1-f_0}t-\sum_{j\ne0}m(jk)\,f_j\e^{\mathrm ijk\Theta},\qquad m(\xi):=\xi\coth(\xi t).$ The partial fractions $\coth x=1/x+\sum_{n\ge1}2x/(x^2+n^2\pi^2)$ and
$\int_\R(1-\cos\xi y)\frac a2\e^{-a|y|}\dd y=\xi^2/(\xi^2+a^2)$, with $a=n\pi/t$, give
$m(\xi)=1/t+\int_\R(1-\cos\xi y)K_t(y)\dd y$, where $0\le m(\xi)-1/t\le|\xi|$. Since
$\sum_j|j||f_j|<\infty$, we may interchange sum and integral and obtain $-\Psi_s(\Theta,t)=\frac{1-f(\Theta)}t+\int_\R K_t(y)\Big(\frac{f(\Theta+y)+f(\Theta-y)}2-f(\Theta)\Big)\dd y .$ At the trough, $f=0$, $f$ is even about $\pi/k$, and $-\Psi_s=|\widehat\nabla\Psi|=G_\alpha(t)^{1/2}$ by
\eqref{eq:bern-flat}; this is the identity in \eqref{eq:troughid}. The inequality follows because $K_t>0$ is decreasing in $|y|$ and
$f\ge0$ is $2\pi/k$-periodic.
\end{proof}

By Lemma~\ref{lem:compare}(b) and Proposition~\ref{prop:heights}, the inequality
$G_\alpha(\eta_{\min})^{1/2}\eta_{\min}>1$ says that $\eta_{\min}>x_1(\alpha)$. Part (c) is a quantitative
form of it: a trough close to the conjugate level forces $\Psi$ to be small on the trough level, in
the average weighted by $K_t$, which decays like $\e^{-\pi|y|/t}$.

For $\alpha>\alpha_0$ put $\hat x_3(\alpha):=\ln z_-(\alpha)$ if $(\alpha-\e^\gamma)^2\ge\cc$ and
$\hat x_3(\alpha):=x_3(\alpha)$ otherwise, and
\begin{equation}\label{eq:Falpha}
F_\alpha(t,e):=\int_t^eG_\alpha(x)\dd x-G_\alpha(t)^{1/2}\Big(1-\frac te\Big)\qquad(0<t<e).
\end{equation}

\begin{theorem}\label{thm:stagnation}
Let $\mathcal K_k=\{(\alpha(s),H(s)):s\in\R\}$ be the curve of Theorem~\ref{thm:global}, and let
$0<s_1\le\infty$. Suppose that for every $s\in(0,s_1)$, with $\alpha=\alpha(s)$ and with $\eta_{\min}$,
$\eta_{\max}$ the trough and crest levels of $H(s)$,
\begin{enumerate}[label=\textup{(G\arabic*${}'$)}]
\item[\textup{(G1)}] $G_\alpha$ is log-concave on $[\eta_{\min},\eta_{\max}]$, that is, $g_*'\ge0$ there;
\item[\textup{(G2${}'$)}] $G_\alpha'(\eta_{\min})\le0$, or $F_\alpha\big(\eta_{\min},\hat x_3(\alpha)\big)\ge0$.
\end{enumerate}
Then for every $s\in(0,s_1)$ the free surface of $H(s)$ satisfies
\begin{equation}\label{eq:curvbound}
\varkappa<g_*(\eta)\cos\chi\qquad\text{on }\overline\Gamma_+,
\end{equation}
$H(s)$ satisfies \textup{(N5)}, the azimuthal relative speed $|V-R|$ decreases strictly along every
streamline from the trough to the crest, and
\[
\max_{\overline\Omega}\p_pH(s)=\p_pH(s)(0,0)=G_{\alpha(s)}\big(\eta_{\max}(s)\big)^{-1/2}.
\]
If $s_1=\infty$, then $G_{\alpha(s)}(\eta_{\max}(s))\to0$ as $s\to\infty$, that is, the relative
velocity at the crest tends to zero. Moreover, every limit point $(\bar\alpha,\bar\eta)$ of
$(\alpha(s),\eta_{\max}(s))$ as $s\to\infty$ satisfies
\[
(\bar\alpha-\e^\gamma)^2\ge\cc\qquad\text{and}\qquad\e^{\bar\eta}=z_-(\bar\alpha),
\]
so the crest approaches the level at which the Bernoulli head $G_{\bar\alpha}$ vanishes.
\end{theorem}

\begin{proof}
We write $\varrho_s$ for the function \eqref{eq:rho} of $H(s)$, regarded as a function of
$\Theta\in[0,\pi/k]$. On the free surface $\chi=-\arctan\eta'$ and
$\varkappa=-\eta''(1+\eta'^2)^{-3/2}$, where $\eta=H(s)|_{p=0}$ depends continuously on $s$ in
$C^{3,\beta}$, and $B_{\alpha(s)}(\e^\eta)>0$. Hence $\varrho_s\in C([0,\pi/k])$ depends
continuously on $s>0$. As $s\downarrow0$, $H(s)\to H_0$ and $\alpha(s)\to\alpha_k$, so $\varrho_s$
converges uniformly to the constant $\frac{\e^\gamma(\alpha_k-\e^\gamma)}{\cc}-1>\frac1\gamma,$ which is $g_*(\gamma)$ at $\alpha=\alpha_k$ by \eqref{eq:gstar} and $B_\alpha(\e^\gamma)=\cc$; the
inequality holds because $\alpha_k>\alpha_0=\e^\gamma+\cc\,\e^{-\gamma}(1+1/\gamma)$.

Let $s_\sharp$ be the supremum of all $s'\in(0,s_1]$ such that $\varrho_s>0$ on $[0,\pi/k]$ for
$0<s<s'$. Then $s_\sharp>0$. Suppose that $s_\sharp<s_1$, and let $H:=H(s_\sharp)$ and
$\varrho:=\varrho_{s_\sharp}$. By continuity $\varrho\ge0$ on $\overline\Gamma_+$, and $\varrho$
vanishes at some point $P_0\in\overline\Gamma_+$, for otherwise $\varrho_s>0$ for all $s$ near
$s_\sharp$. By Theorem~\ref{thm:global}(ii), $H\in\mathcal N_k$. By Lemma~\ref{lem:frakA}(b),(c),
$\mathfrak a\ge0$ on $\p D_+$ and $\mathfrak a>0$ on the bottom, so $\mathfrak a>0$ in $D_+$ by the
strong maximum principle.

\emph{Case 1: $P_0\in\Gamma_+$.} Then $\mathfrak a(P_0)=0=\min_{\overline{D_+}}\mathfrak a$, and the
Hopf lemma gives $\p_n\mathfrak a(P_0)<0$. On the other hand, $\varrho$ attains its minimum $0$ over
$\Gamma_+$ at $P_0$, so $\p_t\varrho(P_0)=0$. Along $\Gamma_+$ we have $\p_t\chi=-\varkappa$ and
$\p_t\eta=\sin\chi$, hence $\p_t[g_*(\eta)]=g_*'(\eta)\sin\chi$, and
\[
\p_t\varrho=g_*'\sin\chi\cos\chi+g_*\varkappa\sin\chi-\p_t\varkappa,\qquad
\p_t\mathfrak b=g_*'\sin^3\chi-2g_*\varkappa\sin\chi\cos\chi+\p_t\varkappa\cos\chi+\varkappa^2\sin\chi .
\]
At $P_0$ we have $\varkappa=g_*\cos\chi$, so $\p_t\varrho=0$ gives
$\p_t\varkappa=(g_*'+g_*^2)\sin\chi\cos\chi$. Inserting this into $\p_t\mathfrak b$, we obtain from
Lemma~\ref{lem:frakA}(c) $\p_n\mathfrak a(P_0)=\p_t\mathfrak b(P_0)=g_*'(\eta)\sin\chi\ge0$ by (G1), a contradiction.

\emph{Case 2: $P_0$ is the crest.} Then $g_*(\eta_{\max})=\varkappa_c$, and $\mathfrak a_\Theta=0$ at
the crest by Lemma~\ref{lem:frakA}(d). Since $\mathfrak a$ is odd in $\Theta$, also
$\mathfrak a=\mathfrak a_s=\mathfrak a_{ss}=\mathfrak a_{\Theta\Theta}=0$ on $\Theta=0$. Near the
crest $D$ is symmetric with respect to $\Theta=0$, its boundary is of class $C^{5,\beta_0}$, and the
line $\Theta=0$ contains the normal to $\p D$ at the crest. Serrin's corner lemma
\cite{Serrin1971}, applied to the harmonic function $\mathfrak a>0$ in $D_+$, therefore gives
$\p_d\mathfrak a>0$ or $\p_d^2\mathfrak a>0$ at the crest, for the direction $d=(1,-1)/\sqrt2$, which
enters $D_+$ non-tangentially. Here $\p_d\mathfrak a=0$ and $\p_d^2\mathfrak a=-\mathfrak a_{\Theta s}$,
so $\mathfrak a_{\Theta s}<0$. The Cauchy--Riemann equations $\mathfrak a_s=-\mathfrak b_\Theta$ and
$\mathfrak b_s=\mathfrak a_\Theta$ give $\mathfrak b_{\Theta\Theta}=-\mathfrak a_{\Theta s}>0$ and
$\mathfrak b_s=0$ at the crest. The free surface, parametrised by arclength in the direction $t$, has
acceleration $\varkappa\nu$, and at the crest $t=(-1,0)$ and $\nu=(0,-1)$. Hence the second derivative
of $\mathfrak b$ along the free surface at the crest is $\p_t^2\mathfrak b=\mathfrak b_{\Theta\Theta}-\varkappa_c\mathfrak b_s=\mathfrak b_{\Theta\Theta}>0 .$ On the other hand, differentiate $\mathfrak b=g_*\sin^2\chi+\varkappa\cos\chi$ twice along the free
surface and use $\chi=0$, $\p_t\chi=-\varkappa_c$, $g_*=\varkappa_c$ and, by symmetry, $\p_t\varkappa=0$ at
the crest. This gives $\p_t^2\mathfrak b=2g_*\varkappa_c^2+\varkappa_2-\varkappa_c^3=\varkappa_2+\varkappa_c^3$, where
$\varkappa_2:=\p_t^2\varkappa$ at the crest. Thus $\varkappa_2+\varkappa_c^3>0$. Finally, $\varrho\ge0$ on
$\overline\Gamma_+$, and $\varrho$ is even with respect to the crest and vanishes there, so
$\p_t^2\varrho\ge0$ at the crest. Differentiating the formula for $\p_t\varrho$ in Case~1 once more,
we get at the crest
\[
0\le\p_t^2\varrho=-g_*'(\eta_{\max})\varkappa_c-g_*\varkappa_c^2-\varkappa_2
=-g_*'(\eta_{\max})\varkappa_c-(\varkappa_2+\varkappa_c^3)<-g_*'(\eta_{\max})\varkappa_c .
\]
Since $\varkappa_c>0$, this contradicts (G1).

\emph{Case 3: $P_0$ is the trough.} Then $0=\varrho=g_*(\eta_{\min})-\varkappa_t$, so
$g_*(\eta_{\min})=\varkappa_t<0$, that is, $G_\alpha'(\eta_{\min})>0$. By (G2${}'$),
$F_\alpha(\eta_{\min},\hat x_3)\ge0$. On the other hand, $\mathfrak a>0$ in $D_+$, so
Lemma~\ref{lem:ffsym}(b) gives $F_\alpha(\eta_{\min},\eta_{\max})<0$. Moreover
$\p_eF_\alpha(t,e)=G_\alpha(e)-G_\alpha(t)^{1/2}t/e^2<0$ for $t=\eta_{\min}$ and $\gamma<e<x_3$, because
$G_\alpha(e)<e^{-2}$ there (Lemma~\ref{lem:galpha}(b) and \eqref{eq:galpha}) and
$G_\alpha(\eta_{\min})^{1/2}\eta_{\min}>1$ (Lemma~\ref{lem:compare}(b)). Finally
$\eta_{\max}<\hat x_3\le x_3$: by Proposition~\ref{prop:heights} and Lemma~\ref{lem:surfspeed},
$\eta_{\max}<x_3$ and, if $(\alpha-\e^\gamma)^2\ge\cc$, $\eta_{\max}<\ln z_-$; and $\ln z_-<x_3$
because $\ln z_->\gamma$ and $g_\alpha(\ln z_-)>0$ by \eqref{eq:galpha}, as $G_\alpha(\ln z_-)=0$.
Hence $F_\alpha(\eta_{\min},\hat x_3)<F_\alpha(\eta_{\min},\eta_{\max})<0$, a contradiction.

Hence we obtain $s_\sharp=s_1$. For $s\in(0,s_1)$, Lemma~\ref{lem:frakA}(e) now gives (N5),
\eqref{eq:curvbound} and the formula for $\max\p_pH(s)$. By Remark~\ref{rem:globalinterp},
$|V-R|=|p_0|\e^{-s}/H_p$ at the point $(\Theta,s)=(q,H(q,p))$ of the streamline $\Psi=-p$. Along
it, $s=H(q,p)$ is non-increasing in $q\in[0,\pi/k]$ by (N1), and $1/H_p$ is strictly increasing by
(N2) and (N5); hence $|V-R|$ increases strictly from the crest line to the trough line.

Let $s_1=\infty$. By Theorem~\ref{thm:global}(iv), $G_{\alpha(s)}(\eta_{\max}(s))\to0$; at the crest
the relative velocity is azimuthal, of size $|p_0|\e^{-\eta_{\max}}G_\alpha(\eta_{\max})^{1/2}$. By
Theorem~\ref{thm:global}(iii) and Lemma~\ref{lem:velbound}, $\alpha_0<\alpha(s)\le A_*$ and
$\gamma<\eta_{\max}(s)<x_3(A_*)$, so limit points exist. By \eqref{eq:Bquad},
$B_{\alpha(s)}(\e^{\eta_{\max}(s)})=p_0^2\e^{-2\eta_{\max}(s)}G_{\alpha(s)}(\eta_{\max}(s))\to0$, hence
$B_{\bar\alpha}(\e^{\bar\eta})=0$. Since $B_{\bar\alpha}(y)=(y-\bar\alpha)^2-\big((\bar\alpha-\e^\gamma)^2-\cc\big)$,
this gives $(\bar\alpha-\e^\gamma)^2\ge\cc$ and $\e^{\bar\eta}\in\{z_-(\bar\alpha),z_+(\bar\alpha)\}$. If
$(\bar\alpha-\e^\gamma)^2=\cc$, then $z_-(\bar\alpha)=z_+(\bar\alpha)$. Otherwise
$(\alpha(s_n)-\e^\gamma)^2>\cc$ for large $n$ along a sequence $s_n\to\infty$ realising the limit
point, and Lemma~\ref{lem:surfspeed} gives
$\e^{\eta_{\max}(s_n)}<z_-(\alpha(s_n))\to z_-(\bar\alpha)$, so $\e^{\bar\eta}\le z_-(\bar\alpha)<z_+(\bar\alpha)$.
In both cases $\e^{\bar\eta}=z_-(\bar\alpha)$.
\end{proof}

The simplest sufficient condition for (G2${}'$) is
\begin{enumerate}[label=\textup{(G2)}]
\item $G_\alpha'(\eta_{\min})\le0$, that is, $g_*(\eta_{\min})\ge0$.
\end{enumerate}
The hypotheses (G1) and (G2) hold automatically in an explicit range of parameters. Recall that
$r=\cc\,\e^{-2\gamma}$, and let $\sigma_\alpha(y)=2y^2-3\alpha y+\cc+2\alpha\e^\gamma-\e^{2\gamma}$ be
the quadratic of the proof of Lemma~\ref{lem:surfspeed}, so that $G_\alpha'$ has the sign of
$\sigma_\alpha(\e^x)$.

\begin{corollary}\label{cor:stagnation}
\leavevmode
\begin{enumerate}[label=\textup{(\alph*)}]
\item Let $\alpha>\alpha_0$ and $(\alpha-\e^\gamma)^2\ge\cc$. Then $(\ln G_\alpha)''\le0$ on each interval where
$G_\alpha>0$. In particular, \textup{(G1)} holds for every solution with this value of $\alpha$.
\item Let $\alpha>\alpha_0$ and $\sigma_\alpha(1)\le0$. Then $G_\alpha'<0$ on $(0,\gamma]$. In
particular, \textup{(G2)} holds for every non-laminar solution with this value of $\alpha$.
\item Let $Y:=\e^\gamma\le3/2$, put $\beta(\gamma):=1+(2-3/Y)(1+1/\gamma)$, and suppose that
\begin{equation}\label{eq:region}
r\ge\frac{\gamma^2}{(1+\gamma)^2}\qquad\text{and}\qquad(Y-1)(Y-2)+rY^2\beta(\gamma)\le0 .
\end{equation}
Then \textup{(G1)} and \textup{(G2)} hold for all $s>0$, and the conclusions of
Theorem~\ref{thm:stagnation} hold with $s_1=\infty$. The function $\beta$ is increasing on
$(0,\ln\frac32]$ and vanishes at $\gamma_b\approx0.2971$. Hence, for $0<\gamma\le\gamma_b$,
\eqref{eq:region} reduces to $r\ge\gamma^2/(1+\gamma)^2$. For $\gamma_b<\gamma\le\ln\frac32$ it
defines a non-empty interval of $r$, for instance $0.0533\le r\le3.46$ for $\gamma=0.3$ and
$0.0832\le r\le1/9$ for $\gamma=\ln\frac32$.
\end{enumerate}
\end{corollary}

\begin{proof}
(a) Put $y=\e^x$ and $\delta_\alpha:=(\alpha-\e^\gamma)^2-\cc$, so that
$B_\alpha(y)=(y-\alpha)^2-\delta_\alpha$ by \eqref{eq:Bquad}. Since
$\ln G_\alpha=2x+\ln B_\alpha(\e^x)-\ln p_0^2$, a short computation gives, with $w:=\alpha-y$,
\begin{equation}\label{eq:lnG2}
(\ln G_\alpha)''(x)=\frac{y}{B_\alpha(y)^2}\Big(B_\alpha(y)(4y-2\alpha)-4y(y-\alpha)^2\Big)
=-\frac{2y}{B_\alpha(y)^2}\big(\alpha w^2-2\delta_\alpha w+\alpha\delta_\alpha\big).
\end{equation}
If $\delta_\alpha\ge0$, the quadratic in $w$ has the discriminant
$4\delta_\alpha(\delta_\alpha-\alpha^2)\le0$, because $\alpha^2-\delta_\alpha=\cc+2\alpha\e^\gamma-\e^{2\gamma}>0$
for $\alpha>\alpha_0>\e^\gamma$. Hence $(\ln G_\alpha)''\le0$ wherever $B_\alpha>0$. Every level in
$[\eta_{\min},\eta_{\max}]$ is attained on the free surface, where
$G_\alpha=|\widehat\nabla\Psi|^2>0$.

(b) By the proof of Lemma~\ref{lem:surfspeed}, $\sigma_\alpha$ is convex and $\sigma_\alpha(\e^\gamma)<0$.
If $\sigma_\alpha(1)\le0$, then $\sigma_\alpha<0$ on $(1,\e^\gamma]$, so $G_\alpha'<0$ on $(0,\gamma]$
by \eqref{eq:gstar}. By Proposition~\ref{prop:heights}, $0<x_1(\alpha)<\eta_{\min}<\gamma$.

(c) If $Y\le3/2$, then $\sigma_\alpha(1)=2+\cc-Y^2-(3-2Y)\alpha$ is non-increasing in $\alpha$.
Inserting $\alpha_0=Y+(\cc/Y)(1+1/\gamma)$ and $\cc=rY^2$, we find
$\sigma_{\alpha_0}(1)=(Y-1)(Y-2)+rY^2\beta(\gamma)$. Hence the second condition in \eqref{eq:region}
gives $\sigma_\alpha(1)\le0$ for all $\alpha\ge\alpha_0$, and the first gives
$(\alpha-\e^\gamma)^2\ge\cc$ for all $\alpha>\alpha_0$ (see the remark after
Lemma~\ref{lem:surfspeed}). Since $\alpha(s)>\alpha_0$ for $s>0$ by Theorem~\ref{thm:global}(iii), (a)
and (b) give (G1) and (G2) for all $s>0$. Finally,
$\beta'(\gamma)=3\e^{-\gamma}(1+1/\gamma)+(3\e^{-\gamma}-2)/\gamma^2>0$ for $\e^\gamma\le3/2$, and
$(Y-1)(Y-2)<0$ for $0<\gamma<\ln2$.
\end{proof}

Condition (G2${}'$) can also be checked a priori, without (G2). For $\alpha>\alpha_0$ let
\[
Z(\alpha):=\Big\{t\in\big(x_1(\alpha),\gamma\big):\ G_\alpha'(t)>0\ \text{ and }\ F_\alpha\big(t,\hat x_3(\alpha)\big)<0\Big\}.
\]
Since $x_1(\alpha)<\eta_{\min}<\gamma$ by Proposition~\ref{prop:heights}, (G2${}'$) holds for every
wave with a given $\alpha$ if $Z(\alpha)=\emptyset$. The set $Z(\alpha)$ depends only on $\gamma$, $r$
and $\alpha$, not on $k$.

\begin{theorem}\label{thm:stagCAP}
Let $r_{\min}(\gamma)=\gamma^2/(1+\gamma)^2$ and let $\mathcal P$ be the union
of the following parameter regions:
\begin{center}
\begin{tabular}{ll}\toprule
Logarithmic depth & Flux parameter \\\midrule
$0<\gamma\le263/200$ & $r\ge r_{\min}(\gamma)$ \\
$263/200<\gamma\le27/20$ & $r_{\min}(\gamma)\le r\le2/3$ \\
$263/200<\gamma\le11/8$ & $r_{\min}(\gamma)\le r\le1/2$ \\
$263/200<\gamma\le7/5$ & $r_{\min}(\gamma)\le r\le2/5$ \\\bottomrule
\end{tabular}
\end{center}
If $(\gamma,r)\in\mathcal P$, then $Z(\alpha)=\emptyset$ for every $\alpha>\alpha_0$. Consequently, for every $k\ge1$, (G1) and (G2${}'$) hold along $\mathcal K_k$ for all
$s>0$, and the conclusions of Theorem~\ref{thm:stagnation} hold with $s_1=\infty$: the free surfaces
satisfy \eqref{eq:curvbound}, $\max\p_pH(s)$ is attained at the crest, the relative velocity at the
crest tends to zero, and the crest approaches the level at which the Bernoulli head vanishes.
\end{theorem}

\begin{proof}
Since $r\ge r_{\min}(\gamma)$, one has $(\alpha-e^\gamma)^2\ge\cc$
for every $\alpha>\alpha_0$, and (G1) follows from
Corollary~\ref{cor:stagnation}(a). The direct maximum-point argument of
Theorem~\ref{thm:direct-trough} proves $Z(\alpha)=\varnothing$
on the whole stated parameter set, hence (G2${}'$). Theorem~\ref{thm:stagnation}
therefore applies along the entire branch.
\end{proof}

The trough enters the proof at the first possible failure of the curvature inequality. At that point the local maximum principle alone gives no contradiction. The integral estimate in Lemma~\ref{lem:ffsym} leads to the explicit exclusion proved next.
\subsection{Direct exclusion of low trough levels}\label{app:capZ}\label{sec:trough-direct}
Throughout this argument $\eps=\kappa=0$. To exclude low trough levels, we first locate the maximum of the Bernoulli head. We then integrate the corresponding normalized profile explicitly.
We write $\gamma_*=263/200$ and let $\mathcal P$ be the parameter set in
Theorem~\ref{thm:stagCAP}.

\begin{theorem}\label{thm:direct-trough}
For every $(\gamma,r)\in\mathcal P$ and every $\alpha>\alpha_0$,
the set $Z(\alpha)$ defined in Section~\ref{sec:stagnation} is empty.
\end{theorem}

Distinguish the radial variable \(v\) from the logarithmic level \(x\).  Write
\[
 Y=e^\gamma,\qquad
 \alpha=Y+rY\beta,\qquad
 \beta_0=1+\frac1\gamma,\qquad
 \rho=\frac1r,\qquad \mathsf c=rY^2.
\]
Since \(\alpha_0=Y+rY\beta_0\), our assumptions imply
\begin{equation}
\label{eq:td-normalised-domain}
 \beta>\beta_0>1,
 \qquad 0<\rho\le\beta_0^2<\beta^2.
\end{equation}
It is useful to allow \(\rho=0\) as a limiting case.  Set
\begin{align}
 \mathfrak G(v)
   &:=\frac{v^2}{\gamma^2}
       \bigl[1+2\beta(1-v)+\rho(1-v)^2\bigr],
       &G(x)&:=\mathfrak G(e^{x-\gamma}),
       \label{eq:td-head}\\
 s(v)&:=1+\beta(2-3v)+\rho(2v-1)(v-1),
       &c&:=1+\frac1{\beta+\sqrt{\beta^2-\rho}},
       \label{eq:td-s-c}\\
 e&:=\gamma+\log c,
       &F(t,e)&:=\int_t^eG(x)\,dx
             -\sqrt{G(t)}\left(1-\frac te\right).
       \label{eq:td-F}
\end{align}
Thus \(c=e^{-\gamma}z_-(\alpha)\) and \(e=\widehat x_3(\alpha)\).  In particular \(e>\gamma\).  The required assertion follows if
\begin{equation}
\label{eq:td-candidate}
 F(t,e)>0
 \quad\text{whenever}\quad
 0<t<\gamma,\qquad G'(t)>0,\qquad G(t)t^2>1.
\end{equation}
Indeed, by Lemma~\ref{lem:galpha}, these are precisely the levels at which trough exclusion must be checked.

\begin{lemma}
\label{lem:td-maximum}
Assume that \(\beta>1\), \(\rho\ge0\), and \(\beta^2>\rho\).  There is a unique
\(m\in(1/2,1)\) at which \(\mathfrak G\) is maximal on \((0,c)\).  Put
\[
 t_m:=\gamma+\log m,\qquad M:=\mathfrak G(m).
\]
If a level \(t\) satisfies the three conditions in
\eqref{eq:td-candidate}, then
\begin{equation}
\label{eq:td-maximum-consequences}
 0<t<t_m<\gamma<e,\qquad Mt_m^2>1,
 \qquad F(t,e)>F(t_m,e).
\end{equation}
\end{lemma}

\begin{proof}
Direct differentiation gives $G'(x)=\frac{2v^2}{\gamma^2}s(v),\qquad v=e^{x-\gamma}.$ Moreover, $s(0)=1+2\beta+\rho>0,\qquad s(1/2)=1+\frac\beta2>0,\qquad s(1)=1-\beta<0.$ For \(\rho>0\), the quadratic \(s\) is strictly convex: it has one root
\(m\in(1/2,1)\) and another root greater than \(1\).
The first zero \(c\) of \(\mathfrak G\) is simple and
\(\mathfrak G'(c)<0\), so the second root of \(s\) lies beyond \(c\).
For \(\rho=0\), the same assertions follow from the decreasing affine function
\(s\).  Hence \(G\) increases up to \(t_m\) and decreases from \(t_m\) to \(e\).

A candidate level lies below \(t_m\).  On \([t,t_m]\), the function
\(x\mapsto G(x)x^2\) is strictly increasing, because \(x>0\), \(G>0\), and
\(G'\ge0\).  Thus \(Mt_m^2>1\).  Differentiating \(F\) with respect to its
first argument gives
\begin{equation}
\label{eq:td-F-monotone}
 \partial_tF(t,e)
 =-\sqrt{G(t)}\left(\sqrt{G(t)}-\frac1e\right)
  -\frac{G'(t)}{2\sqrt{G(t)}}\left(1-\frac te\right)<0
\end{equation}
throughout \([t,t_m]\): feasibility gives
\(\sqrt{G(t)}>1/t>1/e\), and this inequality persists as the first argument
increases toward \(t_m\).  Integration proves the last assertion.
\end{proof}

\begin{lemma}
\label{lem:td-area}
In the setting of Lemma~\ref{lem:td-maximum}, let \(a=c/m\).  Then
\(3/2\le a<2\), and
\begin{align}
 I(a)&:=\frac1M\int_{t_m}^eG(x)\,dx
       =\frac{(a-1)(-2a^2+3a+7)}{12},
       \label{eq:td-I}\\
 0&<I(a)<\log a,\qquad e=t_m+\log a,
       \label{eq:td-I-range}\\
 F(t_m,e)&=
 \sqrt M\left(\sqrt M\,I(a)-\frac{\log a}{t_m+\log a}\right).
       \label{eq:td-F-maximum}
\end{align}
The identities extend continuously to \(a=2\), the double-root limit.
\end{lemma}

\begin{proof}
For \(\rho>0\), factor $\mathfrak G(v)=\frac{\rho}{\gamma^2}v^2(c-v)(d-v), \qquad d>c>1.$ With \(a=c/m\) and \(b=d/m\), the stationary-point equation is $\frac1{a-1}+\frac1{b-1}=2.$ Since \(b>a>1\), it follows that \(3/2<a<2\), and $b=\frac{3a-4}{2a-3}.$ Putting \(u=v/m\), and using \(dx=dv/v\), gives $I(a)=\int_1^a \frac{u(a-u)(b-u)}{(a-1)(b-1)}\,du =\frac{(a-1)(-2a^2+3a+7)}{12}.$ For \(\rho=0\), one has \(a=3/2\) and the same calculation is the limit
\(b\to\infty\).  Since \(0<G(x)<M\) for \(t_m<x<e\), its integral divided by
\(M\) lies strictly between \(0\) and \(e-t_m=\log a\).
Equation~\eqref{eq:td-F-maximum} follows directly from \eqref{eq:td-F}.
\end{proof}

We define the area ratio
\begin{equation}
\label{eq:td-L}
 L(a):=\frac{I(a)\log a}{\log a-I(a)}.
\end{equation}
For a candidate, \(\sqrt M\,t_m>1\).  Therefore
\begin{equation}
\label{eq:td-final-reduction}
 \begin{split}
 F(t_m,e)
 &>\sqrt M\left(\frac{I(a)}{t_m}
                   -\frac{\log a}{t_m+\log a}\right)\\
 &=\frac{\sqrt M\,[I(a)\log a-t_m(\log a-I(a))]}
          {t_m(t_m+\log a)}.
 \end{split}
\end{equation}
Consequently, it suffices to show \(t_m<L(a)\).

\begin{lemma}
\label{lem:td-uniform-L}
For \(3/2\le a\le2\),
\begin{equation}
\label{eq:td-uniform-L}
 L(a)>C_0,\qquad C_0:=\frac{83}{80}.
\end{equation}
\end{lemma}

\begin{proof}
From \eqref{eq:td-I}, $I'(a)=\frac{(2-a)(3a+1)}6,\qquad I''(a)=\frac{5-6a}{6}<0,\qquad I'''(a)=-1.$ In particular \(0<I\le5/12<C_0\).  The desired inequality is equivalent to
positivity of $J(a):=\frac{C_0I(a)}{C_0-I(a)}-\log a.$ Its derivative has the sign of $P(a):=aC_0^2 I'(a)-(C_0-I(a))^2.$ A second differentiation gives the exact identity
\begin{equation}
\label{eq:td-P-concavity}
 P''(a)=C_0^2\bigl(2I''(a)+aI'''(a)\bigr)
       +2(C_0-I(a))I''(a)-2(I'(a))^2<0.
\end{equation}
Also $P(3/2)=\frac{33871}{184320}>0, \qquad P(2)=-\left(C_0-\frac5{12}\right)^2<0.$ Strict concavity implies that \(P\) changes sign exactly once, from positive to
negative, on this interval.  Hence \(J\) first increases and then decreases;
its minimum occurs at an endpoint.

For \(0<z<1\), the positive power series for \(2\operatorname{arctanh}z\)
yields the elementary estimate
\begin{equation}
\label{eq:td-log-upper}
 \log\frac{1+z}{1-z}
 <2z+\frac{2z^3}{3(1-z^2)}.
\end{equation}
Indeed, after the linear term, replacing every denominator \(2j+1\) by \(3\)
produces the stated geometric-series upper bound.  Thus $\log\frac32<\frac{73}{180}<\frac{581}{1432} =\frac{C_0I(3/2)}{C_0-I(3/2)},$  $\log2<\frac{25}{36}<\frac{415}{596} =\frac{C_0I(2)}{C_0-I(2)}.$ Both endpoint values of \(J\) are positive, proving the lemma.
\end{proof}

For a candidate maximum with $m>2/3$, elimination of $\beta$ from $s(m)=0$ gives
\begin{equation}\label{eq:td-eliminate-beta}
 \beta=\frac{1+\rho(2m-1)(m-1)}{3m-2},\qquad
 \gamma^2M=\frac{m^3[1-\rho(1-m)^2]}{3m-2}
 \le\frac{m^3}{3m-2}.
\end{equation}
For $0<m<1$, differentiation gives the elementary inequality
\begin{equation}\label{eq:td-log-m}
 \log m<\frac{2(m-1)}{m+1}.
\end{equation}

\begin{proposition}
\label{prop:td-263-200}
Under \eqref{eq:td-normalised-domain}, condition \eqref{eq:td-candidate} holds for
\(0<\gamma\le263/200\).
\end{proposition}

\begin{proof}
We set $\Gamma:=\frac{263}{200},\qquad m_0:=\frac{757}{1000}.$ If a candidate had \(m\ge m_0\), then
\eqref{eq:td-eliminate-beta} and \eqref{eq:td-log-m} would imply $Mt_m^2< \frac{m^3[((\Gamma+2)m+\Gamma-2)]^2} {\Gamma^2(m+1)^2(3m-2)}<1.$ The second inequality follows from
\begin{align*}
 &\Gamma^2(m+1)^2(3m-2)
     -m^3[((\Gamma+2)m+\Gamma-2)]^2
       =(1-m)Q_\Gamma(m),\\
 &Q_\Gamma(m)=(\Gamma+2)^2m^4
 +(3\Gamma^2+4\Gamma-4)m^3
 +\Gamma^2m^2-3\Gamma^2m-2\Gamma^2,
\end{align*}
and the exact endpoint and derivative bounds
\[
 Q_\Gamma(m_0)
 =\frac{443926866382169}{40000000000000000}>0,
 \qquad
 Q_\Gamma'(m)\ge\frac34\Gamma^2+5\Gamma-1>0
 \quad(m\ge1/2).
\]
Thus \(m<m_0\).  We set \(z=(1-m_0)/(1+m_0)=243/1757\).  Keeping the first two
terms of the positive \(\operatorname{arctanh}\) series gives $\log m_0<-2\left(z+\frac{z^3}{3}\right).$ Therefore $t_m<\Gamma-2\left(z+\frac{z^3}{3}\right) =\frac{1124523209059}{1084789018600}<\frac{83}{80}.$ For completeness, the exact margin in the final comparison is $\frac{83}{80} -\frac{1124523209059}{1084789018600} =\frac{1890795477}{2169578037200}>0.$ Lemma~\ref{lem:td-uniform-L}, \eqref{eq:td-final-reduction} and Lemma~\ref{lem:td-maximum} complete the proof.
\end{proof}

For the remaining regions we may assume \(\gamma>\gamma_*>13/10\).  Their
normalised parameters satisfy
\begin{equation}
\label{eq:td-extension-common}
 \rho\ge\frac32,
 \qquad \rho\le\beta_0^2<\left(\frac{23}{13}\right)^2<\frac{16}{5},
 \qquad \gamma\le\frac75.
\end{equation}
We first localise the maximum and the shape parameter uniformly.

\begin{lemma}
\label{lem:td-extension-localisation}
Under \eqref{eq:td-extension-common}, a candidate has \(m<3/4\).
If \(m\le2/3\), then \(t_m<C_0\).  If \(m>2/3\), then
\(3/2<a<5/3\).
\end{lemma}

\begin{proof}
Suppose that \(m\ge3/4\).  Equations~\eqref{eq:td-eliminate-beta} and
\eqref{eq:td-log-m} give
\[
 \gamma^2M\le
 \frac{m^3[1-\tfrac32(1-m)^2]}{3m-2},
 \qquad
 0<\frac{t_m}{\gamma}
 \le1+\frac57\log m
 <\frac{17m-3}{7(m+1)}.
\]
It follows that $Mt_m^2< \frac{m^3[2-3(1-m)^2](17m-3)^2} {98(m+1)^2(3m-2)}<1.$ The last step is justified by
\begin{align*}
 &98(m+1)^2(3m-2)-m^3[2-3(1-m)^2](17m-3)^2
     =(1-m)Q_*(m),\\
 &Q_*(m)=-867m^6+1173m^5+245m^4+401m^3+98m^2-294m-196.
\end{align*}
Indeed, \(Q_*(3/4)=38369/4096>0\), and $Q_*'(m)=-5202m^5+5865m^4+980m^3+1203m^2+196m-294>0$ on \([3/4,1]\): the first two terms together are at least \(663m^4>0\),
and \(1203m^2>294\).  This contradicts feasibility, so \(m<3/4\).

If \(m\le2/3\), then $t_m\le\frac75+\log\frac23 <\frac75-\frac25=1<C_0.$ It remains to consider \(2/3<m<3/4\).  Since
\(\rho(c-1)(d-1)=1\), the root representation gives
\begin{equation}
\label{eq:td-rho-shape}
 \rho=\mathcal R(a,m):=
 \frac{2a-3}{(am-1)((3a-4)m-(2a-3))}.
\end{equation}
For \(m>2/3\) and \(3/2<a<2\), both denominator factors are positive; the
second equals \(m/2+(3m-2)(a-3/2)\).  Moreover
\begin{equation}
\label{eq:td-rho-monotone}
 \partial_a\mathcal R(a,m)=
 \frac{2m(a-1)(2-a)(3m-2)}
 {(am-1)^2((3a-4)m-(2a-3))^2}>0.
\end{equation}
At \(a=5/3\), $\mathcal R(5/3,m) =\frac1{((5/3)m-1)(3m-1)} >\mathcal R(5/3,3/4)=\frac{16}{5}.$ The bound \(\rho<16/5\) therefore implies \(a<5/3\).
\end{proof}

\begin{lemma}
\label{lem:td-linear-L}
For \(3/2\le a\le5/3\),
\begin{equation}
\label{eq:td-linear-L}
 L(a)>a-\frac{19}{40}.
\end{equation}
\end{lemma}

\begin{proof}
We set \(C(a)=a-19/40\).  Applying \eqref{eq:td-log-upper} with \(z=(a-1)/(a+1)\) we obtain $\log a<\ell(a):= \frac{(a-1)(a^2+10a+1)}{6a(a+1)}.$ Since \(C(a)>I(a)\), it is enough to prove
\(C I-(C-I)\ell>0\).  Polynomial multiplication gives
\begin{align*}
 C I-(C-I)\ell
 &=-\frac{(a-1)P_5(a)}{1440a(a+1)},\\
 P_5(a)&=280a^5+66a^4-1943a^3+1256a^2+819a+26.
\end{align*}
Writing \(a=3/2+x\), we have $P_5''(a)=5600x^3+25992x^2+28518x+5707>0.$ Thus \(P_5\) is strictly convex on the required interval.  Its endpoint values
are $P_5(3/2)=-\frac{67}{4},\qquad P_5(5/3)=-\frac{1312}{243}.$ Convexity places \(P_5\) below the chord joining these two negative values.
Hence \(P_5<0\), which proves the claim.
\end{proof}

We now use the particular upper bounds on \(r\).  Let
\begin{equation}
\label{eq:td-three-pairs}
 (\Gamma,R)\in
 \left\{\left(\frac{27}{20},\frac32\right),
         \left(\frac{11}{8},2\right),
         \left(\frac75,\frac52\right)\right\}.
\end{equation}
In the corresponding parameter region, \(\gamma\le\Gamma\) and \(\rho\ge R\).
Define
\begin{equation}
\label{eq:td-A-Gamma}
 A_\Gamma(m)=\Gamma+\frac{19}{40}+\frac{2(m-1)}{m+1}.
\end{equation}
For \(m<1\),
\begin{equation}
\label{eq:td-level-A}
 t_m=\gamma+\log m
 <\Gamma+\frac{2(m-1)}{m+1}
 =A_\Gamma(m)-\frac{19}{40}.
\end{equation}

\begin{lemma}
\label{lem:td-three-comparisons}
For each pair in \eqref{eq:td-three-pairs}, assume
\(2/3<m<3/4\), \(3/2<a<5/3\), and \(\mathcal R(a,m)=\rho\ge R\).
Then \(a\ge A_\Gamma(m)\).
\end{lemma}

\begin{proof}
We write $A=A_\Gamma(m)$, and let $l$ solve $A_\Gamma(l)=3/2$. For $m\le l$ the conclusion follows from $a>3/2$. For $l<m< u=3/4$, direct substitution in \eqref{eq:td-A-Gamma} gives $3/2<A<5/3$. By \eqref{eq:td-rho-monotone} it suffices to show $\mathcal R(A,m)<R$. Clearing the positive denominators reduces this to positivity of the quartic $B(m)=(1+m)^2\{R(Am-1)[(3A-4)m-(2A-3)]-(2A-3)\}.$ Set $z=(m-l)/(u-m)>0$, and expand
$(1+z)^4 B((l+uz)/(1+z))=\sum_{j=0}^4\binom4j b_jz^j$.
The entries below are $\lfloor1000 b_j\rfloor$ for $j=0,\ldots,4$, obtained by polynomial multiplication. Equivalently they follow from
$b_j=\sum_{i=0}^j\binom ji (u-l)^i B^{(i)}(l)/(i!\binom4i)$.
\begin{center}
\begin{tabular}{cccrrrrr}\toprule
$\Gamma$&$R$&$l$&$0$&$1$&$2$&$3$&$4$\\\midrule
$27/20$&$3/2$&$67/93$&128&102&76&53&32\\
$11/8$&$2$&$33/47$&108&76&49&30&21\\
$7/5$&$5/2$&$13/19$&63&33&14&12&37\\\bottomrule
\end{tabular}
\end{center}
Every coefficient is positive, proving all three comparisons.
\end{proof}

\begin{proof}[Completion of Theorem~\ref{thm:direct-trough}]
The unrestricted range follows from Proposition~\ref{prop:td-263-200}.
For any of the other three regions, take a candidate level and apply
Lemma~\ref{lem:td-extension-localisation}.
If \(m\le2/3\), that lemma gives \(t_m<C_0<L(a)\).
If \(m>2/3\), Lemma~\ref{lem:td-three-comparisons},
\eqref{eq:td-level-A}, and Lemma~\ref{lem:td-linear-L} give $t_m<A_\Gamma(m)-\frac{19}{40} \le a-\frac{19}{40}<L(a).$ In both cases, \eqref{eq:td-final-reduction} gives \(F(t_m,e)>0\), and
Lemma~\ref{lem:td-maximum} gives \(F(t,e)>0\).
This proves \(Z(\alpha)=\varnothing\) throughout
\(\mathcal P\).

\end{proof}

\begin{remark}
\label{rem:td-parameterwise}
The unrestricted-depth proof uses only \(\beta>1\), \(\rho\ge0\), and
\(\beta^2\ge\rho\), with the double-root case obtained by the same shape
calculation at \(a=2\).  It does not use the particular bounds
\(\beta>1+1/\gamma\) and \(\rho\le(1+1/\gamma)^2\) except to guarantee these
conditions.  Therefore, for a given \(\alpha>\alpha_0\) and
\(\gamma\le263/200\), the condition
\[
 \delta_\alpha=(\alpha-e^\gamma)^2-\mathsf c\ge0
\]
already implies \(Z(\alpha)=\varnothing\).
The restriction \(r\ge r_{\min}(\gamma)\) ensures that this condition holds
for every \(\alpha>\alpha_0\); it is not intrinsic to the maximum-point
exclusion itself.
\end{remark}

\subsection{A uniform bound for the slope}\label{sec:slope}
Theorem~\ref{thm:stagnation} gives $-\Psi_s\ge G_\alpha(\eta_{\max})^{1/2}$ on $\overline D$. For each wave with nonzero crest speed, this estimate bounds the angle $\chi$ of the free surface away from $\pi/2$. The bound may deteriorate along $\mathcal K_k$. To obtain a uniform estimate, we extend (N5) to a one-parameter family of nodal properties and impose a stronger concavity condition on $G_\alpha$.
As in the proof of Lemma~\ref{lem:robin} we write $Q:=|\widehat\nabla\Psi|=\e^{\mathfrak q}$. For $m\ge0$
let
\[
\mathfrak h_m:=\mathfrak f'\e^{m\mathfrak f},\qquad
T_m:=\operatorname{Re}\mathfrak h_m,\qquad U_m:=\operatorname{Im}\mathfrak h_m .
\]
Since $\e^{\mathfrak f}=-\Psi_s-\mathrm i\Psi_\Theta=Q\e^{\mathrm i\chi}$ and
$\p_\Theta\e^{m\mathfrak f}=m\mathfrak h_m$,
\begin{equation}\label{eq:Pm}
T_m=\frac1m\,\p_\Theta P_m\quad(m\ne0),\qquad P_m:=Q^m\cos(m\chi),
\end{equation}
and $T_0=\mathfrak a=\p_\Theta\mathfrak q$. On $\overline\Gamma_+$ we put
\[
\varrho_m:=g_*(\eta)\cos\big((m+1)\chi\big)-\varkappa\,\frac{\sin((m+1)\chi)}{\sin\chi},
\]
with the value $g_*(\eta)-(m+1)\varkappa$ at the crest and at the trough, so that $\varrho_0=\varrho$.
Finally let
\[
w_m:=\frac{m(4m+5)}{3(m+1)},\qquad
\omega_m(\chi):=m+\cot\chi\,\cot\big((m+1)\chi\big)-(m+1)\cot^2\big((m+1)\chi\big).
\]
The function $m\mapsto w_m$ increases from $w_0=0$ to $\infty$; we have $w_1=3/2$, and $w_m=2$ for
$m=(1+\sqrt{97})/8\approx1.356$.

\begin{lemma}\label{lem:Tm}
Let $(\alpha,H)\in\mathcal V_k$ solve $\Fop(\alpha,H)=0$ with $H\in\mathcal N_k$ and
$\alpha>\alpha_0$, and let $m\ge0$.
\begin{enumerate}[label=\textup{(\alph*)}]
\item $T_m$ is harmonic in $D$, odd with respect to $\Theta=0$ and to $\Theta=\pi/k$, and
$T_m=Q^m\mathfrak a>0$ on the bottom $(0,\pi/k)\times\{0\}$.
\item On $\overline\Gamma_+$,
$T_m=Q^m\varrho_m\sin\chi$ and
$U_m=Q^m\big(g_*(\eta)\sin\chi\,\sin((m+1)\chi)+\varkappa\cos((m+1)\chi)\big)$. On $\Gamma_+$,
$\p_nT_m=\p_tU_m$.
\item On the lines $\Theta=0$ and $\Theta=\pi/k$,
$\p_\Theta T_m=Q^m(\mathfrak a_\Theta-m\mathfrak b^2)$. At the crest and at the trough,
$\p_\Theta T_m=Q^m\varkappa\varrho_m$.
\item If $m>0$, $0<\chi<\pi/2$ and $(m+1)\chi<\pi$, then $\omega_m(\chi)<w_m$. If moreover
$(m+1)\chi<\pi/2$, then $\omega_m(\chi)>m$.
\item Let $m>0$ and $Q_c:=Q(0,\eta_{\max})=G_\alpha(\eta_{\max})^{1/2}$. If $T_m\ge0$ in $D_+$, then
$P_m\ge Q_c^m$ on $\overline D$. In particular $m|\chi|<\pi/2$ and $Q\ge Q_c$ on $\overline D$, and if
$m\ge1$, then also $-\Psi_s\ge Q_c$ on $\overline D$. Moreover $-\Psi_s$ is non-decreasing in $s$ on
the trough line.
\end{enumerate}
\end{lemma}

\begin{proof}
(a) $\mathfrak h_m$ is holomorphic in $D$, so $T_m$ is harmonic. Since $\mathfrak q$ and
$\mathfrak b$ are even and $\chi$ and $\mathfrak a$ are odd with respect to both lines,
$T_m=Q^m(\mathfrak a\cos(m\chi)-\mathfrak b\sin(m\chi))$ is odd. On the bottom $\chi=0$, so
$T_m=Q^m\mathfrak a>0$ by Lemma~\ref{lem:frakA}(b).

(b) By Lemma~\ref{lem:frakA}(c),
$Q^{-m}T_m=\mathfrak a\cos(m\chi)-\mathfrak b\sin(m\chi)
=g_*\sin\chi\cos((m+1)\chi)-\varkappa\sin((m+1)\chi)=\varrho_m\sin\chi$. From $Q^{-m}U_m=\mathfrak a\sin(m\chi)+\mathfrak b\cos(m\chi)$ we obtain the formula for $U_m$. The last
claim follows from the Cauchy--Riemann equations for $T_m+\mathrm iU_m$, as in
Lemma~\ref{lem:frakA}(c).

(c) On these lines $\mathfrak a=\chi=0$ and $\mathfrak b_\Theta=0$, so $\mathfrak f'=\mathrm i\mathfrak b$,
$\mathfrak f''=\mathfrak a_\Theta$ and
$\p_\Theta T_m=\operatorname{Re}\mathfrak h_m'=\operatorname{Re}\big((\mathfrak f''+m\mathfrak f'^2)\e^{m\mathfrak f}\big)
=Q^m(\mathfrak a_\Theta-m\mathfrak b^2)$. At the crest and at the trough
$\mathfrak a_\Theta=\varkappa\varrho=\varkappa(g_*-\varkappa)$ by Lemma~\ref{lem:frakA}(d) and $\mathfrak b=\varkappa$
by Lemma~\ref{lem:frakA}(c), so $\mathfrak a_\Theta-m\mathfrak b^2=\varkappa\varrho_m$.

(d) Let $n:=m+1$, $\varphi:=n\chi\in(0,\pi)$ and $\ell(x):=x\cot x$. Then
$\omega_m(\chi)-m=\cot\varphi\,(\cot\chi-n\cot\varphi)$, and $w_m-m=m(m+2)/(3(m+1))>0$. If
$\varphi\ge\pi/2$, then $\cot\varphi\le0<\cot\chi-n\cot\varphi$, so $\omega_m(\chi)\le m<w_m$. Let
$\varphi<\pi/2$. Since $\cot\chi=\ell(\chi)/\chi$ and $n\cot\varphi=\ell(\varphi)/\chi$, $\omega_m(\chi)-m=\frac{\ell(\varphi)\,\big(\ell(\chi)-\ell(\varphi)\big)}{\chi\varphi}.$ The function $\ell$ is positive and strictly decreasing on $(0,\pi/2)$, so $\omega_m(\chi)>m$. For the
upper bound we use the expansion $1-\ell(x)=\sum_{j\ge1}c_jx^{2j}$, $|x|<\pi$, whose coefficients
$c_j=2^{2j}|B_{2j}|/(2j)!$ are positive. Since
$\varphi^{2j}-\chi^{2j}\le j\varphi^{2j-2}(\varphi^2-\chi^2)$ for
$0<\chi<\varphi$, with strict inequality for $j\ge2$,
\[
\ell(\chi)-\ell(\varphi)=\sum_{j\ge1}c_j\big(\varphi^{2j}-\chi^{2j}\big)
<(\varphi^2-\chi^2)\sum_{j\ge1}jc_j\varphi^{2j-2}=-(\varphi^2-\chi^2)\frac{\ell'(\varphi)}{2\varphi}.
\]
Moreover $-\ell\ell'\le2\varphi/3$ on $(0,\pi/2)$: this inequality is equivalent to $E(\varphi):=4\sin^3\varphi-3\cos\varphi\,(2\varphi-\sin2\varphi)\ge0,$ and $E(0)=0$, $E'(\varphi)=3\sin\varphi\,(2\varphi-\sin2\varphi)>0$. Hence $\omega_m(\chi)-m<\frac{\varphi^2-\chi^2}{3\chi\varphi}=\frac{n^2-1}{3n}=w_m-m .$ (e) By \eqref{eq:Pm}, $P_m$ is non-decreasing in $\Theta$ on every horizontal segment in
$\overline{D_+}$. By (N1), $\eta$ is non-increasing on $[0,\pi/k]$, so for every
$(\Theta_0,s_0)\in\overline{D_+}$ the segment $[0,\Theta_0]\times\{s_0\}$ lies in $\overline{D_+}$
and $P_m(\Theta_0,s_0)\ge P_m(0,s_0)$. On the crest line $\chi=0$, and by (c) and the
Cauchy--Riemann equations $\mathfrak b_s=\mathfrak a_\Theta\ge m\mathfrak b^2\ge0$; since $\mathfrak b=0$
on the bottom, $\mathfrak b\ge0$ and $\p_s\mathfrak q=-\mathfrak b\le0$ there. Hence $Q$ is
non-increasing on the crest line, and $P_m(0,s_0)=Q(0,s_0)^m\ge Q_c^m$. Thus $P_m\ge Q_c^m>0$ on
$\overline{D_+}$, and by symmetry on $\overline D$. In particular $Q^m\ge P_m\ge Q_c^m$ and
$\cos(m\chi)>0$ on $\overline D$; since $\chi$ is continuous, $\chi=0$ on the bottom and $\overline D$
is connected, $m|\chi|<\pi/2$ on $\overline D$. If $m\ge1$, the function $m\ln\cos\chi-\ln\cos(m\chi)$
vanishes at $\chi=0$, and for $m|\chi|<\pi/2$ its derivative $m\big(\tan(m\chi)-\tan\chi\big)$ has the
sign of $\chi$; hence $(-\Psi_s)^m=Q^m\cos^m\chi\ge P_m\ge Q_c^m$. On the trough line, (c) gives $\mathfrak b_s=\mathfrak a_\Theta\le m\mathfrak b^2$,
and $\mathfrak b=0$ at the bottom. If $\mathfrak b(\pi/k,s_1)>0$ for some $s_1$, let $s_0<s_1$ be the
largest zero of $\mathfrak b(\pi/k,\cdot)$ below $s_1$; on $[s_0,s_1]$ we have
$\mathfrak b_s\le L\mathfrak b$ with $L:=m\max\mathfrak b$, so $\e^{-Ls}\mathfrak b$ is non-increasing
there and $\mathfrak b(\pi/k,s_1)\le0$, a contradiction. Hence $\mathfrak b\le0$ and
$\p_s\mathfrak q\ge0$ on the trough line, where $-\Psi_s=Q$.
\end{proof}

\begin{theorem}\label{thm:slope}
Let $\mathcal K_k$ be the curve of Theorem~\ref{thm:global}, let $0<s_1\le\infty$ and $m>0$, and
suppose that for every $s\in(0,s_1)$ the solution $H(s)$ satisfies \textup{(G2${}'$)} and
\begin{enumerate}[label=\textup{(G3${}_m$)}]
\item $g_*'\ge w_m\,g_*^2$ on $[\eta_{\min},\eta_{\max}]$.
\end{enumerate}
Then for every $s\in(0,s_1)$, $\varrho_m>0$ on $\overline\Gamma_+$, $T_m>0$ in $D_+$, the
conclusions of Theorem~\ref{thm:stagnation} hold, and $m|\chi|<\pi/2$ on $\overline D$. In particular,
if $m>1$, then
\[
0<\chi<\frac\pi{2m}\quad\text{on }\Gamma_+ .
\]
If $s_1=\infty$, the relative velocity at the crest tends to zero along $\mathcal K_k$; if moreover
$m>1$, the angle between the free surface and the circles $R=\mathrm{const}$ stays below $\pi/(2m)$.
\end{theorem}

\begin{proof}
Since $w_m>0$, (G3${}_m$) implies (G1), so Theorem~\ref{thm:stagnation} applies: $\mathfrak a>0$ in
$D_+$ for all $s\in(0,s_1)$. As in Case~3 of its proof, Lemma~\ref{lem:ffsym}(b) then gives
$F_\alpha(\eta_{\min},\hat x_3)<F_\alpha(\eta_{\min},\eta_{\max})<0$, and (G2${}'$) forces
$G_\alpha'(\eta_{\min})\le0$, that is, $g_*(\eta_{\min})\ge0$, for all $s\in(0,s_1)$.

We now repeat the continuation argument of Theorem~\ref{thm:stagnation} with $\varrho_m$ in place
of $\varrho$. The function $\varrho_{m,s}$ of $H(s)$ depends continuously on $s>0$ in
$C([0,\pi/k])$, and it converges uniformly to $g_*(\gamma)>1/\gamma$ as $s\downarrow0$, because
$\chi$ and $\varkappa$ tend to zero uniformly. Let $s_\sharp$ be the supremum of all $s'\in(0,s_1]$ such
that $\varrho_{m,s}>0$ on $[0,\pi/k]$ for $0<s<s'$, and suppose that $s_\sharp<s_1$. For
$H:=H(s_\sharp)$, $\varrho_m\ge0$ on $\overline\Gamma_+$ with a zero at some $P_0\in\overline\Gamma_+$.
By Lemma~\ref{lem:Tm}(a),(b), $T_m\ge0$ on $\p D_+$ and $T_m>0$ on the bottom, so $T_m>0$ in $D_+$.
On $\overline\Gamma_+$ we have $0\le\chi<\pi/2$, and $m\chi<\pi/2$ by Lemma~\ref{lem:Tm}(e); since
$(m+1)/(2\max\{m,1\})\le1$, it follows that $0\le(m+1)\chi<\pi$ on $\overline\Gamma_+$.

\emph{Case 1: $P_0\in\Gamma_+$.} Then $T_m(P_0)=0=\min_{\overline{D_+}}T_m$, so
$\p_nT_m(P_0)<0$ by the Hopf lemma, and $\p_tT_m(P_0)=0$. We write $\varphi:=(m+1)\chi\in(0,\pi)$.
Along $\Gamma_+$ we have $\p_t\chi=-\varkappa$, $\p_t\eta=\sin\chi$ and $\p_t\mathfrak q=-g_*\sin\chi$. At
$P_0$, $T_m=0$ gives $\varkappa=g_*\sin\chi\cot\varphi$, and $\p_tT_m=0$ gives
\[
\p_t\varkappa\,\sin\varphi=g_*'\sin^2\chi\cos\varphi-g_*\varkappa\cos\chi\cos\varphi
+(m+1)\big(g_*\varkappa\sin\chi\sin\varphi+\varkappa^2\cos\varphi\big).
\]
Differentiating the formula for $U_m$ in Lemma~\ref{lem:Tm}(b) and inserting these two relations,
we obtain
\[
\p_nT_m(P_0)=\p_tU_m(P_0)=Q^m\,\frac{\sin^2\chi}{\sin\varphi}\Big(g_*'(\eta)-\omega_m(\chi)\,g_*(\eta)^2\Big)
\ge Q^m\,\frac{\sin^2\chi}{\sin\varphi}\Big(g_*'-w_mg_*^2\Big)\ge0
\]
by Lemma~\ref{lem:Tm}(d) and (G3${}_m$), a contradiction.

\emph{Case 2: $P_0$ is the crest.} Then $g_*(\eta_{\max})=(m+1)\varkappa_c$, and $\p_\Theta T_m=0$ at
the crest by Lemma~\ref{lem:Tm}(c). Since $T_m>0$ in $D_+$ is odd in $\Theta$, Serrin's corner lemma
gives $\p_t^2U_m>0$ at the crest, exactly as in Case~2 of the proof of
Theorem~\ref{thm:stagnation}. At the crest $\chi=0$, $\p_t\chi=-\varkappa_c$, $\p_t\varkappa=0$,
$\p_t\eta=0$, $\p_t^2\eta=-\varkappa_c$, $\p_t\mathfrak q=0$ and $\p_t^2\mathfrak q=g_*\varkappa_c$. Expanding
the formulas of Lemma~\ref{lem:Tm}(b) and the definition of $\varrho_m$ to second order along the
free surface and using $g_*=(m+1)\varkappa_c$, we find, with $\varkappa_2:=\p_t^2\varkappa$ at the crest,
\begin{gather*}
\p_t^2U_m=Q_c^m\big(\varkappa_2+(m+1)(2m+1)\varkappa_c^3\big),\\
\p_t^2\varrho_m=-g_*'\varkappa_c-(m+1)\varkappa_2-(m+1)\Big((m+1)^2-\frac{m(m+2)}3\Big)\varkappa_c^3 .
\end{gather*}
Thus $\varkappa_2>-(m+1)(2m+1)\varkappa_c^3$. Since $\varrho_m\ge0$ on $\overline\Gamma_+$ is even with
respect to the crest and vanishes there, $\p_t^2\varrho_m\ge0$, and therefore
\[
0\le\p_t^2\varrho_m<-g_*'\varkappa_c+\frac{m(m+1)(4m+5)}3\,\varkappa_c^3
=\varkappa_c\big(w_mg_*^2-g_*'\big)\Big|_{\eta_{\max}},
\]
which contradicts (G3${}_m$) because $\varkappa_c>0$.

\emph{Case 3: $P_0$ is the trough.} Then $g_*(\eta_{\min})=(m+1)\varkappa_t<0$, which contradicts
$g_*(\eta_{\min})\ge0$.

Hence we obtain $s_\sharp=s_1$, and Lemma~\ref{lem:Tm}(e) gives the bound for $\chi$. The last statement
follows from Theorem~\ref{thm:stagnation}.
\end{proof}

Since $g_*'-2pg_*^2=\big((1-p)G_\alpha'^2-G_\alpha G_\alpha''\big)/(2G_\alpha^2)$, condition
(G3${}_m$) says that $G_\alpha^{w_m/2}$ is concave on $[\eta_{\min},\eta_{\max}]$; (G1) is the limit
$m\to0$. In the range of Theorem~\ref{thm:stagCAP}, (G3${}_m$) holds with an explicit $m$; this gives
a uniform bound on the slope when $m>1$, that is, when $w_m>3/2$.

\begin{corollary}\label{cor:slope}
Let $r\ge\gamma^2/(1+\gamma)^2$, and suppose that $Z(\alpha)=\emptyset$
for every $\alpha>\alpha_0$, as is the case for $(\gamma,r)\in\mathcal P$.
Then for every $k\ge1$ and every $s>0$:
\begin{enumerate}[label=\textup{(\alph*)}]
\item $G_{\alpha(s)}'(\eta_{\min}(s))\le0$.
\item If
$\alpha_0\le\dfrac{631}{25}\,\delta_{\alpha_0}^{1/2}$,
there is $\theta_*=\theta_*(\gamma,r)<\pi/2$, independent of $k$ and $s$,
such that $0<\chi<\theta_*$ on $\Gamma_+$.
\item If $\alpha_0\le6\delta_{\alpha_0}^{1/2}$, equivalently
\begin{equation}\label{eq:r2}
35\beta_0^2r^2-2(\beta_0+18)r-1\ge0,
\qquad \beta_0:=1+1/\gamma,
\end{equation}
then \textup{(b)} holds with
$\theta_*=4\pi/(1+\sqrt{97})\approx1.1583$, that is, $66.4^\circ$.
\end{enumerate}
\end{corollary}

\begin{proof}
The lower bound on $r$ gives \textup{(G1)} by
Corollary~\ref{cor:stagnation}(a), and $Z(\alpha)=\emptyset$ gives
\textup{(G2${}'$)} by its definition. The first paragraph of the proof of
Theorem~\ref{thm:slope} therefore proves (a).

Under either of the additional hypotheses in (b) and (c),
$\delta_{\alpha_0}>0$. For $0<x<\ln z_-$ set
$y=\e^x$, $A=y/(z_--y)$ and $B=y/(z_+-y)$. Then $g_*=(A+B)/2-1,\qquad g_*'=(A+A^2+B+B^2)/2, \qquad B=\frac{b_\alpha A}{A+b_\alpha+1},$ where
$b_\alpha=z_-/(z_+-z_-)
=(\alpha\delta_\alpha^{-1/2}-1)/2$.
The function $b_\alpha$ is decreasing, since
\[
\frac{\mathrm d}{\mathrm d\alpha}\frac{\alpha^2}{\delta_\alpha}
=-\frac{2\alpha\bigl(\e^\gamma(\alpha-\e^\gamma)+\cc\bigr)}
{\delta_\alpha^2}<0.
\]
Thus $0<b_\alpha\le b_{\alpha_0}$. On the surface range we have
$A\ge B>0$ and $A+B\ge2$ by (a). At $A+B=2$ the condition
\textup{(G3${}_m$)} is automatic, because $g_*=0<g_*'$.
It remains to bound $\Lambda=g_*'/g_*^2$ on the set $\Sigma$ with
$0\le b\le b_{\alpha_0}$, $A>0$, $B=bA/(A+b+1)$ and $A+B>2$.

For (b) the hypothesis is $b_{\alpha_0}\le303/25$. We set $F(A,B):=A^2-6AB+B^2+16(A+B)-12, \qquad \Lambda-\frac32=\frac{F(A,B)}{2(A+B-2)^2}.$ First $F>0$ whenever $0\le B<4+\sqrt{26}/2$ and
$A\ge\max\{B,2-B\}$. Indeed, for $B\le1$ the minimum in $A$ occurs at
$A=2-B$ and equals $8((B-1)^2+2)$; for $1<B<4$ the minimum occurs at
$A=B$ and equals $4(8B-B^2-3)$; and for $4\le B<4+\sqrt{26}/2$ we have
$F(A,B)\ge F(3B-8,B)=4(16B-2B^2-19)>0$.
If $A\le12$, then
$B\le909/157<6<4+\sqrt{26}/2$, so this proves $F>0$.

If $A\ge12$, define
\begin{align*}
P(A,b)&:=(A+b+1)^2F\left(A,\frac{bA}{A+b+1}\right)\\
&=A^4+(18-4b)A^3+(21+44b-4b^2)A^2
 +(32b^2+24b-8)A-12(b+1)^2.
\end{align*}
For $A\ge12$ and $b\ge0$, $\partial_bP =2b(32A-4A^2-12)-4A^2(A-11)+24A-24<0.$ Consequently $P(A,b)\ge P(A,303/25)$. With $u=5A-104\ge-44$,
\begin{align*}
15625P(A,303/25)
&=25u^4+6590u^3+412869u^2+104272u+1185344\\
&\ge122909u^2+104272u+1185344>0.
\end{align*}
Here $6590u^3\ge-289960u^2$, and the last quadratic has positive leading
coefficient and discriminant $-571885132800$. Thus $\Lambda>3/2$ on
$\Sigma$. Also $\Lambda\to2$ as $A\to\infty$, uniformly for
$0\le b\le303/25$, while $\Lambda\to\infty$ as $A+B\downarrow2$.
Continuity on the remaining compact set gives
$W:=\inf_\Sigma\Lambda\in(3/2,2]$.
We choose $m>1$ with $w_m=W$. Theorem~\ref{thm:slope}, with
$s_1=\infty$, proves (b) with $\theta_*=\pi/(2m)$.

For (c), $4(g_*'-2g_*^2)=2\bigl((5-2B)A+5B-4\bigr).$ If $B\le1$, the expression in parentheses is at least
$2(B-1)^2+4>0$, using $A\ge2-B$.
If $1<B\le5/2$, it is at least $-2(B^2-5B+2)>0$, using $A\ge B$.
The hypothesis in (c) gives $B<b_\alpha\le b_{\alpha_0}\le5/2$.
Hence \textup{(G3${}_m$)} holds with $w_m=2$, namely
$m=(1+\sqrt{97})/8$, which proves the stated angle bound.
Finally, $\alpha_0=\e^\gamma(1+r\beta_0)$ and
$\delta_{\alpha_0}=\e^{2\gamma}(r^2\beta_0^2-r)$ give
\eqref{eq:r2} upon squaring.
\end{proof}
Condition (b) is $r\ge r_\sharp(\gamma)$, where $r_\sharp$ is the positive root of
\[
397536\beta_0^2r^2-(1250\beta_0+398161)r-625=0.
\]
This follows by squaring the positive quantities in (b). The gap $r_{\min}\le r<r_\sharp$ is relevant to the limiting-wave theorem: crest stagnation there does not by itself bound the slope.
\subsection{Limiting waves and the Stokes corner}\label{sec:corner}
By Theorem~\ref{thm:stagnation}, the relative velocity at the crest tends to zero along
$\mathcal K_k$. In the range of Corollary~\ref{cor:slope}(b) the profiles are moreover uniformly
Lipschitz, so the waves have limits along subsequences. We show that every such limit is a weak
solution with a stagnation point at the crest, where the free surface has the corner of Stokes:
it meets the circles $R=\mathrm{const}$ at the angles $\pm30^\circ$, so that the angle of the fluid
domain at the crest is $120^\circ$. The proof combines a bound for the relative speed with the
monotonicity formula of Weiss \cite{Weiss1999}, in the form used by V\u{a}rv\u{a}ruc\u{a} and Weiss
\cite{VarvarucaWeiss2011} for the Stokes conjecture; the Bernoulli head $G_\alpha$ enters only
through its simple zero at the stagnation level.

For $\alpha>\alpha_0$ with $\delta_\alpha>0$ we put $L_\alpha:=\ln z_-(\alpha)$ and
\[
K_\alpha:=\max\Big\{\sup_{x_1(\alpha)<x<L_\alpha}\frac{G_\alpha(x)}{L_\alpha-x},\ \frac1{x_1(\alpha)^2L_\alpha}\Big\},
\]
which is finite because $G_\alpha$ is positive on $[x_1(\alpha),L_\alpha)$ and has a simple zero at
$L_\alpha$. It depends continuously on $\alpha$.

\begin{lemma}\label{lem:speedbound}
Let $(\alpha,H)\in\mathcal V_k$ solve $\Fop(\alpha,H)=0$ with $H\in\mathcal N_k$, $\alpha>\alpha_0$
and $\delta_\alpha>0$. Then, on $\overline D$,
\[
|\widehat\nabla\Psi|^2\le K_\alpha\,(L_\alpha-s),\qquad \Psi\le K_\alpha^{1/2}\,(L_\alpha-s)^{3/2}.
\]
\end{lemma}

\begin{proof}
Since $\Psi$ is harmonic and $L_\alpha-s$ is affine, $|\widehat\nabla\Psi|^2-K_\alpha(L_\alpha-s)$ is
subharmonic in $D$. On the free surface it equals $G_\alpha(\eta)-K_\alpha(L_\alpha-\eta)\le0$,
because $x_1(\alpha)<\eta<L_\alpha$ by Proposition~\ref{prop:heights} and
Lemma~\ref{lem:surfspeed}. On the bottom, $|\widehat\nabla\Psi|^2<\eta_{\min}^{-2}<x_1(\alpha)^{-2}\le
K_\alpha L_\alpha$ by Lemma~\ref{lem:compare}(c). The maximum principle gives the first inequality.
For $(\Theta,s)\in\overline D$ we have $s\le\eta(\Theta)<L_\alpha$ and $\Psi(\Theta,\eta(\Theta))=0$, so
$\Psi(\Theta,s)\le\int_s^{\eta(\Theta)}|\widehat\nabla\Psi(\Theta,t)|\dd t\le(\eta(\Theta)-s)\,
K_\alpha^{1/2}(L_\alpha-s)^{1/2}$.
\end{proof}

\begin{theorem}\label{thm:corner}
Let $(\gamma,r)\in\mathcal P$ with $r\ge r_\sharp(\gamma)$, where $\mathcal P$ is the set of
Theorem~\ref{thm:stagCAP}, let $k\ge1$, let $\mathcal K_k$ be the curve of
Theorem~\ref{thm:global}, and let $s_n\to\infty$. We write $\eta_n$ and $\Psi_n$ for the profile and the
stream function of $H(s_n)$, with $\Psi_n:=0$ above the free surface. After passing to a
subsequence, $\alpha(s_n)\to\alpha_\infty$, $\eta_n\to\eta_\infty$ uniformly and $\Psi_n\to\Psi_\infty$
uniformly on $[0,2\pi/k]\times[0,\infty)$, and the following hold, with $L:=L_{\alpha_\infty}$,
$x_c:=(0,L)$ and $D_\infty:=\{(\Theta,s):0<s<\eta_\infty(\Theta)\}$.
\begin{enumerate}[label=\textup{(\alph*)}]
\item $\alpha_0<\alpha_\infty\le A_*$ and $\delta_{\alpha_\infty}>0$. The profile $\eta_\infty$ is even,
$2\pi/k$-periodic, non-increasing on $[0,\pi/k]$ and Lipschitz with constant $\tan\theta_*$, where
$\theta_*$ is the angle of Corollary~\ref{cor:slope}(b); $\eta_\infty(0)=L$ and
$x_1(\alpha_\infty)<\eta_\infty(\pi/k)<\gamma$.
\item $\Psi_\infty$ is Lipschitz, harmonic and positive in $D_\infty$, equal to $1$ on $s=0$ and to
$0$ on and above the graph of $\eta_\infty$, and $|\widehat\nabla\Psi_\infty|^2\le K(L-s)$ in $D_\infty$,
where $K:=\sup_nK_{\alpha(s_n)}<\infty$. For every $\xi\in C^1_c(\R^2;\R^2)$ that vanishes near
$\{s=0\}$,
\begin{equation}\label{eq:domvar}
\int_{D_\infty}\Big(|\widehat\nabla\Psi_\infty|^2\,\mathrm{div}\,\xi-2\,\widehat\nabla\Psi_\infty\cdot
(D\xi)\,\widehat\nabla\Psi_\infty+\mathrm{div}\big(G_{\alpha_\infty}(s)\,\xi\big)\Big)\dd\Theta\dd s=0 .
\end{equation}
\item \emph{(Stokes corner.)} $\eta_\infty(\Theta)<L$ for $0<|\Theta|\le\pi/k$, and $\eta_\infty$ has
one-sided derivatives $\eta_\infty'(0\pm)=\mp1/\sqrt3$ at the crest. Moreover, with
$\mu:=-G_{\alpha_\infty}'(L)>0$ and $\phi(y)\in(-\pi,\pi]$ the angle between $y\in\R^2$ and $(0,-1)$,
\[
\frac{\Psi_\infty(x_c+\rho y)}{\rho^{3/2}}\longrightarrow\frac{\sqrt{2\mu}}3\,|y|^{3/2}\cos\Big(\frac32\,\phi(y)\Big)
\,\mathbf 1_{\{|\phi(y)|<\pi/3\}}
\qquad\text{locally uniformly in }y,\ \text{as }\rho\downarrow0 .
\]
\item \emph{(Regularity away from the crests.)} $\eta_n\to\eta_\infty$ in $C^j$ on every compact subset
of $\R\setminus\frac{2\pi}k\mathbb Z$, for every $j$. Thus $\eta_\infty$ is of class $C^\infty$ there,
$\Psi_\infty$ is smooth up to the free surface away from the crests, and
$|\widehat\nabla\Psi_\infty|^2=G_{\alpha_\infty}(\eta_\infty)$ holds there in the classical sense.
Moreover $\widehat\nabla\Psi_\infty\ne0$ on $\overline{D_\infty}$ except at the crests: the relative
flow of the limiting wave stagnates only there.
\end{enumerate}
\end{theorem}
\begin{figure}[htbp]\centering
\includegraphics[width=.91\linewidth]{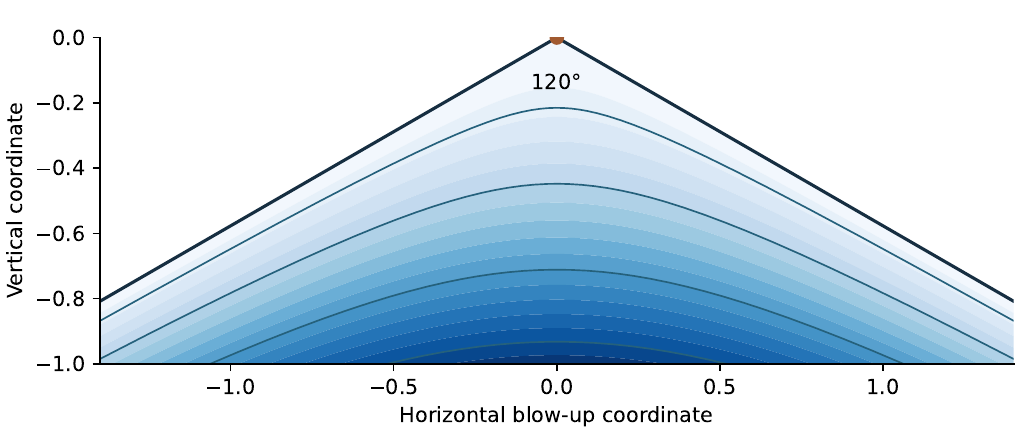}
\caption{Local blow-up at the limiting irrotational crest under the hypotheses of the preceding theorem. The two rays have slopes $\pm1/\sqrt3$ and enclose a fluid angle of $120^\circ$. Contours in the fluid show the limiting stream function after division by $\sqrt{2\mu}/3$. This is the asymptotic crest geometry.}\label{fig:atlas-corner}
\end{figure}

Thus every limiting wave has a stagnation point of Stokes type at the crest: the free surface has a
symmetric corner whose sides meet the circles $R=\mathrm{const}$ at $\pm30^\circ$, and the fluid
occupies an angle of $120^\circ$ there. Away from the crests the limiting wave is a smooth classical
solution, so the crests are its only singular points, as for the Stokes wave of greatest height.
Identity \eqref{eq:domvar} is the weak form of the boundary conditions $\Psi_\infty=0$ and
$|\widehat\nabla\Psi_\infty|^2=G_{\alpha_\infty}(\eta_\infty)$ on the whole free surface.

\begin{proof}
\emph{Step 1: compactness.} By Theorems~\ref{thm:global}, \ref{thm:stagnation} and
\ref{thm:stagCAP} and Corollary~\ref{cor:slope}(b), for every $s>0$ we have
$\alpha_0<\alpha(s)\le A_*$, $x_1(\alpha(s))<\eta_{\min}(s)<\gamma<\eta_{\max}(s)<L_{\alpha(s)}$ and
$|\eta'|=\tan\chi\le\tan\theta_*$ by \eqref{eq:qchi}, and every limit point $(\bar\alpha,\bar\eta)$
of $(\alpha(s),\eta_{\max}(s))$ as $s\to\infty$ satisfies $\e^{\bar\eta}=z_-(\bar\alpha)$. Since
$r>\gamma^2/(1+\gamma)^2$, we have $\alpha_0-\e^\gamma=r\e^\gamma\beta_0>(r\e^{2\gamma})^{1/2}=\cc^{1/2}$,
so $\delta_\alpha\ge\delta_{\alpha_0}>0$ for $\alpha\ge\alpha_0$. By the Arzel\`a--Ascoli theorem we may
assume that $\alpha(s_n)\to\alpha_\infty$ and $\eta_n\to\eta_\infty$ uniformly, and (a) follows, except for
the strict inequalities $\alpha_\infty>\alpha_0$ and $x_1(\alpha_\infty)<\eta_\infty(\pi/k)<\gamma$, which we prove
in Step~8.

The speed constants are uniformly bounded even before we know that $\alpha_\infty>\alpha_0$. Indeed, set
\[
d_*:=x_1(A_*)>0,\qquad
L_*:=\max_{\alpha_0\le\alpha\le A_*}L_\alpha,
\qquad
C_G:=\max_{\substack{\alpha_0\le\alpha\le A_*\\d_*\le x\le L_*}}
|G_\alpha'(x)|.
\]
These constants are finite, since $\delta_{\alpha_0}>0$.
For $x_1(\alpha)<x<L_\alpha$, the mean value theorem and
$G_\alpha(L_\alpha)=0$ give
$G_\alpha(x)/(L_\alpha-x)\le C_G$. Also
$x_1(\alpha)\ge d_*$ and $L_\alpha>\gamma$. Consequently $K_\alpha\le\max\{C_G,(d_*^2\gamma)^{-1}\} \qquad(\alpha_0<\alpha\le A_*).$ By Lemma~\ref{lem:speedbound} the $\Psi_n$ are uniformly Lipschitz on $\overline{D_n}$; since they vanish
on the free surface, their extensions by zero are uniformly Lipschitz as well, and we may assume
that $\Psi_n\to\Psi_\infty$ uniformly. Every compact subset of $D_\infty$ lies in the fluid domains
$D_n$ of $H(s_n)$ for large $n$, so $\Psi_\infty$ is harmonic in $D_\infty$, the derivatives converge
locally uniformly there, and the bound for $|\widehat\nabla\Psi_\infty|^2$ follows from
Lemma~\ref{lem:speedbound}. The boundary values are clear, and $\Psi_\infty>0$ in $D_\infty$ by the
strong maximum principle.

\emph{Step 2: the identity \eqref{eq:domvar}.} For a classical solution and $\xi$ vanishing near the
bottom, Rellich's identity
$\mathrm{div}\big(|\widehat\nabla\Psi|^2\xi-2(\widehat\nabla\Psi\cdot\xi)\widehat\nabla\Psi\big)
=|\widehat\nabla\Psi|^2\mathrm{div}\,\xi-2\widehat\nabla\Psi\cdot(D\xi)\widehat\nabla\Psi$ and the
divergence theorem give
\[
\int_D\big(|\widehat\nabla\Psi|^2\mathrm{div}\,\xi-2\widehat\nabla\Psi\cdot(D\xi)\widehat\nabla\Psi\big)
=-\int_\Gamma G_\alpha(\eta)\,\xi\cdot\nu\dd\sigma=-\int_D\mathrm{div}\big(G_\alpha(s)\xi\big),
\]
because $\widehat\nabla\Psi=-|\widehat\nabla\Psi|\nu$ and $|\widehat\nabla\Psi|^2=G_\alpha(\eta)$ on
$\Gamma$. Thus \eqref{eq:domvar} holds for $\Psi_n$, $D_n$ and $\alpha(s_n)$. To pass to the limit
it suffices to show that $\widehat\nabla\Psi_n\to\widehat\nabla\Psi_\infty$ in $L^2_{\rm loc}$. Let
$\varphi\in C^\infty_c(\R^2)$ vanish near $s=0$. Since $\Psi_n$ is harmonic in $D_n$ and vanishes on
the free surface,
$\int|\widehat\nabla\Psi_n|^2\varphi=-\int\Psi_n\widehat\nabla\Psi_n\cdot\widehat\nabla\varphi$, and
the right side converges to $-\int\Psi_\infty\widehat\nabla\Psi_\infty\cdot\widehat\nabla\varphi$.
The same identity holds for $\Psi_\infty$: the function $(\Psi_\infty-\varepsilon)_+\varphi$ is
Lipschitz with compact support in $D_\infty$, so
$\int\widehat\nabla\Psi_\infty\cdot\widehat\nabla\big((\Psi_\infty-\varepsilon)_+\varphi\big)=0$, and we
let $\varepsilon\downarrow0$. Hence $\int|\widehat\nabla\Psi_n|^2\varphi\to\int|\widehat\nabla\Psi_\infty|^2
\varphi$, which together with the weak-$*$ convergence of $\widehat\nabla\Psi_n$ gives the
$L^2_{\rm loc}$ convergence. The same argument gives
$\int_{B_r(x)}|\widehat\nabla\Psi_\infty|^2=\int_{\p B_r(x)}\Psi_\infty\,\p_\nu\Psi_\infty$ for every
ball with $r<\operatorname{dist}(x,\{s=0\})$ and almost every $r$.

\emph{Step 3: growth at the crest.} By Step~1 and the proof of Lemma~\ref{lem:speedbound},
$\Psi_\infty\le K^{1/2}(L-s)^{3/2}\le K^{1/2}|x-x_c|^{3/2}$ in $D_\infty$. We claim that
\begin{equation}\label{eq:limsupslope}
\limsup_{\Theta\downarrow0}\frac{L-\eta_\infty(\Theta)}{\Theta}\ge\frac1{\sqrt3}.
\end{equation}
Otherwise there are $0\le\tau<\tau'<1/\sqrt3$ and $\rho_0>0$ such that $L-\eta_\infty(\Theta)\le
\tau|\Theta|$ for $|\Theta|\le\rho_0$. Let $\omega:=2\arctan(1/\tau')>2\pi/3$, $\nu:=\pi/\omega<3/2$,
and use polar coordinates $(\rho,\phi)$ about $x_c$, with $\phi$ measured from the downward vertical.
For small $\rho_1>0$ the closed sector $\overline S:=\{0\le\rho\le\rho_1,\ |\phi|\le\omega/2\}$
satisfies $\overline S\setminus\{x_c\}\subset D_\infty$, because $|\phi|\le\arctan(1/\tau')<
\arctan(1/\tau)$ there. The function $h:=\rho^\nu\cos(\nu\phi)$ is harmonic in $S$, vanishes on its
straight sides, and is at most $\rho_1^\nu$ on its arc, where $\Psi_\infty\ge m>0$. By the maximum
principle $\Psi_\infty\ge m\rho_1^{-\nu}h$ in $S$, so $\Psi_\infty(x_c-(0,\rho))\ge m\rho_1^{-\nu}\rho^\nu$,
which contradicts $\Psi_\infty\le K^{1/2}\rho^{3/2}$ for small $\rho$.

\emph{Step 4: a monotonicity formula.} Let $u(y):=\Psi_\infty(x_c+y)$, let $\chi_\Omega$ be the
characteristic function of $\Omega:=D_\infty-x_c$, and let $g(y_2):=G_{\alpha_\infty}(L+y_2)$. Then
$g(y_2)=-\mu y_2+e(y_2)$ with $|e(y_2)|\le C_1y_2^2$ and $|e'(y_2)|\le C_1|y_2|$ for $|y_2|\le L$. For
$0<r<r_0:=\min\{L,\pi/k\}/2$ put
\[
W(r):=r^{-3}\int_{B_r}\big(|\nabla u|^2+g(y_2)\chi_\Omega\big)\dd y-\frac32\,r^{-4}\int_{\p B_r}u^2\dd\mathcal H^1 .
\]
Inserting $\xi(y):=y\,\varphi_\varepsilon(|y|)$ into \eqref{eq:domvar}, where $\varphi_\varepsilon$ is
smooth and decreases to $\mathbf 1_{[0,r]}$, and using $\mathrm{div}(g\,y)=3g+(y_2g'-g)$, we obtain
for almost every $r$
\[
r\int_{\p B_r}\big(|\nabla u|^2+g\chi_\Omega\big)=2r\int_{\p B_r}(\p_\nu u)^2
+\int_{B_r}\big(3g+y_2g'-g\big)\chi_\Omega .
\]
Together with $\int_{B_r}|\nabla u|^2=\int_{\p B_r}u\,\p_\nu u$ (Step~2) and
$\frac{\mathrm d}{\mathrm dr}\int_{\p B_r}u^2=\frac1r\int_{\p B_r}u^2+2\int_{\p B_r}u\,\p_\nu u$, a
direct computation, as in \cite{Weiss1999,VarvarucaWeiss2011}, gives
\begin{equation}\label{eq:weiss}
W'(r)=2r^{-3}\int_{\p B_r}\Big(\p_\nu u-\frac{3u}{2r}\Big)^2\dd\mathcal H^1
+r^{-4}\int_{B_r}\big(y_2g'(y_2)-g(y_2)\big)\chi_\Omega\dd y .
\end{equation}
Since $|y_2g'-g|=|y_2e'-e|\le2C_1r^2$ on $B_r$, the last term is at least $-2\pi C_1$, so
$W(r)+2\pi C_1r$ is non-decreasing. By Steps~1 and~3, $W$ is bounded, so $W(0+):=\lim_{r\downarrow0}W(r)$
exists.

\emph{Step 5: blow-ups.} For $\rho>0$ let $u_\rho(y):=\rho^{-3/2}u(\rho y)$ and
$\zeta_\rho(t):=\rho^{-1}(\eta_\infty(\rho t)-L)$. By Steps~1 and~3, $u_\rho\le K^{1/2}|y|^{3/2}$ and
$|\nabla u_\rho(y)|^2\le K|y|$, and the $\zeta_\rho$ are even, vanish at $0$, are non-increasing on
$[0,\infty)$ and Lipschitz with constant $\tan\theta_*$. Hence every sequence $\rho_j\downarrow0$ has a
subsequence along which $u_{\rho_j}\to u_0$ and $\zeta_{\rho_j}\to\zeta_0$ locally uniformly. The limit
$u_0\ge0$ is harmonic in the connected set $\Omega_0:=\{y_2<\zeta_0(y_1)\}$ and vanishes outside it;
as in Step~2, $\nabla u_{\rho_j}\to\nabla u_0$ in $L^2_{\rm loc}$, and since $\rho^{-1}g(\rho y_2)\to
-\mu y_2$ together with its derivative, $u_0$ satisfies \eqref{eq:domvar} with $D_\infty$,
$\Psi_\infty$ and $G_{\alpha_\infty}(s)$ replaced by $\Omega_0$, $u_0$ and $-\mu y_2$. By scaling,
$W(\rho r)$ is the energy of $u_\rho$ at radius $r$, computed with $\rho^{-1}g(\rho y_2)$; integrating
\eqref{eq:weiss} from $\rho_jr_1$ to $\rho_jr_2$ and letting $j\to\infty$ we get
$\int_{B_{r_2}\setminus B_{r_1}}|y|^{-3}\big(\p_\nu u_0-\frac{3u_0}{2|y|}\big)^2\dd y=0$, so $u_0$ is
homogeneous of degree $3/2$, and
\[
W(0+)=\int_{B_1}\big(|\nabla u_0|^2-\mu y_2\chi_{\Omega_0}\big)\dd y-\frac32\int_{\p B_1}u_0^2\dd\mathcal H^1 .
\]
If $u_0\equiv0$, \eqref{eq:domvar} for $u_0$ reduces to $\int\chi_{\Omega_0}\mathrm{div}(-\mu y_2\xi)=0$;
taking $\xi$ supported in $\{y_2<0\}$ we see that $\chi_{\Omega_0}$ is constant on $\{y_2<0\}$. Since
$-\tan\theta_*|y_1|\le\zeta_0(y_1)\le0$, this forces $\zeta_0\equiv0$, and then
$W(0+)=-\mu\int_{B_1\cap\{y_2<0\}}y_2\dd y=2\mu/3$. If $u_0\not\equiv0$, then $u_0>0$ in $\Omega_0$, so
$\Omega_0=\{u_0>0\}$ is a cone and $\zeta_0(t)=-a|t|$ for some $a\ge0$. Writing
$u_0=|y|^{3/2}f(\phi)$ we get $f''+\frac94f=0$, $f>0$ on $(-\phi_0,\phi_0)$ and $f(\pm\phi_0)=0$, where
$\tan\phi_0=1/a$; hence $\phi_0=\pi/3$, $a=1/\sqrt3$ and $f=c\cos(3\phi/2)$, and \eqref{eq:domvar} near
the two rays, where $u_0$ is smooth, gives $|\nabla u_0|^2=-\mu y_2$ there, that is, $\frac94c^2=\frac\mu2$.
In this case $\int_{B_1}|\nabla u_0|^2=\int_{\p B_1}u_0\p_\nu u_0=\frac32\int_{\p B_1}u_0^2$, so
$W(0+)=-\mu\int_{B_1\cap\Omega_0}y_2\dd y=\mu/\sqrt3$.

\emph{Step 6: conclusion.} By \eqref{eq:limsupslope} there are $\rho_j\downarrow0$ with
$\zeta_{\rho_j}(1)\to-\limsup_{\Theta\downarrow0}(L-\eta_\infty(\Theta))/\Theta\le-1/\sqrt3$. A blow-up
along a subsequence has $\zeta_0(1)<0$, so it is not flat, and Step~5 gives $W(0+)=\mu/\sqrt3$. Since
$\mu/\sqrt3\ne2\mu/3$, no blow-up is flat, so every blow-up is the one found in Step~5 for
$u_0\not\equiv0$. As this limit does not depend on the sequence, $u_\rho$ converges to it and
$\zeta_\rho(t)\to-|t|/\sqrt3$, locally uniformly as $\rho\downarrow0$, which is (c) at the crest. The convergence $\zeta_\rho(1)\to-1/\sqrt3$ gives
$\eta_\infty(\Theta)<L$ for every sufficiently small $\Theta>0$.
Since $\eta_\infty$ is non-increasing on $[0,\pi/k]$, the same strict
inequality holds for all $0<\Theta\le\pi/k$, and symmetry gives it for
$0<|\Theta|\le\pi/k$.

\emph{Step 7: regularity away from the crests.} Let $0<\Theta_0<\pi/k$. By (c) there is $\varepsilon>0$
with $\eta_\infty\le L-3\varepsilon$ on $[\Theta_0,2\pi/k-\Theta_0]$, so $\eta_n\le L-2\varepsilon$
there for large $n$. We write $Q_n:=|\widehat\nabla\Psi_n|$. By (N5), $\p_\Theta\ln Q_n=\mathfrak a\ge0$
in $D_+$, and $Q_n$ is non-increasing along the crest line (the proof of Lemma~\ref{lem:Tm}(e) with
$m=0$). Since the free surface is non-increasing on $[0,\pi/k]$, every point
$(\Theta,s)\in\overline{D_n}$ with $0\le\Theta\le\pi/k$ and $s\le L-\varepsilon$ can be joined to
$(0,s)$ by a horizontal segment in $\overline{D_n}$, and therefore
$Q_n(\Theta,s)\ge Q_n(0,s)\ge Q_n(0,L-\varepsilon)$. As $n\to\infty$,
$Q_n(0,L-\varepsilon)\to|\widehat\nabla\Psi_\infty(0,L-\varepsilon)|>0$, because $\p_s\Psi_\infty$ is
harmonic and non-positive in $D_\infty$, hence negative there. By symmetry and
Lemma~\ref{lem:speedbound}, $c\le Q_n\le C$ on $\overline{D_n}\cap\{s\le L-\varepsilon\}$ for large
$n$, with $c,C>0$ independent of $n$. In this set the functions $\mathfrak q_n=\ln Q_n$ are therefore
uniformly bounded, $|\chi_n|\le\theta_*$, and on the free surface over
$[\Theta_0,2\pi/k-\Theta_0]$ the boundary values $\frac12\ln G_{\alpha(s_n)}(\eta_n)$ of $\mathfrak q_n$
are uniformly Lipschitz. The estimates in the proof of Lemma~\ref{lem:apriori} (the boundary Hölder
estimate, the interior gradient bound for $\chi_n$, and the boundary Schauder estimates) are local,
and applied in this set they give bounds for $\eta_n$ in $C^{j,\beta}([\Theta_1,2\pi/k-\Theta_1])$ for
every $j$ and every $\Theta_1>\Theta_0$, independent of $n$. Since $\Theta_0$ was arbitrary, the first
part of (d) follows from the Arzelà--Ascoli theorem and the convergence of $\Psi_n$ and $\alpha(s_n)$.
The lower bound $Q_n\ge c$ passes to the limit on $\overline{D_\infty}\cap\{s\le L-\varepsilon\}$, and by
(c) the only points of $\overline{D_\infty}$ with $s=L$ are the crests.

\emph{Step 8: the trough.} Step~1 gives $\alpha_\infty\ge\alpha_0$ and $t:=\eta_\infty(\pi/k)\le\gamma$. Let
$t_n:=\eta_{\min}(s_n)\to t$. By Step~7, the second derivatives of $\Psi_n$ are bounded near the trough
uniformly in $n$; since $\Psi_n$ and $\p_\Theta\Psi_n$ vanish at the trough, $\Psi_n(\pi/k+y,t_n)\le
Cy^2$. Hence we may pass to the limit in the trough identity \eqref{eq:troughid} for $H(s_n)$ and get $G_{\alpha_\infty}(t)^{1/2}\,t-1=t\int_\R K_t(y)\,\Psi_\infty\Big(\frac\pi k+y,t\Big)\dd y>0,$ because $\Psi_\infty(\cdot,t)\ge0$ and $\Psi_\infty>0$ in $D_\infty$, which contains points of the line
$s=t$ near the crest. Thus $g_{\alpha_\infty}(t)<0$, whereas $g_{\alpha_\infty}(L)>0$ because
$G_{\alpha_\infty}(L)=0$. The proof of Proposition~\ref{prop:heights}, with $t$ and $L$ in place of
$\eta_{\min}$ and $\eta_{\max}$, gives $\alpha_\infty\ne\alpha_0$, so $\alpha_\infty>\alpha_0$, and
Lemma~\ref{lem:galpha}(b) then gives $x_1(\alpha_\infty)<t<\gamma$.
\end{proof}

\begin{corollary}
\label{cor:limiting-separation}
Under the hypotheses of Theorem~\ref{thm:corner}, every limiting profile
satisfies
\[
\eta_\infty'(\Theta)<0\qquad(0<\Theta<\pi/k).
\]
There are constants $\varepsilon_\alpha,\varepsilon_h>0$, depending on
the fixed branch, such that, for every $s>0$,
\[
\alpha(s)\ge\alpha_0+\varepsilon_\alpha,
\qquad
\eta_{\min}(s)\ge x_1(\alpha(s))+\varepsilon_h.
\]
No explicit lower bound for either separation constant is claimed.
\end{corollary}

\begin{proof}
Let $v=-\partial_\Theta\Psi_\infty$ in
$D_{\infty,+}=D_\infty\cap\{0<\Theta<\pi/k\}$.
Interior convergence gives $v\ge0$, and $v$ is harmonic.
It is not identically zero: otherwise $\Psi_\infty$ would depend only on
$s$, hence be affine in $s$, and its zero boundary value would force the
profile to be constant, which contradicts Theorem~\ref{thm:corner}(a),(c).
Thus $v>0$ in $D_{\infty,+}$ by the strong maximum principle.
At every regular point of the free surface with $0<\Theta<\pi/k$,
the calculation in Lemma~\ref{lem:robin} applies and gives $\partial_n v=(g_*(\eta_\infty)\cos\chi-\varkappa)v.$ If $v$ vanished there, Hopf's lemma would give $\partial_n v<0$,
contradicting this identity. Hence $v>0$ on that surface arc.
Since $\Psi_{\infty,s}<0$ there by Theorem~\ref{thm:corner}(d) and
the uniform angle bound, differentiation of the boundary condition gives
$\eta_\infty'=v/\Psi_{\infty,s}<0$.

Suppose that no $\varepsilon_\alpha>0$ existed. There would be a sequence
$s_n>0$ with $\alpha(s_n)\to\alpha_0$. The local bifurcation limit is
$\alpha_k>\alpha_0$, so $s_n$ cannot tend to $0$. Continuity and
Theorem~\ref{thm:global}(iii) exclude a subsequence converging to a finite
positive parameter. A subsequence therefore tends to $+\infty$, and
Theorem~\ref{thm:corner}(a) gives a further subsequence with limit
$\alpha_\infty>\alpha_0$, a contradiction.

The height separation is proved in the same way. The function
$\eta_{\min}(s)-x_1(\alpha(s))$ is continuous and positive for $s>0$;
its limit as $s\downarrow0$ is $\gamma-x_1(\alpha_k)>0$.
Any sequence on which it tends to zero must consequently escape to
$+\infty$. Theorem~\ref{thm:corner}(a), and the continuity of $x_1$ at
$\alpha_\infty>\alpha_0$, exclude that possibility.
\end{proof}

The argument identifies the blow-up at each limiting crest. Uniqueness of the limiting wave along the entire global branch remains open. A double zero of the Bernoulli head would require a different nondegeneracy argument.

\section{Global branches with constant relative vorticity}\label{sec:kappaglobal}
\subsection{Laminar comparison and the nodal class}
We now extend Sections~\ref{sec:BT}--\ref{sec:globalthm} to every admissible triple
$(\gamma,\kappa,w)$. Throughout this section $\eps=0$, $k\ge1$, and
\[
\varpi:=\frac{2\kappa}{p_0},\qquad\text{so that}\qquad\widehat\Delta\Psi=-\varpi\,\e^{2s}
\]
by \eqref{eq:flatpoisson}. By \eqref{eq:laminar} and \eqref{eq:admissible}, $p_0$, $W_s$ and $w$ have
the same sign, so $\varpi$ has the sign of $\nu=\kappa/w$; $p_0$ may have either sign. The sets
$\mathcal V_k$ and $\mathcal N_k$, the flat variables of Section~\ref{sec:comparison} and the
Bernoulli head $G_\alpha$ are defined as before, and \eqref{eq:flat} and \eqref{eq:bern-flat} hold for
every $\kappa$, by Lemma~\ref{lem:logpolar} and \eqref{eq:bern-phys}. What changes is that $\Psi$ is
no longer harmonic. We use instead the harmonic function $\Phi:=\Psi+\tfrac14\varpi\e^{2s}$ and,
in place of the holomorphic function of the proof of Lemma~\ref{lem:apriori}, a solution of a
$\bar\partial$-equation. We shall see that the global theorem holds for every admissible triple. The
sign of $\varpi$ decides where the flow can degenerate, and $|\kappa|\le1$ decides whether the
conjugate flows have the structure of Lemma~\ref{lem:galpha}.

\emph{Laminar flows of depth $x$.} For $x>0$ let
\begin{gather*}
\mu(x):=\frac{\kappa(\e^{2x}-1)/2-p_0}{x},\qquad W_s(x):=\kappa\e^{x}-\mu(x)\e^{-x},\\
\bar\Psi_x(s):=\frac{\kappa(\e^{2x}-\e^{2s})/2-\mu(x)(x-s)}{p_0}.
\end{gather*}
This is the laminar flow \eqref{eq:laminar} with depth $x$, flux $p_0$ and vorticity $2\kappa$:
\[
\bar\Psi_x''=-\varpi\e^{2s},\qquad\bar\Psi_x(0)=1,\qquad\bar\Psi_x(x)=0,\qquad
\bar\Psi_x'(x)=-\frac{\e^{x}W_s(x)}{p_0},\qquad \bar\Psi_x'(0)=\frac{\mu(x)-\kappa}{p_0},
\]
and $(\mu(\gamma),W_s(\gamma),\bar\Psi_\gamma)=(\mu,W_s,\Psi_0)$. We call the depth $x$
\emph{admissible} if $\bar\Psi_x'<0$ on $[0,x]$. Then $\bar\Psi_x(\bar H_x(p))=-p$ defines a
$q$-independent solution $\bar H_x$ of the interior equation in \eqref{eq:DJ} with $\bar H_x(-1)=0$,
$\bar H_x(0)=x$ and $\p_p\bar H_x>0$. Conversely, if $H(p)$ solves
$H_{pp}+\varpi\e^{2H}H_p^3=0$ with $H(-1)=0$ and $H_p>0$, then $\Psi:=-H^{-1}$ solves
$\Psi''=-\varpi\e^{2s}$, $\Psi(0)=1$, $\Psi(H(0))=0$, so $H=\bar H_x$ with $x=H(0)$. For $\kappa=0$,
$\bar\Psi_x=1-s/x$ and $\bar H_x=x(1+p)$. Linearising $\Fop(\alpha,\cdot)$ at $\bar H_x$ gives the
operator \eqref{eq:L} with $(\gamma,\cc,h_1,H_0)$ replaced by $(x,W_s(x)^2,\p_p\bar H_x(0),\bar H_x)$,
because $p_0^2\e^{-2x}(\p_p\bar H_x(0))^{-2}=W_s(x)^2$. By Proposition~\ref{prop:dispersion}, its
critical values are
\begin{equation}\label{eq:kalphax}
\begin{aligned}
\alpha_j(x)&:=\e^{x}-2\kappa W_s(x)+\e^{-x}W_s(x)^2\big(1+j\coth(jx)\big)\quad(j\ge1),\\
\alpha_0(x)&:=\e^{x}-2\kappa W_s(x)+\e^{-x}W_s(x)^2\big(1+x^{-1}\big).
\end{aligned}
\end{equation}
For $\kappa=0$ these are the functions $\alpha_k(x)$, $\alpha_0(x)$ of Section~\ref{sec:global}, and
$\alpha_j(\gamma)=\alpha_j$ for every $\kappa$.

\begin{lemma}\label{lem:kstructure}
Let $(\gamma,\kappa,w)$ be admissible. Lemmas~\ref{lem:fredholm-global}--\ref{lem:noreturn} hold,
with the laminar flow $\bar x(1+p)$ of Lemmas~\ref{lem:Nclosed} and~\ref{lem:laminarlimit} replaced by
$\bar H_{\bar x}$, and $\alpha_k(\bar x)$ given by \eqref{eq:kalphax}. In particular,
Theorem~\ref{thm:BT} applies at $(\alpha_k,H_0)$ as in Section~\ref{sec:BT}.
\end{lemma}

\begin{proof}
\emph{Lemmas~\ref{lem:fredholm-global} and~\ref{lem:proper}.} The vorticity adds to $\mathcal A_H$ the
lower-order terms $\varpi\e^{2H}H_p^2(3\psi_p+2H_p\psi)$. By \eqref{eq:Bform}, the boundary operator
keeps the form $\mathcal B_{\alpha,H}$, now with $b_p=-2\cc h_1^2(1+H_q^2)H_p^{-3}\e^{-2(H-\gamma)}<0$.
The proofs use only the ellipticity, the obliqueness and the isomorphism $\Mop$ of
Lemma~\ref{lem:fredholm}, so they apply verbatim. Hypotheses (H1) and (H2) hold by
Lemma~\ref{lem:range}.

\emph{Lemmas~\ref{lem:Nopen} and~\ref{lem:Nlocal}.} The first is independent of the equation. For the
second, $\varphi_k=\sigma_k(p)\cos kq$ with $\sigma_k=H_{0,p}\sinh(kH_0)$
(Proposition~\ref{prop:dispersion}), so that $\sigma_k>0$ on $(-1,0]$, $\sigma_k(-1)=0$ and
$\sigma_k'(-1)=kH_{0,p}(-1)^2>0$; the proof applies verbatim.

\emph{Lemma~\ref{lem:Nclosed}.} Differentiating the interior equation in $q$ now gives
$\mathcal A_Hw+\varpi\e^{2H}H_p^2(3w_p+2H_pw)=0$ for $w=H_q$, with a zeroth-order term of no definite
sign. We apply the argument to $V:=w/H_p$ instead, which is $\Psi_\Theta$ at corresponding points by
\eqref{eq:flat}. For a function $u(\Theta,s)=U(q,p)$, the chain rule and $p=-\Psi$ give
$\widehat\Delta u=U_{qq}-2H_qH_p^{-1}U_{qp}+(1+H_q^2)H_p^{-2}U_{pp}-(\widehat\Delta\Psi)U_p$. Since
$\Psi_\Theta$ is harmonic by \eqref{eq:flatpoisson}, $(1+H_q^2)V_{pp}-2H_qH_pV_{qp}+H_p^2V_{qq}+\varpi\e^{2H}H_p^2V_p=0\qquad\text{in }\Omega,$ an equation with the principal part of $\mathcal A_H$ and without zeroth-order term. On $p=0$, $V$
satisfies the oblique condition $b_pH_pV_p+b_qH_pV_q+(b_pH_{pp}+b_qH_{qp}+b_0H_p)V=0$, and $V=0$ on
$p=-1$. Moreover $V$ has the sign of $w$ and vanishes where $w$ does. At the points where the proof
of Lemma~\ref{lem:Nclosed} examines $w_q$, $w_p$ or $w_{qp}$, namely zeros of $w$ on the boundary and
corners where in addition $w_q=w_p=0$, these derivatives of $V$ equal those of $w$ divided by $H_p$.
The proof therefore applies to $V$. In the laminar case $H_q\equiv0$, so $H_{pp}+\varpi\e^{2H}H_p^3=0$,
$H(-1)=0$, $H_p>0$, and $H=\bar H_{\bar x}$ with $\bar x=H(0)$.

\emph{Lemma~\ref{lem:laminarlimit}.} The limit $\bar v$ solves the problem linearised at
$\bar H_{\bar x}$, which is \eqref{eq:L} at depth $\bar x$, as noted above. As in the proof of
Proposition~\ref{prop:dispersion}, its kernel in the odd $2\pi/k$-periodic functions is spanned by
$\p_p\bar H_{\bar x}\sinh(jk\bar H_{\bar x})\sin(jkq)$ for the single $j$, if any, with
$\bar\alpha=\alpha_{jk}(\bar x)$. The factor in front of $\sin(jkq)$ is positive on $(-1,0]$, so the sign
argument of Lemma~\ref{lem:laminarlimit} gives $j=1$.

\emph{Lemma~\ref{lem:noreturn}} uses only Theorems~\ref{thm:BT}, \ref{thm:local} and
Lemma~\ref{lem:Nlocal}.
\end{proof}

\emph{Conjugate flows.} The laminar flow of depth $x$ satisfies \eqref{eq:bern-flat} with
$|\widehat\nabla\Psi|=|\bar\Psi_x'(x)|$ if and only if $g_\alpha(x)=0$, where now
\begin{equation}\label{eq:kgalpha}
g_\alpha(x):=W_s(x)^2+2\alpha\e^x-\e^{2x}-2\lambda(\alpha)=W_s(x)^2-p_0^2\e^{-2x}G_\alpha(x).
\end{equation}
For $\kappa=0$ this is \eqref{eq:galpha}. Let
\[
x_1(\alpha):=\min\{x>0:\ g_\alpha(x)=0\},\qquad \Upsilon(x):=\frac{x\,(9x^2+12x+28)}{3(x+1)(3x^2+4x+4)} .
\]

\begin{lemma}\label{lem:kgalpha}
\leavevmode
\begin{enumerate}[label=\textup{(\alph*)}]
\item $g_\alpha(\gamma)=0$ for every $\alpha$, $g_\alpha'(x)=2\e^x\big(\alpha-\alpha_0(x)\big)$, and
$\alpha_0(\gamma)=\alpha_0$. As $x\to0$, $\alpha_0(x)\to\infty$, and $g_\alpha(x)\to\infty$ uniformly
for $\alpha$ in bounded sets. In particular $x_1(\alpha)\in(0,\gamma]$ is well defined, and
$x_*(A):=\inf_{|\alpha|\le A}x_1(\alpha)>0$ for every $A>0$.
\item $\alpha_0''(x)\ge\e^x\big(1-\kappa^2\Upsilon(x)\big)$ for $x>0$, and $0<\Upsilon<1$. If
$|\kappa|\le1$, then $\alpha_0$ is strictly convex on $(0,\infty)$, $\alpha_0(x)\to\infty$ and
$g_\alpha(x)\to-\infty$ as $x\to\infty$, and Lemma~\ref{lem:galpha} holds verbatim.
\item If $\varpi<0$, then $W_s(x)/p_0$ is strictly decreasing on $(0,\infty)$ and has a single zero
$x_c>\gamma$; the depth $x$ is admissible if and only if $x<x_c$.
\item If $\kappa/w<0$, then $\alpha_0''(x)>0$ on $(0,x_c)$,
independently of $|\kappa|$.
\end{enumerate}
\end{lemma}

\begin{proof}
(a) From $(x\mu(x))'=\kappa\e^{2x}$ one gets $W_s'(x)=2\kappa\e^x-(1+x^{-1})W_s(x)$, and the formula
for $g_\alpha'$ follows. $g_\alpha(\gamma)=0$ is \eqref{eq:lambda}. As $x\downarrow0$, direct expansion gives
\[
 W_s(x)=\frac{p_0}{x}+O(1),\qquad
 \alpha_0(x)=\frac{p_0^2}{x^3}+O(x^{-2}),\qquad
 g_\alpha(x)=\frac{p_0^2}{x^2}+O(x^{-1}),
\]
where the last error is uniform for $\alpha$ in bounded sets. Since $p_0\ne0$,
these estimates give the asserted limits and the uniform lower bound for
$x_1(\alpha)$.

(b) Since $\kappa\e^x-W_s(x)=\mu(x)\e^{-x}$, $\alpha_0(x)=(1-\kappa^2)\e^x+\mu(x)^2\e^{-3x}+x^{-1}\e^{-x}W_s(x)^2 .$ Set $P(x):=(x+1)(3x^2+4x+4),\qquad N(x):=9x^2+16x+12.$ Twice differentiating $\alpha_0$, using
$W_s'=2\kappa\e^x-(1+x^{-1})W_s$, gives the explicit identity
\begin{equation}\label{eq:kalpha-second}
\alpha_0''(x)=\e^x+\frac{4(2x+3)}x\kappa^2\e^x
 -\frac{2N(x)}{x^2}\kappa W_s(x)
 +\frac{3P(x)}{x^3}\e^{-x}W_s(x)^2.
\end{equation}
Completing the square yields
\begin{equation}\label{eq:kalpha-square}
\alpha_0''(x)=\e^x\bigl(1-\kappa^2\Upsilon(x)\bigr)
 +\frac{3P(x)}{x^3}\e^{-x}
 \left(W_s(x)-\frac{\kappa\e^x xN(x)}{3P(x)}\right)^2.
\end{equation}
Indeed, the remaining coefficient is determined by
$12(2x+3)P(x)-N(x)^2=-x^2(9x^2+12x+28)$.
This proves the lower bound without a computer-algebra certificate. Clearly $\Upsilon>0$, and $1-\Upsilon=(9x^2-4x+12)/\big(3(x+1)(3x^2+4x+4)\big)>0$. Let $|\kappa|\le1$. Then
$\alpha_0''>0$. As $x\to\infty$, $W_s(x)=\kappa\e^x(1+o(1))$ if $\kappa\ne0$, so the expression for $\alpha_0$ at the beginning of (b) gives
$\alpha_0(x)\to\infty$, both for $|\kappa|<1$ and for $|\kappa|=1$; for $\kappa=0$ this is
Lemma~\ref{lem:galpha}(a). Hence $g_\alpha'\le-2\e^x$ for large $x$, and $g_\alpha\to-\infty$. The proof of
Lemma~\ref{lem:galpha}(b) uses only these properties, (a), and $\p_\alpha g_\alpha(x)=2(\e^x-\e^\gamma)$,
which holds for every $\kappa$.

(c) By the symmetry $(p_0,\kappa,w)\mapsto(-p_0,-\kappa,-w)$ we may assume $p_0>0>\kappa$. Then
$W_s(x)=\kappa\big(\e^x-\frac{\sinh x}x\big)+p_0\frac{\e^{-x}}x$, the function $\e^x-\frac{\sinh x}x$ is
increasing (its derivative is $\e^x-\sum_{n\ge1}\frac{2n\,x^{2n-1}}{(2n+1)!}>0$) and $\e^{-x}/x$ is
decreasing. So $W_s$ decreases strictly from $+\infty$ to $-\infty$, and $W_s(\gamma)=W_s>0$. For
$x<x_c$ we have $\bar\Psi_x'(x)<0$, and $\bar\Psi_x'$ is increasing because $\bar\Psi_x''>0$; hence
$\bar\Psi_x'<0$ on $[0,x]$. For $x\ge x_c$, $\bar\Psi_x'(x)\ge0$.

(d) For $0<x<x_c$, $W_s(x)$ has the sign of $p_0$, whereas $\kappa$ has the
opposite sign. Hence $\kappa W_s(x)<0$, and every term on the right-hand
side of \eqref{eq:kalpha-second} is non-negative, with the first strictly
positive.
\end{proof}

For $\kappa=0$, (b) contains Lemma~\ref{lem:galpha}(a). For $|\kappa|>1$ the bound in (b) is negative for
large $x$; then $\alpha_0(x)\to-\infty$ and $g_\alpha(x)\to+\infty$ as $x\to\infty$, and
$\alpha_0$ need not be quasi-convex.

\begin{lemma}\label{lem:kcompare}
Let $(\alpha,H)\in\mathcal V_k$ be a non-laminar solution of $\Fop(\alpha,H)=0$, and assume that
$\varpi\ge0$, or that $\varpi<0$ and $H\in\mathcal N_k$. We write $\hat x:=\eta_{\max}$ and
$\check x:=\eta_{\min}$. Then:
\begin{enumerate}[label=\textup{(\alph*)}]
\item $\bar\Psi_{\hat x}>0$ on $[0,\hat x)$, $\bar\Psi_{\hat x}'(\hat x)<0$, and
$|\widehat\nabla\Psi|<|\bar\Psi_{\hat x}'(\hat x)|$ at every point $(\Theta,\eta_{\max})$ of the free surface;
\item $\bar\Psi_{\check x}<0$ on $(\check x,\hat x]$, $\bar\Psi_{\check x}'(\check x)<0$, and
$|\widehat\nabla\Psi|>|\bar\Psi_{\check x}'(\check x)|$ at every point $(\Theta,\eta_{\min})$;
\item $\bar\Psi_{\check x}'(0)<\Psi_s<\bar\Psi_{\hat x}'(0)$ on the bottom $s=0$.
\end{enumerate}
If $\varpi<0$, then moreover $\eta_{\max}<x_c$. For $\kappa=0$ this is Lemma~\ref{lem:compare}.
\end{lemma}

\begin{proof}
As in the proof of Lemma~\ref{lem:compare}, $\eta_{\min}<\eta_{\max}$. For $x>0$ the function
$Z_x:=\Psi-\bar\Psi_x$ is harmonic in $D$, because $\widehat\Delta\bar\Psi_x=-\varpi\e^{2s}$; it vanishes on
$s=0$ and equals $-\bar\Psi_x(\eta(\Theta))$ on the free surface.

\emph{Signs of $\bar\Psi_{\hat x}$ and $\bar\Psi_{\check x}$.} If $\varpi\ge0$, then $\bar\Psi_x$ is
concave with $\bar\Psi_x(0)=1$ and $\bar\Psi_x(x)=0$. Hence $\bar\Psi_x(s)\ge1-s/x>0$ for $s<x$,
$\bar\Psi_x'(x)\le-1/x$, and $\bar\Psi_x(s)\le-(s-x)/x<0$ for $s>x$; this gives the first two
assertions of (a) and of (b). Let $\varpi<0$ and $H\in\mathcal N_k$. By (N1) and \eqref{eq:flat},
$\Psi(\cdot,s)$ is even and strictly decreasing on $(0,\delta)$ for $0<s<\eta_{\max}$ and some
$\delta>0$, so $\Psi_{\Theta\Theta}(0,s)\le0$. The function $f(s):=\Psi(0,s)$ on the crest line
therefore satisfies $f''=-\varpi\e^{2s}-\Psi_{\Theta\Theta}\ge\bar\Psi_{\hat x}''$, $f(0)=1$,
$f(\hat x)=0$ and $f'<0$. Thus $f-\bar\Psi_{\hat x}$ is convex and vanishes at $0$ and $\hat x$, so
$\bar\Psi_{\hat x}\ge f>0$ on $[0,\hat x)$ and $\bar\Psi_{\hat x}'(\hat x)\le f'(\hat x)<0$. Since
$\bar\Psi_{\hat x}'$ is increasing, $\hat x$ is admissible, and $\eta_{\max}<x_c$ by
Lemma~\ref{lem:kgalpha}(c). Next, $\bar\Psi_{\check x}-\bar\Psi_{\hat x}$ is linear and vanishes at
$s=0$, so $\bar\Psi_{\check x}(s)=\bar\Psi_{\hat x}(s)+\beta s$ with
$\beta=-\bar\Psi_{\hat x}(\check x)/\check x<0$. The convex function $\bar\Psi_{\check x}$ vanishes at
$\check x$ and equals $\beta\hat x<0$ at $\hat x$, so it is negative on $(\check x,\hat x]$. If
$\bar\Psi_{\check x}'(\check x)\ge0$, convexity would give $\bar\Psi_{\check x}\ge0$ to the right of
$\check x$; hence $\bar\Psi_{\check x}'(\check x)<0$.

\emph{Maximum principle.} Since $\eta(\Theta)\in[\check x,\hat x]$, the first part gives
$Z_{\hat x}\le0$ and $Z_{\check x}\ge0$ on the free surface, with strict inequalities where
$\check x<\eta<\hat x$. The strong maximum principle gives $Z_{\hat x}<0<Z_{\check x}$ in $D$. As in
the proof of Lemma~\ref{lem:compare}, Hopf's lemma at a crest gives $\Psi_s>\bar\Psi_{\hat x}'(\hat x)$,
at a trough $\Psi_s<\bar\Psi_{\check x}'(\check x)$, and on the bottom (c). At a crest
$|\widehat\nabla\Psi|=|\Psi_s|$ and $\bar\Psi_{\hat x}'(\hat x)<\Psi_s<0$; at a trough
$\Psi_s<\bar\Psi_{\check x}'(\check x)<0$. This proves (a) and (b).
\end{proof}

\begin{proposition}\label{prop:kheights}
Let $(\alpha,H)$ be as in Lemma~\ref{lem:kcompare}. Then:
\begin{enumerate}[label=\textup{(\alph*)}]
\item $g_\alpha(\eta_{\max})>0>g_\alpha(\eta_{\min})$. In particular $\eta_{\max}\ne\gamma\ne\eta_{\min}$
and $\eta_{\min}>x_1(\alpha)$.
\item If $\eta_{\min}<\gamma<\eta_{\max}$, then $\alpha>\alpha_-:=\min_{0<x\le\gamma}\alpha_0(x)$.
\item If $|\kappa|\le1$, the conclusions of Proposition~\ref{prop:heights} hold; there
$x_1(\alpha)$ and $x_3(\alpha)$ are the zeros of \eqref{eq:kgalpha} given by Lemma~\ref{lem:galpha},
which applies by Lemma~\ref{lem:kgalpha}(b).
\end{enumerate}
\end{proposition}

\begin{proof}
(a) At a crest, \eqref{eq:bern-flat} and Lemma~\ref{lem:kcompare}(a) give
$p_0^2\e^{-2\hat x}G_\alpha(\hat x)<p_0^2\e^{-2\hat x}\bar\Psi_{\hat x}'(\hat x)^2=W_s(\hat x)^2$, that is,
$g_\alpha(\hat x)>0$; likewise (b) of that lemma gives $g_\alpha(\check x)<0$. Since $g_\alpha(\gamma)=0$,
$\eta_{\max},\eta_{\min}\ne\gamma$, and $g_\alpha>0$ on $(0,x_1(\alpha))$ gives $\eta_{\min}>x_1(\alpha)$.
(b) Since $g_\alpha(\eta_{\min})<0=g_\alpha(\gamma)$, $g_\alpha'>0$ somewhere in $(\eta_{\min},\gamma)$, so
$\alpha>\alpha_0(x)\ge\alpha_-$ there. (c) The proof of Proposition~\ref{prop:heights} uses only (a) and
the properties of $g_\alpha$ in Lemma~\ref{lem:kgalpha}(b).
\end{proof}

\begin{proposition}
\label{prop:kalpha-adverse}
Let $(\alpha,H)$ satisfy the hypotheses of Lemma~\ref{lem:kcompare}, with
$H\in\mathcal N_k$, $\kappa/w\le0$, and
$\eta_{\min}<\gamma<\eta_{\max}$. Then
\[
 \alpha>\alpha_0
   >\e^\gamma+\cc\e^{-\gamma}>0.
\]
\end{proposition}
\begin{proof}
If $\kappa=0$, take $I=(0,\infty)$; then $\alpha_0$ is strictly
convex on $I$ by Lemma~\ref{lem:galpha}. If $\kappa/w<0$, take
$I=(0,x_c)$; Lemma~\ref{lem:kcompare} places
$[\eta_{\min},\eta_{\max}]$ in $I$, and
\eqref{eq:kalpha-second} proves strict convexity there.
Suppose, to obtain a contradiction, that
$\alpha\le\alpha_0(\gamma)$. If $\alpha_0'(\gamma)\le0$, strict convexity
gives $\alpha_0(x)>\alpha_0(\gamma)\ge\alpha$ for $0<x<\gamma$.
Thus $g_\alpha'=2\e^x(\alpha-\alpha_0(x))<0$ there, and
$g_\alpha(\eta_{\min})>g_\alpha(\gamma)=0$, which contradicts
Proposition~\ref{prop:kheights}(a). If $\alpha_0'(\gamma)>0$, strict
convexity instead gives $\alpha_0(x)>\alpha$ for $x>\gamma$ in $I$.
Then $g_\alpha(\eta_{\max})<g_\alpha(\gamma)=0$, the other contradiction.
Finally, writing $\nu=\kappa/w$, $\alpha_0-\e^\gamma-\cc\e^{-\gamma} =\e^\gamma w^2(\gamma^{-1}-2\nu)>0.$ \end{proof}

For $\kappa/w>0$ and $|\kappa|>1$ the conjugate flows give no upper bound for the crest, since
$g_\alpha(x)\to\infty$ as $x\to\infty$. The flux through the crest line gives one instead. Let
$\phi(u):=1-u+u\ln u$ for $0<u\le1$ and $\phi(0):=1$; $\phi$ decreases from $1$ to $0$.

\begin{lemma}\label{lem:kcrest}
Let $(\alpha,H)\in\mathcal V_k$ be a non-laminar solution of $\Fop(\alpha,H)=0$ with $H\in\mathcal N_k$, and
let $Q_c:=G_\alpha(\eta_{\max})^{1/2}$.
\begin{enumerate}[label=\textup{(\alph*)}]
\item For $0\le\sigma\le\eta_{\max}$,
\[
-\Psi_s(0,\sigma)\ge Q_c-\tfrac12\varpi\big(\e^{2\eta_{\max}}-\e^{2\sigma}\big).
\]
In particular $-\Psi_s\ge Q_c$ on the crest line if $\varpi\le0$.
\item Let $\varpi>0$, $y:=\e^{\eta_{\max}}$ and $u:=1-B_\alpha(y)^{1/2}/(|\kappa|y)$, with $B_\alpha$ as in
\eqref{eq:Bquad}. Then $u<1$ and
\begin{equation}\label{eq:kcrest}
\tfrac14\varpi\big(y^2\phi(u^+)-1-2\ln y\big)\le1,\qquad u^+:=\max\{u,0\}.
\end{equation}
Consequently, for every $A>0$ there is $X(A)$ such that $\eta_{\max}\le X(A)$ whenever $|\alpha|\le A$.
\end{enumerate}
\end{lemma}

\begin{proof}
(a) By (N1) and \eqref{eq:flat}, $\Psi(\cdot,\sigma)$ is even and decreasing on $(0,\delta)$ for some
$\delta>0$, so $\Psi_{\Theta\Theta}(0,\sigma)\le0$ for $0<\sigma<\eta_{\max}$, as in the proof of
Lemma~\ref{lem:kcompare}. Hence $h:=\Psi(0,\cdot)-\bar\Psi_{\eta_{\max}}$ satisfies
$h''=-\Psi_{\Theta\Theta}(0,\cdot)\ge0$, so $h'(\sigma)\le h'(\eta_{\max})$, that is,
$-\Psi_s(0,\sigma)\ge-\Psi_s(0,\eta_{\max})+\bar\Psi_{\eta_{\max}}'(\eta_{\max})-\bar\Psi_{\eta_{\max}}'(\sigma)$. At the crest
$\Psi_\Theta=0$ and $-\Psi_s=Q_c$ by \eqref{eq:bern-flat}, and $\bar\Psi_x''=-\varpi\e^{2s}$.

(b) By \eqref{eq:Bquad}, $Q_c=yB_\alpha(y)^{1/2}/|p_0|>0$, and $|p_0|\varpi=2|\kappa|$; hence
$a:=\frac12\varpi y^2-Q_c=\frac12\varpi y^2u$, and $u<1$. By (a), $-\Psi_s(0,\sigma)\ge\ell(\sigma):=\frac12\varpi\e^{2\sigma}-a$,
and $-\Psi_s>0$. Since $\Psi(0,0)=1$ and $\Psi(0,\eta_{\max})=0$, $1=\int_0^{\eta_{\max}}-\Psi_s(0,\sigma)\dd\sigma\ge\int_0^{\eta_{\max}}\max\{\ell(\sigma),0\}\dd\sigma .$ If $a\le0$, then $\ell(\sigma)\ge\frac12\varpi\e^{2\sigma}$, and the right-hand side is at least
$\frac14\varpi(y^2-1)$. If $a>0$, then $\ell$ vanishes at $\sigma_a:=\frac12\ln(uy^2)<\eta_{\max}$. If
$\sigma_a\le0$, then $\ell\ge0$ on $[0,\eta_{\max}]$, and the right-hand side is
$\frac14\varpi(y^2-1)-a\ln y\ge\frac14\varpi(y^2-1-2\ln y)$, because $uy^2\le1$. If $\sigma_a>0$, the
right-hand side is at least $\int_{\sigma_a}^{\eta_{\max}}\ell\dd\sigma=\frac14\varpi y^2\phi(u)$. In all three cases
\eqref{eq:kcrest} follows, because $0\le\phi\le1$.
Finally, $B_\alpha(y)/y^2=1-2\alpha/y+(\cc+2\alpha\e^\gamma-\e^{2\gamma})/y^2\to1$ as $y\to\infty$, uniformly for
$|\alpha|\le A$. Hence $u^+\to(1-1/|\kappa|)^+<1$, and $\phi(u^+)$ stays bounded away from zero for large $y$,
so the left-hand side of \eqref{eq:kcrest} tends to infinity as $y\to\infty$, uniformly for $|\alpha|\le A$.
\end{proof}

\subsection{Uniform estimates and continuation}
\begin{lemma}\label{lem:kvelbound}
Let $A>0$. There is $C_1=C_1(A)$ such that every non-laminar solution $(\alpha,H)\in\mathcal V_k$ of
$\Fop(\alpha,H)=0$ with $H\in\mathcal N_k$, $|\alpha|\le A$ and $\eta_{\max}\le A$ satisfies
$|\widehat\nabla\Psi|\le C_1$ in $\overline D$, $H_p\ge1/C_1$ and $|H_q|\le C_1H_p$ in $\overline\Omega$.
Moreover, with $\Gamma_\alpha:=\max_{[\eta_{\min},\eta_{\max}]}G_\alpha$,
\begin{equation}\label{eq:kvel}
|\widehat\nabla\Psi|^2\le2\Gamma_\alpha+2\varpi^2\e^{4\eta_{\max}}\quad\text{in }\overline D,
\end{equation}
and if $\varpi\ge0$, then $|\widehat\nabla\Psi|^2\le\Gamma_\alpha$ in $\overline D$, with equality only on
the free surface.
\end{lemma}

\begin{proof}
By Proposition~\ref{prop:kheights}(a) and Lemma~\ref{lem:kgalpha}(a), $\eta_{\min}\ge x_*(A)>0$. The
function $\Phi=\Psi+\frac14\varpi\e^{2s}$ is harmonic, so $|\widehat\nabla\Phi|^2$ is subharmonic and
attains its maximum on $\p D$. On the bottom, $|\widehat\nabla\Phi|=|\Psi_s+\frac12\varpi|$, and
Lemma~\ref{lem:kcompare}(c) gives $|\Psi_s|<|\bar\Psi_{\check x}'(0)|$. Since
$\bar\Psi_{\check x}'(\check x)-\bar\Psi_{\check x}'(0)=-\frac12\varpi(\e^{2\check x}-1)$ and
$\bar\Psi_{\check x}'<0$ on $[0,\check x]$ (Lemma~\ref{lem:kcompare}(b),(c)), Lemma~\ref{lem:kcompare}(b) gives
$|\bar\Psi_{\check x}'(0)|\le|\widehat\nabla\Psi(\pi/k,\eta_{\min})|+\frac12\varpi^-(\e^{2\check x}-1)\le
\Gamma_\alpha^{1/2}+\frac12\varpi^-(\e^{2\check x}-1)$, where $\varpi^\pm:=\max\{\pm\varpi,0\}$. On the free surface
$|\widehat\nabla\Phi|\le G_\alpha(\eta)^{1/2}+\frac12|\varpi|\e^{2\eta}$. Hence
$|\widehat\nabla\Phi|\le\Gamma_\alpha^{1/2}+\frac12|\varpi|\e^{2\eta_{\max}}$ in $\overline D$. From $|\widehat\nabla\Psi|\le|\widehat\nabla\Phi|+\frac12|\varpi|\e^{2s}$ we obtain \eqref{eq:kvel}. Since
$\Gamma_\alpha$ is bounded for $|\alpha|\le A$ and $x_*(A)\le\eta_{\min}<\eta_{\max}\le A$, this gives
$C_1$, and the bounds on $H$ follow from \eqref{eq:flat}.

If $\varpi\ge0$, then $|\widehat\nabla\Psi|^2$ is subharmonic, because
$\widehat\Delta|\widehat\nabla\Psi|^2=2|\widehat\nabla^2\Psi|^2-4\varpi\e^{2s}\Psi_s\ge0$. On the bottom,
$|\Psi_s|<|\bar\Psi_{\check x}'(0)|\le|\bar\Psi_{\check x}'(\check x)|$ by concavity, and
$|\bar\Psi_{\check x}'(\check x)|<|\widehat\nabla\Psi|$ at the trough. By the strong maximum principle the
maximum is therefore attained only on the free surface, where $|\widehat\nabla\Psi|^2=G_\alpha(\eta)$.
\end{proof}

\begin{lemma}\label{lem:kapriori}
Let $M>0$. There is $C_2=C_2(M)$ such that every non-laminar solution $(\alpha,H)\in\mathcal V_k$ of
$\Fop(\alpha,H)=0$ with $H\in\mathcal N_k$, $|\alpha|\le M$ and $\max_{\overline\Omega}H_p\le M$ satisfies
$\|H\|_{C^{3,\beta}(\overline\Omega)}\le C_2$ and $\min_{\overline\Omega}H_p\ge1/C_2$.
\end{lemma}

\begin{proof}
Constants $C$ depend only on $M$. Since $H(q,-1)=0$, $\eta=H(\cdot,0)\le\max H_p\le M$, so
Lemma~\ref{lem:kvelbound} applies with $A=M$, and $x_*\le\eta_{\min}$ with $x_*=x_*(M)>0$. Let
$F:=\Psi_s+\mathrm i\Psi_\Theta$ and $\zeta=\Theta+\mathrm is$, $\p_{\bar\zeta}=\frac12(\p_\Theta+\mathrm i\p_s)$. Then
$M^{-1}\le-\Psi_s\le|F|\le C_1$, so $\mathfrak f:=\log(-F)=\mathfrak q+\mathrm i\chi$ is defined, with
$|\mathfrak q|\le C$, $|\chi|\le\chi_0<\pi/2$, and \eqref{eq:qchi} holds. Now
$\p_{\bar\zeta}F=\frac{\mathrm i}2\widehat\Delta\Psi=-\frac{\mathrm i}2\varpi\e^{2s}$, so
\begin{equation}\label{eq:kdbar}
\p_{\bar\zeta}\mathfrak f=\rho:=-\frac{\mathrm i\varpi\e^{2s}}{2F},\qquad |\rho|\le C .
\end{equation}
For $\varpi=0$, $\mathfrak f$ is holomorphic, and this is the proof of Lemma~\ref{lem:apriori}. In
general we subtract a Cauchy transform. For a bounded open set $U\subset\C$ and $g\in L^\infty(U)$ let
$T_Ug(\zeta):=-\frac1\pi\iint_Ug(z)(z-\zeta)^{-1}\dd A(z)$. Then $\p_{\bar\zeta}T_Ug=g$ in $U$ in the
sense of distributions, $\|T_Ug\|_{C^{\beta_1}(\C)}\le C(\beta_1,U)\|g\|_{L^\infty}$ for every
$\beta_1<1$, and $T_U$ maps $C^{j,\beta}$-functions with compact support in $U$ into
$C^{j+1,\beta}_{\rm loc}(\C)$ \cite[Ch.~I]{Vekua1962}.
These constants are uniform for the varying domains used below: they
are all contained in a fixed rectangle depending only on $M$ and $k$.
For example, splitting the difference of the two kernels into disks
of radius $2r$ and dyadic annuli gives $|T_Ug(z)-T_Ug(z')| \le C\|g\|_{L^\infty}\,r\bigl(1+\log(R/r)\bigr), \qquad r=|z-z'|\le R/2,$ when $U$ lies in a fixed disk of radius $R$. This implies every
$C^{\beta_1}$ bound with $\beta_1<1$, independently of the shape of $U$.

\emph{Hölder bound.} Extend $\mathfrak f$ to $\widetilde D$ by
$\mathfrak f(\Theta,-s):=\overline{\mathfrak f(\Theta,s)}$, as in the proof of
Lemma~\ref{lem:apriori}. Since $\chi=0$ on $s=0$, the extension is continuous, and
$\p_{\bar\zeta}\mathfrak f=\tilde\rho$ in $\widetilde D$ in the sense of distributions, where
$\tilde\rho(\Theta,-s)=\overline{\rho(\Theta,s)}$, so $|\tilde\rho|\le C$. Let
$U:=\widetilde D\cap\{|\Theta|<3\pi/k\}$. Then $\mathfrak h:=\mathfrak f-T_U\tilde\rho$ is holomorphic
in $U$ by Weyl's lemma. On $s=\pm\eta(\Theta)$, $\operatorname{Re}\mathfrak h=\frac12\ln G_\alpha(\eta)-
\operatorname{Re}T_U\tilde\rho$, which is Hölder continuous with exponent $\frac12$ and norm at most
$C$. The Hölder estimate in the proof of Lemma~\ref{lem:apriori}, applied to $\mathfrak h$ on $U$,
bounds $\mathfrak h$ in $C^{\beta_0}(\overline{\widetilde D}\cap\{|\Theta|\le2\pi/k\})$ for some
$\beta_0\in(0,\frac12]$. Hence $\|\mathfrak f\|_{C^{\beta_0}}\le C$, and $\eta'=-\tan\chi(\Theta,\eta(\Theta))$
gives $\|\eta\|_{C^{1,\beta_0}}\le C$.

\emph{Bootstrap.} Near the bottom, $\Psi-1$ vanishes on $s=0$ and solves
$\widehat\Delta(\Psi-1)=-\varpi\e^{2s}$ in $\{0<s<x_*\}\subset D$, and $|\Psi|\le1$; the boundary
Schauder estimates \cite[Section~6.2]{GT2001} bound $\Psi$ in $C^{5,\beta_0}(\{0\le s\le\frac34x_*\})$.
Let $D':=D\cap\{s>\frac12x_*,\ |\Theta|<3\pi/k\}$. Suppose that $\|\eta\|_{C^{j+1,\beta_0}}\le C$ and
$\|\mathfrak f\|_{C^{j,\beta_0}(\overline{D'})}\le C$ for some $j\ge0$. Then $\rho=\frac{\mathrm i}2\varpi\e^{2s}\e^{-\mathfrak f}$ is
bounded in $C^{j,\beta_0}(\overline{D'})$. Since the free surface and the line $s=\frac12x_*$ are
$C^{j+1,\beta_0}$-curves, there is $\rho_e\in C^{j,\beta_0}(\C)$ with norm at most $C$, support in a fixed
bounded open set $U'$, and $\rho_e=\rho$ in $D'':=D'\cap\{|\Theta|<5\pi/(2k)\}$. Then
$\mathfrak h:=\mathfrak f-T_{U'}\rho_e$ is holomorphic in $D''$. Its real part has
boundary values $\frac12\ln G_\alpha(\eta)-\operatorname{Re}T_{U'}\rho_e$ on the free surface and is
bounded in $C^{j+1,\beta_0}$ on the line $s=\frac12x_*$ by the bottom estimate. Local boundary Schauder
estimates for the Dirichlet problem near the top and bottom boundaries,
away from the artificial vertical edges
(\cite[Corollary~8.36]{GT2001} for $j=0$, \cite[Section~6.2]{GT2001} for $j\ge1$) bound
$\operatorname{Re}\mathfrak h$ in $C^{j+1,\beta_0}$ on $D''\cap\{|\Theta|\le2\pi/k\}$, and by the
Cauchy--Riemann equations so is $\mathfrak h$. Hence
$\|\mathfrak f\|_{C^{j+1,\beta_0}}\le C$ there, and $\eta'=-\tan\chi(\Theta,\eta(\Theta))$ gives
$\|\eta\|_{C^{j+2,\beta_0}}\le C$. Periodicity supplies the same
bounds on translated patches before the next bootstrap step. After four steps,
$\Psi_s=-\e^{\mathfrak q}\cos\chi$ and $\Psi_\Theta=-\e^{\mathfrak q}\sin\chi$ are bounded in
$C^{4,\beta_0}$ near the free surface, and with the bottom estimate
$\|\Psi\|_{C^{5,\beta_0}(\overline D)}\le C$. As in the proof of Lemma~\ref{lem:apriori}, this bounds
$\|H\|_{C^4(\overline\Omega)}$, and $H_p\ge1/C_1$ by Lemma~\ref{lem:kvelbound}.
\end{proof}

\begin{lemma}\label{lem:kflowforce}
Let $(\alpha,H)\in\mathcal V_k$ solve $\Fop(\alpha,H)=0$.
\begin{enumerate}[label=\textup{(\alph*)}]
\item The flow force
$\mathcal S(\Theta):=\int_0^{\eta(\Theta)}\big(\Psi_s^2-\Psi_\Theta^2-2\varpi\e^{2s}\Psi\big)\dd s+\int_\gamma^{\eta(\Theta)}G_\alpha(x)\dd x$
is independent of $\Theta$.
\item If $H\in\mathcal N_k$ is not laminar, then, with $\Gamma_\alpha$ as in Lemma~\ref{lem:kvelbound},
\begin{equation}\label{eq:kffineq}
\frac1{\eta_{\max}}-\varpi^+\big(\e^{2\eta_{\max}}-1\big)+\int_{\eta_{\min}}^{\eta_{\max}}G_\alpha\dd x\le
\eta_{\min}\big(2\Gamma_\alpha+2\varpi^2\e^{4\eta_{\max}}\big)+\varpi^-\big(\e^{2\eta_{\min}}-1\big).
\end{equation}
\item If $H\in\mathcal N_k$, then $v=-\Psi_\Theta$ is positive in $D_+\cup\Gamma_+$, vanishes on the
rest of $\p D_+$, and, with the notation of Lemma~\ref{lem:robin} and $Q:=|\widehat\nabla\Psi|$,
\begin{equation}\label{eq:krobin}
\p_nv=\Big(g_*(\eta)\cos\chi-\varkappa+\frac{\varpi\e^{2\eta}}Q\Big)v\qquad\text{on }\Gamma_+ .
\end{equation}
For $\varphi$ as in Lemma~\ref{lem:robin},
\begin{equation}\label{eq:keigineq}
\int_0^{\pi/k}g_*(\eta)\varphi_\Gamma^2\dd\Theta+\varpi\int_0^{\pi/k}\frac{\e^{2\eta}\varphi_\Gamma^2}{|\Psi_s|}\dd\Theta
\le\int_{D_+}|\widehat\nabla\varphi|^2\dd\Theta\dd s+\frac\pi2\int_0^{\pi/k}\big|(\varphi_\Gamma^2)'\big|\dd\Theta,
\end{equation}
where $\Psi_s$ is evaluated on the free surface.
\end{enumerate}
\end{lemma}

\begin{proof}
(a) By \eqref{eq:flatpoisson},
$\p_\Theta\big(\Psi_s^2-\Psi_\Theta^2-2\varpi\e^{2s}\Psi\big)=2\p_s(\Psi_\Theta\Psi_s)$. Since $\Psi=0$ on the
free surface, the extra term does not contribute to the boundary terms, and the computation of
Lemma~\ref{lem:flowforce} applies.

(b) On the crest line $\int_0^{\eta_{\max}}\Psi_s^2\ge1/\eta_{\max}$, and
$-2\varpi\e^{2s}\Psi\ge-2\varpi^+\e^{2s}$ because $0\le\Psi\le1$, so
$\mathcal S(0)\ge1/\eta_{\max}-\varpi^+(\e^{2\eta_{\max}}-1)+\int_\gamma^{\eta_{\max}}G_\alpha$. On the trough line,
$-2\varpi\e^{2s}\Psi\le2\varpi^-\e^{2s}$, and \eqref{eq:kvel} bounds $\Psi_s^2$. Hence
$\mathcal S(\pi/k)\le\eta_{\min}(2\Gamma_\alpha+2\varpi^2\e^{4\eta_{\max}})+\varpi^-(\e^{2\eta_{\min}}-1)
+\int_\gamma^{\eta_{\min}}G_\alpha$, and (a) gives \eqref{eq:kffineq}. Every level in $[\eta_{\min},\eta_{\max}]$ is
attained on the free surface, so $G_\alpha>0$ on this interval.

(c) $\Psi_\Theta$ is harmonic, and the positivity and the boundary values of $v$ follow from (N1) as in
Lemma~\ref{lem:robin}. The derivation of \eqref{eq:robin1} there uses only \eqref{eq:bern-flat} and
the kinematic condition, so it holds for every $\kappa$. To compute $\p_s(Q^2)$ we use
$\mathfrak f=\mathfrak q+\mathrm i\chi$ and \eqref{eq:kdbar}. For a positively oriented orthonormal
frame $(t,\nu)$ with $t=\e^{\mathrm i\theta}$ we have $\p_t+\mathrm i\p_\nu=2\e^{-\mathrm i\theta}\p_{\bar\zeta}$.
On $\Gamma_+$, $t=-\e^{-\mathrm i\chi}$ and $F=-Q\e^{\mathrm i\chi}$ by \eqref{eq:qchi}, so
$\p_t\mathfrak f+\mathrm i\p_\nu\mathfrak f=-\mathrm i\varpi\e^{2\eta}/Q$. Hence $\p_\nu\chi=\p_t\mathfrak q$
as before, while $\p_\nu\mathfrak q=-\p_t\chi-\varpi\e^{2\eta}/Q=\varkappa-\varpi\e^{2\eta}/Q$. The
computation in the proof of Lemma~\ref{lem:robin} then gives
$\p_s\mathfrak q=-g_*(\eta)\sin^2\chi-(\varkappa-\varpi\e^{2\eta}/Q)\cos\chi$, and \eqref{eq:robin1}
gives \eqref{eq:krobin}. The proof of \eqref{eq:eigineq} applies verbatim and yields
$\int_{D_+}|\widehat\nabla\varphi|^2\ge\int_{\Gamma_+}(g_*\cos\chi-\varkappa+\varpi\e^{2\eta}/Q)\varphi^2\dd\ell$.
Since $\dd\ell=\dd\Theta/\cos\chi$ and $Q\cos\chi=|\Psi_s|$ on $\Gamma_+$ by \eqref{eq:qchi}, this is
\eqref{eq:keigineq}.
\end{proof}

\begin{proposition}\label{prop:kalphabound}
There is $A_*=A_*(\gamma,\kappa,w,k)$ such that $\alpha\le A_*$ for every solution
$(\alpha,H)\in\mathcal V_k$ of $\Fop(\alpha,H)=0$ with $H\in\mathcal N_k$ and
$\eta_{\min}<\gamma<\eta_{\max}$.
\end{proposition}

\begin{proof}
We follow the proof of Proposition~\ref{prop:alphabound} and use its notation $Y$, $y$, $B_\alpha$,
$J$, $a_0$, $a_1$, $E_0$. Put
\[
t_0:=\frac{9(Y-1)^2(Y+2)}{80Y^3},\qquad
b_1:=\max_{t_0\le x\le\gamma}|\bar\Psi_x'(0)|+\tfrac12\varpi^+\big((Y+1)^2-1\big),\qquad K:=16p_0^2b_1^2,
\]
$C_5:=\varpi^-|p_0|\big(1+Y(\pi/k+\gamma)/|W_s|\big)$, and let $E_0$ be defined as there, with this
$t_0$. We may assume
$\alpha>\alpha_1:=\max\{Y+2\cc^{1/2},\,Y+\cc,\,2Y+2,\,2K/Y,\,8|\kappa||W_s|\}$. The first part of the
proof of Lemma~\ref{lem:surfspeed} uses only $\eta_{\min}<\gamma<\eta_{\max}$, so
$\e^\eta<z_-(\alpha)\le Y+1$ on the free surface, $\eta_{\max}\le\ln(Y+1)$, $\alpha-y\ge\alpha/2$ for
$y=\e^\eta$, and \eqref{eq:Bbounds} holds. Also $\Gamma_\alpha\le a_1\alpha+a_0$ and
$\int_{\eta_{\min}}^{\eta_{\max}}G_\alpha\ge\int_{\eta_{\min}}^\gamma G_\alpha\ge\alpha J(\eta_{\min})$, as before.

\emph{Deep troughs: $\eta_{\min}\le t_0$.} Inequality \eqref{eq:kffineq} gives
$\alpha\big(J(\eta_{\min})-2a_1\eta_{\min}\big)\le t_0\big(2a_0+2\varpi^2(Y+1)^4\big)+|\varpi|(Y+1)^2$.
As in Proposition~\ref{prop:alphabound}, $J(t)-2a_1t\ge J(0)-\frac{20Y^3t}{27p_0^2}\ge\frac12J(0)$ for
$t\le t_0$, which bounds $\alpha$.

\emph{Shallow troughs: $\eta_{\min}>t_0$.} The divergence theorem in $D_+$ and
$\widehat\Delta\Psi=-\varpi\e^{2s}$ give
$\int_{\Gamma_+}Q\dd\ell=\int_0^{\pi/k}|\Psi_s(\Theta,0)|\dd\Theta+\varpi\iint_{D_+}\e^{2s}\le\pi b_1/k$, by
Lemma~\ref{lem:kcompare}(c). As in Proposition~\ref{prop:alphabound}, the set
$\mathsf S:=\{\Theta\in[0,\pi/k]:\e^{\eta(\Theta)}\ge Y-K/\alpha\}$ has a complement of length at most
$\pi/(4k)$. Let $\varphi=\sin(k\Theta)\min\{s/t_0,1\}$, so that $\varphi_\Gamma=\sin k\Theta$. We claim
that
\begin{equation}\label{eq:kgflat}
\int_0^{\pi/k}g_*^\flat\,\varphi_\Gamma^2\dd\Theta\le E_0+\pi+C_5,
\end{equation}
where $g_*^\flat:=g_*(\eta)$ if $\varpi\ge0$ or $\eta\le\gamma$, and
$g_*^\flat:=y\alpha/(4B_\alpha(y))-1$ if $\varpi<0$ and $\eta>\gamma$. For $\varpi\ge0$, \eqref{eq:kgflat}
follows from \eqref{eq:keigineq}, whose second term is non-negative. Let $\varpi<0$, and split
$\Gamma_+$ into $\Gamma_\le:=\Gamma_+\cap\{\eta\le\gamma\}$ and $\Gamma_>:=\Gamma_+\cap\{\eta>\gamma\}$; since
$\eta$ is decreasing, these are arcs. On $\Gamma_\le$, $B_\alpha(y)\ge\cc$ gives
$Q=yB_\alpha(y)^{1/2}/|p_0|\ge y|W_s|/|p_0|$, so $\e^{2\eta}/Q\le Y|p_0|/|W_s|$, and $\Gamma_\le$ has length at
most $\pi/k+\gamma$. On $\Gamma_>$, $\e^{2\eta}/Q=y|p_0|/B_\alpha(y)^{1/2}$ and
$\dd\ell\le(1+|\eta'|)\dd\Theta$. Since $B_\alpha(y)=(z_--y)(z_+-y)\ge(z_--y)(z_+-z_-)$ for $Y<y<z_-$,
\[
\int_{\Gamma_>}\frac{y\,|\eta'|}{B_\alpha(y)^{1/2}}\dd\Theta=\int_Y^{\e^{\eta_{\max}}}\frac{\dd y}{B_\alpha(y)^{1/2}}
\le2\Big(\frac{z_--Y}{z_+-z_-}\Big)^{1/2}\le1,
\]
because $z_--Y\le\cc/(\alpha-Y)$ and $z_+-z_-\ge\sqrt3(\alpha-Y)\ge4\cc/(\alpha-Y)$ for $\alpha\ge Y+2\cc^{1/2}$.
Finally, on $\Gamma_>$ we have $B_\alpha\le\cc$ and $g_*\ge y\alpha/(2B_\alpha)-1$. From $\alpha\ge8|\kappa||W_s|=4|\varpi p_0|\cc^{1/2}$ we obtain $|\varpi||p_0|y/B_\alpha^{1/2}\le y\alpha/(4B_\alpha)$.
Inserting these bounds into \eqref{eq:keigineq}, we obtain \eqref{eq:kgflat}.

On $\mathsf S\cap\{\eta\le\gamma\}$ we have $g_*^\flat\ge Y\alpha/(4(\cc+2K))-1$ as in
Proposition~\ref{prop:alphabound}; on $\{\eta>\gamma\}$, $B_\alpha\le\cc$ and $y\ge Y$ give the same
bound, and $g_*^\flat>-1$ everywhere. As there,
$\frac\pi{4k}\big(\frac{Y\alpha}{4(\cc+2K)}-2\big)\le E_0+\pi+C_5$, which bounds $\alpha$. We may take
\[
A_*:=\max\Big\{\alpha_1,\ \frac{2t_0(2a_0+2\varpi^2(Y+1)^4)+2|\varpi|(Y+1)^2}{J(0)},\
\frac{4(\cc+2K)}Y\Big(2+\frac{4k}\pi(E_0+\pi+C_5)\Big)\Big\}.\qedhere
\]
\end{proof}

\begin{theorem}\label{thm:globalkappa}
Let $(\gamma,\kappa,w)$ be admissible and $k\ge1$, and let
$\mathcal K_k=\{(\alpha(s),H(s)):s\in\R\}$ be the curve given by Theorem~\ref{thm:BT} at
$(\alpha_k,H_0)$ in the mathematical parameter space $\alpha\in\R$. Then $\mathcal K_k\subset\mathcal V_k$ consists of solutions of $\Fop(\alpha,H)=0$, and
the following hold.
\begin{enumerate}[label=\textup{(\roman*)}]
\item \emph{(Local structure and symmetry.)} As in Theorem~\ref{thm:global}(i).
\item \emph{(Monotone profiles.)} $H(s)\in\mathcal N_k$ for every $s>0$.
\item \emph{(Bounds.)} For every $s>0$,
\[
\alpha_-<\alpha(s)\le A_*,\qquad x_1(\alpha(s))<\eta_{\min}(s)<\gamma<\eta_{\max}(s),
\]
with $\alpha_-$ and $A_*$ as in Propositions~\ref{prop:kheights} and~\ref{prop:kalphabound}.
In addition,
\begin{equation}\label{eq:kpositive-class}
 |\kappa|\le1\ \text{or}\ \kappa/w\le0
 \quad\Longrightarrow\quad \alpha(s)>\alpha_0>0.
\end{equation}
Moreover:
\begin{enumerate}[label=\textup{(\alph*)}]
\item if $|\kappa|\le1$, then $\eta_{\max}(s)<x_3(\alpha(s))$;
\item if $\kappa/w<0$, then $\eta_{\max}(s)<x_c$;
\item $\eta_{\max}$ and the relative velocity are bounded along $\{(\alpha(s),H(s)):s>0\}$.
\end{enumerate}
\item \emph{(Global alternative.)} $\max_{\overline\Omega}\p_pH(s)\to\infty$ as $s\to\infty$.
\item \emph{(No loop.)} $\mathcal K_k$ is not a closed loop.
\end{enumerate}
For $\kappa=0$, (iii)(a) and (iii)(c) apply, and we recover Theorem~\ref{thm:global}, with a larger
constant $A_*$.
\end{theorem}
Indeed $W=V-R=p_0e^{-s}/H_p$ in the height coordinates. Since $R=e^s\ge1$, evaluating at a maximizer of $H_p$ gives $\min|V-R|\le |p_0|/\max H_p\to0$ along the whole positive continuation. The assertion concerns this curve of graph waves; secondary branches of its connected component are outside the argument.

\begin{proof}
By Lemma~\ref{lem:kstructure}, Theorem~\ref{thm:BT} applies, and (i) and (v) follow as in the proof of
Theorem~\ref{thm:global}.

(ii), (iii) Let $s_*\in(0,\infty]$ be as in the proof of Theorem~\ref{thm:global}. For $0<s<s_*$ the
wave $H(s)\in\mathcal N_k$ is not laminar, so Proposition~\ref{prop:kheights}(a) gives
$\eta_{\min}(s)\ne\gamma\ne\eta_{\max}(s)$. These functions are continuous in $s$. By (C2) and Theorem~\ref{thm:local},
$H(s)=H_0+s(\hat\varphi_k+o(1))$ in $\Xs_k$ as $s\downarrow0$, and $\hat\varphi_k(q,0)=\cos kq$, so
$\eta_{\min}(s)<\gamma<\eta_{\max}(s)$ for small $s>0$. Hence
\begin{equation}\label{eq:kcross}
\eta_{\min}(s)<\gamma<\eta_{\max}(s)\qquad(0<s<s_*).
\end{equation}
If $s_*<\infty$, then as in the proof of Theorem~\ref{thm:global}, $H(s_*)=\bar H_{\bar x}$ for some
$\bar x$; letting $s\uparrow s_*$ in \eqref{eq:kcross} gives $\bar x=\gamma$, so $H(s_*)=H_0$ and
$\alpha(s_*)=\alpha_k$, contradicting Lemma~\ref{lem:noreturn}. Hence we obtain $s_*=\infty$, which is (ii), and
\eqref{eq:kcross} holds for all $s>0$. Propositions~\ref{prop:kheights}(a),(b)
and~\ref{prop:kalphabound} give the displayed bounds. (a) For $|\kappa|\le1$,
Proposition~\ref{prop:kheights}(c) gives $\alpha(s)\ne\alpha_0$, so $\alpha(s)>\alpha_0$ by continuity,
and then \eqref{eq:heights}. The formula for $\alpha_0$ in
Lemma~\ref{lem:kgalpha}(b) also gives $\alpha_0>0$ in this case.
For $\kappa/w\le0$, Proposition~\ref{prop:kalpha-adverse} proves the
additional lower bound in \eqref{eq:kpositive-class}, without requiring
$|\kappa|\le1$. (b) is part of Lemma~\ref{lem:kcompare}. (c) If $|\kappa|\le1$, then
$\eta_{\max}(s)<x_3(A_*)$ by (a) and Lemma~\ref{lem:kgalpha}(b). If $\kappa/w<0$, then $\eta_{\max}(s)<x_c$ by
(b). If $\kappa/w>0$, then $\eta_{\max}(s)\le X(A')$ with $A':=\max\{A_*,|\alpha_-|\}$, by
Lemma~\ref{lem:kcrest}(b). Lemma~\ref{lem:kvelbound} then bounds the relative velocity.

(iv) By (v), (C5)(I) holds. Suppose that $\max\p_pH(s_n)\le M$ along a sequence $s_n\to\infty$, and
let $M':=\max\{M,A_*,|\alpha_-|\}$. By (iii), $|\alpha(s_n)|\le M'$, and Lemma~\ref{lem:kapriori}
with $M'$ in place of $M$ puts all $(\alpha(s_n),H(s_n))$ into a fixed set $Q_j$, a contradiction.
\end{proof}

\begin{remark}\label{rem:kglobalinterp}
The global continuation parameter is $\alpha\in\R$. Under
\eqref{eq:kpositive-class}, the entire curve remains in the physical
range of inward gravity, $\alpha>0$. More generally, the same is true
whenever the explicit lower bound $\alpha_-$ is positive. For
$\kappa/w>0$ and $|\kappa|>1$, our general estimates do not establish
$\alpha>0$ along the whole curve. If one restricts to the physical
component $\alpha>0$, that component may reach $\alpha=0$ before the
mathematical global curve degenerates.

\leavevmode
\begin{enumerate}[label=(\alph*)]
\item As in Remark~\ref{rem:globalinterp}, $V-R=-p_0\e^{-s}\Psi_s=p_0\e^{-s}/H_p$, so (iv) gives
$\min_{\overline D}|V-R|\to0$: the no-stagnation condition \eqref{eq:nostag} degenerates.
\item If $\kappa/w<0$, then $\widehat\Delta\Psi_s=-2\varpi\e^{2s}>0$, so $\Psi_s$ attains its maximum
over $\overline D$ on $\p D$. On the bottom,
$-\Psi_s>-\bar\Psi_{\eta_{\max}}'(0)\ge1/\eta_{\max}>1/x_c$ by Lemmas~\ref{lem:kcompare}(c) and~\ref{lem:kcrestline}(a).
Hence, as for $\kappa=0$, the degeneration takes place on the free surface: either the relative speed
tends to zero somewhere on it, or it becomes radial somewhere.
\item For $\kappa/w>0$, $\Psi_s$ is superharmonic and its maximum need not be on the surface. Lemma~\ref{lem:kcrestline} gives a bed bound while the comparison flow at the crest depth has strictly signed bed velocity.
\end{enumerate}
\end{remark}

\subsection{Where the flow degenerates for constant vorticity}\label{sec:kstag}
For $\kappa\ne0$ the function $\mathfrak f=\log(-\Psi_s-\mathrm i\Psi_\Theta)$ solves \eqref{eq:kdbar}
and is not holomorphic. Then $\mathfrak a=\p_\Theta\mathfrak q$ and the functions $T_m$ of
Section~\ref{sec:slope} are not harmonic, with one exception: $\Psi_s+\mathrm i\Psi_\Theta$ differs
from the holomorphic function $\Phi_s+\mathrm i\Phi_\Theta$ by $-\frac12\varpi\e^{2s}$, which does not
depend on $\Theta$, so $T_1=-\Psi_{s\Theta}$ is harmonic for every $\kappa$. We now propagate the
positivity of $T_1$ along $\mathcal K_k$. This locates the degeneration in Theorem~\ref{thm:globalkappa}(iv)
on the radial line beneath the crest, the crest line $\Theta=0$, and at the crest itself if $\kappa/w\le0$.

Throughout, $(\gamma,\kappa,w)$ is admissible, $\eps=0$, and we use the notation of
Lemmas~\ref{lem:robin} and~\ref{lem:kflowforce}; $B_\alpha$ and $g_*$ are given by \eqref{eq:Bquad}
and \eqref{eq:gstar}, which hold for every $\kappa$ with $\cc=W_s^2$. Where $G_\alpha(x)>0$ we put
\[
\varpi_*(x):=\frac{\varpi\,\e^{2x}}{G_\alpha(x)^{1/2}}=\frac{2\kappa\sgn(w)\,\e^{x}}{B_\alpha(\e^{x})^{1/2}} ,
\]
so that $\varpi_*(\eta)=\varpi\e^{2\eta}/Q$ is the vorticity term in \eqref{eq:krobin}. It has the sign
of $\kappa/w$, vanishes for $\kappa=0$, and $\varpi_*'=(2+g_*)\varpi_*$. On $\overline\Gamma_+$ let
\[
\tilde\varrho:=g_*(\eta)\cos2\chi-2\varkappa\cos\chi+\varpi_*(\eta)\cos\chi ,
\]
which for $\kappa=0$ is the function $\varrho_1$ of Section~\ref{sec:slope}, and for $x>0$ with
$G_\alpha(x)>0$ and $0\le\theta<\pi/2$ let
\[
\mathcal E(x,\theta):=g_*'(x)-\omega_1(\theta)\,g_*(x)^2
+\varpi_*(x)\Big(2\cos\theta-\frac{g_*(x)}{2\cos\theta}\Big),\qquad
\omega_1(\theta)=2-\frac1{2\cos^2\theta} ;
\]
here $\omega_1$ is the function $\omega_m$ of Section~\ref{sec:slope} with $m=1$. We consider the conditions
\begin{enumerate}[label=\textup{(V\arabic*)}]
\item $g_*(\eta_{\min})+\varpi_*(\eta_{\min})\ge0$;
\item $\mathcal E(x,\theta)\ge0$ for $\eta_{\min}\le x\le\eta_{\max}$ and $0\le\theta<\pi/2$.
\end{enumerate}
Every level in $[\eta_{\min},\eta_{\max}]$ is attained on the free surface, where $G_\alpha=Q^2>0$, so
(V1) and (V2) make sense. Since $\omega_1$ decreases from $\omega_1(0)=\frac32=w_1$, for $\kappa=0$
condition (V1) is (G2) and (V2) is (G3${}_1$).

\begin{lemma}\label{lem:T1k}
Let $(\alpha,H)\in\mathcal V_k$ solve $\Fop(\alpha,H)=0$ with $H\in\mathcal N_k$, and let
$T_1:=-\Psi_{s\Theta}$ and $U_1:=-\Psi_{\Theta\Theta}$.
\begin{enumerate}[label=\textup{(\alph*)}]
\item $T_1+\mathrm iU_1$ is holomorphic in $D$ as a function of $\zeta=\Theta+\mathrm is$, and
$T_1,U_1\in C^{3,\beta_0}(\overline D)$. $T_1$ is odd and $U_1$ is even with respect to $\Theta=0$
and to $\Theta=\pi/k$. On the bottom, $U_1=0$, $\p_sT_1=0$, and $T_1=-H_{qp}/H_p^2>0$ for
$0<\Theta<\pi/k$.
\item On $\overline\Gamma_+$,
\[
T_1=Q\,\tilde\varrho\,\sin\chi,\qquad
U_1=Q\big(g_*(\eta)\sin\chi\sin2\chi+\varkappa\cos2\chi+\varpi_*(\eta)\sin^2\chi\big),
\]
and on $\Gamma_+$, $\p_nT_1=\p_tU_1$.
\item At the crest and at the trough, $\p_\Theta T_1=Q\varkappa\tilde\varrho$.
\item If $T_1\ge0$ in $D_+$, then $-\Psi_s$ is non-decreasing in $\Theta$ on every horizontal
segment in $\overline{D_+}$; moreover $U_1\ge0$ on the crest line $\{0\}\times[0,\eta_{\max}]$, $U_1\le0$
on the trough line $\{\pi/k\}\times[0,\eta_{\min}]$, and on both lines
$\p_s(-\Psi_s)=\varpi\e^{2s}-U_1$.
\item If $T_1>0$ in $D_+$, then $\min_{\overline D}(-\Psi_s)$ is attained only on the crest lines
$\Theta\in(2\pi/k)\mathbb Z$; that is, $\max_{\overline\Omega}H_p$ is attained only on the lines
$q\in(2\pi/k)\mathbb Z$. If moreover $\varpi\le0$, it is attained only at the crests, and
\[
\max_{\overline\Omega}H_p=\frac1{\min_{\overline D}(-\Psi_s)}=G_\alpha(\eta_{\max})^{-1/2}.
\]
\end{enumerate}
For $\kappa=0$, $T_1+\mathrm iU_1$ is the function $\mathfrak h_1$ of Section~\ref{sec:slope}.
\end{lemma}

\begin{proof}
(a) Since $\Phi=\Psi+\frac14\varpi\e^{2s}$ is harmonic, $\Phi_s+\mathrm i\Phi_\Theta=\mathrm i(\Phi_\Theta-\mathrm i\Phi_s)$
is holomorphic, and $\Psi_s+\mathrm i\Psi_\Theta=\Phi_s+\mathrm i\Phi_\Theta-\frac12\varpi\e^{2s}$. Hence
$T_1+\mathrm iU_1=-\p_\Theta(\Phi_s+\mathrm i\Phi_\Theta)$ is holomorphic. By the bootstrap in the proof
of Lemma~\ref{lem:kapriori}, $\Psi\in C^{5,\beta_0}(\overline D)$, $\mathfrak f\in C^{4,\beta_0}$ near the free
surface and $\eta\in C^{5,\beta_0}$. The parities follow from those of $\Psi$. On the bottom $\Psi=1$,
so $U_1=0$ there, and the Cauchy--Riemann equation $\p_sT_1=-\p_\Theta U_1$ gives $\p_sT_1=0$.
The bottom point $(\Theta,0)$ corresponds to $(q,p)=(\Theta,-1)$, where $-\Psi_s=1/H_p$ by
\eqref{eq:flat}; hence $T_1=\p_q(1/H_p)=-H_{qp}/H_p^2$ there, which is positive by (N2).

(b) On $\overline\Gamma_+$ we have $-\Psi_s-\mathrm i\Psi_\Theta=\e^{\mathfrak f}=Q\e^{\mathrm i\chi}$, so
$T_1+\mathrm iU_1=\p_\Theta\e^{\mathfrak f}=Q\e^{\mathrm i\chi}(\p_\Theta\mathfrak q+\mathrm i\p_\Theta\chi)$.
In the proof of Lemma~\ref{lem:kflowforce}(c) we found $\p_t\mathfrak q=-g_*\sin\chi$, $\p_t\chi=-\varkappa$,
$\p_\nu\chi=\p_t\mathfrak q$ and $\p_\nu\mathfrak q=\varkappa-\varpi_*(\eta)$. Since
$(1,0)=-\cos\chi\,t-\sin\chi\,\nu$,
\[
\p_\Theta\mathfrak q=\sin\chi\,\big(g_*\cos\chi-\varkappa+\varpi_*\big),\qquad
\p_\Theta\chi=\varkappa\cos\chi+g_*\sin^2\chi ,
\]
and multiplying out gives the formulas for $T_1$ and $U_1$. In the positively oriented frame
$(t,\nu)$ the Cauchy--Riemann equations for $T_1+\mathrm iU_1$ give $\p_\nu T_1=-\p_tU_1$, and $n=-\nu$.

(c) At the crest and at the trough $\chi=0$, $t=(-1,0)$ and $\p_t\chi=-\varkappa$. By (b),
$\p_tT_1=Q\tilde\varrho\cos\chi\,\p_t\chi=-Q\varkappa\tilde\varrho$ there, and $\p_\Theta T_1=-\p_tT_1$.

(d) $T_1=\p_\Theta(-\Psi_s)$. By (N1), $\eta$ is non-increasing on $[0,\pi/k]$, so for
$(\Theta_0,s_0)\in\overline{D_+}$ the segment $[0,\Theta_0]\times\{s_0\}$ lies in $\overline{D_+}$, which gives
the first claim. On the crest line $T_1=0$ and $T_1\ge0$ to its right, so $\p_\Theta T_1\ge0$ there;
by the Cauchy--Riemann equations $\p_sU_1=\p_\Theta T_1\ge0$, and $U_1=0$ at the bottom, so $U_1\ge0$
on the crest line. In the same way $U_1\le0$ on the trough line. On both lines $\Psi_\Theta=0$. From \eqref{eq:flatpoisson} we obtain $-\Psi_{ss}=\Psi_{\Theta\Theta}+\varpi\e^{2s}=\varpi\e^{2s}-U_1$.

(e) Let $(\Theta_0,s_0)\in\overline{D_+}$ with $\Theta_0>0$. The open segment $(0,\Theta_0)\times\{s_0\}$
lies in $D_+$ if $s_0>0$, because $\eta$ is strictly decreasing on $[0,\pi/k]$ by (N1), and on the
bottom if $s_0=0$; in both cases $T_1>0$ on it, by (a). Hence $-\Psi_s(\Theta_0,s_0)>-\Psi_s(0,s_0)$,
and by evenness and periodicity the minimum over $\overline D$ is attained only on the crest lines. By
\eqref{eq:flat}, $\max H_p=1/\min(-\Psi_s)$, at corresponding points. Let $\varpi\le0$. By the Hopf lemma,
$\p_\Theta T_1>0$ on $\{0\}\times(0,\eta_{\max})$, so $U_1>0$ on $\{0\}\times(0,\eta_{\max}]$ by the proof
of (d), and $\p_s(-\Psi_s)=\varpi\e^{2s}-U_1<0$ there. Thus $-\Psi_s$ decreases strictly along the crest
line, and its minimum is attained at the crest, where $\Psi_\Theta=0$ and $-\Psi_s=Q=G_\alpha(\eta_{\max})^{1/2}$
by \eqref{eq:bern-flat}. For $\kappa=0$, $T_1+\mathrm iU_1=\p_\Theta\e^{\mathfrak f}=\mathfrak f'\e^{\mathfrak f}=\mathfrak h_1$.
\end{proof}

\begin{theorem}\label{thm:kstag}
Let $\mathcal K_k=\{(\alpha(s),H(s)):s\in\R\}$ be the curve of Theorem~\ref{thm:globalkappa}, and let
$0<s_1\le\infty$. Suppose that \textup{(V1)} and \textup{(V2)} hold for every $s\in(0,s_1)$, with
$\alpha=\alpha(s)$ and with $\eta_{\min}$, $\eta_{\max}$ the trough and crest levels of $H(s)$. Then for
every $s\in(0,s_1)$:
\begin{enumerate}[label=\textup{(\roman*)}]
\item $\tilde\varrho>0$ on $\overline\Gamma_+$ and $T_1>0$ in $D_+$;
\item on every circle $R=\mathrm{const}$ that meets the fluid, the azimuthal relative speed $|V-R|$ is
smallest exactly on the crest lines $\Theta\in(2\pi/k)\mathbb Z$, and $\max_{\overline\Omega}\p_pH(s)$ is
attained only on the lines $q\in(2\pi/k)\mathbb Z$;
\item if $\kappa/w\le0$, then
$\max_{\overline\Omega}\p_pH(s)=\p_pH(s)(0,0)=G_{\alpha(s)}\big(\eta_{\max}(s)\big)^{-1/2}$.
\end{enumerate}
If $s_1=\infty$, then $\min_{0\le\sigma\le\eta_{\max}(s)}(-\Psi_s)(0,\sigma)\to0$ as $s\to\infty$: the relative
velocity tends to zero at some point of the crest line. If $s_1=\infty$ and $\kappa/w\le0$, then
$G_{\alpha(s)}(\eta_{\max}(s))\to0$, that is, the relative velocity at the crest tends to zero, and every
limit point $(\bar\alpha,\bar\eta)$ of $(\alpha(s),\eta_{\max}(s))$ as $s\to\infty$ satisfies
\[
\bar\alpha>\e^\gamma,\qquad(\bar\alpha-\e^\gamma)^2\ge\cc,\qquad \e^{\bar\eta}=z_-(\bar\alpha).
\]
\end{theorem}

\begin{proof}
We write $\tilde\varrho_s$ for the function $\tilde\varrho$ of $H(s)$, regarded as a function of
$\Theta\in[0,\pi/k]$. As in the proof of Theorem~\ref{thm:stagnation}, $\tilde\varrho_s$ depends
continuously on $s>0$ in $C([0,\pi/k])$; here we use that $B_{\alpha(s)}(\e^\eta)>0$ on the free surface.
As $s\downarrow0$, $\eta\to\gamma$ and $\chi,\varkappa\to0$ uniformly and $\alpha(s)\to\alpha_k$, so
$\tilde\varrho_s$ converges uniformly to $g_*(\gamma)+\varpi_*(\gamma)$, evaluated at $\alpha=\alpha_k$. Since
$B_\alpha(\e^\gamma)=\cc=W_s^2$ and $W_s$ has the sign of $w$, \eqref{eq:kalphax} gives
\[
g_*(\gamma)+\varpi_*(\gamma)=\frac{\e^\gamma(\alpha_k-\e^\gamma)}{W_s^2}-1+\frac{2\kappa\e^\gamma}{W_s}
=k\coth(k\gamma)>0 .
\]
Let $s_\sharp$ be the supremum of all $s'\in(0,s_1]$ such that $\tilde\varrho_s>0$ on $[0,\pi/k]$ for
$0<s<s'$. Then $s_\sharp>0$. Suppose that $s_\sharp<s_1$, let $H:=H(s_\sharp)$, and let
$\tilde\varrho:=\tilde\varrho_{s_\sharp}$. Then $\tilde\varrho\ge0$ on $\overline\Gamma_+$, with a zero at some
$P_0\in\overline\Gamma_+$. By Theorem~\ref{thm:globalkappa}(ii), $H\in\mathcal N_k$, and by
Lemma~\ref{lem:T1k}(a),(b), $T_1\ge0$ on $\p D_+$ and $T_1>0$ on the bottom. Hence $T_1>0$ in $D_+$.
Along $\Gamma_+$ we have $\p_t\chi=-\varkappa$, $\p_t\eta=\sin\chi$, $\p_t\mathfrak q=-g_*\sin\chi$ and
$\p_t[\varpi_*(\eta)]=(2+g_*)\varpi_*\sin\chi$.

The two boundary formulas in Lemma~\ref{lem:T1k}(b) imply the useful
identity
\begin{equation}\label{eq:kU-rho}
 \frac{U_1}{Q}
 =\frac{g_*}{2\cos\chi}+\frac{\varpi_*}{2}
  -\frac{\cos2\chi}{2\cos\chi}\,\tilde\varrho.
\end{equation}

\emph{Case 1: $P_0\in\Gamma_+$.}
Here $T_1(P_0)=0$ and $T_1>0$ in $D_+$, so Hopf's lemma gives
$\p_nT_1(P_0)<0$. Since $\tilde\varrho$ has a minimum at $P_0$,
$\tilde\varrho=\p_t\tilde\varrho=0$ there. Differentiate
\eqref{eq:kU-rho}, using
\[
 \frac{\p_tQ}{Q}=-g_*\sin\chi,\quad
 \p_tg_*=g_*'\sin\chi,\quad
 \p_t\cos\chi=\varkappa\sin\chi,\quad
 \p_t\varpi_*=(2+g_*)\varpi_*\sin\chi.
\]
At $P_0$ the relation $\tilde\varrho=0$ gives
$2\varkappa\cos\chi=g_*\cos2\chi+\varpi_*\cos\chi$.
Substitution yields directly $\p_nT_1(P_0)=\p_tU_1(P_0) =\frac{Q\sin\chi}{2\cos\chi}\, \mathcal E(\eta(P_0),\chi(P_0))\ge0,$ contradicting Hopf's lemma.

\emph{Case 2: $P_0$ is the crest.}
At this point
$2\varkappa_c=g_*(\eta_{\max})+\varpi_*(\eta_{\max})$ and
$\p_\Theta T_1=0$. Since $T_1$ is positive in $D_+$ and odd in
$\Theta$, the same application of Serrin's corner lemma as in
Theorem~\ref{thm:stagnation} gives $\p_t^2U_1>0$ at the crest.
The even function $\tilde\varrho$ is non-negative and vanishes there,
so $\p_t\tilde\varrho=0$ and $\p_t^2\tilde\varrho\ge0$.
At the crest,
\[
 \p_t^2\eta=-\varkappa_c,\quad
 \frac{\p_t^2Q}{Q_c}=g_*\varkappa_c,\quad
 \p_t^2(\cos\chi)=-\varkappa_c^2,\quad
 \p_t^2\varpi_*=-(2+g_*)\varpi_*\varkappa_c.
\]
Twice differentiating \eqref{eq:kU-rho} therefore gives the compact
identity $\p_t^2\tilde\varrho+\frac{2}{Q_c}\p_t^2U_1 =-\varkappa_c\,\mathcal E(\eta_{\max},0).$ Its left side is strictly positive, whereas its right side is
non-positive by (V2) and $\varkappa_c>0$. This is a contradiction.

\emph{Case 3: $P_0$ is the trough.} Then $0=\tilde\varrho=g_*(\eta_{\min})-2\varkappa_t+\varpi_*(\eta_{\min})$, so
$g_*(\eta_{\min})+\varpi_*(\eta_{\min})=2\varkappa_t<0$ by (N3), contradicting (V1).

Hence we obtain $s_\sharp=s_1$, and (i) holds. By Remark~\ref{rem:kglobalinterp}(a), $|V-R|=|p_0|\e^{-s}(-\Psi_s)$,
so Lemma~\ref{lem:T1k}(e) gives (ii) and (iii).

Let $s_1=\infty$. By Theorem~\ref{thm:globalkappa}(iv) and \eqref{eq:flat},
$\min_{\overline D}(-\Psi_s)=1/\max\p_pH(s)\to0$, and by (ii) this minimum is attained on the crest line. Let
$\kappa/w\le0$. By (iii), $G_{\alpha(s)}(\eta_{\max}(s))\to0$. By Theorem~\ref{thm:globalkappa}(iii),
$\alpha(s)\in(\alpha_-,A_*]$ and $\eta_{\max}(s)$ is bounded, since $\kappa/w<0$ or $\kappa=0$; so limit points
exist. By \eqref{eq:Bquad}, $B_{\alpha(s)}(\e^{\eta_{\max}(s)})=p_0^2\e^{-2\eta_{\max}}G_{\alpha(s)}(\eta_{\max})\to0$, hence
$B_{\bar\alpha}(\e^{\bar\eta})=0$ with $\bar\eta\ge\gamma$. If $\bar\alpha\le\e^\gamma$, then
$B_{\bar\alpha}(y)=\cc+(y-\e^\gamma)(y+\e^\gamma-2\bar\alpha)\ge\cc$ for $y\ge\e^\gamma$; hence $\bar\alpha>\e^\gamma$. The
rest of the argument is the end of the proof of Theorem~\ref{thm:stagnation}: the first part of the proof
of Lemma~\ref{lem:surfspeed} uses only $\alpha>\e^\gamma$, $(\alpha-\e^\gamma)^2\ge\cc$ and
$\eta_{\min}<\gamma$.
\end{proof}

\begin{remark}
The proof only uses (V2) on the attained boundary pairs
$(x,\theta)=(\eta(\Theta),\chi(\Theta))$, including the crest.
Thus the theorem also holds if (V2) is replaced by
$\mathcal E(\eta(\Theta),\chi(\Theta))\ge0$ on
$\overline\Gamma_+$ for every wave in the indicated part of the curve.
The rectangular condition (V2) is a convenient stronger criterion.
\end{remark}

For $\kappa=0$, Theorem~\ref{thm:kstag} is the case $m=1$ of Theorem~\ref{thm:slope}, with (G2) in place of
(G2${}'$). For $\kappa/w>0$ the theorem leaves open where on the crest line the degeneration takes place.
Lemma~\ref{lem:kcrest}(a) and the next lemma bound the relative speed there, without (V1) and (V2).

\begin{lemma}\label{lem:kcrestline}
Let $(\alpha,H)\in\mathcal V_k$ be a non-laminar solution of $\Fop(\alpha,H)=0$ with $H\in\mathcal N_k$, and
let $Q_c:=G_\alpha(\eta_{\max})^{1/2}$.
\begin{enumerate}[label=\textup{(\alph*)}]
\item On the bottom, $-\Psi_s>\beta(\eta_{\max})$, where
\[
\beta(x):=-\bar\Psi_x'(0)=\frac1x-\frac\varpi2\Big(\frac{\e^{2x}-1}{2x}-1\Big).
\]
If $\varpi\le0$, then $\beta(x)\ge1/x$. If $\varpi>0$, then $\beta$ is strictly decreasing and has a single
zero $x_b>\gamma$, the depth at which the laminar flow with flux $p_0$ and vorticity $2\kappa$ has zero
relative velocity at the bed $R=1$. Consequently, if the relative speed at the bed tends to zero along a
sequence of such solutions, then $\liminf\eta_{\max}\ge x_b$.
\item Let $\varpi\le0$, $\eta_{\min}<\gamma$, $\alpha>\e^\gamma+\cc\e^{-\gamma}$ and $(\alpha-\e^\gamma)^2\ge\cc$; the
first condition holds along $\mathcal K_k$ by Theorem~\ref{thm:globalkappa}(iii). Then the conclusions of
Lemma~\ref{lem:surfspeed} hold: $\e^\eta<z_-(\alpha)$ on the free surface, and $Q>Q_c$ at every point of the
free surface where $\eta<\eta_{\max}$.
\end{enumerate}
\end{lemma}

\begin{proof}
(a) The inequality is Lemma~\ref{lem:kcompare}(c), and the formula for $\beta$ follows from the
definitions of $\mu(x)$ and $\bar\Psi_x$. The function $(\e^{2x}-1)/(2x)=\sum_{j\ge0}(2x)^j/(j+1)!$ is
increasing and at least $1$; this gives the monotonicity and the bound for $\varpi\le0$. Moreover
$\beta(x)\to-\infty$ as $x\to\infty$ if $\varpi>0$, and $\beta(\gamma)=-\Psi_0'(0)>0$ because the laminar flow
$\Psi_0$ is admissible. If the relative speed at the bed tends to zero, then so does
$\max\{\beta(\eta_{\max}),0\}$, which gives the last claim.

(b) The first part of the proof of Lemma~\ref{lem:surfspeed} applies verbatim, since it uses only
$\eta_{\min}<\gamma<\e^\gamma<\alpha$, $(\alpha-\e^\gamma)^2\ge\cc$ and $B_\alpha(\e^\eta)>0$ on the free surface. Its second part uses
$\alpha>\e^\gamma+\cc\e^{-\gamma}$ and $G_\alpha(\eta_{\max})<G_\alpha(\eta_{\min})$. By Lemma~\ref{lem:kcompare}(a),(b),
the latter follows if $x\mapsto|\bar\Psi_x'(x)|=\e^xW_s(x)/p_0$ is strictly decreasing on
$[\eta_{\min},\eta_{\max}]$. Its derivative is $\e^x\big(\varpi\e^x-W_s(x)/(xp_0)\big)$, by the formula for $W_s'$
in the proof of Lemma~\ref{lem:kgalpha}, and this is negative because $\varpi\le0$ and $W_s(x)/p_0>0$ for
$x<x_c$ (Lemma~\ref{lem:kgalpha}(c) and Lemma~\ref{lem:kcompare}; for $\kappa=0$, $W_s(x)/p_0=\e^{-x}/x$).
\end{proof}

For $\kappa/w>0$, azimuthal stagnation at the bed can occur only if
the limiting crest height reaches at least $x_b$.
For $\varpi>0$, the estimate along the crest line is an estimate for the
normalised quantity $-\Psi_s$:
\[
 -\Psi_s(0,\sigma)\ge Q_c-
 \frac{\varpi}{2}(\e^{2\eta_{\max}}-\e^{2\sigma}).
\]
Since the radial velocity vanishes on that line, it also implies the
physical relative-speed estimate
\[
 |W(\e^\sigma,0)|\ge |W(\e^{\eta_{\max}},0)|
   -|\kappa|(\e^{2\eta_{\max}}-1),
 \qquad 0\le\sigma\le\eta_{\max}.
\]

If $\kappa/w\le0$, Proposition~\ref{prop:kalpha-adverse} gives
$\alpha>\e^\gamma+\cc\e^{-\gamma}$ along $\mathcal K_k$.
Consequently Lemma~\ref{lem:kcrestline}(b) applies whenever
$(\alpha-\e^\gamma)^2\ge\cc$, without a restriction on $|\kappa|$.

Together with Remark~\ref{rem:kglobalinterp}(b), this says that
the waves with $(\alpha-e^\gamma)^2\ge\cc$ can degenerate through crest
stagnation or a nearly radial free surface, without a bound on $|\kappa|$.

Conditions (V1) and (V2) hold a priori for thin layers.

\begin{corollary}\label{cor:kthin}
Let $\alpha_\flat>\e^\gamma+\cc^{1/2}$ be a lower bound for $\alpha$ along $\{(\alpha(s),H(s)):s>0\}$; by
Theorem~\ref{thm:globalkappa}(iii) one may take $\alpha_\flat=\alpha_-$
when this satisfies the displayed strict lower bound. Under
\eqref{eq:kpositive-class}, one may instead take $\alpha_\flat=\alpha_0$,
again provided that $\alpha_0>\e^\gamma+\cc^{1/2}$.
Let $z_\flat:=z_-(\alpha_\flat)$, $b_\flat:=z_\flat/\big(z_+(\alpha_\flat)-z_\flat\big)$ and $A_\flat:=1/(z_\flat-1)$, and suppose that
\begin{equation}\label{eq:kthin}
b_\flat\le\tfrac12\qquad\text{and}\qquad
A_\flat^{1/2}\ge2|\kappa|\,b_\flat^{1/2}+\big(4\kappa^2b_\flat+2\big)^{1/2}.
\end{equation}
Then \textup{(V1)} and \textup{(V2)} hold for every $s>0$, and Theorem~\ref{thm:kstag} applies with $s_1=\infty$.
\end{corollary}

\begin{proof}
Let $s>0$ and $\alpha=\alpha(s)\ge\alpha_\flat$, and write $\delta_\alpha=(\alpha-\e^\gamma)^2-\cc>0$. Since
$z_-(\alpha)=\alpha-\delta_\alpha^{1/2}$ is decreasing for $\alpha>\e^\gamma+\cc^{1/2}$, $z_-(\alpha)\le z_\flat$; and
$b_\alpha:=z_-/(z_+-z_-)\le b_\flat$, as in the proof of Corollary~\ref{cor:slope}. By the first part of the proof of
Lemma~\ref{lem:surfspeed}, every level $y=\e^x$ with $x\in[\eta_{\min},\eta_{\max}]$ satisfies $1<y<z_-(\alpha)$. With
$A:=y/(z_--y)$ and $B:=y/(z_+-y)$ we have, as in the proof of Corollary~\ref{cor:slope}, $g_*=\frac{A+B}2-1,\qquad g_*'-2g_*^2=\frac{(5-2B)A+5B-4}2,\qquad |\varpi_*|=2|\kappa|(AB)^{1/2},$ and $A>1/(z_--1)\ge A_\flat$, $0<B<b_\alpha\le b_\flat$. The function $u\mapsto u^2/2-2|\kappa|b_\flat^{1/2}u-1$ is
non-negative for $u\ge2|\kappa|b_\flat^{1/2}+(4\kappa^2b_\flat+2)^{1/2}$, so \eqref{eq:kthin} with $u=A^{1/2}$ gives
$g_*\ge\frac A2-1\ge2|\kappa|(Ab_\flat)^{1/2}\ge|\varpi_*|$. This implies (V1). For (V2) we write
\[
\mathcal E(x,\theta)=\big(g_*'-2g_*^2\big)+\frac{g_*}{2\cos\theta}\Big(\frac{g_*}{\cos\theta}-\varpi_*\Big)+2\varpi_*\cos\theta .
\]
The first term is at least $2A-2$ because $B\le\frac12$, the second is non-negative because
$g_*\ge|\varpi_*|$, and the third is at least $-2|\varpi_*|\ge-(A-2)$. Hence $\mathcal E\ge A>0$.
\end{proof}

For $\kappa\ne0$ the logarithmic velocity solves \eqref{eq:kdbar}, whose right-hand side can become unbounded as the speed vanishes. The irrotational corner argument therefore leaves the rotational limiting profile open.

\section{Tidal equilibria and resonant bifurcation}\label{sec:forced}
\subsection{The tidal branch off the even spectrum}\label{sec:tidal}

We now include the tidal forcing. To identify the branch that persists away from resonance, we first specify its symmetry class. Let
\[
\Sigma_{\rm e}:=\{\alpha_{2m}:\ m\ge0\}=\{\alpha_0,\alpha_2,\alpha_4,\dots\}.
\]

\begin{theorem}\label{thm:tidal}
Let $I\subset\R\setminus\Sigma_{\rm e}$ be a compact interval and fix the sign in $\Fop^\pm$.
There exist $\eps_0>0$, $r_0>0$ and a real-analytic map
\[
H_T^\pm:I\times(-\eps_0,\eps_0)\to\Xs_\pi\cap\Uo
\]
with the following properties.
\begin{enumerate}[label=\textup{(\roman*)}]
\item $\Fop^\pm(\alpha,\eps;H_T^\pm(\alpha,\eps))=0$ and $H_T^\pm(\alpha,0)=H_0$.
\item If $H\in\Xs_\pi$, $\|H-H_0\|_{\Xs}<r_0$ and $\Fop^\pm(\alpha,\eps;H)=0$ with
$(\alpha,\eps)\in I\times(-\eps_0,\eps_0)$, then $H=H_T^\pm(\alpha,\eps)$.
\item $H_T^\pm(\alpha,\eps)$ is invariant under $\rho$ and $\iota$, and
$H_T^-(\alpha,\eps)=\tau_{\pi/2}H_T^+(\alpha,\eps)$.
\item $H_T^\pm(\alpha,\eps)=H_0+\eps H_{T,1}^\pm(\alpha)+O(\eps^2)$ in $\Xs$, uniformly on $I$, where
\begin{equation}\label{eq:HT1}
H_{T,1}^\pm(\alpha)=\frac{\e^{\gamma}}{4(\alpha-\alpha_0)}\,\hat\varphi_0
\pm\frac{3\e^{\gamma}}{4(\alpha-\alpha_2)}\,\hat\varphi_2,\qquad
\hat\varphi_0:=\frac{\sigma_0}{h_1\gamma},\quad \hat\varphi_2:=\frac{\sigma_2(p)\cos2q}{h_1\sinh(2\gamma)},
\end{equation}
with $\sigma_0,\sigma_2$ as in Proposition~\ref{prop:dispersion}, so that $\hat\varphi_0(\cdot,0)=1$ and
$\hat\varphi_2(q,0)=\cos2q$. For $\kappa=0$, $\hat\varphi_0=1+p$ and
$\hat\varphi_2=\sinh(2\gamma(1+p))\cos2q/\sinh(2\gamma)$.
\item For each $j\ge1$, $\p_\eps^jH_T^\pm(\alpha,0)$ is a trigonometric polynomial in $q$ of degree
at most $2j$, with coefficients smooth in $p$.
\end{enumerate}
In particular, the free surface of $H_T^+$ is
\begin{equation}\label{eq:S1}
S(\Theta)=S_0+\eps\Big[\frac{S_0^2}{4(\alpha-\alpha_0)}+\frac{3S_0^2}{4(\alpha-\alpha_2)}\cos2\Theta\Big]+O(\eps^2).
\end{equation}
\end{theorem}

\begin{proof}
On the even $\pi$-periodic subspace, Lemma~\ref{lem:range}(d) gives an isomorphism for every $\alpha\in I$. The analytic implicit function theorem and a finite cover of $I$ give uniform radii and uniqueness. The quarter-turn conjugacy gives (iii). Differentiation at zero forcing gives $\Lop H_{T,1}^\pm=(0,2e^{2\gamma}P^\pm)$. On mode $m$ the normalised harmonic solution has boundary residual $2e^\gamma(\alpha-\alpha_m)\cos mq$, giving \eqref{eq:HT1} and, after exponentiation, \eqref{eq:S1}. At order $j$, products of lower coefficients of total order $j$ have Fourier degree at most $2j$ by induction. A forcing factor uses one order and adds at most two to the degree. The inverse preserves modes, proving (v).
\end{proof}

\begin{remark}
The local uniqueness of the tidal equilibrium also holds among all
$\pi$-periodic perturbations, without prescribing evenness. Indeed,
the linearised operator on that full periodic space is invertible away
from $\Sigma_{\rm e}$: the sine modes give the same excluded spectrum.
The analytic implicit function theorem applies there, and uniqueness
under reflection then forces the equilibrium to be even.
\end{remark}

\begin{remark}\label{rem:airy}
In dimensional form the $\cos2\Theta$ part of \eqref{eq:S1} reads
\[
\eta^{(2)}(\theta,t)=\frac34\,\frac{GMa^2}{gd^3}\,S_0^2\,\frac{\alpha}{\alpha-\alpha_2}\,\cos2(\theta-\Omega_ft)+O(\eps^2).
\]
The dynamic factor $\alpha/(\alpha-\alpha_2)$ changes sign across $\alpha_2$. Its denominator includes the harmonic extension of the angular disturbance through the layer via $2\coth(2\gamma)$ in the dispersion relation. The thin-layer limit is evaluated in Section~\ref{sec:earth}.

At the mean resonance $\alpha_0$ the laminar Bernoulli constant is stationary with respect to the depth $h$ at fixed $(\alpha,p_0,\kappa)$. Indeed, for $\lambda_{\rm lam}(h)=W_s(h)^2/2+\alpha e^h-e^{2h}/2$, differentiation gives $\lambda_{\rm lam}'(h)=e^h[\alpha-\alpha_0(h)]$, using \eqref{eq:kalphax}. A nondegenerate fold requires $\alpha_0'(\gamma)\ne0$.
\end{remark}

\subsection{Odd modes and reflection symmetry}\label{sec:odd}

For every $n\ge1$ we define, with $T_n=\coth(n\gamma)$ and $r=w^2$,
\begin{gather*}
\delta_n:=\frac{\e^{\gamma}N_n}{4r\,(n\gamma T_n-1)},\qquad
\eta_1:=\frac{3\e^{\gamma}M_1}{8r},\qquad \eta_n:=0\ \ (n\ge2),\\
N_n:=\gamma(1-4\kappa^2)+2\kappa w\,(2n\gamma T_n+3\gamma+1)-r\big(n^2\gamma T_n^2+n\gamma T_n+2nT_n-(n^2-3)\gamma+4\big),\\
M_1:=r\,(3T_1^3+5T_1^2+2T_1+4)-T_1(1-4\kappa^2)-2\kappa w\,(3T_1^2+3T_1+1).
\end{gather*}
Throughout this subsection $n\ge1$ is odd. We fix a compact interval $I$ with $\alpha_n$ in its
interior and $I\cap\Sigma_{\rm e}=\emptyset$; this is possible since $\alpha_n\notin\Sigma_{\rm e}$.
Let $H_T^\pm$ be given by Theorem~\ref{thm:tidal}.

\begin{theorem}\label{thm:odd}
There exist $\eps_1>0$, $\bar s>0$, a neighbourhood $\mathcal W$ of $(\alpha_n,H_0)$ in $\R\times\Xs$, and
real-analytic functions
\[
\alpha_n^\pm:(-\eps_1,\eps_1)\to I,\qquad
A^\pm:(-\bar s^2,\bar s^2)\times(-\eps_1,\eps_1)\to I
\]
with $\alpha_n^\pm(0)=\alpha_n$ and $A^\pm(0,\eps)=\alpha_n^\pm(\eps)$, such that the following hold for
$|\eps|<\eps_1$.
\begin{enumerate}[label=\textup{(\roman*)}]
\item \emph{(Pitchfork.)} There is an analytic family $H^\pm(s,\eps)\in\Xs$, $|s|<\bar s$, with
\[
\Fop^\pm\big(A^\pm(s^2,\eps),\eps;H^\pm(s,\eps)\big)=0 .
\]
It satisfies
\[
H^\pm(0,\eps)=H_T^\pm(\alpha_n^\pm(\eps),\eps),\qquad
\mathsf a_n\big[H^\pm(s,\eps)\big]=s,\qquad H^\pm(-s,\eps)=\rho H^\pm(s,\eps).
\]
Every solution of $\Fop^\pm(\alpha,\eps;H)=0$ with $(\alpha,H)\in\mathcal W$ is either the tidal
equilibrium $H=H_T^\pm(\alpha,\eps)$ or of the form $(A^\pm(s^2,\eps),H^\pm(s,\eps))$.
\item \emph{(Direction.)} $A^\pm(\sigma,\eps)=\alpha_n^\pm(\eps)+\big(\alpha_n^{(2)}+O(\eps)\big)\sigma+O(\sigma^2)$,
with $\alpha_n^{(2)}$ as in Theorem~\ref{thm:direction}.
\item \emph{(Shift and splitting.)} $\alpha_n^\pm(\eps)=\alpha_n+\eps\,(\delta_n\pm\eta_n)+O(\eps^2)$.
\item \emph{(Higher odd modes.)} If $n\ge3$, then $\alpha_n^+(\eps)-\alpha_n^-(\eps)=O(\eps^n)$.
\end{enumerate}
\end{theorem}
\begin{remark}
The estimate $\alpha_n^+(\eps)-\alpha_n^-(\eps)=O(\eps^n)$ says that
splitting cannot occur before order $n$; it does not assert that the
coefficient at that order is non-zero.  For $n=1$, first-order separation
requires $\eta_1\ne0$.  This condition can fail even for $\kappa=0$: it
fails when
\[
r=\frac{T_1}{3T_1^3+5T_1^2+2T_1+4}.
\]
The completeness assertion of Theorem~\ref{thm:odd} is made in each of
the two specified reflection subspaces.  Without an additional
non-degeneracy argument it is not a classification of all nearby waves
in the space without a reflection condition.
\end{remark}

By Lemma~\ref{lem:sym}(d), applying $\tau_{-\pi/2}$ to the solutions of $\Fop^-$ yields solutions of
the original problem $\Fop=\Fop^+$ that are symmetric with respect to the axis perpendicular
to the satellite. Their profiles are, to leading order, $\sin n\Theta$. Thus the single
$O(2)$-pitchfork at $\alpha_n$ of the unforced problem gives rise to two pitchforks:
\begin{itemize}
\item one at $\alpha_n^+(\eps)$, with profiles $\cos n\Theta$ (symmetric about the satellite
axis);
\item one at $\alpha_n^-(\eps)$, with profiles $\sin n\Theta$ (symmetric about the
perpendicular axis).
\end{itemize}
For $n=1$ and $\eta_1\ne0$, the two critical values differ by $2\eta_1\eps$ at first order. For $n\ge3$
they agree to order $\eps^{n-1}$.

\begin{proof}[Proof of Theorem~\ref{thm:odd}, (i)--(iii)]
We treat $\Fop=\Fop^+$; the case $\Fop^-$ is identical. We set $\mathcal G(\alpha,\eps;w):=\Fop\big(\alpha,\eps;H_T(\alpha,\eps)+w\big),\qquad w\in\Xs\ \text{small}.$ Then $\mathcal G$ is analytic and $\mathcal G(\alpha,\eps;0)=0$. Since $\Fop$ commutes with
$\rho$ and $\rho H_T=H_T$, we have $\mathcal G(\alpha,\eps;\rho w)=\rho\,\mathcal G(\alpha,\eps;w)$.
Let $\Mop(\alpha,\eps):=D_w\mathcal G(\alpha,\eps;0)$; then $\Mop(\alpha,0)=\Lop(\alpha)$.

By Proposition~\ref{prop:dispersion}, $\ker\Lop(\alpha_n)=\spn\{\hat\varphi_n\}$, and
$\rho\hat\varphi_n=-\hat\varphi_n$ since $n$ is odd. Let $\Zs_n=\{\psi\in\Xs:\ \mathsf a_n[\psi]=0\},\qquad y_n^*:=(0,\cos nq)\in\Ys,$ and let $Q_n$ be the projection of $\Ys$ onto $\spn\{y_n^*\}$ along
$\mathrm R(\Lop(\alpha_n))=\ker\ell_n$, that is, $Q_nf=\ell_n(f)\,y_n^*/\ell_n(y_n^*)$. The
operator $\rho$ maps $\Zs_n$ and $\ker\ell_n$ into themselves, $\rho y_n^*=-y_n^*$, and
$\ell_n\circ\rho=-\ell_n$. Hence $Q_n$ commutes with $\rho$.

\emph{Lyapunov--Schmidt reduction.} Since $(I-Q_n)\Lop(\alpha_n):\Zs_n\to\ker\ell_n$ is an
isomorphism, the implicit function theorem gives a unique analytic $v(s,\alpha,\eps)\in\Zs_n$,
defined near $(0,\alpha_n,0)$, with $(I-Q_n)\,\mathcal G\big(\alpha,\eps;s\hat\varphi_n+v\big)=0,\qquad v(0,\alpha,\eps)=0.$ The last identity holds because $w=0$ solves $\mathcal G=0$. By equivariance and uniqueness,
$v(-s,\alpha,\eps)=\rho\,v(s,\alpha,\eps)$. We define the reduced function
\[
\mathfrak g(s,\alpha,\eps):=\frac{\ell_n\big(\mathcal G(\alpha,\eps;s\hat\varphi_n+v(s,\alpha,\eps))\big)}{\ell_n\big(\p_\alpha\Lop(\alpha_n)\hat\varphi_n\big)} .
\]
Then $\mathfrak g(-s,\alpha,\eps)=-\mathfrak g(s,\alpha,\eps)$ and $\mathfrak g(0,\alpha,\eps)=0$, so
$\mathfrak g(s,\alpha,\eps)=s\,\Gamma(s^2,\alpha,\eps)$ with $\Gamma$ analytic.

We have
\[
\Gamma(0,\alpha,\eps)=\frac{\ell_n\big(\Mop(\alpha,\eps)(\hat\varphi_n+\p_sv(0,\alpha,\eps))\big)}{\ell_n\big(\p_\alpha\Lop(\alpha_n)\hat\varphi_n\big)},
\qquad \p_sv(0,\alpha_n,0)=0 .
\]
Since $\ell_n$ annihilates $\mathrm R(\Lop(\alpha_n))$, it follows that
$\Gamma(0,\alpha_n,0)=0$ and $\p_\alpha\Gamma(0,\alpha_n,0)=1$. By the implicit function theorem,
$\Gamma(\sigma,\alpha,\eps)=0$ if and only if $\alpha=A(\sigma,\eps)$ with $A$ analytic and
$A(0,0)=\alpha_n$.

We set $\alpha_n^+(\eps):=A(0,\eps)$ and $H^+(s,\eps):=H_T\big(A(s^2,\eps),\eps\big)+s\hat\varphi_n+v\big(s,A(s^2,\eps),\eps\big).$ Then $\mathsf a_n[H^+(s,\eps)]=s$, because $\mathsf a_n[v]=0$ and $H_T\in\Xs_\pi$ has no odd
Fourier modes. Moreover $\rho H^+(s,\eps)=H^+(-s,\eps)$. Local uniqueness follows from the
equivalence between $\mathcal G=0$ and the system ``$(I-Q_n)\mathcal G=0$, $\mathfrak g=0$''.
This proves (i).

For $\eps=0$ the curve coincides with that of Theorem~\ref{thm:local}, by uniqueness and the
identical amplitude normalisation. Hence we obtain $A(\sigma,0)=\alpha_n+\alpha_n^{(2)}\sigma+O(\sigma^2)$,
which gives (ii).

For (iii), implicit differentiation gives
\[
(\alpha_n^+)'(0)=-\p_\eps\Gamma(0,\alpha_n,0)
=-\frac{\ell_n\big(\p_\eps\Mop(\alpha_n,0)\hat\varphi_n\big)}{\ell_n\big(\p_\alpha\Lop(\alpha_n)\hat\varphi_n\big)} ,
\]
where
\[
\p_\eps\Mop(\alpha,0)=D_H\p_\eps\Fop(\alpha,0;H_0)+D^2_H\Fop(\alpha,0;H_0)\big[H_{T,1}^+(\alpha),\,\cdot\,\big].
\]
All terms are explicit through \eqref{eq:Bform}, \eqref{eq:HT1} and $\hat\varphi_n$. The evaluation
(Lemma~\ref{lem:shiftcomp}) gives $\delta_n+\eta_n$. For $\Fop^-$ the sign of the $\cos2q$-part of
the forcing and of $H_{T,1}$ is reversed. Only the coupling $\cos2q\cdot\cos q\ni\cos q$
(present for $n=1$ alone) is odd in that sign; this gives $\delta_n-\eta_n$.
\end{proof}

For part (iv) we use the bookkeeping of Fourier modes, in the space of $2\pi$-periodic
functions \emph{without} parity.

\begin{lemma}\label{lem:band}
Let $\Xs^\natural$ and $\Ys^\natural$ be defined as $\Xs$ and $\Ys$ but without the parity
restriction, and complexify them. Say that a bounded operator $K$ has \emph{bandwidth} $\le b$
if it maps the Fourier mode $\e^{\mathrm imq}$ into $\spn\{\e^{\mathrm i(m+l)q}:|l|\le b\}$ for every $m$.
We write $\Mop(\alpha,\eps)=D_H\Fop(\alpha,\eps;H_T(\alpha,\eps))$ on $\Xs^\natural$. Then
$\Mop(\alpha,\eps)=\sum_{j\ge0}\eps^j\Mop_j(\alpha)$, with convergence in operator norm for small
$|\eps|$, and $\Mop_j(\alpha)$ has bandwidth $\le2j$; in particular $\Mop_0=\Lop(\alpha)$ is
diagonal.
\end{lemma}

\begin{proof}
$D_H\Fop(\alpha,\eps;H)$ is a linear differential operator. Its coefficients are analytic
functions of $H$ and its derivatives up to order two, and of $\eps\Ptide(q)$. Expanding
$H_T=H_0+\sum_j\eps^jH_{T,j}$ and using Theorem~\ref{thm:tidal}(v), the $\eps^j$-coefficient of
each coefficient function is a trigonometric polynomial of degree $\le2j$. Multiplication by
such a function has bandwidth $\le2j$.
\end{proof}

\begin{proof}[Proof of Theorem~\ref{thm:odd}, (iv)]
Let $n\ge3$ be odd. On $\Xs^\natural$ the kernel of $\Lop(\alpha_n)$ is
\[
K=\operatorname{span}\{\widehat\sigma_n(p)\e^{\mathrm inq},
                         \widehat\sigma_n(p)\e^{-\mathrm inq}\},
\qquad
\widehat\sigma_n(p)=\frac{H_{0,p}(p)\sinh(nH_0(p))}{h_1\sinh(n\gamma)}.
\]
Let $P$ be the projection onto $K$ along the space of functions $\psi$ whose
$\e^{\pm\mathrm inq}$-coefficients at $p=0$ vanish. Let $Q$ be the projection onto
$\spn\{(0,\e^{\pm\mathrm inq})\}$ along $\mathrm R(\Lop(\alpha_n))$. Both projections act only on
the Fourier modes $\pm n$, so they commute with the decomposition into Fourier modes.

For $(\alpha,\eps)$ near $(\alpha_n,0)$, the operator $(I-Q)\Mop(I-P)$ is invertible from
$(I-P)\Xs^\natural$ onto $(I-Q)\Ys^\natural$. Hence $\Mop(\alpha,\eps)$ has a non-trivial kernel if
and only if the $2\times2$ Schur complement $S(\alpha,\eps)=Q\Mop P-Q\Mop(I-P)\big[(I-Q)\Mop(I-P)\big]^{-1}(I-Q)\Mop P$ is singular. Expanding the inverse in a Neumann series around the diagonal operator
$(I-Q)\Mop_0(I-P)$ and using Lemma~\ref{lem:band}, we see that the $\eps^j$-coefficient of $S$
has bandwidth $\le2j$.

The reflection $\iota$ commutes with $\Mop$ and interchanges the modes $\e^{\pm\mathrm inq}$.
Hence, in the basis of $K$ given above,
\[
S=\begin{pmatrix}A&B\\B&A\end{pmatrix}.
\]
The coefficient $B$ maps the mode $-n$ to the mode $+n$, a shift of $2n$, so $B=O(\eps^n)$.

The eigenvectors of $S$ are $\cos nq$ (eigenvalue $A+B$) and $\sin nq$ (eigenvalue $A-B$).
Therefore $\alpha_n^+(\eps)$ and $\alpha_n^-(\eps)$ are the zeros of $A+B$ and $A-B$ near
$\alpha_n$; by Lemma~\ref{lem:sym}(d), $\alpha_n^-$ is the critical value in the $\sin nq$ sector of
$\Mop^+$. Since $\p_\alpha A(\alpha_n,0)\ne0$ by Lemma~\ref{lem:range}(b), the two zeros differ by
$O(B)=O(\eps^n)$.
\end{proof}

\subsection{First-order critical shifts}

\begin{lemma}\label{lem:shiftcomp}
Let $n\ge1$ and let $H_{T,1}^\pm$ be given by \eqref{eq:HT1}, where for $n=2$ only its mean part
$\frac{\e^{\gamma}}{4(\alpha_2-\alpha_0)}\hat\varphi_0$ is kept. Then
\[
-\frac{\ell_n\Big(D_H\p_\eps\Fop^\pm(\alpha_n,0;H_0)\hat\varphi_n+D^2_H\Fop(\alpha_n,0;H_0)\big[H_{T,1}^\pm(\alpha_n),\hat\varphi_n\big]\Big)}{\ell_n\big(\p_\alpha\Lop(\alpha_n)\hat\varphi_n\big)}
=\delta_n\pm\eta_n ,
\]
with $\delta_n$ and $\eta_n$ as defined in Section~\ref{sec:odd}.
\end{lemma}

\begin{proof}[Direct computation of the shifts]
For a laminar flow of logarithmic depth $h$ at fixed $(p_0,\kappa)$,
write $w(h)=W_s(h)\e^{-h}$ and $\mathcal A_n(h)=\e^h\bigl[1+w(h)^2(1+n\coth(nh))-2\kappa w(h)\bigr].$ Differentiating the flux relation gives $w'(h)=2\kappa-(2+h^{-1})w(h).$ For the constant tidal profile $\Ptide=1/4$, the centrifugal term
in the surface condition is multiplied by $1+\eps/2$.  The
nearby laminar depth satisfies, at $\alpha=\alpha_n$, $h'(0)=\frac{\e^\gamma}{4(\alpha_n-\alpha_0)}.$ Its mode-$n$ critical value is exactly
$\alpha=\mathcal A_n(h)+(\eps/2)\e^h$.  Therefore
\begin{equation}\label{eq:delta-direct}
\delta_n=\frac{\e^\gamma}{2}
 +\frac{\e^\gamma\mathcal A_n'(\gamma)}
        {4(\alpha_n-\alpha_0)}.
\end{equation}
Here the derivative with respect to $\alpha$ of the laminar
Bernoulli equation vanishes at $h=\gamma$, since its two terms
$2\alpha\e^h$ and $2\alpha\e^\gamma$ cancel there.  Thus the
variation of the critical value does not add another term to
$h'(0)$.  With $L=n\coth(n\gamma)$, direct differentiation gives
\[
\e^{-\gamma}\mathcal A_n'(\gamma)
=1-4\kappa^2+2\kappa w(3+2L+\gamma^{-1})
 +r\bigl[n^2-L^2-(3+2\gamma^{-1})(1+L)\bigr].
\]
Using $\alpha_n-\alpha_0=\e^\gamma r(L-\gamma^{-1})$
in \eqref{eq:delta-direct} gives the displayed formula for $N_n$.

The nonconstant part of the forcing can contribute to the
mode-$n$ projection at first order only when $n=1$, because
$\cos2q\cos nq$ has a mode-$n$ component only in that case.
For $n=1$ put $T=T_1$, and let $E_2$ denote the coefficient
$E_{2k}$ of Lemma~\ref{lem:flat2} with $k=1$.  The coefficient
of $t u\cos q$ in the quadratic surface Bernoulli residual,
normalised by $\e^{2\gamma}$, for the displacement
$t\cos q+u\cos2q$ is $Q_{12}=\frac{r(3T^3+5T^2+2T+2) -2\kappa w(3T^2+3T+1)-T(1-4\kappa^2)}{T} =\frac{4rE_2+2r}{T}.$ This is evaluated at $\alpha=\alpha_1$ and follows by expanding
the two surface conditions to second order, using the harmonic-extension multipliers
$T$ and $2\coth(2\gamma)=T+T^{-1}$.
The mode-$2$ coefficient of the first-order tide at $\alpha_1$
is $-3T/(4r)$ for $\Ptide^+$.  The explicit product of the
forcing and $t\cos q$ contributes $-3/2$ to the same residual.
Hence the anisotropic shift satisfies $\frac{2\eta_1}{\e^\gamma} -\frac{3T}{4r}Q_{12}-\frac32=0, \qquad \eta_1=\frac{3\e^\gamma}{2}(1+E_2).$ Inserting the formula for $E_2$ gives the stated $M_1$.
For $\Ptide^-$ both anisotropic contributions change sign.
For $n=2$ the resonant tide is omitted in the range equation,
and the same constant-profile argument gives $\delta_2$.
\end{proof}

\subsection{Phase selection and the principal resonance}\label{sec:n2}

The tidal equilibrium of Theorem~\ref{thm:tidal} has no analytic $O(\eps)$
continuation through $\alpha=\alpha_2$, because the coefficient of $\cos2q$
in \eqref{eq:HT1} blows up. Resonant solutions instead occur on the
$O(\eps^{1/3})$ amplitude scale. We
therefore perform the reduction at the \emph{unforced} state. Since $\Ptide^+$ is
$\pi$-periodic, we may work in $\Xs_\pi$, where $\ker\Lop(\alpha_2)=\spn\{\hat\varphi_2\}$. Every
solution of $\Fop=0$ in $\Xs$ near $(\alpha_2,H_0,0)$ lies in $\Xs_\pi$: the kernel of $\Lop(\alpha_2)$
in $\Xs$ is the same, so the reductions in $\Xs$ and in $\Xs_\pi$ coincide.

Let $\Zs_2=\{\psi\in\Xs_\pi:\ \mathsf a_2[\psi]=0\}$ and $y_2^*=(0,\cos2q)$, and let $Q_2$ be the
projection of $\Ys_\pi$ onto $\spn\{y_2^*\}$ along $\ker\ell_2$. As in Section~\ref{sec:odd}, there is a
unique analytic $v(a,\alpha,\eps)\in\Zs_2$, defined near $(0,\alpha_2,0)$, with
$(I-Q_2)\Fop(\alpha,\eps;H_0+a\hat\varphi_2+v)=0$. We set
\begin{equation}\label{eq:g2}
g(a,\alpha,\eps):=\frac{\ell_2\big(\Fop(\alpha,\eps;H_0+a\hat\varphi_2+v(a,\alpha,\eps))\big)}{\ell_2\big(\p_\alpha\Lop(\alpha_2)\hat\varphi_2\big)} .
\end{equation}
Near $(\alpha_2,H_0,0)$, the solutions in the even space of $\Fop(\alpha,\eps;H)=0$ are exactly
$H=H_0+a\hat\varphi_2+v(a,\alpha,\eps)$ with $g(a,\alpha,\eps)=0$. Here $a=\mathsf a_2[H-H_0]$ is the
$\cos2\Theta$-coefficient of $\ln(S/S_0)$, and $a>0$ means high water under the satellite.

\begin{proposition}
\label{prop:n2-fullspace}
Let $\Xs^\natural$ denote the space $\Xs$ without its evenness
condition.  There is a neighbourhood of $(\alpha_2,0,H_0)$ such that
every solution of $\Fop(\alpha,\eps;H)=0$ in $\Xs^\natural$ with
$\eps\ne0$ is even and $\pi$-periodic.  Consequently the scalar
reduction \eqref{eq:g2}, including its amplitude normalisation, describes
all nearby solutions with $\eps\ne0$ in the full periodic space.
\end{proposition}

\begin{proof}
Put
\[
e_c=\widehat\sigma_2(p)\cos2q,\qquad
e_s=\widehat\sigma_2(p)\sin2q,\qquad
\widehat\sigma_2(p)=\frac{H_{0,p}(p)\sinh(2H_0(p))}{h_1\sinh(2\gamma)}.
\]
For $u\in\Xs^\natural$, let $a[u]$ and $b[u]$ be the cosine and sine
coefficients of $u(q,0)$ at frequency $2$, each normalised by $1/\pi$.
Thus $P u=a[u]e_c+b[u]e_s$ is an $O(2)$-equivariant projection.
Let $\ell_c=\ell_2$, and define $\ell_s$ by replacing $\cos2q$ by
$\sin2q$ in \eqref{eq:ell}.  With $y_c=(0,\cos2q)$ and
$y_s=(0,\sin2q)$, the projection $Qy=\frac{\ell_c(y)}{\ell_c(y_c)}y_c+ \frac{\ell_s(y)}{\ell_s(y_s)}y_s$ is also $O(2)$-equivariant: the two denominators agree, and the sine
and cosine functionals transform as the corresponding Fourier
coefficients.  The same is true of the normalising constants
\[
d:=\ell_c(\p_\alpha\Lop(\alpha_2)e_c)
   =\ell_s(\p_\alpha\Lop(\alpha_2)e_s)
   =\pi h_1\cc^{-1}\e^\gamma\sinh(2\gamma)>0.
\]
The range equation has a unique analytic solution
$v(a,b,\alpha,\eps)\in\ker P$.  Define $(G_c,G_s)=d^{-1}(\ell_c,\ell_s) \Fop\bigl(\alpha,\eps;H_0+ae_c+be_s+v(a,b,\alpha,\eps)\bigr).$ At $\eps=0$ this reduced map is $O(2)$-equivariant, with the action
on $(a,b)$ given by rotation through twice the physical angle.
Reflection equivariance and analyticity therefore give $(G_c,G_s)(a,b,\alpha,0) =(a,b)\,A(a^2+b^2,\alpha)$ for a real-analytic scalar $A$.  For arbitrary small $\eps$, reflection
$q\mapsto-q$ still gives $G_s=bK(a,b^2,\alpha,\eps)$ with $K$ analytic.  At $\eps=0$, $K=A$, so analytic division by
$\eps$ yields another analytic function $B$ such that $G_c=aK+\eps B.$ The forcing calculation in Theorem~\ref{thm:n2red} gives
$B(0,0,\alpha_2,0)=-3\e^\gamma/4\ne0$.  Shrink the neighbourhood so that
$B$ never vanishes.  If $\eps\ne0$ and $b\ne0$, the equation $G_s=0$ implies $K=0$, whereas $G_c=\eps B\ne0$, a contradiction.
Thus every zero has $b=0$.

Uniqueness in the range equation then gives
$\iota v(a,0,\alpha,\eps)=v(a,0,\alpha,\eps)$.  The same uniqueness,
using the symmetry $\tau_\pi$ which fixes both kernel vectors, gives
$\tau_\pi v=v$.  Hence the whole solution is even and $\pi$-periodic.
At $b=0$, the full range equation restricts to the original scalar
range equation and $G_c=g$ with exactly the normalisation
\eqref{eq:g2}.
\end{proof}

\begin{theorem}\label{thm:n2red}
The function $g$ is real-analytic near $(0,\alpha_2,0)$, and the following hold.
\begin{enumerate}[label=\textup{(\alph*)}]
\item $g(0,\alpha,0)=0$ and $g(-a,\alpha,0)=-g(a,\alpha,0)$.
\item $g_a(0,\alpha_2,0)=0$, $g_{a\alpha}(0,\alpha_2,0)=1$ and $g_{aaa}(0,\alpha_2,0)=-6\alpha_2^{(2)}$.
\item $g_\eps(0,\alpha_2,0)=-\beta$ with $\beta:=\tfrac34\e^{\gamma}$; $g_{a\eps}(0,\alpha_2,0)=-\delta_2$, where
$\delta_2$ is given by the formula for $\delta_n$ in Section~\ref{sec:odd} with $n=2$; and
$g_{\alpha\eps}(0,\alpha_2,0)=0$.
\end{enumerate}
Consequently, in the region $|a|\le C\eps^{1/3}$, $|\alpha-\alpha_2|\le C\eps^{2/3}$ one has
\begin{equation}\label{eq:normalform}
g(a,\alpha,\eps)=(\alpha-\alpha_2-\delta_2\eps)\,a-\alpha_2^{(2)}a^3-\beta\eps+O(\eps^{5/3}).
\end{equation}
\end{theorem}
\begin{figure}[htbp]\centering
\includegraphics[width=.91\linewidth]{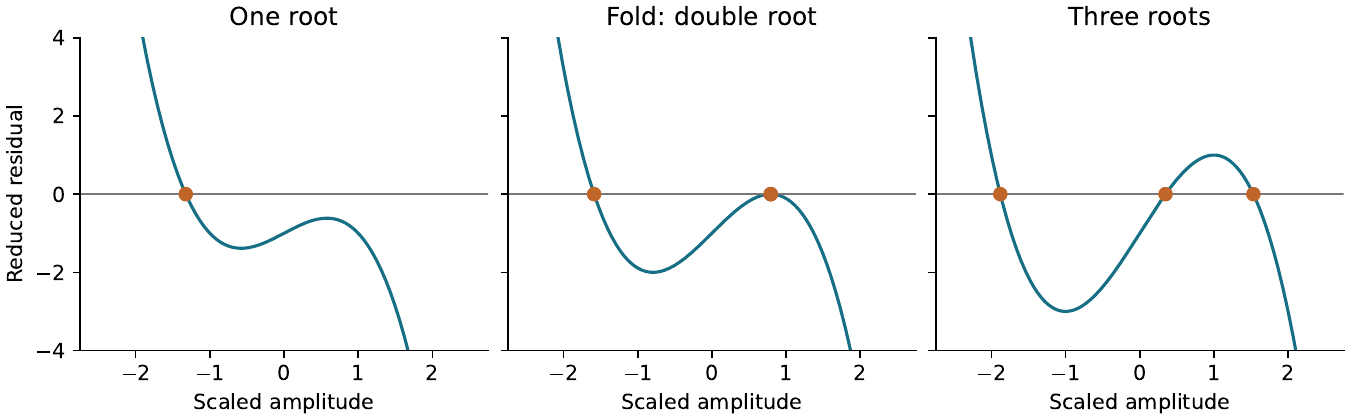}
\caption{Root count in the leading normal form \eqref{eq:normalform} for $c=\alpha_2^{(2)}>0$. Set $a=(\beta\eps/c)^{1/3}X$ and $M=(\alpha-\alpha_2-\delta_2\eps)/(c^{1/3}(\beta\eps)^{2/3})$. The panels plot $MX-X^3-1$ for $M=1,3/2^{2/3},3$. Dots are distinct roots; at the fold the positive root is double. The plotted threshold is for the scaled cubic.}\label{fig:atlas-roots}
\end{figure}

\begin{proof}
(a) For $\eps=0$, $a=0$ gives the solution $H_0$, so $v(0,\alpha,0)=0$ and $g(0,\alpha,0)=0$. The
translation $\tau_{\pi/2}$ acts on $\Xs_\pi$, commutes with $\Fop(\alpha,0;\cdot)$, fixes $H_0$, and
satisfies $\tau_{\pi/2}\hat\varphi_2=-\hat\varphi_2$, $\tau_{\pi/2}\Zs_2=\Zs_2$, $\ell_2\circ\tau_{\pi/2}=-\ell_2$
and $\tau_{\pi/2}y_2^*=-y_2^*$. Hence $v(-a,\alpha,0)=\tau_{\pi/2}v(a,\alpha,0)$ and $g$ is odd in $a$
when $\eps=0$.

(b) The identities $g_a=0$ and $g_{a\alpha}=1$ follow as in the proof of
Theorem~\ref{thm:odd}. For $\eps=0$ the zero set of $g/a$ is the curve
$\alpha=\alpha_2+\alpha_2^{(2)}a^2+O(a^4)$ of Theorem~\ref{thm:local}, with the same amplitude. Writing $g=a\,\Gamma(a^2,\alpha)$ we obtain $\Gamma_\sigma=-\alpha_2^{(2)}\Gamma_\alpha=-\alpha_2^{(2)}$ at the
origin, i.e.\ $g_{aaa}=6\Gamma_\sigma=-6\alpha_2^{(2)}$.

(c) We have $\p_\eps\Fop(\alpha,0;H_0)=(0,-2\e^{2\gamma}\Ptide^+)$. Since $\ell_2$ only sees the
$\cos2q$-component, and $\ell_2(\p_\alpha\Lop\hat\varphi_2)=\pi h_1\cc^{-1}\e^{\gamma}\sinh2\gamma$,
\[
g_\eps(0,\alpha_2,0)=\frac{\tfrac{h_1}{2\cc}\int_\T(-2\e^{2\gamma})\tfrac34\cos^2(2q)\sinh(2\gamma)\dd q}{\pi h_1\cc^{-1}\e^{\gamma}\sinh2\gamma}=-\tfrac34\e^{\gamma}.
\]
Next, $v_\eps:=\p_\eps v(0,\alpha_2,0)\in\Zs_2$ solves
$(I-Q_2)\big[\Lop(\alpha_2)v_\eps+\p_\eps\Fop(\alpha_2,0;H_0)\big]=0$. Because $y_2^*=(0,\cos2q)$,
the projection $(I-Q_2)$ removes the $\cos2q$-part of $\p_\eps\Fop$ entirely. Therefore $v_\eps$ is the
mean tide $\frac{\e^{\gamma}}{4(\alpha_2-\alpha_0)}\hat\varphi_0$, i.e.\ the $m=0$ part of \eqref{eq:HT1}.

Then $g_{a\eps}(0,\alpha_2,0)$ equals $\ell_2\big(D_H\p_\eps\Fop\,\hat\varphi_2+D_H^2\Fop[v_\eps,\hat\varphi_2]\big)$
divided by the normalisation; the term $\ell_2(\Lop\,\p_a\p_\eps v)$ vanishes. This is the same
expression as in the proof of Theorem~\ref{thm:odd}(iii), with $n=2$ and with only the mean part
of the tide; the $\cos2q$-part of $\Ptide^+$ multiplied by $\hat\varphi_2$ produces only the modes
$0$ and $4$. Lemma~\ref{lem:shiftcomp} gives the value $-\delta_2$.

Finally, $\p_\alpha\p_\eps g(0,\alpha_2,0)=\ell_2\big(\p_\alpha\Lop(\alpha_2)v_\eps\big)/\ell_2(\p_\alpha\Lop\hat\varphi_2)=0$,
since $\p_\alpha\Lop\,v_\eps=(0,2\e^{\gamma}v_\eps(\cdot,0))$ is constant in $q$.

For \eqref{eq:normalform}, give $a$, $\alpha-\alpha_2$ and $\eps$ the weights $1$, $2$ and $3$. By
(a)--(c), every Taylor monomial of $g$ of weighted degree $\le4$ appears in
\eqref{eq:normalform}; all remaining monomials have weighted degree $\ge5$.
\end{proof}

\begin{theorem}\label{thm:n2}
Assume that $\alpha_2^{(2)}\neq0$. There are a neighbourhood $N$ of $(0,\alpha_2)$ and $\eps_0>0$ such that
for every $\eps\in(0,\eps_0)$ the following hold.
\begin{enumerate}[label=\textup{(\roman*)}]
\item \emph{(Structure.)} The zero set of $g(\cdot,\cdot,\eps)$ in $N$ is the disjoint union of
two connected analytic curves $\mathcal C_1^\eps$ and $\mathcal C_2^\eps$.
\begin{itemize}
\item $\mathcal C_1^\eps$ is a graph $a=a_1(\alpha)$ over the $\alpha$-interval of $N$.
\item $\mathcal C_2^\eps$ has exactly one fold point $(a_{\rm SN}(\eps),\alpha_{\rm SN}(\eps))$,
which is non-degenerate.
\item $\mathcal C_2^\eps$ meets $\{\alpha=\alpha'\}$ in two points if
$\alpha_2^{(2)}(\alpha'-\alpha_{\rm SN})>0$, and in no point if $\alpha_2^{(2)}(\alpha'-\alpha_{\rm SN})<0$.
\end{itemize}
Hence the equation has three solutions on one side of $\alpha_{\rm SN}(\eps)$ and one on the
other side.
\item \emph{(Fold.)} As $\eps\downarrow0$,
\[
\alpha_{\rm SN}(\eps)=\alpha_2+3\big(\alpha_2^{(2)}\big)^{1/3}\Big(\frac{\beta\eps}{2}\Big)^{2/3}+\delta_2\eps+O(\eps^{4/3}),\qquad
a_{\rm SN}(\eps)=\Big(\frac{\beta\eps}{2\alpha_2^{(2)}}\Big)^{1/3}+O(\eps),
\]
with $\beta=\tfrac34\e^{\gamma}$ and real cube roots.
\item \emph{(Away from the fold.)} Let $\alpha_2^{(2)}>0$ and $\theta>0$ with
$\theta\to0$ and $\eps/\theta^3\to0$.
\begin{itemize}
\item If $\alpha=\alpha_2+\theta^2$, the three solutions satisfy
$a=\beta\eps\theta^{-2}(1+o(1))$ (the direct tide, on $\mathcal C_2^\eps$) and
$a=\pm\theta(\alpha_2^{(2)})^{-1/2}(1+o(1))$. Here the sign $+$ lies on $\mathcal C_2^\eps$ and
the sign $-$ on $\mathcal C_1^\eps$.
\item If $\alpha=\alpha_2-\theta^2$, the unique solution satisfies
$a=-\beta\eps\theta^{-2}(1+o(1))$ (the inverted tide, on $\mathcal C_1^\eps$).
\end{itemize}
For $\alpha$ in a compact subset of the $\alpha$-interval of $N$ avoiding $\alpha_2$ and $\eps$ small, the solution with
$a\to0$ is the tidal equilibrium $H_T^+(\alpha,\eps)$ of Theorem~\ref{thm:tidal}.
\end{enumerate}
\end{theorem}

Thus, for $\alpha_2^{(2)}>0$ (which, by Theorem~\ref{thm:kge2}, holds
whenever $\gamma<1.655$ if $\kappa=0$, and for every admissible $\kappa$ whenever $\gamma<1.357$), we have
the following picture.
\begin{itemize}
\item As $\alpha$ increases past $\alpha_2$, the inverted tide on $\mathcal C_1^\eps$ continues locally to a larger-amplitude $\pi$-periodic wave with a \emph{trough} under the satellite.
\item The direct tide exists only for $\alpha>\alpha_{\rm SN}(\eps)$. As $\alpha$ decreases it grows
and, at the fold, merges with the larger-amplitude local wave with a \emph{crest} under the satellite.
\item The resonant tide saturates at amplitude $O(\eps^{1/3})$, instead of diverging like
$\eps/|\alpha-\alpha_2|$. The three-solution count holds in the local neighbourhood of the theorem.
It is a stationary multiplicity result; dynamical stability and hysteresis are not established here.
\end{itemize}
The case $\alpha_2^{(2)}<0$ is the mirror image, with the fold on the side $\alpha<\alpha_2$.

\begin{proof}
(i) Write $\mu=\alpha-\alpha_2$, $c=\alpha_2^{(2)}\ne0$, and regard
$g$ as a function of $(a,\mu,\eps)$.  At the origin, $g_{a\mu}=1,\qquad g_{aaa}=-6c,\qquad g_\eps=-\beta<0.$ The equation $g_a=0$ can therefore be solved uniquely as
$\mu=m(a,\eps)$, with $m(a,0)=3ca^2+O(a^4)$.  On this critical
set put $h(a,\eps)=g(a,m(a,\eps),\eps)$.  Then $h(a,0)=2ca^3+O(a^5),\qquad h_\eps(0,0)=-\beta.$ A second application of the implicit function theorem solves $h=0$
as $\eps=E(a)=\frac{2c}{\beta}a^3+O(a^4),\qquad E'(a)=\frac{6c}{\beta}a^2+O(a^3).$ After shrinking the neighbourhood, $E'$ has the sign of $c$ for
$a\ne0$.  Thus each sufficiently small $\eps>0$ has exactly one
fold point, and its amplitude has the sign of $c$.  At that point
$g_\mu=a(1+o(1))\ne0$ and $g_{aa}=-6ca(1+o(1))\ne0$, so the
fold is non-degenerate.

For completeness choose a small rectangle
$|a|<a_*$, $|\mu|<\mu_*$ with $0<\mu_*\ll |c|a_*^2$.
On its two amplitude boundaries $g$ has the fixed signs of
$-ca^3$, uniformly for sufficiently small $\eps$.  Also
$g_{aaa}$ has the constant sign of $-c$, so Rolle's theorem allows
at most three zeros at each fixed $(\mu,\eps)$.  At
$\mu=-\operatorname{sgn}(c)\mu_*$ the derivative $g_a$ has constant
sign and there is one zero.  At
$\mu=\operatorname{sgn}(c)\mu_*$ the three simple unforced zeros
persist for small $\eps$.  Zeros cannot leave through an amplitude
boundary, and their number can change only at the unique fold.
There are therefore three zeros exactly when $c(\alpha-\alpha_{\rm SN}(\eps))>0,$ one when this product is negative, and two distinct zeros at the
fold, one of which is double.

Finally, $g(0,\mu,\eps)<0$ throughout the rectangle for
$\eps>0$, so no branch crosses $a=0$.  The branch containing the
fold has $\operatorname{sgn}a=\operatorname{sgn}c$; the other
branch has the opposite sign and is a graph over the whole
parameter interval.  Away from the fold both are analytic graphs
by the implicit function theorem.  At the fold $g_\mu\ne0$ gives
an analytic parametrisation by $a$.  These facts prove the stated
two-component local diagram, without a recognition theorem or a
universal-unfolding argument.

(ii) Put $\eta=\eps^{1/3}$, $a=\eta A$, $\alpha=\alpha_2+\eta^2\Lambda$, and $G(A,\Lambda,\eta):=\eta^{-3}g(\eta A,\alpha_2+\eta^2\Lambda,\eta^3).$ By the weighted-degree argument in the proof of Theorem~\ref{thm:n2red}, $G$ is analytic near
$\eta=0$ and $G(A,\Lambda,\eta)=\Lambda A-\alpha_2^{(2)}A^3-\beta-\delta_2\eta A+O(\eta^2).$ At $\eta=0$ the system $G=G_A=0$ has the unique real solution
$A_0=(\beta/(2\alpha_2^{(2)}))^{1/3}$, $\Lambda_0=3\alpha_2^{(2)}A_0^2$. The Jacobian of $(G,G_A)$
with respect to $(A,\Lambda)$ at that point is
$\left(\begin{smallmatrix}0&A_0\\-6\alpha_2^{(2)}A_0&1\end{smallmatrix}\right)$, with determinant
$6\alpha_2^{(2)}A_0^2\neq0$. The implicit function theorem gives an analytic fold curve
$(A(\eta),\Lambda(\eta))$. Differentiating $G(A(\eta),\Lambda(\eta),\eta)=0$ at $\eta=0$, we obtain
$A_0\Lambda'(0)-\delta_2A_0=0$, i.e.\ $\Lambda'(0)=\delta_2$.
Differentiating $G_A(A(\eta),\Lambda(\eta),\eta)=0$ also, we obtain $-6\alpha_2^{(2)}A_0A'(0)+\Lambda'(0)-\delta_2=0.$ Thus $A'(0)=0$ and $a_{\rm SN}=\eta A_0+O(\eta^3)$,
which gives the improved amplitude remainder $O(\eps)$.
Since
$\alpha_{\rm SN}=\alpha_2+\eta^2\Lambda(\eta)$ and $3\alpha_2^{(2)}A_0^2=3(\alpha_2^{(2)})^{1/3}(\beta/2)^{2/3}$,
the formula follows. The point is non-degenerate because $G_\Lambda=A_0\neq0$ and
$G_{AA}=-6\alpha_2^{(2)}A_0\neq0$.

(iii) Put $a=\theta\hat A$, $\alpha=\alpha_2\pm\theta^2$, $\eps=\nu\theta^3$, and $K(\hat A,\theta,\nu):=\theta^{-3}g(\theta\hat A,\alpha_2\pm\theta^2,\nu\theta^3).$ As before, $K$ is analytic, and $K(\hat A,0,0)=\pm\hat A-\alpha_2^{(2)}\hat A^3$. Moreover
$K(\hat A,\theta,0)$ is odd in $\hat A$, so $K(0,\theta,0)=0$.

In the case $+$ the three roots $\hat A=0,\pm(\alpha_2^{(2)})^{-1/2}$ of $K(\cdot,0,0)$ are simple
and persist for small $(\theta,\nu)$. Since $K(\cdot,\theta,0)$ is odd and
$K_\nu(0,0,0)=g_\eps(0,\alpha_2,0)=-\beta$, we have
$K(\hat A,\theta,\nu)=\hat A\big(1+O(\theta+\hat A^2)\big)-\beta\nu+O\big(\nu(\theta+|\hat A|+\nu)\big)$,
so the root near $0$ satisfies $\hat A=\beta\nu(1+O(\theta+\nu))$.
In the case $-$, $\hat A=0$ is the only real root of $K(\cdot,0,0)$, and the count in (i)
shows that there is no other solution in $N$. Undoing the scaling gives the asymptotics.

The solutions with $a\to0$ converge to $H_0$ as $\eps\to0$. For $\alpha$ in a compact set avoiding
$\alpha_2$, they therefore coincide with $H_T^+(\alpha,\eps)$ by the uniqueness in
Theorem~\ref{thm:tidal}(ii).

It remains to assign the branches. Since $g(0,\alpha,0)=0$ and
$\p_\eps g(0,\alpha,0)=-\beta+O(|\alpha-\alpha_2|)$, after shrinking $N$ we have $g(0,\alpha,\eps)<0$ for
$0<\eps<\eps_0$ and all $(0,\alpha)\in N$. Hence no solution has $a=0$, and $a$ has a constant sign
along each of the connected curves $\mathcal C_1^\eps$ and $\mathcal C_2^\eps$.
\begin{itemize}
\item $\mathcal C_2^\eps$ contains the fold point, where $a_{\rm SN}>0$, so $\mathcal C_2^\eps\subset\{a>0\}$.
\item $\mathcal C_1^\eps$ contains the inverted tide ($a<0$), so $\mathcal C_1^\eps\subset\{a<0\}$.
\end{itemize}
This gives the assignment stated in (iii).
\end{proof}
\begin{corollary}
For $k=2$, $\kappa=0$ and $r=16537/125000$, the cubic coefficient has exactly one zero on $[1.655,1.656]$. It is simple and lies between $1.655331417$ and $1.655331418$.
\end{corollary}
\begin{proof}
Theorem~\ref{thm:kge2} and its subcritical witness give opposite endpoint signs. Put $F(x)=P_2(0,\sqrt r)$ at $\gamma=x$. Substituting $T=\coth(2x)$ in \eqref{eq:Pk} and differentiating with $T'=2(1-T^2)$ yields $-5.367<F'(x)<-5.364$ on the whole interval by the 200-bit interval evaluation in the supplement (\texttt{degenerate\_depth.py}). The same evaluation gives opposite strict signs at the two narrower endpoints. Equation~\eqref{eq:alpha2} has positive prefactor there. The local diagram at this depth requires higher-order terms.
\end{proof}

For an even mode $n=2m\ge4$, the same bandwidth argument permits the first amplitude-independent resonant forcing at order $\eps^m$. Determining its coefficient is a separate problem.

\subsection{Weak forcing along the global curves}\label{sec:persist}
Let $(\gamma,\kappa,w)$ be admissible, and let $\mathcal K_k$ be the global curves of
Theorem~\ref{thm:globalkappa} (Theorem~\ref{thm:global} for $\kappa=0$). Theorems~\ref{thm:odd}
and~\ref{thm:n2} describe the forced waves near the laminar flows. Away from them, weak forcing
perturbs the global curves regularly, except near secondary bifurcation points. Let $k\ge1$ and put $\Xs^{(k)}:=\Xs$, $\Ys^{(k)}:=\Ys$ if $k$ is odd, and
$\Xs^{(k)}:=\Xs_\pi$, $\Ys^{(k)}:=\Ys_\pi$ if $k$ is even. Since $\Ptide^\pm$ is $\pi$-periodic,
$\Fop^\pm(\alpha,\eps;\cdot)$ maps $\Xs^{(k)}\cap\Uo$ into $\Ys^{(k)}$, and $\mathcal K_k\subset\R\times
\Xs^{(k)}$. We call a point $(\alpha,H)\in\mathcal K_k$ \emph{regular} if
$D_{(\alpha,H)}\Fop(\alpha,0;H):\R\times\Xs^{(k)}\to\Ys^{(k)}$ is onto.

\begin{proposition}\label{prop:persist}
Let $k\ge1$, and let $\mathcal C$ be a compact connected subset of $\mathcal K_k$ that consists of
regular points. Then there are $\eps_2>0$ and a neighbourhood $\mathcal W'$ of $\mathcal C$ in
$\R\times\Xs^{(k)}$ such that, for $|\eps|<\eps_2$ and either sign, the zeros of
$\Fop^\pm(\cdot,\eps;\cdot)$ in $\mathcal W'$ form a real-analytic curve, which depends analytically on
$\eps$ and equals $\mathcal K_k\cap\mathcal W'$ for $\eps=0$.
\end{proposition}

\begin{proof}
The derivative $D_H\Fop(\alpha,0;H):\Xs^{(k)}\to\Ys^{(k)}$ is Fredholm of index zero at every point of
$\mathcal K_k$, by the proof of Lemma~\ref{lem:fredholm-global} (Lemma~\ref{lem:kstructure} for
$\kappa\ne0$), so $D_{(\alpha,H)}\Fop$ has index one.
At a regular point it is onto, so its kernel is one-dimensional. By the implicit function theorem
the zeros of $\Fop(\cdot,0;\cdot)$ near such a point form a real-analytic curve, which contains the
piece of $\mathcal K_k$ through the point and therefore coincides with it there; the kernel is its
tangent line. The implicit function theorem, applied to $(\alpha,H,\eps)\mapsto\Fop^\pm(\alpha,\eps;H)$
on a complement of this tangent line, shows that each point of $\mathcal C$ has a neighbourhood
$\mathcal W_0$ and an $\eps_0>0$ such that for $|\eps|<\eps_0$ the zeros in $\mathcal W_0$ form an
analytic curve depending analytically on $\eps$. By compactness finitely many such neighbourhoods cover $\mathcal C$, and by
the uniqueness part of the implicit function theorem the local curves agree on overlaps.
\end{proof}

By local uniqueness these continuations agree on overlaps with regular portions of the branches in Theorems~\ref{thm:odd} and~\ref{thm:n2}. The size of the admissible forcing depends on the compact segment under consideration.

The neighbourhood of a single-mode reduction also depends on the spectral gaps. For example, $\alpha_2-\alpha_1=S_0r\tanh\gamma$. The fold scale in Theorem~\ref{thm:n2} can be comparable to this gap in a thin layer.
\section{The canal limit and computed waves}\label{sec:earth}

The results of Sections~\ref{sec:linear}, \ref{sec:unforced} and
\ref{sec:tidal}--\ref{sec:n2} hold for every admissible constant vorticity. To illustrate these results, we use the astronomical parameters tabulated by NASA
\cite{NASAearth,NASAmoon}. We take $\Omega_f$ to be the Moon's sidereal angular velocity,
$2.662\times10^{-6}\,$s$^{-1}$, and let the ocean co-rotate with the Earth, so that
$\omega_v=7.292\times10^{-5}\,$s$^{-1}$ and $\mu=0$ in \eqref{eq:laminar}. The relative flow is then
the solid-body rotation $W_0=\kappa R$. Thus $w=\kappa$, the triple $(\gamma,\kappa,\kappa)$ is
admissible for every $\gamma>0$, and by Proposition~\ref{prop:dispersion}
\begin{equation}\label{eq:alphamu0}
\alpha_k=S_0\big[1+\kappa^2(kT_k-1)\big]\quad(k\ge1),\qquad \alpha_0=S_0\big[1+\kappa^2(\gamma^{-1}-1)\big].
\end{equation}

\begin{figure}[htbp]\centering
\includegraphics[width=.95\linewidth]{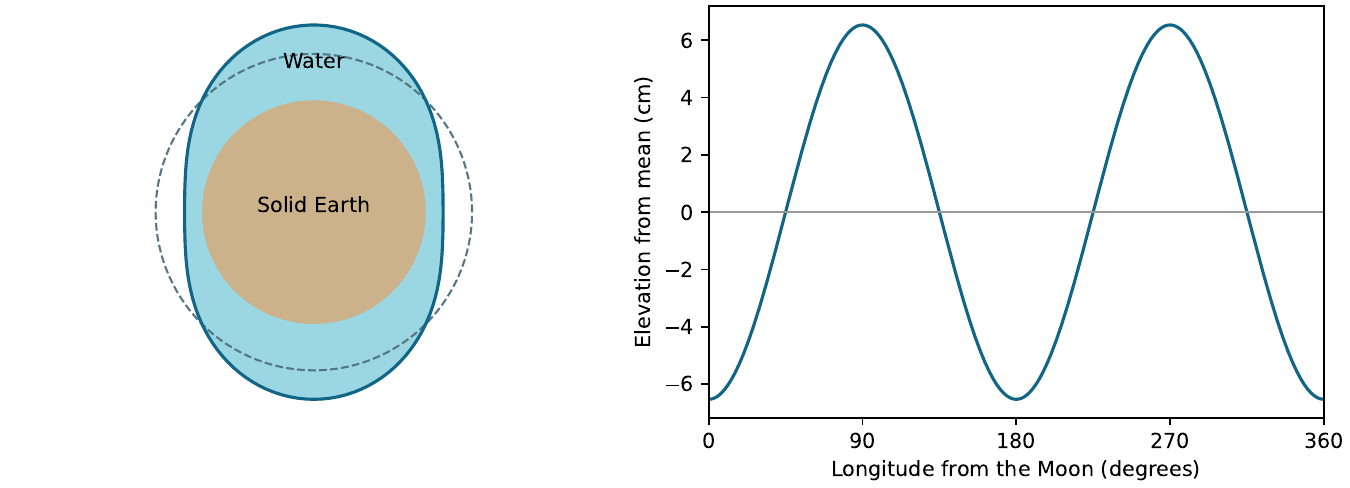}
\caption{Earth--Moon canal at 4 km depth. Left: solid bed (tan), water (blue), mean surface (dashed) and computed free surface (solid). The water depth and its tidal variation are separately exaggerated. Right: the same solution in physical elevation units. The sublunar and antipodal points are troughs.}\label{fig:canal}
\end{figure}

\begin{example}\label{ex:earth}
With $a=6371\,$km one has
\[
\kappa=26.40,\qquad W_s=\kappa S_0,\qquad \alpha=2.17\times10^5,\qquad
\eps=0.0122\approx M/(M_\oplus+M).
\]
In the thin-layer limit, $\alpha_2\simeq\kappa^2/\gamma$ by \eqref{eq:alphamu0}. Resonance occurs when
$gD=\big((\omega_v-\Omega_f)a\big)^2$, i.e.\ at the depth $D_{\rm res}\approx20.4\,$km.

For an ocean of depth $D=4\,$km one finds $\alpha/\alpha_2=0.196$. The dynamic factor of
Remark~\ref{rem:airy} is then $-0.243$, so the $\cos2\Theta$-tide is inverted, with amplitude
about $6.5\,$cm (the equilibrium tide is $26.8\,$cm).
\end{example}
Figure~\ref{fig:canal} shows the forced free surface at fixed laminar flux and Bernoulli constant. Fourier collocation on the reflection-symmetric quarter circle gives a second-mode coefficient of $-6.53286$ cm at both 12 and 24 nodes; the linear formula gives $-6.53299$ cm with the same rounded inputs.

For $\mu=0$ and $D\ll a$, one has $\gamma=\log(1+D/a)\sim D/a$. The singularity in \eqref{eq:S1} then reduces to Airy's matching of the long-wave speed $\sqrt{gD}$ and the relative speed $(\omega_v-\Omega_f)a$ \cite{Airy1842,Lamb1932,Tokieda2013}. The $4$ km example lies on the inverted side. Its $6.5$ cm amplitude is the linear response of the equatorial canal. The local bifurcation theorems give an asymptotic small-forcing regime without an explicit threshold at the physical value $\eps=0.0122$.

We next apply the unforced continuation theorem at the same depth, flux and relative vorticity, taking $\eps=0$. For $x\le\gamma$, $\mu(x)\le0$ and $W_s(x)\ge\kappa e^x$, so the completed-square formula in Lemma~\ref{lem:kgalpha}(b) gives $\alpha_0(x)\ge(1-\kappa^2)e^\gamma+\kappa^2/\gamma>1.1\times10^6$. Taking $\alpha_\flat=10^6$ in Corollary~\ref{cor:kthin}, we obtain $A_\flat^{1/2}>31$ while the right-hand side of \eqref{eq:kthin} is below $1.5$. Hence the unforced branch lies in positive gravity and degenerates on a crest line. At this depth $\eta_{\max}<\log z_\flat<0.001<x_b$, so Lemma~\ref{lem:kcrestline} excludes the bed. The same estimates give $\limsup Q_c\le3.12$, about $0.2\%$ of the normalised laminar crest speed. Whether the crest speed actually tends to zero for these parameters remains open.

\subsection{Numerical comparisons}\label{sec:numerics}
We use two independent calculations of the unforced coefficient to illustrate the sign comparison in Section~\ref{sec:unforced}. One solves the conformal-strip Bernoulli equation by Fourier collocation \cite{Babenko1987,Toland1996}. The other uses Fourier--Chebyshev collocation in the height variables of \cite[(2.23)]{LWpreprint}. At $(\gamma,\lambda)=(0.3,1.15)$ they give $22.5072$ and $22.5066$, respectively, for the quadratic gravity coefficient, compared with $22.5060$ from \eqref{eq:alpha2}. At $(0.2,1.4)$ the corresponding values are $2015.146$, $2015.13$ and $2015.148$. The numerical values are extrapolations of $(\alpha-\alpha_1)/a^2$ to zero amplitude.

\begin{figure}[htbp]\centering
\includegraphics[width=\textwidth]{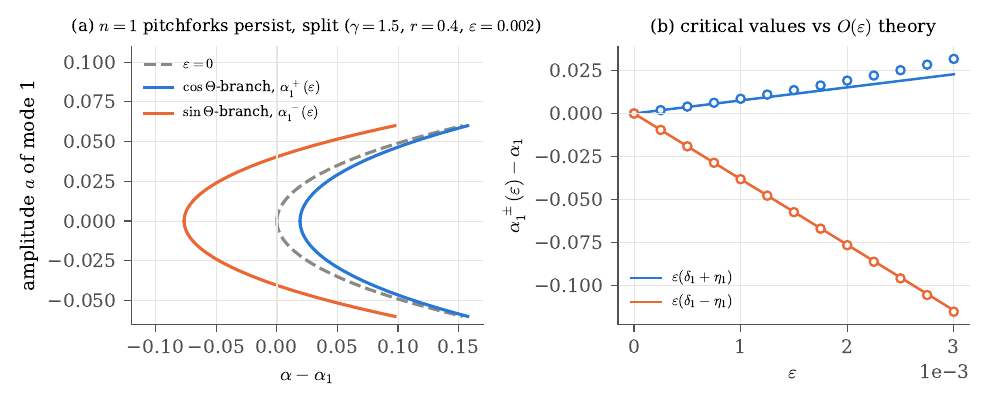}
\caption{Theorem~\ref{thm:odd} for $n=1$ at $\gamma=1.5$, $\kappa=0$, $w=-\sqrt{0.4}$ (the irrotational solver depends only on $w^2$). (a) The $O(2)$-pitchfork at
$\eps=0$ (dashed) splits into the $\cos\Theta$-pitchfork at $\alpha_1^+(\eps)$ and the
$\sin\Theta$-pitchfork at $\alpha_1^-(\eps)$; here $\eps=2\times10^{-3}$. (b) $\alpha_1^\pm(\eps)-\alpha_1$:
conformal solver (circles) against $\eps(\delta_1\pm\eta_1)$ (lines).}\label{fig:n1}
\end{figure}

\begin{figure}[htbp]\centering
\includegraphics[width=\textwidth]{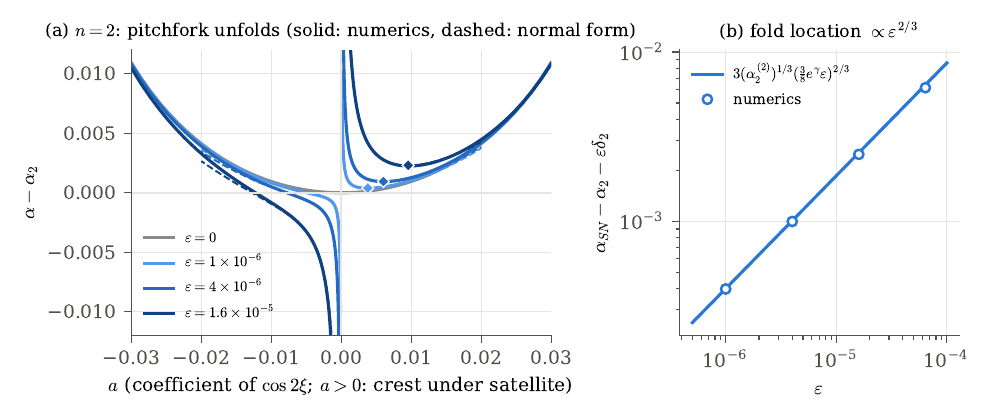}
\caption{Theorem~\ref{thm:n2} at the parameters of \cite[Example~6.2]{LWpreprint}. (a) Solution
curves for several $\eps$: conformal solver (solid) and the normal form \eqref{eq:normalform}
(dashed); the folds are marked. (b) Fold location versus $\eps$, with slope $2/3$.}\label{fig:n2}
\end{figure}

Figure~\ref{fig:n1} resolves the two reflection branches at the first odd mode. At the principal resonance, Figure~\ref{fig:n2} shows how the folded branch approaches the predicted $\eps^{2/3}$ law. For the parameters of \cite[Example~6.2]{LWpreprint} and $\eps=10^{-6}$ the computed fold is at $\alpha=1.5197652$, compared with $1.5197670$ from the asymptotic formula.
The conformal solver records Fourier coefficients on its native angle $\xi$. If $X(\xi)=\Theta$ and $Y(\xi)=\log R$, the physical amplitude is
\[
 a_k=\frac1\pi\int_0^{2\pi}(Y(\xi)-\gamma)\cos(kX(\xi))X'(\xi)\,d\xi .
\]
The native coefficient uses $\cos(k\xi)$ and unit weight instead. The two agree to first order at the laminar state. The fold gravity is independent of this change of amplitude. As a separate check at $\eps=10^{-6}$, physical-angle collocation at $N=8,12,16,24$ gives $\alpha_{\rm SN}=1.519765192908$ and $a_2\simeq0.003801$. On $8N+1$ independent angles the largest Bernoulli residual among these runs was below $3\times10^{-15}$. The supplement includes the amplitude minimization and residual checks.

\begin{figure}[htbp]\centering
\includegraphics[width=\textwidth]{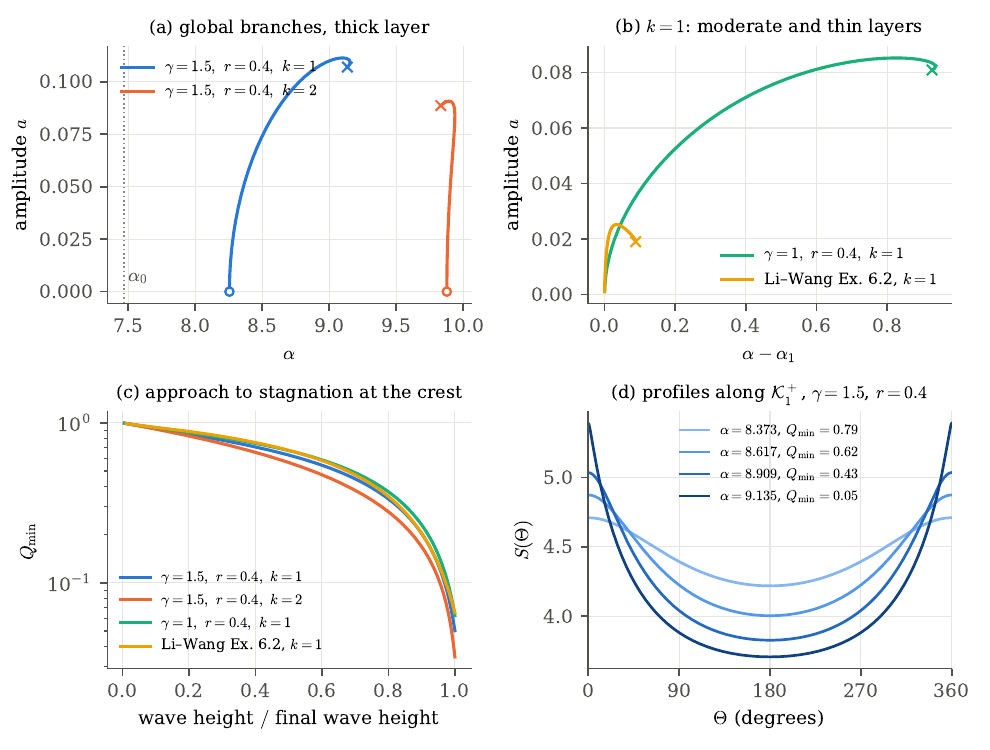}
\caption{Global branches of Theorem~\ref{thm:global}. (a) $\alpha$ and the amplitude $a$
($\cos k\xi$-coefficient of $\ln R$) for $\gamma=1.5$, $r=0.4$, $k=1,2$; circles mark $\alpha_k$,
crosses the last computed point. (b) The same for $k=1$ at $\gamma=1$, $r=0.4$ and at the
parameters of \cite[Example~6.2]{LWpreprint}, against $\alpha-\alpha_1$. (c) $Q_{\min}$, attained at the
crest, against the wave height divided by its final value. (d) Free surfaces along
$\mathcal K_1^+$ for $\gamma=1.5$, $r=0.4$.}\label{fig:global}
\end{figure}

Figure~\ref{fig:global} follows the unforced irrotational curve. The last computed state on each curve is regular and has positive crest speed. The plotted depth $\gamma=1.5$ exceeds the explicit depth range of Theorem~\ref{thm:stagCAP}. In these computations the profile steepens near the crest as that speed decreases.

\begin{figure}[htbp]\centering
\includegraphics[width=.93\linewidth]{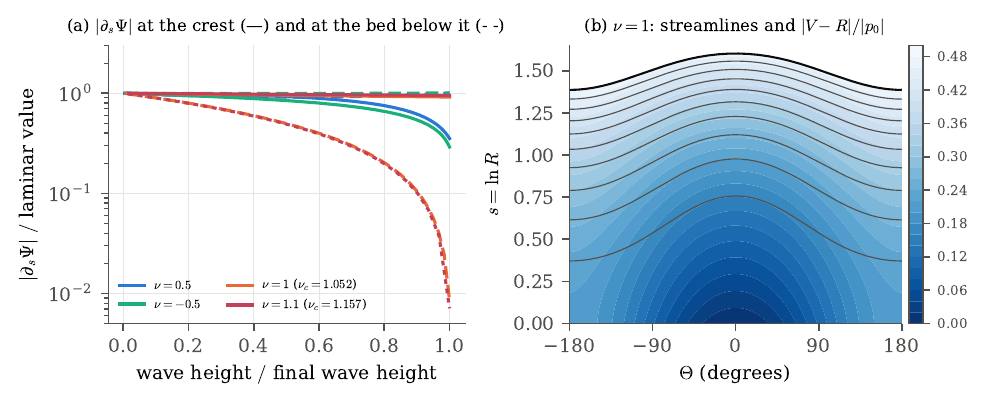}
\caption{Computed rotational branches. Left: normalised crest and bed values of $|\Psi_s|$ for $(\gamma,\kappa,w)=(1,0.4,0.8)$, $(1.5,-0.3162,0.6325)$, $(1.5,0.6325,0.6325)$ and $(1,0.66,0.6)$, as labelled in the plot. Right: streamlines and the normalised azimuthal speed in a nearly critical bed flow at $\gamma=1.5$, $\kappa/w=1$.}\label{fig:kglobal}
\end{figure}

For general vorticity, $\Phi=\Psi+\kappa e^{2s}/(2p_0)$ is harmonic. In conformal coordinates its surface normal derivative follows from the strip Dirichlet-to-Neumann map, so the same Fourier discretisation gives the azimuthal velocity in the layer. Figure~\ref{fig:kglobal} compares unforced profiles whose smallest azimuthal speed occurs at the surface crest or at the bed underneath it. The plotted values are from finite-amplitude continuation and do not determine the limiting speed.
\appendix\section{Coefficient estimates for the discriminant}
\subsection{Proof of Theorem~\ref{thm:k1-vorticity}}\label{app:details-k1-vorticity}
\begin{proof}
The terms of degrees four and three in $\nu$ cancel in
$\operatorname{Disc}(\nu)=A_1(\nu)^2-4a_0A_2(\nu)$, giving $\operatorname{Disc}(\nu)=d_2\nu^2+d_1\nu+D,$ where $d_2=144X^2-16a_0a_1-4a_0b_2,\qquad d_1=-24a_1X-4a_0b_1.$ Substitution of the coefficients and $\nu_c=(1+t)/2$ yields
\begin{align*}
d_2&=48(xt-1)Q_2,\\
\operatorname{Disc}(\nu_c)&=(xt-1)Q_D,\\
\operatorname{Disc}'(\nu_c)&=-8(xt-1)Q_S.
\end{align*}
By Lemmas~\ref{lem:k1-aux} and~\ref{lem:k1-Qsigns}, $d_2<0,\qquad\operatorname{Disc}(\nu_c)<0, \qquad\operatorname{Disc}'(\nu_c)>0.$ For $\nu\le\nu_c$, concavity gives
$\operatorname{Disc}'(\nu)\ge\operatorname{Disc}'(\nu_c)>0$, and therefore
$\operatorname{Disc}(\nu)\le\operatorname{Disc}(\nu_c)<0$.
Since $a_0>0$, this implies
$A_2(\nu)>A_1(\nu)^2/(4a_0)\ge0$, so $P_1>0$ for every real $r$.
Equation~\eqref{eq:k1-alpha} now gives $\alpha_1^{(2)}>0$.
\end{proof}
\subsection{Proof of Theorem~\ref{thm:mode-full}}\label{app:details-mode-full}
\begin{proof}
For $0<x\le1$, Lemma~\ref{lem:mode-hyperbolic} implies $0<h<\frac13,\quad 0<z<\frac13,\quad h>\frac{x^2}{4}, \quad T>\frac1x,\quad Th<x,\quad d>\frac23.$ First, $q_2<k\left(\frac23+\frac13-2\right)-x=-k-x<0.$ Each coefficient in \eqref{eq:mode-L} is negative. More precisely,
\[
\begin{aligned}
 T(-16h+24z-48)+32x
 &<-40T+32x<-\frac{40}{x}+32x\le-\frac8x,\\
 8h+28z-56&<\frac83+\frac{28}{3}-56=-44,\\
 24(Th-2x)&<-24x.
\end{aligned}
\]
For the cubic coefficient in \eqref{eq:mode-C},
\[
\begin{split}
 d(12z-24)-12z+8h+8
 &<-20d+8h+8\\
 &<-\frac{40}{3}+\frac83+8=-\frac83.
\end{split}
\]
For its quadratic coefficient,
\[
\begin{aligned}
 T(-28h+20z-40)+40x
 &<T\left(-7x^2-\frac{100}{3}\right)+40x\\
 &<-\frac{100}{3x}+33x\le-\frac1{3x}.
\end{aligned}
\]
The factor multiplying $T$ is negative, so replacing $T$ by its lower bound $1/x$ gives a valid upper bound. Finally, $16z-7h-32<\frac{16}{3}-32=-\frac{80}{3}, \qquad12(Th-2x)<-12x.$ These estimates prove $q_2,L,C<0$ for $0<x\le1$.

We next consider $1\le x\le7/6$. In this interval,
\begin{equation}\label{eq:mode-extended-bounds}
 \frac14<h<\frac{49}{108},\qquad
 \frac49<d<\frac34,\qquad z<\frac13,\qquad T>\frac65.
\end{equation}
The bounds on $h$ follow from its strict increase and
\eqref{eq:mode-hyperbolic-bounds}. For $d$, use
\[
 \sinh1>1+\frac16=\frac76,
 \qquad
 \sinh\frac76<\frac76+
 \frac{(7/6)^3}{6(1-(7/6)^2/20)}
 =\frac{17521}{12078}<\frac32.
\]
For the last bound in \eqref{eq:mode-extended-bounds}, geometric estimates on
Taylor remainders give
\[
 e<\frac83+\frac{1/24}{1-1/5}=\frac{87}{32}<\frac{11}{4},
 \qquad
 e^{1/3}<\frac43+\frac{1/18}{1-1/9}
 =\frac{67}{48}<\frac75.
\]
Thus $e^{7/3}=e^2e^{1/3}<847/80<11$, and
$T=\coth x>6/5$ for $x\le7/6$.

Using \eqref{eq:mode-extended-bounds}, $q_2<k\left(\frac{49}{54}+\frac13-2\right)-x =-\frac{41}{54}k-x<0.$ The three coefficients making up $L$ are negative, since
\[
\begin{aligned}
 T(-16h+24z-48)+32x&<-40+\frac{112}{3}=-\frac83,\\
 8h+28z-56&<\frac{98}{27}+\frac{28}{3}-56<0,\\
 Th-2x&<-x.
\end{aligned}
\]
To control the cubic coefficient in $C$, observe that the convex quadratic
$12d^2-12d+8$ has maximum $23/4$ on $[4/9,3/4]$. Hence
\[
\begin{split}
 d(12z-24)-12z+8h+8
 &=h(12d^2-12d+8)-24d+8\\
 &<\frac{49}{108}\frac{23}{4}-\frac{32}{3}+8
 =-\frac{25}{432}<0.
\end{split}
\]
For the quadratic coefficient,
\[
\begin{split}
 T(-28h+20z-40)+40x
 &<\frac65\left(-7+\frac{20}{3}-40\right)
       +40\frac76\\
 &=-\frac{26}{15}<0.
\end{split}
\]
The other coefficients remain negative because
$16z-7h-32<16/3-32<0$ and $Th-2x<-x$.
We finally consider $7/6\le x\le6/5$. On this interval,
\begin{equation}\label{eq:mode-final-bounds}
 \frac25<h<\frac9{20},\qquad
 \frac37<d<\frac12,\qquad z<\frac13,\qquad T>\frac{19}{16}.
\end{equation}
Indeed, $h$ is increasing and
$h(7/6)>(7/6)(6/5)-1=2/5$. The Taylor lower bound $e^{12/5}>\sum_{j=0}^{5}\frac{(12/5)^j}{j!} =\frac{166093}{15625}>\frac{53}{5}$ implies $T(6/5)<29/24$ and $h(6/5)<9/20$. For $d$, use
\[
 \sinh\frac76>\frac76+\frac{(7/6)^3}{6}
 =\frac{1855}{1296}>\sqrt2,
 \qquad
 \sinh\frac65<\frac65+
 \frac{(6/5)^3}{6(1-(6/5)^2/20)}
 =\frac{219}{145}<\sqrt{\frac73}.
\]
For the lower bound on $T$, the geometric Taylor estimate $e^{2/5}<\frac75+\frac{2/25}{1-2/15} =\frac{97}{65}<\frac32$ together with $e<11/4$ yields
$e^{12/5}<363/32<35/3$, hence $T>19/16$.

Now $q_2<-23k/30-x<0$. The quadratic coefficient of $L$ is bounded above by
$-40+192/5=-8/5$, and its other coefficients are negative as before.
For the cubic coefficient of $C$, \eqref{eq:mode-final-bounds} gives
\[
\begin{split}
 d(12z-24)-12z+8h+8
 &=h\left[12(d-1/2)^2+5\right]-24d+8\\
 &<\frac9{20}\frac{248}{49}-\frac{72}{7}+8
 =-\frac2{245}<0.
\end{split}
\]
For its quadratic coefficient,
\[
\begin{split}
 T(-28h+20z-40)+40x
 &<\frac{19}{16}\left(-\frac{56}{5}+\frac{20}{3}-40\right)+48\\
 &=-\frac{293}{60}<0.
\end{split}
\]
The last two coefficients are again negative. We have therefore proved
$q_2,L,C<0$ for $0<x\le6/5$. Since $k,h,b>0$,
\eqref{eq:mode-full-factor} yields $\Delta_k(\nu)<0$.
As $a_0>0$, we obtain
$A_2>A_1^2/(4a_0)\ge0$ and hence $P_k>0$ for every real $r$.
Equation~\eqref{eq:all-mode-prefactor} proves the theorem.
\end{proof}
\clearpage\section*{Acknowledgements and disclosures}
\markboth{Acknowledgements and disclosures}{}
\phantomsection
\addcontentsline{toc}{section}{Acknowledgements and disclosures}

The original manuscript reports assistance from Claude (Anthropic) in symbolic
computations, numerical verification and manuscript preparation. The present revision
used ChatGPT (OpenAI) for mathematical review, derivation and integration of direct
estimates, and editorial preparation. The source package contains exact
algebra checks and the finite-mode interval verifier.

\end{document}